\documentclass[final,3p,11pt, times]{article}
\usepackage{authblk}
\usepackage{amsmath}
\usepackage{amssymb,latexsym}
\usepackage[german,english]{babel}
\usepackage{url}
\usepackage{graphicx}
\usepackage{gastex}
\usepackage{longtable}
\usepackage{lscape}
\usepackage{tabularx}
\usepackage{multicol}
\usepackage{verbatim}
\usepackage{multirow}
\usepackage{graphicx}
\usepackage{epsfig,graphics,graphicx}
\usepackage{geometry}
\usepackage[usenames]{xcolor}
\usepackage{authblk}
\usepackage{verbatim}
\usepackage{multirow}
\usepackage{ragged2e}
\usepackage{xcolor}
\usepackage[T1]{fontenc}
\usepackage{authblk}
\usepackage{lipsum}
\usepackage{hyperref}
\hypersetup{colorlinks=pink,linkcolor=blue,filecolor=red,urlcolor=blue,}
\usepackage{tikz}

\usepackage{lineno}
\usepackage{tikz}

\hbadness=\maxdimen

\newtheorem{rem}{{\bf Remark}}
\newtheorem{deff}{Definition}[subsection]
\newtheorem{lem}{{\bf Lemma}}[subsection]
\newtheorem{thm}{{\bf Theorem}}[section]
\newtheorem{prop}{{\bf Proposition}}[section]
\newtheorem{cor}{{\bf Corollary}}[section]

\newtheorem{conj}{{\bf Conjecture}}[section]

\begin{document}
	
\title{{Hamiltonicity of doubly semi-equivelar maps on the torus}}
	
\author{Yogendra Singh, Anand Kumar Tiwari and Seema Kushwaha}

\affil{\small Department of Applied Science, Indian Institute of Information Technology, Allahabad 211015, India. E-mail: $\{$rss2018503,anand,seemak$\}$@iiita.ac.in}
	
\maketitle
	
\hrule
	
\begin{abstract} The well-known twenty types of 2-uniform tilings of the plane give rise infinitely many doubly semi-equivelar maps on the torus. In this article, we show that every such doubly semi-equivelar map on the torus contains a Hamiltonian cycle. As a consequence, we establish the Nash-Williams conjecture for the graphs associated with these doubly semi-equivelar maps by showing that these graphs are either 3-connected or 4-connected.

\end{abstract}
	
\smallskip
	
	
	
\hrule

\section{Introduction} \label{s1}
Let $G$ denote a simple, connected graph. A cycle in $G$ is called a Hamiltonian cycle if it covers all the vertices of the graph. The existence of a Hamiltonian cycle is one of the important properties of a graph and the problem to find a Hamiltonian cycle is NP-complete.

Let $F$ denote a closed surface. An embedding of $G$ into $F$ is called a map $M$ on $F$ if the closure of each component of $F\setminus G$ is a $p$-gonal 2-disk $D_p$,$(p \geq 3)$, also called face, and the non-empty intersection of any such two disks is either a vertex or an edge. A map on $F$ is called a triangulation if each disk is triangular, $i.e.$, $p=3$. A cyclic sequence of faces around a vertex $v$ is called the face-cycle at $v$. In a map, if the face-cycle at each vertex is same, the map is called semi-equivelar. 
 
A graph $G$ is $k$-connected if the removal of $k-1$ vertices leaves a non-trivial connected graph. The connectivity $\kappa(G)$ is the maximum non-negative integer $k$ for which $G$ is $k$-connected. A map is called  $k$-connected (resp. Hamiltonian)  if the underlying graph is $k$-connected (resp. Hamiltonian). 
  

The connectivity of a graph plays an important role to discuss Hamiltonicity.  In the case of maps, significant progress has been made using the connectivity of the underlying graphs. It goes back to 1931 when Whitney \cite{whitney(1931)} initiated the theory of Hamiltonian cycle in maps and he used the connectivity of the underlying graphs. He proved that every triangulation of the sphere is Hamiltonian if it is 4-connected. Tutte \cite{tutte(1956)} generalized this result for an arbitrary 4-connected map on the sphere. Grunbaum \cite{grunbaum(1970)} conjectured that every 4-connected map on the Projective plane is Hamiltonian, which was proved by Thomas and Yu \cite{ThYu(1994)}. Brunet et al.  \cite{BrNaNe(1999)} showed that every 5-connected triangulation on the Klein bottle is Hamiltonian. For the torus, Grunbaum \cite{grunbaum(1970)} and Nash-Williams \cite{NaWi(1973)} posed the following conjecture.


\begin{conj} \normalfont \label{c1}
	Every 4-connected map on the torus is Hamiltonian.	
\end{conj}\label{conj}

The above conjecture is yet to prove, however, several interesting related results have been established so far. Altshuler \cite{altshuler(1972)} established the conjecture for every 6-connected map on the torus. Brunet and Richter \cite{BrRi(1995)} established the conjecture for every 5-connected triangulation on the torus, which was further improved by Thomas and Yu \cite{ThYu(1997)} for any 5-connected map. Kawarabyyashi and Ozeki \cite{KAOZ(2011)} proved the conjecture for any 4-connected triangulation on the torus. The current progress on connectivity and Hamiltonicity of the underlying on the torus can be found in  \cite{KaOz(2016)}, \cite{OzZa(2018)}, \cite{ThYuZa(2005)},  \cite{NaOz(2012)}.

Since the plane is the universal cover of the torus, it is natural to explore the maps on the torus that are associated with the $k$-uniform tilings. In this regard, the eleven 1-uniform tilings, also known as Archimedean tilings, provide semi-equivelar maps on the torus, see \cite{tu(2017)}, \cite{BrKu(2008)}, \cite{MaUp(2018)}. Altshuler \cite{altshuler(1972)} showed that the semi-equivelar maps on the torus of types $[4^4]$ and $[3^6]$ are Hamiltonian. Bouwer and Chernoff \cite{BoCh(1988)} proved that every semi-equivelar map on the torus of type $[6^3]$ is Hamiltonian. Recently, Maity and Upadhyay \cite{MaUp(2020)} showed that every semi-equivelar map on the torus of the remaining eight types $[3.12^2]$, $[4.8^2]$, $[4.6
.12]$, $[3^4.6]$, $[3^3.4^2]$, $[3^2.4.3.4]$, $[3.6.3.6]$ and $[3^4.6]$ are Hamiltonian. Thus, we have the following proposition. 

\begin{prop} \label{p1}
	All the semi-equivelar maps (corresponding to the 1-uniform tilings) on the torus are Hamiltonian.	
\end{prop}

A map with exactly two distinct face-cycles under certain conditions is called a doubly semi-equivelar map or DSEM (defined precisely in Section \ref{s2}). 
There are exactly twenty distinct 2-uniform tilings of the plane \cite{kro(1969)} and every 2-uniform tilling gives doubly semi-equivelar maps on the torus.  In fact, one can construct infinitely many doubly semi-equivelar maps on the torus by taking the quotient of 2-uniform tilings.  
In this article, we show that every doubly semi-equivelar map on the torus corresponding to a 2-uniform tilling  is Hamiltonian. Here, we state our main results. 

\begin{thm} \label{t1}
	All the doubly semi-equivelar maps (corresponding to the 2-uniform tilings) on the torus are Hamiltonian.
\end{thm}

\begin{thm} \label{t2}
	Let $M$ be a doubly semi-equivelar map corresponding to a 2-uniform tilling $T$. 
	\begin{enumerate}
	\item 	If $T$ is one of the following: $[3^3.4^2:4^4]_1$, $[3^3.4^2:4^4]_2$, $[3^6: 3^3.4^2]_1$, $[3^6: 3^3.4^2]_2$, $[3^3.4^2: 3^2.4.3.4]_1$, $[3^3.4^2: 3^2.4.3.4]_2$, $[3^6: 3^2.4.3.4]$, $[3.4^2.6: 3.6.3.6]_1$, $[3.4^2.6: 3.6.3.6]_2$, $[3^2.6^2: 3.6.3.6]$, $[3^6: 3^2.6^2],$ $[3^4.6: 3^2.6^2]$, $[3^2.4.3.4: 3.4.6.4]$, $[3^6: 3^4.6]_1$, $[3^6: 3^4.6]_2$, $[3^3.4^2: 3.4.6.4]$, $[3^6: 3^2.4.12]$ or $[3.4^2.6: 3.4.6.4]$, then $M$ is 4-connected.
	
	\item If $T$ is $[3.4.3.12:3.12^2]$ or $[3.4.6.4: 4.6.12]$, then $M$ is 3-connected.
	\end{enumerate}
\end{thm}

From Theorem \ref{t1} and Theorem \ref*{t2}, we establish Conjecture \ref{c1} for the doubly semi-equivelar maps. In other words:

\begin{cor} 
	All 4-connected doubly semi-equivelar maps (corresponding to 2-uniform tilings) on the torus are Hamiltonian.  
\end{cor}

This article is organized as follows. In Section \ref{s2}, we give some basic definitions 
and notations. In Section \ref{s2.5}, we discuss the connectivity of DSEMs (corresponding to the 2-uniform tilings) on the torus.  In Section \ref{s3}, we describe Hamiltonicity of such DSEMs. To explore a Hamiltonian cycle in a DSEM $M$, we construct a planar representation, called an $M(i,j,k)$ representation, which is obtained by cutting $M$ along two non-homologous cycles at a vertex. Clearly, finding a Hamiltonian cycle in $M$ is equivalent to find a Hamiltonian cycle in its planar representation. We give a detailed description for the DSEMs corresponding to $[3^6:3^3.4^2]_1$ and $[3^6: 3^3.4^2]_2$ and avoid a similar discussion for the remaining types. In support of our argument, we provide a diagram to describe a Hamiltonian cycle (indicated by thick black cycle) in a given DSEM (wherever it is required). 


\section{Definitions and notations} \label{s2}
Let $G$ be the underlying graph of a map $M$ with the vertex set $V(G)$ and the edge set $E(G)$. Let $u$-$v$ denote the edge joining $u, v \in V(G)$. The notation $P = P(u_1, \ldots, u_n)$ denotes a path $u_1$-$u_2$- $\cdots$ -$u_n$. A path with a single vertex $u$ is referred to as a path of length zero, and is denoted as $P = P(u)$. In $P$, the vertices $u_1$ and $u_2$ are called the boundary vertices, and $u_i$, for $2 \leq i \leq n-1$, are called the inner vertices. A path $P_1$ is an extension of another path $P_2$ if $P_2$ is a proper sub graph of $P_1$. A path $P = P(u_1, \ldots, u_n)$ is called a cycle, if it is closed, that is, $u_1 = u_n$. We denote a cycle by $C = C_n(u_1, \ldots, u_n)$. A cycle $C$ is contractible if it bounds a 2-disk $D_p$, otherwise called non-contractible. The notations $G_1 \cup G_2$ and $G_1 \cap G_2$ denote the usual union and intersection of two graphs $G_1$ and $G_2$. For $U \subset V (G)$, $G - U$ is the graph obtained from $G$ by deleting the vertices of $U$. For the graph theory related terms, we refer to Bondy and Murty \cite{BoMu(2008)}. 

\smallskip

The face-sequence of a vertex $v$  in a map $M$, denoted as $f_{seq}(v)$, is a finite sequence $(p_1^{n_1}, \ldots, p_k^{n_k})$ such that  $n_i$ numbers of $D_{p_i}$, for $1 \leq i \leq k$,  incident at $v$ in the given order. A map $M$ is called semi-equivelar of type $[p_1^{n_1}. \ldots. p_k^{n_k}]$ if the face-sequence of each vertex is $(p_1^{n_1}, \ldots, p_k^{n_k})$, see \cite{tu(2017)}.

\smallskip

Let $v$ be a vertex in a map $M$ and $K_v$ denote the set of all the faces $D_{p_i}$, $i \in \{1,2,\ldots,l\}$, incident at $v$. The geometric career $|K_v|$ of $K_v$ is the union of all these faces. Note that $|K_v|$ is a 2-disk with boundary cycle $C_m$, $m \geq 3$. We call  $C_m$  the {\it link of $v$}, and denote it as ${\rm lk}(v)$. Let ${\rm lk}(v) = C_k(v_1, v_2, \ldots, v_k)$. The face-sequence of the link of $v$ is given by $f_{seq} ({\rm lk}(v)) = (f_{seq} (v_1), f_{seq} (v_2), \ldots$, $f_{seq} (v_k))$. Since ${\rm lk}(v)$ is a cycle, for any $v_i, v_j \in {\rm lk}(v)$, there are two possible paths between them in ${\rm lk}(v)$, one is along the clockwise and the other along the anticlockwise direction.

\smallskip
Let $M$ be a map with precisely two distinct face-sequences $f_1$ and $f_2$. Then $M$ is called a doubly semi-equivelar map (or DSEM) if every vertex with the face sequence $f_1$ has the face-sequence of its link as $(f_{11}, f_{12}, \ldots,f_{1k_1})$ and every vertex with face sequence $f_2$ has the face-sequence of its link as $(f_{21}, f_{22},\ldots,f_{2k_2})$. Observe that $f_{ij}$ is either $f_1$ or $f_2$ for each $i,j$. We say $M$ is of type $[f_{1}^{(f_{11},f_{12},\ldots,f_{1k_1})} :f_2^{(f_{21},f_{22},\ldots,f_{2k_2})}]$. Clearly, if $f_1 = f_2$, then $M$ turns into a semi-equivelar map. As discussed earlier, the 2-uniform tilings of the plane induce the same type doubly semi-equivelar maps on the torus. We denote the types of such DSEMs by the same notation as used for the respective tilings.


\section{Connectivity of doubly semi-equivelar maps} \label{s2.5}

The Nash-Williams conjecture states about the Hamiltonicity of a map with respect to the $k$-connectivity of the underlying graph. To establish the conjecture in our context, here, we determine the $k$-connectivity of the graphs associated with certain DSEMs.

\begin{lem}\label{l1}
	
	Let $M$ be a DSEM (corresponding to a 2-uniform tiling of the plane) on the torus. Then removal of any two distinct vertices leaves the map connected.
	
\end{lem}

\noindent{\bf Proof:} Let $V(G)$ denote the vertex set of the underlying graph $G$ of $M$. Note that $G$ is connected. 
Suppose, if possible, there is a vertex $z_1 \in V(G)$ such that $G' = G-z_1$ is disconnected. Then, there exist distinct vertices $u, v \in V(G)$ such that every path between $u$ and $v$ passes through $z_1$. Choose such a path $P = u$-$a_1$-$a_2$- $\cdots$ -$a_{n-1}$-$z_1$-$a_{n+1}$- $\cdots$ -$v$. 
Let ${\rm lk}(z_{1}) = C_i(w_1, w_2, \ldots, w_i)$, $P_1 = u$-$a_1$-$a_2$- $\cdots$ -$a_{n-1}$ and $P_2 = a_{n+1}$-$a_{n+2}$- $\cdots$ -$v$. Note that $a_{n-1}, a_{n+1} \in{\rm lk}(z_{1})$, there is a path $P_3$ between $a_{n-1}$ and $a_{n+1}$ in ${\rm lk}(z_{1})$. 
Now, observe that $\hat{P} = P_1 \cup P_2 \cup P_3$ is a connected sub graph of $G'$ such that $u, v \in V(\hat{P})$ and $z_1 \notin V(\hat{P})$. Therefore $G'$ is connected.

Suppose, if possible, there is a vertex $z_2 \in V(G)$ such that $G''= G'-z_2$ is disconnected. Then, there exist distinct vertices $u', v' \in G''$ such that $z_2$ is a vertex of any path connecting $u'$ and $v'$, say $P' = u'$-$b_1$-$b_2$- $\cdots$ -$b_{n-1}$-$b_n$-$b_{n+1}$- $\cdots$ -$v'$. Then, either $z_2 \in {\rm lk}(z_{1})$ or $z_2 \notin {\rm lk}(z_{1})$. 

First, suppose $z_2 \notin {\rm lk}(z_1)$, then $G''$ is connected by a similar argument used for $G'$. Now suppose, $z_2 \in {\rm lk}(z_1)$. Then $z_2 = w_k$ and $z_2$-$w_{k-1}$, $z_2$-$w_{k+1}$ are edges for some $1 \leq k \leq i$. Since $G'$ is connected, we get paths $P_1'$ from $u'$ to $w_{k+1}$ and $P_2'$ from $w_{k-1}$ to $v'$ which are not passing through $z_2$. Also, we get a path $P_3'$ from $w_{k-1}$ to $w_{k+1}$ in ${\rm lk}(z_2)$. Note that $P_3'$ is a path in $G$ not passing through $z_1$. Thus we get a connected sub graph $\hat{P'} = P_1' \cup P_2' \cup P_3'$ such that $u',v' \in V(\hat{P'})$ and $z_2 \notin V(\hat{P'})$. Hence $G''$ is connected. \hfill$\Box$

\smallskip
The proof of the next theorem follows from Lemma \ref{l1}.

\begin{thm} \label{t3}
	Let $M$ be a doubly semi-equivelar map of type $Y$, where $Y \in \{[3.4.3.12:3.12^2]$, $[3.4.6.4: 4.6.12]\}$. Then $M$ is 3-connected.
	
\end{thm}

\begin{thm} \label{t4}
	Let $M$ be a doubly semi-equivelar map of type $X$, where $X \in \{[3^3.4^2:4^4]_1$, $[3^3.4^2:4^4]_2$, $[3^6: 3^3.4^2]_1$, $[3^6: 3^3.4^2]_2$, $[3^3.4^2: 3^2.4.3.4]_1$, $[3^3.4^2: 3^2.4.3.4]_2$, $[3^6: 3^2.4.3.4]$, $[3.4^2.6: 3.6.3.6]_1$, $[3.4^2.6: 3.6.3.6]_2$, $[3^2.6^2: 3.6.3.6]$, $[3^6: 3^2.6^2],$ $[3^4.6: 3^2.6^2]$, $[3^2.4.3.4: 3.4.6.4]$, $[3^6: 3^4.6]_1$, $[3^6: 3^4.6]_2$, $[3^3.4^2: 3.4.6.4]$, $[3^6: 3^2.4.12]$ or $[3.4^2.6: 3.4.6.4]\}$. Then $M$ is 4-connected.
	
\end{thm}

\noindent{\bf Proof:} 
Let $V(G)$ denote the vertex set of the underlying graph $G$ of $M$, which is connected. Note that the degree of every vertex of $M$ is grater or equal to 4. Let $z_1, z_2 \in V(G)$. Then by Lemma \ref{l1}, $G'= G-\{z_1\}$ and $G''= G-\{z_1, z_2\}$ are connected. Suppose, if possible, there is a vertex $z_3 \in V(G)$ such that $G''' = G - \{z_1, z_2, z_3\}$ is disconnected. Then we get $u, v \in V(G)$ such that there is no path between $u$ and $v$ in $G'''$. Let ${\rm lk}(z_{1}) = C_i(w_1, w_2, \cdots, w_i)$, ${\rm lk}(z_{2}) = C_j(x_1, x_2, \cdots, x_j)$ and ${\rm lk}(z_{3}) = C_k(y_1, y_2, \cdots, y_k)$.
Now, as per the positions of $z_1, z_2, z_3$, we have the following cases:

1. $z_a \notin {\rm lk}(z_{b})$, where $a \neq b$ and $a, b \in \{1,2,3\}$.

2. $z_a \in {\rm lk}(z_{b})$ such that $z_c \notin {\rm lk}(z_{a})$ and $ z_c \notin {\rm lk}(z_{b})$, where $a \neq b \neq c \neq a$ and $a,b, c \in \{1,2,3\}$.

3. $z_2, z_3 \in {\rm lk}(z_{1})$ and $z_2 \in {\rm lk}(z_{3})$.

4. $z_a, z_b \in {\rm lk}(z_{c})$ such that $z_a \notin {\rm lk}(z_{b})$, where $a \neq b \neq c \neq a$ and $a,b, c \in \{1,2,3\}$. \smallskip

In Case 1 and 2, $G'''$ remains connected by Lemma \ref{l1}. 
So we discuss the remaining cases.

\noindent{\bf Case 3:} Since $z_2, z_3 \in {\rm lk}(z_{1})$, $z_2= w_h$, $z_3 = w_l$  for some $ 1 \leq h, l \leq i$. Note that, $w_{h-1}, w_{h+1}$ are adjacent to $z_2$ and $w_{l-1}, w_{l+1}$ are adjacent to $z_3$. Hence, $w_{h-1}, w_{h+1} \in {\rm lk}(z_{2})$, $w_{l-1}, w_{l+1} \in {\rm lk}(z_{3})$ and $w_{h-1} = x_m, w_{h+1} = x_n$, $w_{l-1} = y_r, w_{l+1} = y_s$ for some  $1 \leq m, n \leq j$, $1 \leq r, s \leq k$.  Let $R = x_{m}$-$x_{m+1}$- $\cdots$ -$x_{n-1}$-$x_{n}$ be the path in ${\rm lk}(z_{2})$ from $x_m$ to $x_n$ such that $z_1 \notin V(R)$ and $S = y_r$-$y_{r+1}$- $\cdots$ -$y_{s-1}$-$y_{s}$ be the path in ${\rm lk}(z_{3})$ from $y_r$ to $y_s$ such that $z_1 \notin V(S)$.

\textbf{Subcase 3.1:} If $z_2, z_3$ are not adjacent such that $w_{l-1}, w_{l+1} \notin \{w_{h-1}, w_{h+1}\}$. Then $A_1 = w_{l+1}$-$w_{l+2}$- $\cdots$ -$w_{h-2}$-$w_{h-1}$ and $A_2 = w_{h+1}$-$w_{h+2}$- $\cdots$ -$w_{l-2}$-$w_{l-1}$ are paths in ${\rm lk}(z_{1})$ such that $z_2, z_3 \notin V(A_1) \cup V(A_2)$. Since $z_2 \in {\rm lk}(z_{3})$ and $z_3 \in {\rm lk}(z_{2})$, we have $z_2 = y_p, z_3 = x_q$ for some $1 \leq p \leq k$, $1 \leq q \leq j$ and $z_2$-$y_{p+1}$, $z_2$-$y_{p-1}$, $z_3$-$x_{q-1}$, $z_3$-$x_{q+1}$ are edges.  Consider the path $R = x_{m}$-$x_{m+1}$- $\cdots$ -$x_{q-1}$-$x_q(=z_3)$-$x_{q+1}$- $\cdots$ -$x_{n-1}$-$x_{n}$ in ${\rm lk}(z_{2})$, then we get sub paths $A_3 = x_{q+1}$-$x_{q+2}$- $\cdots$ -$x_n(=w_{h+1})$ and $A_4 = x_{q-1}$-$x_{q-2}$- $\cdots$ -$x_m(=w_{h-1})$ of $R$. Similarly, the path $S = y_r$-$y_{r+1}$- $\cdots$ -$y_{p-1}$-$y_p(=z_2)$-$y_{p+1}$- $\cdots$ -$y_{s-1}$-$y_{s}$ in ${\rm lk}(z_{3})$ gives sub paths $A_5 = y_{p+1}$-$y_{p+2}$- $\cdots$ -$y_s(=w_{l+1})$ and $A_6 = y_{p-1}$-$y_{p-2}$- $\cdots$ -$y_r(=w_{l-1})$. In this manner, we get connected sub graphs $W_1 = A_3 \cup A_2 \cup A_6$ and $W_2 = A_5 \cup A_1 \cup A_4$ of $G'''$. 
If $V(W_1) \cap V(W_2) \neq \phi$, then $W_1\cup W_2$ is connected, so $G'''$ is connected. If $V(W_1) \cap V(W_2) = \phi$, then the links of all the vertices of $G'''$ remain complete except for the vertices of $W_1$ and $W_2$ as they belong to the links of $z_1, z_2$ or $z_3$. Note that there exists either a vertex or an edge in $W_1$ contained in a face $F_1$, where $F_1$ is a face in $M$ that does not contain the vertices $z_1$, $z_2$ or $z_3$. Similarly, for a vertex or an edge in $W_2,$ we get a face $H_1.$
Corresponding to the faces $F_1, H_1$, we get faces $F_2, H_2$ such that $F_1 \cap F_2 \neq \phi, H_1 \cap H_2 \neq \phi$, otherwise links of some vertices of these faces are not complete. 
In this manner, we get  finite sequences of faces  $F_1, F_2, \cdots, F_e$ and $H_1, H_2, \cdots, H_f$ corresponding to $W_1$ and $W_2$, respectively. Since $M$ is a map on the finite number of vertices, after some stage, we get $F_u \cap H_v \neq \phi$, for some $1 \leq u \leq e$, $1 \leq v \leq f$. Therefore, every vertex of $W_1$ is connected to every vertex of $W_2$. Hence, $G'''$ is connected.

If $z_2, z_3$ are not adjacent such that either $w_{l-1} = w_{h+1}$ or $w_{l+1} = w_{h-1}$, then proceeding similarly, we get connected sub graphs $W'_1$ and $W'_2$ of $G'''$. By using the similar argument as above, $G'''$ is connected.

\textbf{Subcase 3.2:} Let $z_2, z_3$ are adjacent. As $z_2, z_3 \in {\rm lk}(z_{1})$, $z_2 = w_h, z_3 = w_{h+1}$ for some $1 \leq h \leq i$ and $z_2$-$w_{h-1}$, $z_2$-$z_3(=w_{h+1})$, $z_3$-$w_{h+2}$ are edges. Clearly, there is a path $B_1=w_{h+2}$-$w_{h+3}$- $\cdots$ -$w_{h-2}$-$w_{h-1}$ between $w_{h+2}$ and $w_{h-1}$ such that $z_1, z_2, z_3 \notin V(B_1)$. Since $z_3 \in {\rm lk}(z_{1})$ and $z_3$-$w_{h+2}$ is an edge, so $w_{h+2} \in {\rm lk}(z_{3})$. Now $z_2(=w_h)$-$z_3(=w_{h+1})$ is an edge, so $z_2 \in {\rm lk}(z_{3})$. If $z_2=y_p$ for some $1 \leq p \leq k$ in ${\rm lk}(z_{3})$, then $B_2=y_{p+1}$-$y_{p+2}$- $\cdots$ -$w_{h+2}$ is a path between $y_{p+1}$-$w_{h+2}$ such that $z_1,z_2, z_3 \notin V(B_2)$. As $w_{h-1}$-$z_2$ and $y_{p+1}$-$z_2$ are edges, this gives $w_{h-1}, y_{p+1} \in {\rm lk}(z_{2})$. Hence, we have a path $B_3$ between $w_{h-1}$ and $y_{p+1}$. Thus, from above, we get a connected sub graph $W= B_1 \cup B_3 \cup B_2$ of $G'''$. 
Since $G''$ is connected, every vertices of $G'''$ is connected with some vertices of $W$, therefore $G'''$ is connected.

\textbf{Case 4:} 
Without loss of generality, suppose that $z_2, z_3 \in {\rm lk}(z_{1})$ and $z_2 \notin {\rm lk}(z_{3})$. As discussed in Case 3, $z_2, z_3 \in {\rm lk}(z_{1})$ gives $z_2= w_h$, $z_3 = w_l$ for some $ 1 \leq h, l \leq i$,  $w_{h-1}, w_{h+1} \in {\rm lk}(z_{2})$, $w_{l-1}, w_{l+1} \in {\rm lk}(z_{3})$ and $w_{h-1} = x_m, w_{h+1} = x_n$, $w_{l-1} = y_r, w_{l+1} = y_s$ for some  $1 \leq m, n \leq j$, $1 \leq r, s \leq k$.  Also, we get a path $R = x_{m}$-$x_{m+1}$- $\cdots$ -$x_{n-1}$-$x_{n}$ in ${\rm lk}(z_{2})$ from $x_m$ to $x_n$ such that $z_1 \notin V(R)$ and  a path $S = y_r$-$y_{r+1}$- $\cdots$ -$y_{s-1}$-$y_{s}$ in ${\rm lk}(z_{3})$ from $y_r$ to $y_s$ such that $z_1 \notin V(S)$.

If $z_2, z_3$ are such that $w_{l-1}, w_{l+1} \notin \{w_{h-1}, w_{h+1}\}$. Then $A_1 = w_{l+1}$-$w_{l+2}$- $\cdots$ -$w_{h-2}$-$w_{h-1}$ and $A_2 = w_{h+1}$-$w_{h+2}$- $\cdots$ -$w_{l-2}$-$w_{l-1}$ are paths in ${\rm lk}(z_{1})$ such that $z_2, z_3 \notin V(A_1) \cup V(A_2)$. Clearly, $Z= R \cup A_2 \cup S \cup A_3$ is a connected sub graph of $G'''$. Therefore $G'''$ is connected.

If $z_2, z_3$ are such that either $w_{l-1} = w_{h+1}$ or $w_{l+1} = w_{h-1}$, then proceeding similarly, we get a connected sub graph $Z'$ of $G'''$. By using the similar argument as above, $G'''$ is connected.

Hence, by the Cases 1-4, $M$ is 4-connected. \hfill$\Box$

\vspace{.2cm}

\noindent{\bf Proof of Theorem \ref{t2}:} The proof follows from Theorem \ref{t3} and Theorem \ref{t4}. \hfill$\Box$

\section{Hamiltonicity of doubly semi-equivelar maps on the torus} \label{s3}

Throughout this section by a DSEM, we mean a doubly semi-equivelar map on the torus. Let $M$ be a DSEM of type $[f_1: f_2]$ with vertex set $V(M)$. For $j=1,2$, let $V_{f_j}$ denote the set of vertices with face-sequence $f_j$. The notations $|V(M)|$, $|V_{f_j}|$ denote the cardinality of the vertex set $V(M)$, $V_{f_j}$ respectively.

\subsection{DSEMs of type $[3^6:3^3.4^2]_1$ and $[3^6:3^3.4^2]_2$} \label{s3.1}

Let $M^r$ be a DSEM of type $[3^6:3^3.4^2]_r$ with the vertex set $V(M^r)$, where $r \in \{1,2\}$. For $M^1$, we see that the number of triangular faces is $4 |V_{(3^6)}|$ or $2 |V_{(3^3,4^2)}|$. Thus, if $M^1$ exists, then $2|V_{(3^6)}| = |V_{(3^3,4^2)}|$. Similarly, for the existence of $M^2$, we get $|V_{(3^6)}| = |V_{(3^3,4^2)}|$. Let $u$ be a vertex with face-sequence $(3^6)$ or $(3^3,4^2)$. We denote their respective links by ${\rm lk}(u) = C_6(u_1,u_2,u_3,u_4,u_5,u_6)$ or ${\rm lk}(u) = C_7(\boldsymbol{u_1},u_2,\boldsymbol{u_3},u_4,u_5,u_6,u_7)$ (labeling of vertices in the links is considered anti-clockwise). The bold appearance of some $u_i$ means $u$ is not adjacent with $u_i$ by an edge.

Considering Figures 4.1.1 and 4.1.2 as the planar drawings of the maps $M^1$ and $M^2$ respectively, we define certain types of paths as follows.

\begin{deff} \label{d4.1.1} \normalfont
A path $P_{1} = P( \ldots, y_{i-1},y_{i},y_{i+1}, \ldots)$ in $M^r$ is of type $A_{1}$ if, either
$(i)$ every vertex of $P_1$ has the face-sequence $(3^6)$ or $(ii)$ each vertex of $P_1$ has the face-sequence $(3^3,4^2)$ such that all the triangles (square) which incident on its inner vertices lie on the one side  of the path $P_1$. See thick black paths in Figure 4.1.1 and Figure 4.1.2.
\end{deff}

\begin{picture}(0,0)(-35,11.75)
\setlength{\unitlength}{5.25mm}

\drawpolygon(1,0)(.5,1)(-.5,1)(-1,0)(-.5,-1)(.5,-1)
\drawpolygon(3,0)(2.5,1)(1.5,1)(1,0)(1.5,-1)(2.5,-1)
\drawpolygon(5,0)(4.5,1)(3.5,1)(3,0)(3.5,-1)(4.5,-1)

\drawline[AHnb=0](4.5,1)(4.5,2)
\drawline[AHnb=0](4.5,-1)(4.5,-2)

\drawline[AHnb=0](-.5,1)(-.5,2)
\drawline[AHnb=0](-.5,-1)(-.5,-2)

\drawline[AHnb=0](.5,-1)(.5,-2)
\drawline[AHnb=0](1.5,-1)(1.5,-2)
\drawline[AHnb=0](2.5,-1)(2.5,-2)
\drawline[AHnb=0](3.5,-1)(3.5,-2)

\drawline[AHnb=0](.5,1)(.5,2)
\drawline[AHnb=0](1.5,1)(1.5,2)

\drawline[AHnb=0](2.5,1)(2.5,2)
\drawline[AHnb=0](3.5,1)(3.5,2)


\drawline[AHnb=0](-1.11,0)(5.11,0)
\drawline[AHnb=0](-1,-2)(5,-2)
\drawline[AHnb=0](-1,2)(5,2)
\drawline[AHnb=0](-1,1)(5,1)


\drawline[AHnb=0](.5,1)(-.5,-1)
\drawline[AHnb=0](-.5,1)(.5,-1)

\drawline[AHnb=0](1.5,-1)(2.5,1)
\drawline[AHnb=0](2.5,-1)(1.5,1)

\drawline[AHnb=0](3.5,-1)(4.5,1)
\drawline[AHnb=0](4.5,-1)(3.5,1)


\drawline[AHnb=0](.5,2)(.75,2.5)
\drawline[AHnb=0](.5,2)(.25,2.5)

\drawline[AHnb=0](-.5,2)(-.75,2.5)
\drawline[AHnb=0](-.5,2)(-.25,2.5)

\drawline[AHnb=0](1.5,2)(1.75,2.5)
\drawline[AHnb=0](1.5,2)(1.25,2.5)

\drawline[AHnb=0](3.5,2)(3.75,2.5)
\drawline[AHnb=0](3.5,2)(3.25,2.5)

\drawline[AHnb=0](4.5,2)(4.75,2.5)
\drawline[AHnb=0](4.5,2)(4.25,2.5)

\drawline[AHnb=0](2.5,2)(2.25,2.5)
\drawline[AHnb=0](2.5,2)(2.75,2.5)

\drawline[AHnb=0](.5,-2)(.75,-2.5)
\drawline[AHnb=0](.5,-2)(.25,-2.5)

\drawline[AHnb=0](-.5,-2)(-.75,-2.5)
\drawline[AHnb=0](-.5,-2)(-.25,-2.5)

\drawline[AHnb=0](1.5,-2)(1.75,-2.5)
\drawline[AHnb=0](1.5,-2)(1.25,-2.5)

\drawline[AHnb=0](3.5,-2)(3.75,-2.5)
\drawline[AHnb=0](3.5,-2)(3.25,-2.5)

\drawline[AHnb=0](4.5,-2)(4.75,-2.5)
\drawline[AHnb=0](4.5,-2)(4.25,-2.5)

\drawline[AHnb=0](2.5,-2)(2.25,-2.5)
\drawline[AHnb=0](2.5,-2)(2.75,-2.5)


\drawpolygon[fillcolor=black](-1,2)(5,2)(5,2.1)(-1,2.1)

\drawpolygon[fillcolor=black](-1,1)(5,1)(5,1.1)(-1,1.1)

\drawpolygon[fillcolor=black](-1,-1)(5,-1)(5,-1.1)(-1,-1.1)

\drawpolygon[fillcolor=black](-1.12,0)(5.1,0)(5.1,.1)(-1.12,.1)

\drawpolygon[fillcolor=black](-1,-2)(5,-2)(5,-2.1)(-1,-2.1)

\drawpolygon[fillcolor=green](.75,-2.5)(.5,-2)(.5,-1)(-.5,1)(-.5,2)(-.75,2.5)(-.65,2.5)(-.4,2)(-.4,1)(.6,-1)(.6,-2)(.85,-2.5)

\drawpolygon[fillcolor=green](1.75,-2.5)(1.5,-2)(1.5,-1)(.5,1)(.5,2)(.25,2.5)(.35,2.5)(.6,2)(.6,1)(1.6,-1)(1.6,-2)(1.85,-2.5)

\drawpolygon[fillcolor=green](2.75,-2.5)(2.5,-2)(2.5,-1)(1.5,1)(1.5,2)(1.25,2.5)(1.35,2.5)(1.6,2)(1.6,1)(2.6,-1)(2.6,-2)(2.85,-2.5)

\drawpolygon[fillcolor=green](3.75,-2.5)(3.5,-2)(3.5,-1)(2.5,1)(2.5,2)(2.25,2.5)(2.35,2.5)(2.6,2)(2.6,1)(3.6,-1)(3.6,-2)(3.85,-2.5)

\drawpolygon[fillcolor=green](4.75,-2.5)(4.5,-2)(4.5,-1)(3.5,1)(3.5,2)(3.25,2.5)(3.35,2.5)(3.6,2)(3.6,1)(4.6,-1)(4.6,-2)(4.85,-2.5)

\put(1,2.65) {\scriptsize {\tiny $A_2$}}
\put(5.5,-1.1) {\scriptsize {\tiny $A_1$}}
\put(5.5,-.1) {\scriptsize {\tiny $A_1$}}

\put(-.7,-3.4){\scriptsize {\tiny {\bf Figure 4.1.1:} $[3^6 : 3^3.4^2]_1$}}
\end{picture}

\begin{picture}(0,0)(-90,4)
\setlength{\unitlength}{5mm}


\drawline[AHnb=0](-1,-2)(-1,-3)
\drawline[AHnb=0](0,-2)(0,-3)
\drawline[AHnb=0](1,-2)(1,-3)
\drawline[AHnb=0](2,-2)(2,-3)
\drawline[AHnb=0](3,-2)(3,-3)
\drawline[AHnb=0](4,-2)(4,-3)
\drawline[AHnb=0](5,-2)(5,-3)

\drawline[AHnb=0](-.5,1)(-.5,2)
\drawline[AHnb=0](.5,1)(.5,2)
\drawline[AHnb=0](1.5,1)(1.5,2)
\drawline[AHnb=0](2.5,1)(2.5,2)
\drawline[AHnb=0](3.5,1)(3.5,2)
\drawline[AHnb=0](4.5,1)(4.5,2)


\drawline[AHnb=0](-1,-3)(5,-3)
\drawline[AHnb=0](-1,-2)(5,-2)
\drawline[AHnb=0](-1,1)(5,1)
\drawline[AHnb=0](-1,2)(5,2)
\drawline[AHnb=0](-1.11,0)(5.11,0)


\drawline[AHnb=0](-.75,.5)(-.5,1)
\drawline[AHnb=0](-1,-2)(.5,1)
\drawline[AHnb=0](0,-2)(1.5,1)
\drawline[AHnb=0](1,-2)(2.5,1)
\drawline[AHnb=0](2,-2)(3.5,1)
\drawline[AHnb=0](3,-2)(4.5,1)
\drawline[AHnb=0](4,-2)(4.75,-.5)

\drawline[AHnb=0](-.75,-.5)(0,-2)
\drawline[AHnb=0](-.5,1)(1,-2)
\drawline[AHnb=0](.5,1)(2,-2)
\drawline[AHnb=0](1.5,1)(3,-2)
\drawline[AHnb=0](2.5,1)(4,-2)
\drawline[AHnb=0](3.5,1)(5,-2)
\drawline[AHnb=0](4.5,1)(4.75,0.5)


\drawline[AHnb=0](.5,2)(.75,2.5)
\drawline[AHnb=0](.5,2)(.25,2.5)

\drawline[AHnb=0](-.5,2)(-.75,2.5)
\drawline[AHnb=0](-.5,2)(-.25,2.5)

\drawline[AHnb=0](1.5,2)(1.75,2.5)
\drawline[AHnb=0](1.5,2)(1.25,2.5)

\drawline[AHnb=0](3.5,2)(3.75,2.5)
\drawline[AHnb=0](3.5,2)(3.25,2.5)

\drawline[AHnb=0](4.5,2)(4.75,2.5)
\drawline[AHnb=0](4.5,2)(4.25,2.5)

\drawline[AHnb=0](2.5,2)(2.25,2.5)
\drawline[AHnb=0](2.5,2)(2.75,2.5)

\drawline[AHnb=0](-1,-3)(-.75,-3.5)
\drawline[AHnb=0](-1,-3)(-1.25,-3.5)

\drawline[AHnb=0](1,-3)(1.25,-3.5)
\drawline[AHnb=0](1,-3)(.75,-3.5)

\drawline[AHnb=0](0,-3)(-.25,-3.5)
\drawline[AHnb=0](0,-3)(.25,-3.5)

\drawline[AHnb=0](2,-3)(2.25,-3.5)
\drawline[AHnb=0](2,-3)(1.75,-3.5)

\drawline[AHnb=0](4,-3)(4.25,-3.5)
\drawline[AHnb=0](4,-3)(3.75,-3.5)

\drawline[AHnb=0](5,-3)(5.25,-3.5)
\drawline[AHnb=0](5,-3)(4.75,-3.5)

\drawline[AHnb=0](3,-3)(2.75,-3.5)
\drawline[AHnb=0](3,-3)(3.25,-3.5)


\drawpolygon[fillcolor=black](-1,2)(5,2)(5,2.1)(-1,2.1)

\drawpolygon[fillcolor=black](-1,1)(5,1)(5,1.1)(-1,1.1)

\drawpolygon[fillcolor=black](-1,0)(5,0)(5,.1)(-1,.1)

\drawpolygon[fillcolor=black](-1,-1)(5,-1)(5,-1.1)(-1,-1.1)

\drawpolygon[fillcolor=black](-1.2,-2)(5.2,-2)(5.2,-2.1)(-1.2,-2.1)
\drawpolygon[fillcolor=black](-1.2,-3)(5.2,-3)(5.2,-3.1)(-1.2,-3.1)

\drawpolygon[fillcolor=green](4.25,-3.5)(4,-3)(4,-2)(2.5,1)(2.5,2)(2.25,2.5)(2.35,2.5)(2.6,2)(2.6,1)(4.1,-2)(4.1,-3)(4.35,-3.5)

\drawpolygon[fillcolor=green](5.25,-3.5)(5,-3)(5,-2)(3.5,1)(3.5,2)(3.25,2.5)(3.35,2.5)(3.6,2)(3.6,1)(5.1,-2)(5.1,-3)(5.35,-3.5)

\drawpolygon[fillcolor=green](3.25,-3.5)(3,-3)(3,-2)(1.5,1)(1.5,2)(1.25,2.5)(1.35,2.5)(1.6,2)(1.6,1)(3.1,-2)(3.1,-3)(3.35,-3.5)

\drawpolygon[fillcolor=green](2.25,-3.5)(2,-3)(2,-2)(.5,1)(.5,2)(.25,2.5)(.35,2.5)(.6,2)(.6,1)(2.1,-2)(2.1,-3)(2.35,-3.5)

\drawpolygon[fillcolor=green](1.25,-3.5)(1,-3)(1,-2)(-.5,1)(-.5,2)(-.75,2.5)(-.65,2.5)(-.4,2)(-.4,1)(1.1,-2)(1.1,-3)(1.35,-3.5)

\put(1,2.65) {\scriptsize {\tiny $A_2$}}
\put(5.5,-2.1) {\scriptsize {\tiny $A_1$}}
\put(5.5,-1.1) {\scriptsize {\tiny $A_1$}}

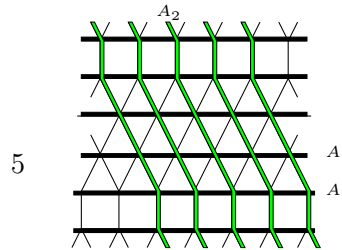
\put(-1,-4.3){\scriptsize {\tiny {\bf Figure 4.1.2:} $[3^6 : 3^3.4^2]_2$}}
\end{picture}

\vspace{2cm}

\begin{deff} \label{d4.1.2} \normalfont A path $P_{2} = P(\ldots, z_{i-1},z_{i},z_{i+1}, \ldots)$ in $M^1$, such that $z_{i-1}, z_i, z_{i+1}$ are inner vertices of $P_2$ or an extended path of $P_2$, is of type $A_{2}$ (shown by green colored paths in Figure 4.1.1), if either of the following three conditions follows for each vertex of the path.

	\begin{enumerate}

		\item  If ${\rm lk}(z_{i})=C_7(\boldsymbol{m},z_{i-1},\boldsymbol{n},o,p,z_{i+1},q)$ and  ${\rm lk}(z_{i-1})=C_7(\boldsymbol{q},m,r,z_{i-2},n,\boldsymbol{o},z_{i})$,  then ${\rm lk}(z_{i+1})=C_6(s,q,z_i,p,t,z_{i+2})$.
		
		\item If ${\rm lk}(z_{i})=C_7(\boldsymbol{m},z_{i+1},\boldsymbol{n},o,p,z_{i-1},q)$ and  ${\rm lk}(z_{i-1})=C_6(r,z_{i-2},s,q,z_{i},p)$, then ${\rm lk}(z_{i+1})=C_7(\boldsymbol{o},z_{i},\boldsymbol{q},m,t,z_{i+2},n)$.

		\item  If ${\rm lk}(z_{i})=C_6(z_{i-1},m,n,z_{i+1},o,p)$ and ${\rm lk}(z_{i-1})=C_7(\boldsymbol{r},z_{i-2},\boldsymbol{s},t,m,z_{i},p)$, then ${\rm lk}(z_{i+1})=C_7(\boldsymbol{u},z_{i+2},\boldsymbol{v},w,o,z_{i},n)$.
		
	\end{enumerate}
	
\end{deff}


\begin{deff} \label{d4.1.3} \normalfont A path $P_{3} = P(\ldots, z_{i-1},z_{i},z_{i+1}, \ldots)$ in $M^2$, such that $z_{i-1}, z_{i}, z_{i+1}$ are inner vertices of $P_{3}$ or an extended path of $P_{3}$, is of type $A_{2}$ (shown by green colored  paths in Figure 4.1.2), if either of the following four conditions follows for each vertex of the path.
	
	\begin{enumerate} \item  If ${\rm lk}(z_{i})=C_7(\boldsymbol{m},z_{i-1},\boldsymbol{n},o,p,z_{i+1},q)$ and ${\rm lk}(z_{i-1})=C_7(r,z_{i-2},n,\boldsymbol{o},z_i,\boldsymbol{q},m)$, then ${\rm lk}(z_{i+1})=C_6(z_{i},p,s,z_{i+2},t,q), {\rm lk}(z_{i+2})=C_6(z_{i+1},s,u,z_{i+3},v,t)$.

	\item If ${\rm lk}(z_{i})=C_7(\boldsymbol{m},z_{i+1},\boldsymbol{n},o,p,z_{i-1},q)$ and ${\rm lk}(z_{i-1})=C_6(z_{i-2},s,q,z_{i},p,r)$, then ${\rm lk}(z_{i+1})=C_7(\boldsymbol{o},z_{i}, \boldsymbol{q},m,t,z_{i+2},n)$, ${\rm lk}(z_{i+2})=C_6(z_{i+1},t,u,z_{i+3},v,n)$.
		
	\item If  ${\rm lk}(z_{i})=C_6(z_{i-1},m,n,z_{i+1},o,p)$ and ${\rm lk}(z_{i-1})=C_7(\boldsymbol{q},z_{i-2}, \boldsymbol{r},s,m,z_{i},p)$, then ${\rm lk}(z_{i+1})=C_6(z_{i},n,t,z_{i+2},u,o)$, ${\rm lk}(z_{i+2})=C_7(\boldsymbol{v},z_{i+3}, \boldsymbol{w},x,u,z_{i+1},t)$.
		
	\item  If ${\rm lk}(z_{i})=C_6(z_{i-1},m,n,z_{i+1},o,p)$ and ${\rm lk}(z_{i-1})=C_6(z_{i-2},r,m,z_{i},p,q)$, then  
	${\rm lk}(z_{i+1})=C_7(\boldsymbol{s},z_{i+2},\boldsymbol{t},u,o,z_{i},n)$, ${\rm lk}(z_{i+2})=C_7(\boldsymbol{u}, z_{i+1},\boldsymbol{n},s,v,z_{i+3},t)$.
	\end{enumerate}
\end{deff}

\smallskip



\begin{rem}
	Recall that $M^r$ is a map with a finite vertex set. Hence, for a given  maximal path $P=$ $(v_1$-$v_2$- $\cdots$- $v_{r-1}$-$v_r)$ of type $A_{\alpha}$, $\alpha \in \{1,2\}$, there is an edge $e = v_1$-$v_r$ in $M^r$ such that  $P\cup e$ is a cycle $Q$. Depending on $A_\alpha$, the corresponding cycle is called a cycle of type $A_{\alpha}$. Then we see that such cycles are non-contractible.
\end{rem}

\begin{lem} \label{l4.1.1}
	The cycle $Q$ of type $A_{\alpha}$, for $\alpha \in \{1,2\}$, is non-contractible.
\end{lem}

\noindent{\bf Proof.}  On the contrary, suppose $Q$ is a contractible cycle of type $A_{1}$. Then $Q$ is the boundary of a 2-disk, say $D$. Let $v$, $e$, and $f$ denote the number of vertices, edges, and faces of $D$ respectively. Let $m$ and $n$ represent the inner and boundary vertices of $D$, respectively. Here, $m= m_1+m_2$, in which  $m_1$ and $m_2$ represent the number of vertices with face-sequences of $(3^6)$ and $(3^3,4^2)$, respectively. Then we have two cases depending on $Q$. In the first case, if the incident faces on $Q$ are quadrangles,  then $v=n+m_{1}+m_{2}$, $e= 3n/2 + 6m_{1}/2 + 5m_{2}/2$, and  $f=2n/4 + 6m_{1}/3 + 3m_{2}/3 + 2m_{2}/4$. In the second case, if the incident faces on $Q$ are triangles, then $v= n + m_{1}+ m_{2}$, $e= 4n/2 + 6m_{1}/2 + 5m_{2}/2$, and  $f= 3n/3 + 6m_{1}/3 + 2m_{2}/4 +3m_{2}/3$. For both cases, we obtain the Euler characteristic $\chi(D) = v-e+f=0$. We get a contradiction, as the Euler characteristic of a 2-disk is 1. Hence, $Q$ is not contractible. Similarly, if $Q$ is of type $A_2$, then $Q$ is non-contractible.  \hfill $\Box$

\vspace{.2cm}



Let $Q$ be a cycle of type $A_1$ with the vertex set $V(Q)$. Let $S_Q$ denote the set of all the faces incident at $v$ for all $v \in V(Q)$. 
Then the geometric carrier $|S_{Q}|$ is a cylinder, as $Q$ is non-contractible. Let $\partial |S_{Q}|$ denote the boundary of $|S_{Q}|$. 

Let $Q_1$ and $Q_2$ be  cycles of the same type in $M^r$. We say that cycles $Q_{1}$ and $Q_{2}$ are homologous if there is a cylinder whose boundary is $\{Q_1, Q_2\}$.

\begin{lem} \label{l4.1.2} Let $Q$ be a cycle of type $A_1$ such that $\partial |S_{Q}| = \{Q_{1}, Q_{2}\}$. Then $Q$, $Q_1$, and $Q_2$ are of the same type with equal length.
\end{lem}

\noindent{\bf Proof.}  Let $Q_n \in \{Q_1, Q_2\}$. Without loss of generality, let $Q_n = Q_1$, consider the faces that are incident with both $Q$ and $Q_1$. Now, depending on $Q$, we see that if the faces incident with both $Q$ and $Q_1$ are quadrangles (respectively  triangles), then the faces on the other side of $Q_1$ must be triangles (respectively triangles or quadrangles). This shows that $Q_1$ is of type $A_1$. Similarly, $Q_2$ is also of type $A_1$.

 

Let $Q = C_l(v_1, \ldots, v_l)$, $Q_{1} = C_{m}(u_1, \ldots, u_{m})$, and $Q_2 = C_{n}(w_1, \ldots, w_{n})$. Now we show $l=m=n$ by contradiction. For this, it is enough to deal with the case, when $l < m < n$. Following the definition of type $A_1$, the face-sequences of $v_{1},v_{2}, \ldots, v_{l-1},v_{l}$ are same throughout the cycle $Q$. Since $l < m$, ${\rm lk}(v_{l})$ contains the vertices $u_{l},u_{l+i}$, and $w_{l}$ for some $i > 0$. Then, the face-sequences of $v_{l}$ and $v_{l-1}$ are not same. This is not possible. Therefore, $l = m = n$.    \hfill$\Box$ 

\vspace{.2cm}

{\bf $M(i,j,k)$-representation of $M^r$:} Let $u \in V(M^r)$ and $Q_{\alpha}$ be cycles of type $A_\alpha$ through $u$, where $ \alpha \in \{1,2\}$. Let $Q_{1} = C_i(u_1,u_2, \ldots, u_i)$. We cut $M^r$ first along the cycle $Q_{1}$. This gives a cylinder, say $R_1$, bounded by identical cycle $Q_{1}$. We say that a cycle is  horizontal (resp. vertical) if it is $Q_{1}$ or homologous to $Q_{1}$ (resp. $Q_2$ or homologous to $Q_2$). In $R_{1}$, starting from the vertex $u$, make another cut along the path $P \subset Q_{2}$, until it reaches $Q_{1}$ again for the first time. As a result, we get the unfold torus (a planer representation), say $R_2$. The idea of $M(i,j,k)$-representation is taken from \cite{MaUp(2018)} and a detailed description of such representation for the given type DSEMs is avalible in \cite{YoTi(2020)}.

\smallskip
Without loss of generality, suppose that the quadrangular faces are incident on $Q_1$. In $R_2$, let there are $j$ cycles, as shown in
Fig. 4.1.3 and Fig. 4.1.4, which are homologous to $Q_1$ along $P$. Since length$(Q_1) = i$ and the number of horizontal cycles along $P$ is $j$. So, we denote $R_{2}$ by $(i,j)$ representation. To reconstruct the map $M^r$ from its $(i, j)$ representation, one gets a natural way of identification of the vertical cycle of $R_2$, but identification of the horizontal cycle generally needs some shifting so that the lower horizontal cycle is identified with the upper horizontal cycle. Suppose $k$ is such shifting, $i.e.$, $k+1$ is the starting vertex of the upper horizontal cycle. This gives another representation say $M^r(i, j, k)$ representation of the $(i, j)$ representation of $M^r$. The admissible relations among $i, j, k$ of the $M^r(i,j,k)$ is given in Lemma \ref{l4.1.3}.

\begin{picture}(0,0)(-14,47.75)
\setlength{\unitlength}{6mm}

\drawpolygon(0,0)(7,0)(7,6)(0,6)

\drawline[AHnb=0](0,1)(7,1)
\drawline[AHnb=0](0,2)(7,2)
\drawline[AHnb=0](0,3)(7,3)
\drawline[AHnb=0](0,4)(7,4)
\drawline[AHnb=0](0,5)(7,5)
\drawline[AHnb=0](0,6)(7,6)

\drawline[AHnb=0](1,0)(1,6)
\drawline[AHnb=0](2,0)(2,6)
\drawline[AHnb=0](3,0)(3,6)
\drawline[AHnb=0](4,0)(4,6)
\drawline[AHnb=0](5,0)(5,6)
\drawline[AHnb=0](6,0)(6,6)

\drawline[AHnb=0](0,1)(1,2)
\drawline[AHnb=0](1,1)(2,2)
\drawline[AHnb=0](2,1)(3,2)
\drawline[AHnb=0](3,1)(4,2)
\drawline[AHnb=0](4,1)(5,2)
\drawline[AHnb=0](5,1)(6,2)
\drawline[AHnb=0](6,1)(7,2)

\drawline[AHnb=0](0,2)(1,3)
\drawline[AHnb=0](1,2)(2,3)
\drawline[AHnb=0](2,2)(3,3)
\drawline[AHnb=0](3,2)(4,3)
\drawline[AHnb=0](4,2)(5,3)
\drawline[AHnb=0](5,2)(6,3)
\drawline[AHnb=0](6,2)(7,3)

\drawline[AHnb=0](0,4)(1,5)
\drawline[AHnb=0](1,4)(2,5)
\drawline[AHnb=0](2,4)(3,5)
\drawline[AHnb=0](3,4)(4,5)
\drawline[AHnb=0](4,4)(5,5)
\drawline[AHnb=0](5,4)(6,5)
\drawline[AHnb=0](6,4)(7,5)

\drawline[AHnb=0](0,5)(1,6)
\drawline[AHnb=0](1,5)(2,6)
\drawline[AHnb=0](2,5)(3,6)
\drawline[AHnb=0](3,5)(4,6)
\drawline[AHnb=0](4,5)(5,6)
\drawline[AHnb=0](5,5)(6,6)
\drawline[AHnb=0](6,5)(7,6)

\put(-.7,-.1){\scriptsize {\tiny $a_{11}$}}
\put(.7,-.25){\scriptsize {\tiny $a_{12}$}}
\put(1.7,-.25){\scriptsize {\tiny $a_{13}$}}
\put(4.35,-.25){\scriptsize {\tiny $a_{1(i-1)}$}}
\put(5.9,-.25){\scriptsize {\tiny $a_{1i}$}}
\put(6.9,-.25){\scriptsize {\tiny $a_{11}$}}

\put(-.7,.9){\scriptsize {\tiny $a_{21}$}}
\put(-.7,1.9){\scriptsize {\tiny $a_{31}$}}

\put(7.1,.9){\scriptsize {\tiny $a_{21}$}}
\put(7.1,1.9){\scriptsize {\tiny $a_{31}$}}

\put(-1.25,6.25){\scriptsize {\tiny $a_{1(k+1)}$}}
\put(.4,6.25){\scriptsize {\tiny $a_{1(k+2)}$}}
\put(4.2,6.25){\scriptsize {\tiny $a_{1(k-1)}$}}
\put(5.8,6.25){\scriptsize {\tiny $a_{1k}$}}
\put(6.6,6.25){\scriptsize {\tiny $a_{1(k+1)}$}}

\put(-.8,4.9){\scriptsize {\tiny $a_{j1}$}}
\put(7.15,4.9){\scriptsize {\tiny $a_{j1}$}}


\put(-1.8,-.1){\scriptsize {\tiny $Q_{1}$}}
\put(-1.8,.9){\scriptsize {\tiny $Q_{2}$}}
\put(-1.8,1.9){\scriptsize {\tiny $Q_{3}$}}
\put(-2.2,3.9){\scriptsize {\tiny $Q_{j-1}$}}
\put(-1.8,4.9){\scriptsize {\tiny $Q_{j}$}}
\put(-1.8,5.9){\scriptsize {\tiny $Q_{1}$}}

\put(3,-.35){\scriptsize $\ldots$}
\put(3,6.2){\scriptsize $\ldots$}

\put(-.5,3){\scriptsize $\vdots$}

\put(7.5,3){\scriptsize $\vdots$}

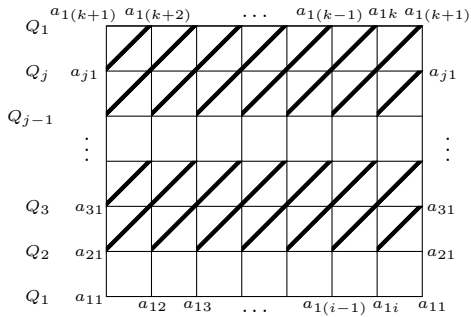
\put(-.8,-1.3){\scriptsize {\tiny {\bf Figure 4.1.3:} $M(i,j,k)$ of DSEM of type $[3^6 : 3^3.4^2]_1$ }}

\drawpolygon[fillcolor=black](0,2)(1,3)(.9,3)(0,2.1)
\drawpolygon[fillcolor=black](0,1)(2,3)(1.9,3)(0,1.1)
\drawpolygon[fillcolor=black](1,1)(3,3)(2.9,3)(1,1.1)
\drawpolygon[fillcolor=black](2,1)(4,3)(3.9,3)(2,1.1)
\drawpolygon[fillcolor=black](3,1)(5,3)(4.9,3)(3,1.1)
\drawpolygon[fillcolor=black](4,1)(6,3)(5.9,3)(4,1.1)
\drawpolygon[fillcolor=black](5,1)(7,3)(6.9,3)(5,1.1)
\drawpolygon[fillcolor=black](6,1)(7,2)(6.9,2)(6,1.1)

\drawpolygon[fillcolor=black](0,5)(1,6)(.9,6)(0,5.1)
\drawpolygon[fillcolor=black](0,4)(2,6)(1.9,6)(0,4.1)
\drawpolygon[fillcolor=black](1,4)(3,6)(2.9,6)(1,4.1)
\drawpolygon[fillcolor=black](2,4)(4,6)(3.9,6)(2,4.1)
\drawpolygon[fillcolor=black](3,4)(5,6)(4.9,6)(3,4.1)
\drawpolygon[fillcolor=black](4,4)(6,6)(5.9,6)(4,4.1)
\drawpolygon[fillcolor=black](5,4)(7,6)(6.9,6)(5,4.1)
\drawpolygon[fillcolor=black](6,4)(7,5)(6.9,5)(6,4.1)

\end{picture}

\begin{picture}(0,0)(-95,45)
\setlength{\unitlength}{6mm}

\drawpolygon(0,0)(7,0)(7,8)(0,8)

\drawline[AHnb=0](0,1)(7,1)
\drawline[AHnb=0](0,2)(7,2)
\drawline[AHnb=0](0,3)(7,3)
\drawline[AHnb=0](0,4)(7,4)
\drawline[AHnb=0](0,5)(7,5)
\drawline[AHnb=0](0,6)(7,6)
\drawline[AHnb=0](0,7)(7,7)
\drawline[AHnb=0](0,8)(7,8)

\drawline[AHnb=0](1,0)(1,8)
\drawline[AHnb=0](2,0)(2,8)
\drawline[AHnb=0](3,0)(3,8)
\drawline[AHnb=0](4,0)(4,8)
\drawline[AHnb=0](5,0)(5,8)
\drawline[AHnb=0](6,0)(6,8)

\drawline[AHnb=0](0,1)(3,4)
\drawline[AHnb=0](1,1)(4,4)
\drawline[AHnb=0](2,1)(5,4)
\drawline[AHnb=0](3,1)(6,4)

\drawline[AHnb=0](0,2)(2,4)
\drawline[AHnb=0](0,3)(1,4)

\drawline[AHnb=0](4,1)(7,4)
\drawline[AHnb=0](5,1)(7,3)
\drawline[AHnb=0](6,1)(7,2)

\drawline[AHnb=0](0,5)(3,8)
\drawline[AHnb=0](1,5)(4,8)
\drawline[AHnb=0](2,5)(5,8)
\drawline[AHnb=0](3,5)(6,8)

\drawline[AHnb=0](0,6)(2,8)
\drawline[AHnb=0](0,7)(1,8)

\drawline[AHnb=0](4,5)(7,8)
\drawline[AHnb=0](5,5)(7,7)
\drawline[AHnb=0](6,5)(7,6)

\drawpolygon[fillcolor=black](0,3)(1,4)(.9,4)(0,3.1)
\drawpolygon[fillcolor=black](0,2)(2,4)(1.9,4)(0,2.1)
\drawpolygon[fillcolor=black](0,1)(3,4)(2.9,4)(0,1.1)
\drawpolygon[fillcolor=black](1,1)(4,4)(3.9,4)(1,1.1)
\drawpolygon[fillcolor=black](2,1)(5,4)(4.9,4)(2,1.1)
\drawpolygon[fillcolor=black](3,1)(6,4)(5.9,4)(3,1.1)
\drawpolygon[fillcolor=black](4,1)(7,4)(6.9,4)(4,1.1)
\drawpolygon[fillcolor=black](5,1)(7,3)(6.9,3)(5,1.1)
\drawpolygon[fillcolor=black](6,1)(7,2)(6.9,2)(6,1.1)

\drawpolygon[fillcolor=black](0,7)(1,8)(.9,8)(0,7.1)
\drawpolygon[fillcolor=black](0,6)(2,8)(1.9,8)(0,6.1)
\drawpolygon[fillcolor=black](0,5)(3,8)(2.9,8)(0,5.1)
\drawpolygon[fillcolor=black](1,5)(4,8)(3.9,8)(1,5.1)
\drawpolygon[fillcolor=black](2,5)(5,8)(4.9,8)(2,5.1)
\drawpolygon[fillcolor=black](3,5)(6,8)(5.9,8)(3,5.1)
\drawpolygon[fillcolor=black](4,5)(7,8)(6.9,8)(4,5.1)
\drawpolygon[fillcolor=black](5,5)(7,7)(6.9,7)(5,5.1)
\drawpolygon[fillcolor=black](6,5)(7,6)(6.9,6)(6,5.1)

\put(-.7,-.1){\scriptsize {\tiny $a_{11}$}}
\put(.7,-.25){\scriptsize {\tiny $a_{12}$}}
\put(1.7,-.25){\scriptsize {\tiny $a_{13}$}}
\put(4.25,-.25){\scriptsize {\tiny $a_{1(i-1)}$}}
\put(5.8,-.25){\scriptsize {\tiny $a_{1i}$}}
\put(7.1,-.1){\scriptsize {\tiny $a_{11}$}}

\put(-.7,.9){\scriptsize {\tiny $a_{21}$}}
\put(-.7,1.9){\scriptsize {\tiny $a_{31}$}}
\put(-.7,2.9){\scriptsize {\tiny $a_{41}$}}
\put(-.7,6.9){\scriptsize {\tiny $a_{j1}$}}
\put(-1.25,8.2){\scriptsize {\tiny $a_{1(k+1)}$}}
\put(.4,8.2){\scriptsize {\tiny $a_{1(k+2)}$}}
\put(4,8.2){\scriptsize {\tiny $a_{1(k-1)}$}}
\put(5.6,8.2){\scriptsize {\tiny $a_{1k}$}}
\put(6.6,8.2){\scriptsize {\tiny $a_{1(k+1)}$}}

\put(7.1,.9){\scriptsize {\tiny $a_{21}$}}
\put(7.1,1.9){\scriptsize {\tiny $a_{31}$}}
\put(7.1,2.9){\scriptsize {\tiny $a_{41}$}}
\put(7.1,6.9){\scriptsize {\tiny $a_{j1}$}}

\put(-1.8,-.1){\scriptsize {\tiny $Q_{1}$}}
\put(-1.8,.9){\scriptsize {\tiny $Q_{2}$}}
\put(-1.8,1.9){\scriptsize {\tiny $Q_{3}$}}
\put(-1.8,2.9){\scriptsize {\tiny $Q_{4}$}}

\put(-2,6){\scriptsize {\tiny $Q_{j-1}$}}

\put(-1.8,7){\scriptsize {\tiny $Q_{j}$}}
\put(-1.8,8){\scriptsize {\tiny $Q_{1}$}}

\put(3,-.35){\scriptsize $\ldots$}
\put(3,8.2){\scriptsize $\ldots$}

\put(-.5,4.35){\scriptsize $\vdots$}
\put(7.5,4.35){\scriptsize $\vdots$}

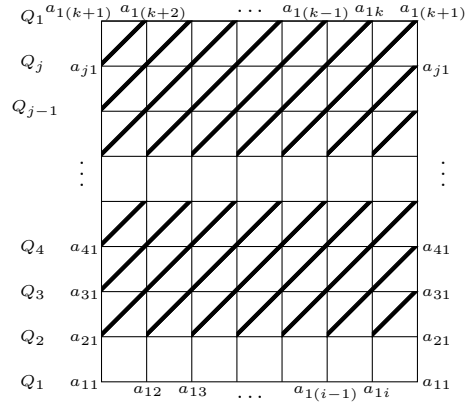
\put(-1.7,-1){\scriptsize {\tiny {\bf Figure 4.1.4:} $M(i,j,k)$ of DSEM of type $[3^6 : 3^3.4^2]_2$ }} 

\end{picture}

\vspace{5cm}
\begin{lem}  \label{l4.1.3} Let $r\in\{1,2\}.$ Let $M^r$ be a DSEM of type $[3^6: 3^3.4^2]_r$. Then $M^r$ admits an $M^r(i,j,k)$ representation, iff the following holds: $(i)$ $i \geq 3$ and $j=(r+2)m$, where $m \in \mathbb{N}$, $(ii)$ $ij \geq 3(r+2) $, $(iii)$ $ 0 \leq k \leq i-1 $.
\end{lem}

\noindent{\bf Proof.} Note that an $M^r(i, j, k)$ of $M^r$ contains $j$ number of $A_{1}$ type disjoint horizontal cycles of length $i$. Since these cycles cover all the vertices of $M^r$, the number of vertices in $M^r$ is $n = ij$. Clearly if $i \leq 2$, $M^r$ is not a map. So $i \geq 3$. If $j=1$ then $M^r$ is not a map and if $j =2$ then $M^r$ has no vertices of face-sequence $(3^6)$. If $j=(r+2)m+1$ or $(r+2)m+2$, then $2|V_{(3^6)}| \neq |V_{(3^3.4^2)}|$ for $M^1$ and $|V_{(3^6)}| \neq |V_{(3^3.4^2)}|$ for $M^2$. So $j = (r+2)m$. Thus $n = ij \geq 3(r+2)$. Since the length of the horizontal cycle is $i$, we get $k \in \{0,1, \ldots, i-1\}$.  The converse part follows directly by constructing $M^r(i,j,k)$ representation for the given values of $i,j$, and $k$. 
\hfill $\Box$



\vspace{.2cm}
In the upcoming subsections, we proceed in a similar way for the remaining types DSEMs. We determine the conditions among $i,j,k$ to construct its $M(i,j,k)$ representation. For this, we consider fixed types of paths on maps, such paths can be defined using the definition of link (as in  Subsection \ref{s3.1}). Since these DSEMs are on the finite vertex set, given every maximal path $P$ of such types there is an edge $e$ such that $P \cup \{e\}$ is a non-contractible, as in Lemma \ref{l4.1.1}.

\subsection{DSEMs of types $[3^3.4^2:4^4]_1$ and  $[3^3.4^2:4^4]_2$} \label{s3.2}

Let $M^r$ be a DSEM of type $[3^3.4^2:4^4]_r$, where $ r \in \{1,2\}$. Then for the existence of $M^1$ and $M^2$ we have $|V_{(3^3,4^2)}| = 2|V_{(4^4)}|$ and $|V_{(3^3,4^2)}| = |V_{(4^4)}|$ respectively. 
Now, consider the following paths in $M^r$ as follows. 

A path $P_{1} = P( \ldots, y_{i-1},y_{i},y_{i+1}, \ldots)$ in $M^r$, say  of type $B_{1}$, indicated by thick black paths. The vertices $y_i$'s have the face-sequence either $(4^4)$ or $(3^3,4^2)$. 

A path $P_2=P( \ldots, y_{i-2},z_{i-2},z_{i-1},y_{i},z_{i},z_{i+1},y_{i+1} \ldots)$ in $M^1$, say of type $B_2$, indicated by green paths. The vertices $y_i$'s and $z_i$'s have the face-sequences $(4^4)$ and $(3^3,4^2)$ respectively.  

A path $P_3=P( \ldots, y_{i-1},z_{i-2},z_{i-1},y_{i},y_{i+1},z_{i},z_{i+1},y_{i+2} \ldots)$ in $M^2$, say of type $B_2$, indicated by green paths. The vertices $y_i$'s and $z_i$'s have the face-sequences $(4^4)$ and $(3^3,4^2)$ respectively. 

\begin{picture}(0,0)(-28,24)
\setlength{\unitlength}{5.25mm}

\drawpolygon(0,0)(6,0)(6,2)(0,2)


\drawline[AHnb=0](1,0)(1,2)
\drawline[AHnb=0](2,0)(2,2)
\drawline[AHnb=0](3,0)(3,2)
\drawline[AHnb=0](4,0)(4,2)
\drawline[AHnb=0](5,0)(5,2)

\drawline[AHnb=0](.5,3)(.5,3.5)
\drawline[AHnb=0](1.5,3)(1.5,3.5)

\drawline[AHnb=0](2.5,3)(2.5,3.5)
\drawline[AHnb=0](3.5,3)(3.5,3.5)

\drawline[AHnb=0](4.5,3)(4.5,3.5)
\drawline[AHnb=0](5.5,3)(5.5,3.5)


\drawline[AHnb=0](6,0)(6.5,0)
\drawline[AHnb=0](0,0)(-.5,0)
\drawline[AHnb=0](6,3)(6.5,3)
\drawline[AHnb=0](0,3)(-.5,3)

\drawline[AHnb=0](-.5,1)(6.5,1)
\drawline[AHnb=0](-.5,2)(6.5,2)

\drawline[AHnb=0](-.5,3)(6.5,3)


\drawline[AHnb=0](0,2)(.5,3)
\drawline[AHnb=0](1,2)(.5,3)

\drawline[AHnb=0](1,2)(1.5,3)
\drawline[AHnb=0](2,2)(1.5,3)

\drawline[AHnb=0](2,2)(2.5,3)
\drawline[AHnb=0](3,2)(2.5,3)

\drawline[AHnb=0](3,2)(3.5,3)
\drawline[AHnb=0](4,2)(3.5,3)

\drawline[AHnb=0](4,2)(4.5,3)
\drawline[AHnb=0](5,2)(4.5,3)

\drawline[AHnb=0](5,2)(5.5,3)
\drawline[AHnb=0](6,2)(5.5,3)


\drawline[AHnb=0](0,0)(.25,-.5)
\drawline[AHnb=0](1,0)(1.25,-.5)
\drawline[AHnb=0](2,0)(2.25,-.5)
\drawline[AHnb=0](3,0)(3.25,-.5)
\drawline[AHnb=0](4,0)(4.25,-.5)
\drawline[AHnb=0](5,0)(5.25,-.5)
\drawline[AHnb=0](6,0)(6.25,-.5)

\drawline[AHnb=0](0,0)(-.25,-.5)
\drawline[AHnb=0](1,0)(.75,-.5)
\drawline[AHnb=0](2,0)(1.75,-.5)
\drawline[AHnb=0](3,0)(2.75,-.5)
\drawline[AHnb=0](4,0)(3.75,-.5)
\drawline[AHnb=0](5,0)(4.75,-.5)
\drawline[AHnb=0](6,0)(5.75,-.5)

\drawline[AHnb=0](0,2)(-.25,2.5)
\drawline[AHnb=0](6,2)(6.25,2.5)

\drawpolygon[fillcolor=green](.25,-.5)(0,0)(0,2)(-.25,2.5)(-.15,2.5)(.1,2)(.1,0)(.35,-.5)

\drawpolygon[fillcolor=green](1.25,-.5)(1,0)(1,2)(.5,3)(.5,3.5)(.6,3.5)(.6,3)(1.1,2)(1.1,0)(1.35,-.5)

\drawpolygon[fillcolor=green](2.25,-.5)(2,0)(2,2)(1.5,3)(1.5,3.5)(1.6,3.5)(1.6,3)(2.1,2)(2.1,0)(2.35,-.5)

\drawpolygon[fillcolor=green](3.25,-.5)(3,0)(3,2)(2.5,3)(2.5,3.5)(2.6,3.5)(2.6,3)(3.1,2)(3.1,0)(3.35,-.5)

\drawpolygon[fillcolor=green](4.25,-.5)(4,0)(4,2)(3.5,3)(3.5,3.5)(3.6,3.5)(3.6,3)(4.1,2)(4.1,0)(4.35,-.5)

\drawpolygon[fillcolor=green](5.25,-.5)(5,0)(5,2)(4.5,3)(4.5,3.5)(4.6,3.5)(4.6,3)(5.1,2)(5.1,0)(5.35,-.5)

\drawpolygon[fillcolor=green](6.25,-.5)(6,0)(6,2)(5.5,3)(5.5,3.5)(5.6,3.5)(5.6,3)(6.1,2)(6.1,0)(6.35,-.5)

\drawpolygon[fillcolor=black](-.5,0)(6.5,0)(6.5,.1)(-.5,.1)

\drawpolygon[fillcolor=black](-.5,1)(6.5,1)(6.5,1.1)(-.5,1.1)

\drawpolygon[fillcolor=black](-.5,2)(6.5,2)(6.5,2.1)(-.5,2.1)

\drawpolygon[fillcolor=black](-.5,3)(6.5,3)(6.5,3.1)(-.5,3.1)

\put(3.4,3.75) {\scriptsize {\tiny $B_2$}}
\put(6.7,-.1) {\scriptsize {\tiny $B_1$}}

\end{picture}

\begin{picture}(60,0)(-85,20)
\setlength{\unitlength}{5mm}

\drawpolygon(0,0)(6,0)(6,3)(0,3)


\drawline[AHnb=0](1,0)(1,3)
\drawline[AHnb=0](2,0)(2,3)
\drawline[AHnb=0](3,0)(3,3)
\drawline[AHnb=0](4,0)(4,3)
\drawline[AHnb=0](5,0)(5,3)

\drawline[AHnb=0](.5,4)(.5,4.5)
\drawline[AHnb=0](1.5,4)(1.5,4.5)

\drawline[AHnb=0](2.5,4)(2.5,4.5)
\drawline[AHnb=0](3.5,4)(3.5,4.5)

\drawline[AHnb=0](4.5,4)(4.5,4.5)
\drawline[AHnb=0](5.5,4)(5.5,4.5)


\drawline[AHnb=0](6,0)(6.5,0)
\drawline[AHnb=0](0,0)(-.5,0)
\drawline[AHnb=0](6,3)(6.5,3)
\drawline[AHnb=0](0,3)(-.5,3)

\drawline[AHnb=0](-.5,1)(6.5,1)
\drawline[AHnb=0](-.5,2)(6.5,2)

\drawline[AHnb=0](-.5,4)(6.5,4)


\drawline[AHnb=0](0,3)(.5,4)
\drawline[AHnb=0](1,3)(.5,4)

\drawline[AHnb=0](1,3)(1.5,4)
\drawline[AHnb=0](2,3)(1.5,4)

\drawline[AHnb=0](2,3)(2.5,4)
\drawline[AHnb=0](3,3)(2.5,4)

\drawline[AHnb=0](3,3)(3.5,4)
\drawline[AHnb=0](4,3)(3.5,4)

\drawline[AHnb=0](4,3)(4.5,4)
\drawline[AHnb=0](5,3)(4.5,4)

\drawline[AHnb=0](5,3)(5.5,4)
\drawline[AHnb=0](6,3)(5.5,4)


\drawline[AHnb=0](0,0)(.25,-.5)
\drawline[AHnb=0](1,0)(1.25,-.5)
\drawline[AHnb=0](2,0)(2.25,-.5)
\drawline[AHnb=0](3,0)(3.25,-.5)
\drawline[AHnb=0](4,0)(4.25,-.5)
\drawline[AHnb=0](5,0)(5.25,-.5)
\drawline[AHnb=0](6,0)(6.25,-.5)

\drawline[AHnb=0](0,0)(-.25,-.5)
\drawline[AHnb=0](1,0)(.75,-.5)
\drawline[AHnb=0](2,0)(1.75,-.5)
\drawline[AHnb=0](3,0)(2.75,-.5)
\drawline[AHnb=0](4,0)(3.75,-.5)
\drawline[AHnb=0](5,0)(4.75,-.5)
\drawline[AHnb=0](6,0)(5.75,-.5)

\drawline[AHnb=0](0,3)(-.25,3.5)
\drawline[AHnb=0](6,3)(6.25,3.5)


\drawpolygon[fillcolor=black](-.5,0)(6.5,0)(6.5,.1)(-.5,.1)

\drawpolygon[fillcolor=black](-.5,1)(6.5,1)(6.5,1.1)(-.5,1.1)

\drawpolygon[fillcolor=black](-.5,2)(6.5,2)(6.5,2.1)(-.5,2.1)

\drawpolygon[fillcolor=black](-.5,3)(6.5,3)(6.5,3.1)(-.5,3.1)

\drawpolygon[fillcolor=black](-.5,4)(6.5,4)(6.5,4.1)(-.5,4.1)

\drawpolygon[fillcolor=green](.25,-.5)(0,0)(0,3)(-.25,3.5)(-.15,3.5)(.1,3)(.1,0)(.35,-.5)

\drawpolygon[fillcolor=green](1.25,-.5)(1,0)(1,3)(.5,4)(.5,4.5)(.6,4.5)(.6,4)(1.1,3)(1.1,0)(1.35,-.5)

\drawpolygon[fillcolor=green](2.25,-.5)(2,0)(2,3)(1.5,4)(1.5,4.5)(1.6,4.5)(1.6,4)(2.1,3)(2.1,0)(2.35,-.5)

\drawpolygon[fillcolor=green](3.25,-.5)(3,0)(3,3)(2.5,4)(2.5,4.5)(2.6,4.5)(2.6,4)(3.1,3)(3.1,0)(3.35,-.5)

\drawpolygon[fillcolor=green](4.25,-.5)(4,0)(4,3)(3.5,4)(3.5,4.5)(3.6,4.5)(3.6,4)(4.1,3)(4.1,0)(4.35,-.5)

\drawpolygon[fillcolor=green](5.25,-.5)(5,0)(5,3)(4.5,4)(4.5,4.5)(4.6,4.5)(4.6,4)(5.1,3)(5.1,0)(5.35,-.5)

\drawpolygon[fillcolor=green](6.25,-.5)(6,0)(6,3)(5.5,4)(5.5,4.5)(5.6,4.5)(5.6,4)(6.1,3)(6.1,0)(6.35,-.5)

\put(3.4,4.75) {\scriptsize {\tiny $B_2$}}
\put(6.7,-.1) {\scriptsize {\tiny $B_1$}}

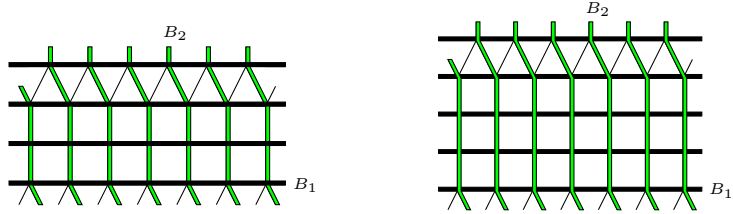
\put(-12,-1.8){\scriptsize {\tiny { Figure 4.2.1 : Paths of type $B_1, B_2$ in DSEM of type $[3^3.4^2:4^4]_1$ and $[3^3. 4^2:4^4]_2$ }} }

\end{picture}

\begin{picture}(0,0)(-14,80.5)
\setlength{\unitlength}{6mm}

\drawpolygon(0,0)(7,0)(7,6)(0,6)

\drawline[AHnb=0](0,1)(7,1)
\drawline[AHnb=0](0,2)(7,2)
\drawline[AHnb=0](0,3)(7,3)
\drawline[AHnb=0](0,4)(7,4)
\drawline[AHnb=0](0,5)(7,5)
\drawline[AHnb=0](0,6)(7,6)

\drawline[AHnb=0](1,0)(1,6)
\drawline[AHnb=0](2,0)(2,6)
\drawline[AHnb=0](3,0)(3,6)
\drawline[AHnb=0](4,0)(4,6)
\drawline[AHnb=0](5,0)(5,6)
\drawline[AHnb=0](6,0)(6,6)

\drawline[AHnb=0](0,2)(1,3)
\drawline[AHnb=0](1,2)(2,3)
\drawline[AHnb=0](2,2)(3,3)
\drawline[AHnb=0](3,2)(4,3)
\drawline[AHnb=0](4,2)(5,3)
\drawline[AHnb=0](5,2)(6,3)
\drawline[AHnb=0](6,2)(7,3)

\drawline[AHnb=0](0,5)(1,6)
\drawline[AHnb=0](1,5)(2,6)
\drawline[AHnb=0](2,5)(3,6)
\drawline[AHnb=0](3,5)(4,6)
\drawline[AHnb=0](4,5)(5,6)
\drawline[AHnb=0](5,5)(6,6)
\drawline[AHnb=0](6,5)(7,6)

\put(-.7,-.1){\scriptsize {\tiny $a_{11}$}}
\put(.7,-.35){\scriptsize {\tiny $a_{12}$}}
\put(1.7,-.35){\scriptsize {\tiny $a_{13}$}}
\put(4.5,-.35){\scriptsize {\tiny $a_{1(i-1)}$}}
\put(5.9,-.35){\scriptsize {\tiny $a_{1i}$}}
\put(6.9,-.35){\scriptsize {\tiny $a_{11}$}}

\put(-.7,.9){\scriptsize {\tiny $a_{21}$}}

\put(7.1,.9){\scriptsize {\tiny $a_{21}$}}

\put(-1.25,6.25){\scriptsize {\tiny $a_{1(k+1)}$}}
\put(.4,6.25){\scriptsize {\tiny $a_{1(k+2)}$}}
\put(4.2,6.25){\scriptsize {\tiny $a_{1(k-1)}$}}
\put(5.8,6.25){\scriptsize {\tiny $a_{1k}$}}
\put(6.6,6.25){\scriptsize {\tiny $a_{1(k+1)}$}}

\put(-.8,4.9){\scriptsize {\tiny $a_{j1}$}}
\put(7.15,4.9){\scriptsize {\tiny $a_{j1}$}}


\put(-1.8,-.1){\scriptsize {\tiny $Q_{1}$}}
\put(-1.8,.9){\scriptsize {\tiny $Q_{2}$}}
\put(-2.2,3.9){\scriptsize {\tiny $Q_{j-1}$}}
\put(-1.8,4.9){\scriptsize {\tiny $Q_{j}$}}
\put(-1.8,5.9){\scriptsize {\tiny $Q_{1}$}}

\put(3,-.35){\scriptsize $\ldots$}
\put(3,6.2){\scriptsize $\ldots$}

\put(-.5,2.5){\scriptsize $\vdots$}

\put(7.5,2.5){\scriptsize $\vdots$}

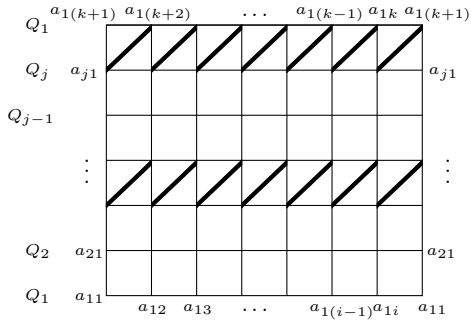
\put(-.8,-1.5){\scriptsize {\tiny {\bf Figure 4.2.2:} $M(i,j,k)$ of DSEM of type $[3^3.4^2 : 4^4]_1$ }}

\drawpolygon[fillcolor=black](0,2.05)(0,1.95)(1,2.9)(1,3)
\drawpolygon[fillcolor=black](1,2.05)(1,1.95)(2,2.9)(2,3)
\drawpolygon[fillcolor=black](2,2.05)(2,1.95)(3,2.9)(3,3)
\drawpolygon[fillcolor=black](3,2.05)(3,1.95)(4,2.9)(4,3)
\drawpolygon[fillcolor=black](4,2.05)(4,1.95)(5,2.9)(5,3)
\drawpolygon[fillcolor=black](5,2.05)(5,1.95)(6,2.9)(6,3)
\drawpolygon[fillcolor=black](6,2.05)(6,1.95)(7,2.9)(7,3)

\drawpolygon[fillcolor=black](0,5.05)(0,4.95)(1,5.9)(1,6)
\drawpolygon[fillcolor=black](1,5.05)(1,4.95)(2,5.9)(2,6)
\drawpolygon[fillcolor=black](2,5.05)(2,4.95)(3,5.9)(3,6)
\drawpolygon[fillcolor=black](3,5.05)(3,4.95)(4,5.9)(4,6)
\drawpolygon[fillcolor=black](4,5.05)(4,4.95)(5,5.9)(5,6)
\drawpolygon[fillcolor=black](5,5.05)(5,4.95)(6,5.9)(6,6)
\drawpolygon[fillcolor=black](6,5.05)(6,4.95)(7,5.9)(7,6)

\end{picture}

\begin{picture}(0,0)(-95,76.25)
\setlength{\unitlength}{6mm}

\drawpolygon(0,0)(7,0)(7,8)(0,8)

\drawline[AHnb=0](0,1)(7,1)
\drawline[AHnb=0](0,2)(7,2)
\drawline[AHnb=0](0,3)(7,3)
\drawline[AHnb=0](0,4)(7,4)
\drawline[AHnb=0](0,5)(7,5)
\drawline[AHnb=0](0,6)(7,6)
\drawline[AHnb=0](0,7)(7,7)
\drawline[AHnb=0](0,8)(7,8)

\drawline[AHnb=0](1,0)(1,8)
\drawline[AHnb=0](2,0)(2,8)
\drawline[AHnb=0](3,0)(3,8)
\drawline[AHnb=0](4,0)(4,8)
\drawline[AHnb=0](5,0)(5,8)
\drawline[AHnb=0](6,0)(6,8)

\drawpolygon[fillcolor=black](0,3.05)(0,2.95)(1,3.9)(1,4)
\drawpolygon[fillcolor=black](1,3.05)(1,2.95)(2,3.9)(2,4)
\drawpolygon[fillcolor=black](2,3.05)(2,2.95)(3,3.9)(3,4)
\drawpolygon[fillcolor=black](3,3.05)(3,2.95)(4,3.9)(4,4)
\drawpolygon[fillcolor=black](4,3.05)(4,2.95)(5,3.9)(5,4)
\drawpolygon[fillcolor=black](5,3.05)(5,2.95)(6,3.9)(6,4)
\drawpolygon[fillcolor=black](6,3.05)(6,2.95)(7,3.9)(7,4)

\drawpolygon[fillcolor=black](0,7.05)(0,6.95)(1,7.9)(1,8)
\drawpolygon[fillcolor=black](1,7.05)(1,6.95)(2,7.9)(2,8)
\drawpolygon[fillcolor=black](2,7.05)(2,6.95)(3,7.9)(3,8)
\drawpolygon[fillcolor=black](3,7.05)(3,6.95)(4,7.9)(4,8)
\drawpolygon[fillcolor=black](4,7.05)(4,6.95)(5,7.9)(5,8)
\drawpolygon[fillcolor=black](5,7.05)(5,6.95)(6,7.9)(6,8)
\drawpolygon[fillcolor=black](6,7.05)(6,6.95)(7,7.9)(7,8)

\put(-.7,-.1){\scriptsize {\tiny $a_{11}$}}
\put(.7,-.35){\scriptsize {\tiny $a_{12}$}}
\put(1.7,-.35){\scriptsize {\tiny $a_{13}$}}
\put(4.4,-.35){\scriptsize {\tiny $a_{1(i-1)}$}}
\put(5.8,-.35){\scriptsize {\tiny $a_{1i}$}}
\put(7.1,-.1){\scriptsize {\tiny $a_{11}$}}

\put(-.7,.9){\scriptsize {\tiny $a_{21}$}}
\put(-.7,1.9){\scriptsize {\tiny $a_{31}$}}
\put(-.7,6.9){\scriptsize {\tiny $a_{j1}$}}
\put(-1.25,8.2){\scriptsize {\tiny $a_{1(k+1)}$}}
\put(.4,8.2){\scriptsize {\tiny $a_{1(k+2)}$}}
\put(4.2,8.2){\scriptsize {\tiny $a_{1(k-1)}$}}
\put(5.6,8.2){\scriptsize {\tiny $a_{1k}$}}
\put(6.6,8.2){\scriptsize {\tiny $a_{1(k+1)}$}}

\put(7.1,.9){\scriptsize {\tiny $a_{21}$}}
\put(7.1,1.9){\scriptsize {\tiny $a_{31}$}}
\put(7.1,6.9){\scriptsize {\tiny $a_{j1}$}}

\put(-1.8,-.1){\scriptsize {\tiny $Q_{1}$}}
\put(-1.8,.9){\scriptsize {\tiny $Q_{2}$}}
\put(-1.8,1.9){\scriptsize {\tiny $Q_{3}$}}

\put(-2,6){\scriptsize {\tiny $Q_{j-1}$}}

\put(-1.8,7){\scriptsize {\tiny $Q_{j}$}}
\put(-1.8,8){\scriptsize {\tiny $Q_{1}$}}

\put(3,-.35){\scriptsize $\ldots$}
\put(3,8.2){\scriptsize $\ldots$}

\put(-.5,4){\scriptsize $\vdots$}

\put(7.5,4){\scriptsize $\vdots$}

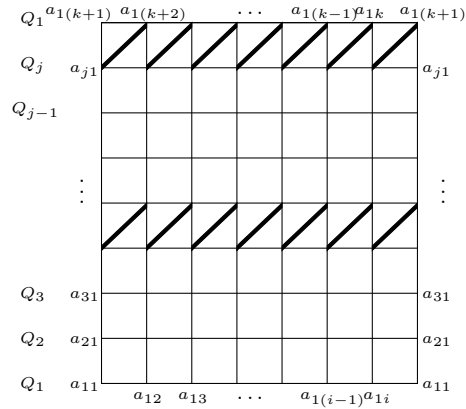
\put(-.7,-1.35){\scriptsize {\tiny {\bf Figure 4.2.3:} $M(i,j,k)$ of DSEM of type $[3^3.4^2 : 4^4]_2$ }} 
\end{picture}

\vspace{8.75cm}

Now, first cutting $M^r$ along a cycle of type $B_1$ and then along a cycle of type $B_2$, we get $M^r(i,j,k)$-representation of $M^r$, shown in Figures 4.2.2 and 4.2.3. To get the relation among $i,j,k$, we proceed similarly as in Lemma \ref{l4.1.3}. This gives the following result.


\begin{lem} \label{l3.2.1}
	Let $r \in \{1,2\}$. The DSEM $M^r$ of type $[3^3.4^2:4^4]_r$ admits an $M^r(i,j,k)$-representation iff the following holds: $(i)$ $i \geq 3$ and $j=(r+2)m$, where $m \in \mathbb{N} $, $(ii)$ $ij \geq 3(r+2) $, $(iii)$ $ 0 \leq k \leq i-1 $.
\end{lem}

%

\subsection{DSEMs of types $[3^3.4^2: 3^2.4.3.4]_1$}\label{s3.3}

Let $M$ be a DSEM of type $[3^3.4^2: 3^2.4.3.4]_1$. Then, for the existence of $M$, we have $2|V_{(3^3,4^2)}| = |V_{(3^2,4,3,4)}|$. Now consider the following path through a vertex in Figure 4.3.1. One can define such a path using the definition of link as in Subsection \ref{s3.1}. 

A path $P_1 = P(\ldots,w_i,z_{j},z_{j+1},z_{j+2},z_{j+3}$, $w_{i+1}, \ldots)$ in $M$, say of type $C_1$, indicated by thick black paths in Figure 4.3.1. The vertices $w_i$'s and $z_j$'s have face-sequences $(3^3,4^2)$ and $(3^2,4,3,4)$ respectively.

Note that, through a vertex $v$ with the face-sequence $(3^2,4,3,4)$, we see two paths of type $C_1$. Now as in the previous section, cutting $M$ along these paths one by one, we get its $M(i,j,k)$ representation. The relations among $i, j, k$ are given by Lemma \ref{l3.3.1}.


\begin{picture}(0,0)(-5,33)
\setlength{\unitlength}{7mm}

\drawline[AHnb=0](0,0)(2,0)(3,.5)(4,.5)(5,0)(7,0)(7.5,.25)
\drawline[AHnb=0](0,-1)(2,-1)(3,-1.5)(4,-1.5)(5,-1)(7,-1)(7.5,-.75)

\drawline[AHnb=0](0,0)(0,-1)
\drawline[AHnb=0](1,0)(1,-1)
\drawline[AHnb=0](2,0)(2,-1)
\drawline[AHnb=0](3,.5)(3,-1.5)
\drawline[AHnb=0](4,.5)(4,-1.5)
\drawline[AHnb=0](5,0)(5,-1)
\drawline[AHnb=0](6,0)(6,-1)
\drawline[AHnb=0](7,0)(7,-1)

\drawline[AHnb=0](3,-.5)(4,-.5)

\drawline[AHnb=0](3,-.5)(2,-1)
\drawline[AHnb=0](3,-.5)(2,0)

\drawline[AHnb=0](4,-.5)(5,-1)
\drawline[AHnb=0](4,-.5)(5,0)

\drawline[AHnb=0](.5,3)(1.5,3)(2.5,2.5)(4.5,2.5)(5.5,3)(6.5,3)(7.5,2.5)(8,2.5)

\drawline[AHnb=0](.5,1)(1.5,1)(2.5,1.5)(4.5,1.5)(5.5,1)(6.5,1)(7.5,1.5)(8,1.5)
\drawline[AHnb=0](.5,3)(.5,1)

\drawline[AHnb=0](1.5,3)(1.5,1)
\drawline[AHnb=0](2.5,2.5)(2.5,1.5)
\drawline[AHnb=0](4.5,2.5)(4.5,1.5)
\drawline[AHnb=0](5.5,3)(5.5,1)
\drawline[AHnb=0](6.5,3)(6.5,1)
\drawline[AHnb=0](7.5,2.5)(7.5,1.5)

\drawline[AHnb=0](3.5,2.5)(3.5,1.5)

\drawline[AHnb=0](.5,2)(1.5,2)
\drawline[AHnb=0](5.5,2)(6.5,2)

\drawline[AHnb=0](2.5,2.5)(1.5,2)
\drawline[AHnb=0](2.5,1.5)(1.5,2)

\drawline[AHnb=0](4.5,2.5)(5.5,2)
\drawline[AHnb=0](4.5,1.5)(5.5,2)

\drawline[AHnb=0](7.5,2.5)(6.5,2)
\drawline[AHnb=0](7.5,1.5)(6.5,2)

\drawline[AHnb=0](0,0)(.5,1)
\drawline[AHnb=0](1,0)(1.5,1)
\drawline[AHnb=0](3,.5)(3.5,1.5)
\drawline[AHnb=0](4,.5)(4.5,1.5)
\drawline[AHnb=0](5,0)(5.5,1)
\drawline[AHnb=0](6,0)(6.5,1)

\drawline[AHnb=0](3,.5)(2.5,1.5)
\drawline[AHnb=0](2,0)(1.5,1)
\drawline[AHnb=0](7,0)(6.5,1)

\drawline[AHnb=0](2.5,2.5)(2.75,3)
\drawline[AHnb=0](3.5,2.5)(3.25,3)

\drawline[AHnb=0](3.5,2.5)(3.75,3)

\drawline[AHnb=0](4.5,2.5)(4.25,3)

\drawline[AHnb=0](5.5,3)(5.75,3.5)
\drawline[AHnb=0](6.5,3)(6.75,3.5)
\drawline[AHnb=0](6.5,3)(6.25,3.5)

\drawline[AHnb=0](5.5,3)(5.25,3.5)

\drawline[AHnb=0](7.75,3)(7.5,2.5)

\drawline[AHnb=0](2,-1)(1.5,-2)
\drawline[AHnb=0](1,-1)(1.5,-2)
\drawline[AHnb=0](1,-1)(.5,-2)
\drawline[AHnb=0](0,-1)(.5,-2)
\drawline[AHnb=0](.5,-2)(1.5,-2)
\drawline[AHnb=0](.5,-2)(.5,-2.5)
\drawline[AHnb=0](1.5,-2)(1.5,-2.5)
\drawline[AHnb=0](1.5,-2)(2.5,-2.5)
\drawline[AHnb=0](6,0)(5.5,1)
\drawline[AHnb=0](4,.5)(3.5,1.5)
\drawline[AHnb=0](4.5,-2.5)(5,-2.75)
\drawline[AHnb=0](2.5,-2.5)(2,-2.75)

\drawline[AHnb=0](3,-1.5)(3.5,-2.5)
\drawline[AHnb=0](4,-1.5)(3.5,-2.5)

\drawline[AHnb=0](2.5,-2.5)(2.5,-2.9)
\drawline[AHnb=0](3.5,-2.5)(3.5,-2.9)
\drawline[AHnb=0](4.5,-2.5)(4.5,-2.9)

\drawline[AHnb=0](3,-1.5)(2.5,-2.5)
\drawline[AHnb=0](2.5,-2.5)(3.5,-2.5)
\drawline[AHnb=0](3.5,-2.5)(4.5,-2.5)
\drawline[AHnb=0](4,-1.5)(4.5,-2.5)

\drawline[AHnb=0](5.5,-2)(4.5,-2.5)
\drawline[AHnb=0](6.5,-2)(7,-2.25)
\drawline[AHnb=0](7.5,1.5)(7.75,1)

\drawline[AHnb=0](5,-1)(5.5,-2)
\drawline[AHnb=0](6,-1)(5.5,-2)
\drawline[AHnb=0](6,-1)(6.5,-2)
\drawline[AHnb=0](7,-1)(6.5,-2)
\drawline[AHnb=0](5.5,-2)(6.5,-2)

\drawline[AHnb=0](1,0)(.5,1)
\drawline[AHnb=0](-.5,1.5)(.5,1)
\drawline[AHnb=0](-.5,1.5)(.5,2)
\drawline[AHnb=0](-.5,2.5)(.5,3)
\drawline[AHnb=0](-.5,2.5)(.5,2)
\drawline[AHnb=0](-.5,1.5)(-.5,2.5)

\drawline[AHnb=0](-.5,.25)(0,0)
\drawline[AHnb=0](-.5,-.25)(0,0)

\drawline[AHnb=0](0,-1)(-.5,-1.25)

\drawline[AHnb=0](-.5,1.5)(-.75,1)
\drawline[AHnb=0](-.5,2.5)(-.75,3)

\drawline[AHnb=0](.5,-2)(0,-2.25)
\drawline[AHnb=0](0,-1)(-.5,-.75)

\drawline[AHnb=0](7,-1)(7.5,-1.25)
\drawline[AHnb=0](7,0)(7.5,-.25)

\drawline[AHnb=0](5.5,-2)(5.5,-2.4)
\drawline[AHnb=0](6.5,-2)(6.5,-2.4)

\drawline[AHnb=0](.5,3)(.25,3.5)
\drawline[AHnb=0](.5,3)(.75,3.5)
\drawline[AHnb=0](1.5,3)(1.25,3.5)
\drawline[AHnb=0](1.5,3)(1.75,3.5)

\drawpolygon[fillcolor=black](-.5,.25)(0,0)(2,0)(3,.5)(4,.5)(5,0)(7,0)(7.5,.25)(7.5,.35)(7,0.1)(5,0.1)(4,.6)(3,.6)(2,0.1)(0,0.1)(-.5,.35)

\drawpolygon[fillcolor=black](-.5,2.5)(.5,3)(1.5,3)(2.5,2.5)(4.5,2.5)(5.5,3)(6.5,3)(7.5,2.5)(8,2.5)(8,2.6)(7.5,2.6)(6.5,3.1)(5.5,3.1)(4.5,2.6)(2.5,2.6)(1.5,3.1)(.5,3.1)(-.5,2.6)

\drawpolygon[fillcolor=black](-.5,1.5)(.5,1)(1.5,1)(2.5,1.5)(4.5,1.5)(5.5,1)(6.5,1)(7.5,1.5)(8,1.5)(8,1.6)(7.5,1.6)(6.5,1.1)(5.5,1.1)(4.5,1.6)(2.5,1.6)(1.5,1.1)(.5,1.1)(-.5,1.6)

\drawpolygon[fillcolor=black](-.5,-1.25)(0,-1)(2,-1)(3,-1.5)(4,-1.5)(5,-1)(7,-1)(7.5,-.75)(7.5,-.65)(7,-0.9)(5,-.9)(4,-1.4)(3,-1.4)(2,-.9)(0,-.9)(-.5,-1.15)

\drawpolygon[fillcolor=black](2.5,-3)(2.5,-2.5)(3,-1.5)(3,.5)(2.5,1.5)(2.5,2.5)(2.75,3)(2.85,3)(2.6,2.5)(2.6,1.5)(3.1,.5)(3.1,-1.5)(2.6,-2.5)(2.6,-3)

\drawpolygon[fillcolor=black](2.5,-3)(2.5,-2.5)(3,-1.5)(3,.5)(2.5,1.5)(2.5,2.5)(2.75,3)(2.85,3)(2.6,2.5)(2.6,1.5)(3.1,.5)(3.1,-1.5)(2.6,-2.5)(2.6,-3)

\drawpolygon[fillcolor=black](4.5,-3)(4.5,-2.5)(4,-1.5)(4,.5)(4.5,1.5)(4.5,2.5)(4.25,3)(4.35,3)(4.6,2.5)(4.6,1.5)(4.1,.5)(4.1,-1.5)(4.6,-2.5)(4.6,-3)

\drawpolygon[fillcolor=black](2.5,-3)(2.5,-2.5)(3,-1.5)(3,.5)(2.5,1.5)(2.5,2.5)(2.75,3)(2.85,3)(2.6,2.5)(2.6,1.5)(3.1,.5)(3.1,-1.5)(2.6,-2.5)(2.6,-3)

\drawpolygon[fillcolor=black](1.5,-2.5)(1.5,-2)(2,-1)(2,.0)(1.5,1)(1.5,3)(1.75,3.5)(1.85,3.5)(1.6,3)(1.6,1)(2.1,0)(2.1,-1)(1.6,-2)(1.6,-2.5)

\drawpolygon[fillcolor=black](.5,-2.5)(.5,-2)(0,-1)(0,.0)(.5,1)(.5,3)(.25,3.5)(.35,3.5)(.6,3)(.6,1)(0.1,0)(0.1,-1)(.6,-2)(.6,-2.5)

\drawpolygon[fillcolor=black](5.5,-2.5)(5.5,-2)(5,-1)(5,.0)(5.5,1)(5.5,3)(5.25,3.5)(5.35,3.5)(5.6,3)(5.6,1)(5.1,0)(5.1,-1)(5.6,-2)(5.6,-2.5)

\drawpolygon[fillcolor=black](6.5,-2.5)(6.5,-2)(7,-1)(7,.0)(6.5,1)(6.5,3)(6.75,3.5)(6.85,3.5)(6.6,3)(6.6,1)(7.1,0)(7.1,-1)(6.6,-2)(6.6,-2.5)

\put(7.8,.3){\scriptsize {\tiny $C_1$}}
\put(2.75,3.45){\scriptsize {\tiny $C_1$}}

\put(1,-4){\scriptsize {\tiny {\bf Figure 4.3.1:} Paths of type $C_1$}}

\end{picture}


\begin{picture}(0,0)(-80,48)
\setlength{\unitlength}{6mm}

\drawpolygon(0,0)(12,0)(12,8)(0,8)

\drawline[AHnb=0](0,0)(0,8)
\drawline[AHnb=0](1,0)(1,8)
\drawline[AHnb=0](2,0)(2,8)
\drawline[AHnb=0](3,0)(3,8)
\drawline[AHnb=0](4,0)(4,8)
\drawline[AHnb=0](5,0)(5,8)
\drawline[AHnb=0](6,0)(6,8)
\drawline[AHnb=0](7,0)(7,8)
\drawline[AHnb=0](8,0)(8,8)
\drawline[AHnb=0](9,0)(9,8)
\drawline[AHnb=0](10,0)(10,8)
\drawline[AHnb=0](11,0)(11,8)


\drawline[AHnb=0](0,1)(12,1)
\drawline[AHnb=0](0,2)(12,2)
\drawline[AHnb=0](0,3)(12,3)
\drawline[AHnb=0](0,4)(12,4)
\drawline[AHnb=0](0,5)(12,5)
\drawline[AHnb=0](0,6)(12,6)
\drawline[AHnb=0](0,7)(12,7)

\drawline[AHnb=0](0,5)(1,5.5)  \drawline[AHnb=0](0,6)(1,5.5)
\drawline[AHnb=0](4,5)(5,5.5)  \drawline[AHnb=0](4,6)(5,5.5)
\drawline[AHnb=0](8,5)(9,5.5)  \drawline[AHnb=0](8,6)(9,5.5)

\drawline[AHnb=0](1,5.5)(2,5.5)
\drawline[AHnb=0](5,5.5)(6,5.5)
\drawline[AHnb=0](9,5.5)(10,5.5)

\drawline[AHnb=0](2,5.5)(3,5)  \drawline[AHnb=0](2,5.5)(3,6)
\drawline[AHnb=0](6,5.5)(7,5)  \drawline[AHnb=0](6,5.5)(7,6)
\drawline[AHnb=0](10,5.5)(11,5)  \drawline[AHnb=0](10,5.5)(11,6)

\drawline[AHnb=0](3.5,5)(3.5,6)
\drawline[AHnb=0](7.5,5)(7.5,6)
\drawline[AHnb=0](11.5,5)(11.5,6)

\drawline[AHnb=0](1,6)(1.5,7)  \drawline[AHnb=0](2,6)(1.5,7)
\drawline[AHnb=0](5,6)(5.5,7)  \drawline[AHnb=0](6,6)(5.5,7)
\drawline[AHnb=0](9,6)(9.5,7)  \drawline[AHnb=0](10,6)(9.5,7)

\drawline[AHnb=0](3.5,6)(3,7)  \drawline[AHnb=0](3.5,6)(4,7)
\drawline[AHnb=0](7.5,6)(7,7)  \drawline[AHnb=0](7.5,6)(8,7)
\drawline[AHnb=0](11.5,6)(11,7)  \drawline[AHnb=0](11.5,6)(12,7)

\drawline[AHnb=0](1.5,7)(1.5,8)
\drawline[AHnb=0](5.5,7)(5.5,8)
\drawline[AHnb=0](9.5,7)(9.5,8)

\drawline[AHnb=0](3,7.5)(4,7.5)
\drawline[AHnb=0](7,7.5)(8,7.5)
\drawline[AHnb=0](11,7.5)(12,7.5)

\drawline[AHnb=0](0,7.5)(1,7)  \drawline[AHnb=0](0,7.5)(1,8)
\drawline[AHnb=0](2,7)(3,7.5)  \drawline[AHnb=0](2,8)(3,7.5)

\drawline[AHnb=0](4,7.5)(5,7)  \drawline[AHnb=0](4,7.5)(5,8)
\drawline[AHnb=0](6,7)(7,7.5)  \drawline[AHnb=0](6,8)(7,7.5)

\drawline[AHnb=0](8,7.5)(9,7)  \drawline[AHnb=0](8,7.5)(9,8)
\drawline[AHnb=0](10,7)(11,7.5)  \drawline[AHnb=0](10,8)(11,7.5)

\drawline[AHnb=0](1,1)(1.5,0)
\drawline[AHnb=0](2,1)(1.5,0)

\drawline[AHnb=0](1,5)(1.5,4)(2,5)
\drawline[AHnb=0](5,5)(5.5,4)(6,5)
\drawline[AHnb=0](9,5)(9.5,4)(10,5)

\drawline[AHnb=0](3,4)(3.5,5)(4,4)
\drawline[AHnb=0](7,4)(7.5,5)(8,4)
\drawline[AHnb=0](11,4)(11.5,5)(12,4)

\drawline[AHnb=0](7,4)(7.5,5)(8,4)
\drawline[AHnb=0](11,4)(11.5,5)(12,4)

\drawline[AHnb=0](5,5)(5.5,4)
\drawline[AHnb=0](6,5)(5.5,4)

\drawline[AHnb=0](5,1)(5.5,0)
\drawline[AHnb=0](6,1)(5.5,0)
\drawline[AHnb=0](9,1)(9.5,0)
\drawline[AHnb=0](10,1)(9.5,0)

\drawline[AHnb=0](0,1)(1,1.5)
\drawline[AHnb=0](0,2)(1,1.5)
\drawline[AHnb=0](4,1)(5,1.5)
\drawline[AHnb=0](4,2)(5,1.5)
\drawline[AHnb=0](8,1)(9,1.5)
\drawline[AHnb=0](8,2)(9,1.5)

\drawline[AHnb=0](1,2)(1.5,3)
\drawline[AHnb=0](2,2)(1.5,3)
\drawline[AHnb=0](5,2)(5.5,3)
\drawline[AHnb=0](6,2)(5.5,3)
\drawline[AHnb=0](9,2)(9.5,3)
\drawline[AHnb=0](10,2)(9.5,3)

\drawline[AHnb=0](1,3)(0,3.5)
\drawline[AHnb=0](1,4)(0,3.5)
\drawline[AHnb=0](5,3)(4,3.5)
\drawline[AHnb=0](5,4)(4,3.5)
\drawline[AHnb=0](9,3)(8,3.5)
\drawline[AHnb=0](9,4)(8,3.5)

\drawline[AHnb=0](3,0)(3.5,1)
\drawline[AHnb=0](4,0)(3.5,1)
\drawline[AHnb=0](7,0)(7.5,1)
\drawline[AHnb=0](8,0)(7.5,1)
\drawline[AHnb=0](11,0)(11.5,1)
\drawline[AHnb=0](12,0)(11.5,1)

\drawline[AHnb=0](2,1.5)(3,1)
\drawline[AHnb=0](2,1.5)(3,2)
\drawline[AHnb=0](6,1.5)(7,1)
\drawline[AHnb=0](6,1.5)(7,2)
\drawline[AHnb=0](10,1.5)(11,1)
\drawline[AHnb=0](10,1.5)(11,2)

\drawline[AHnb=0](3.5,2)(3,3)
\drawline[AHnb=0](3.5,2)(4,3)
\drawline[AHnb=0](7.5,2)(7,3)
\drawline[AHnb=0](7.5,2)(8,3)
\drawline[AHnb=0](11.5,2)(11,3)
\drawline[AHnb=0](11.5,2)(12,3)

\drawline[AHnb=0](3,3.5)(2,4)
\drawline[AHnb=0](3,3.5)(2,3)
\drawline[AHnb=0](7,3.5)(6,4)
\drawline[AHnb=0](7,3.5)(6,3)
\drawline[AHnb=0](11,3.5)(10,4)
\drawline[AHnb=0](11,3.5)(10,3)

\drawline[AHnb=0](7,3.5)(8,3.5)
\drawline[AHnb=0](11,3.5)(12,3.5)

\drawline[AHnb=0](1,1.5)(2,1.5)
\drawline[AHnb=0](3.5,1)(3.5,2)
\drawline[AHnb=0](5,1.5)(6,1.5)
\drawline[AHnb=0](7.5,1)(7.5,2)
\drawline[AHnb=0](9,1.5)(10,1.5)
\drawline[AHnb=0](11.5,1)(11.5,2)

\drawline[AHnb=0](1.5,3)(1.5,4)
\drawline[AHnb=0](3,3.5)(4,3.5)
\drawline[AHnb=0](5.5,3)(5.5,4)
\drawline[AHnb=0](7,3.5)(7,3.5)
\drawline[AHnb=0](9.5,3)(9.5,4)
\drawline[AHnb=0](11,3.5)(11,3.5)

\put(-.4,-.5){\scriptsize {\tiny $a_{11}$}}
\put(.5,-.5){\scriptsize {\tiny $a_{12}$}}
\put(1.3,-.5){\scriptsize {\tiny $a_{13}$}}
\put(2.15,-.5){\scriptsize {\tiny $a_{14}$}}
\put(10.8,-.5){\scriptsize {\tiny $a_{1i}$}}
\put(12,-.5){\scriptsize {\tiny $a_{11}$}}

\put(-.85,1){\scriptsize {\tiny $a_{21}$}}
\put(12.1,1){\scriptsize {\tiny $a_{21}$}}

\put(-.85,2){\scriptsize {\tiny $a_{31}$}}
\put(12.1,1.9){\scriptsize {\tiny $a_{31}$}}

\put(-.9,7){\scriptsize {\tiny $a_{j1}$}}
\put(12.1,7){\scriptsize {\tiny $a_{j1}$}}

\put(-1.25,8.2){\scriptsize {\tiny $a_{1(k+1)}$}}
\put(.5,8.2){\scriptsize {\tiny $a_{1(k+2)}$}}
\put(10.45,8.2){\scriptsize {\tiny $a_{1k}$}}
\put(11.45,8.2){\scriptsize {\tiny $a_{1(k+1)}$}}

\put(-1.5,-.2){\scriptsize {\tiny $Q_1$}}
\put(-1.5,1){\scriptsize {\tiny $Q_2$}}
\put(-1.5,2){\scriptsize {\tiny $Q_3$}}
\put(-1.75,7){\scriptsize {\tiny $Q_{j}$}}
\put(-1.75,8){\scriptsize {\tiny $Q_{1}$}}

\put(-1.2,4.5){\scriptsize $\vdots$}

\put(4.4,-.5){\scriptsize $\ldots$}
\put(6.4,-.5){\scriptsize $\ldots$}
\put(8.4,-.5){\scriptsize $\ldots$}

\put(4.4,8.2){\scriptsize $\ldots$}
\put(6.4,8.2){\scriptsize $\ldots$}
\put(8.4,8.2){\scriptsize $\ldots$}

\drawpolygon[fillcolor=black](0,5)(0,0)(11,0)(11.5,1)(11,1)(10,1.5)(10,1)(9,1)(9,1.5)(8,1)(7,1)(6,1.5)(6,1)(5,1)(5,1.5)(4,1)(3,1)(2,1.5)(2,1)(1,1)(1,1.5)(1,2)(11.5,2)(11,3)(11,3.5)(10,3)(9,3)(8,3.5)(8,3)(7,3)(7,3.5)(6,3)(5,3)(4,3.5)(4,3)(3,3)(3,3.5)(2,3)(1,3)(1,4)(11,4)(11.5,5)(11,5)(11,5.1)(11.6,5.1)(11.1,3.9)(1.1,3.9)(1.1,3.1)(2,3.1)(3.1,3.65)(3.1,3.1)(3.9,3.1)(3.9,3.65)(5,3.1)(6,3.1)(7.1,3.6)(7.1,3.1)(7.9,3.1)(7.9,3.65)(9,3.1)(10,3.1)(11.1,3.6)(11.1,3.1)(11.7,1.9)(1.1,1.9)(1.1,1.1)(1.9,1.1)(1.9,1.65)(3.1,1.1)(3.9,1.1)(5.1,1.6)(5.1,1.1)(5.9,1.1)(5.9,1.6)(7,1.1)(7.9,1.1)(9.1,1.6)(9.1,1.1)(9.9,1.1)(9.9,1.6)(11,1.1)(11.65,1.1)(11.1,-.1)(-.1,-.1)(-.1,5)

\drawpolygon[fillcolor=black](0,5)(0,7.5)(-.1,7.5)(-.1,5)

\drawpolygon[fillcolor=black](-.1,7.5)(1,7)(2,7)(3,7.5)(3,7)(4,7)(4,7.5)(5,7)(6,7)(7,7.5)(7,7)(8,7)(8,7.5)(9,7)(10,7)(11,7.5)(11,7)(11.1,7)(11.1,7.6)(10,7.1)(9,7.1)(7.9,7.65)(7.9,7.1)(7.1,7.1)(7.1,7.65)(6,7.1)(5,7.1)(3.9,7.65)(3.9,7.1)(3.1,7.1)(3.1,7.65)(1.9,7.1)(1,7.1)(-.1,7.6)

\drawpolygon[fillcolor=black](11,7)(11.5,6)(1,6)(1,5.9)(11.65,5.9)(11.1,7)

\drawpolygon[fillcolor=black](11,5)(10,5.5)(10,5)(9,5)(9,5.5)(8,5)(7,5)(6,5.5)(6,5)(5,5)(5,5.5)(4,5)(3,5)(2,5.5)(2,5)(1,5)(1,6)(1.1,6)(1.1,5.1)(1.9,5.1)(1.9,5.65)(3.1,5.1)(3.9,5.1)(5.1,5.6)(5.1,5.1)(5.9,5.1)(5.9,5.6)(7,5.1)(7.9,5.1)(9.1,5.6)(9.1,5.1)(9.9,5.1)(9.9,5.6)(11,5.1)

\put(-1,-1.5){\scriptsize {\tiny {\bf Figure 4.3.2:} $M(i,j=4m+4,k)$ of DSEM of type $[3^3.4^2 : 3^2.4.3.4]_1$}}

\end{picture}

\vspace{5cm}

\begin{picture}(0,0)(-40,48)
\setlength{\unitlength}{5.8mm}

\drawpolygon(0,0)(12,0)(12,6)(0,6)


\drawline[AHnb=0](1,0)(1,6)
\drawline[AHnb=0](2,0)(2,6)
\drawline[AHnb=0](3,0)(3,6)
\drawline[AHnb=0](4,0)(4,6)
\drawline[AHnb=0](5,0)(5,6)
\drawline[AHnb=0](6,0)(6,6)
\drawline[AHnb=0](7,0)(7,6)
\drawline[AHnb=0](8,0)(8,6)
\drawline[AHnb=0](9,0)(9,6)
\drawline[AHnb=0](10,0)(10,6)
\drawline[AHnb=0](11,0)(11,6)


\drawline[AHnb=0](0,1)(12,1)
\drawline[AHnb=0](0,2)(12,2)
\drawline[AHnb=0](0,3)(12,3)
\drawline[AHnb=0](0,4)(12,4)
\drawline[AHnb=0](0,5)(12,5)
\drawline[AHnb=0](0,6)(12,6)
\drawline[AHnb=0](0,5)(12,5)

\drawline[AHnb=0](1,1)(1.5,0)
\drawline[AHnb=0](2,1)(1.5,0)

\drawline[AHnb=0](5,1)(5.5,0)
\drawline[AHnb=0](6,1)(5.5,0)

\drawline[AHnb=0](0,1)(1,1.5)
\drawline[AHnb=0](0,2)(1,1.5)

\drawline[AHnb=0](4,1)(5,1.5)
\drawline[AHnb=0](4,2)(5,1.5)

\drawline[AHnb=0](1,2)(1.5,3)
\drawline[AHnb=0](2,2)(1.5,3)

\drawline[AHnb=0](5,2)(5.5,3)
\drawline[AHnb=0](6,2)(5.5,3)

\drawline[AHnb=0](1,3)(0,3.5)
\drawline[AHnb=0](1,4)(0,3.5)

\drawline[AHnb=0](5,3)(4,3.5)
\drawline[AHnb=0](5,4)(4,3.5)

\drawline[AHnb=0](0,5)(1,5.5)
\drawline[AHnb=0](0,6)(1,5.5)

\drawline[AHnb=0](4,5)(5,5.5)
\drawline[AHnb=0](4,6)(5,5.5)

\drawline[AHnb=0](3,0)(3.5,1)
\drawline[AHnb=0](4,0)(3.5,1)

\drawline[AHnb=0](7,0)(7.5,1)
\drawline[AHnb=0](8,0)(7.5,1)

\drawline[AHnb=0](2,1.5)(3,1)
\drawline[AHnb=0](2,1.5)(3,2)

\drawline[AHnb=0](6,1.5)(7,1)
\drawline[AHnb=0](6,1.5)(7,2)

\drawline[AHnb=0](3.5,2)(3,3)
\drawline[AHnb=0](3.5,2)(4,3)

\drawline[AHnb=0](7.5,2)(7,3)
\drawline[AHnb=0](7.5,2)(8,3)

\drawline[AHnb=0](3,3.5)(2,4)
\drawline[AHnb=0](3,3.5)(2,3)

\drawline[AHnb=0](7,3.5)(6,4)
\drawline[AHnb=0](7,3.5)(6,3)

\drawline[AHnb=0](2,5.5)(3,5)
\drawline[AHnb=0](2,5.5)(3,6)

\drawline[AHnb=0](6,5.5)(7,5)
\drawline[AHnb=0](6,5.5)(7,6)

\drawline[AHnb=0](1,5)(1.5,4)
\drawline[AHnb=0](2,5)(1.5,4)

\drawline[AHnb=0](5,5)(5.5,4)
\drawline[AHnb=0](6,5)(5.5,4)

\drawline[AHnb=0](3,4)(3.5,5)
\drawline[AHnb=0](4,4)(3.5,5)

\drawline[AHnb=0](7,4)(7.5,5)
\drawline[AHnb=0](8,4)(7.5,5)

\drawline[AHnb=0](7,3.5)(8,3.5)

\drawline[AHnb=0](1,1.5)(2,1.5)
\drawline[AHnb=0](3.5,1)(3.5,2)

\drawline[AHnb=0](1.5,3)(1.5,4)
\drawline[AHnb=0](3,3.5)(4,3.5)

\drawline[AHnb=0](1,5.5)(2,5.5)
\drawline[AHnb=0](3.5,5)(3.5,6)

\drawline[AHnb=0](5,1.5)(6,1.5)
\drawline[AHnb=0](7.5,1)(7.5,2)

\drawline[AHnb=0](5.5,3)(5.5,4)
\drawline[AHnb=0](7,3.5)(7,3.5)

\drawline[AHnb=0](5,5.5)(6,5.5)
\drawline[AHnb=0](7.5,5)(7.5,6)

\drawline[AHnb=0](1,1)(1.5,0)
\drawline[AHnb=0](2,1)(1.5,0)

\drawline[AHnb=0](1,5)(1.5,4)(2,5)
\drawline[AHnb=0](5,5)(5.5,4)(6,5)
\drawline[AHnb=0](9,5)(9.5,4)(10,5)

\drawline[AHnb=0](3,4)(3.5,5)(4,4)
\drawline[AHnb=0](7,4)(7.5,5)(8,4)
\drawline[AHnb=0](11,4)(11.5,5)(12,4)

\drawline[AHnb=0](7,4)(7.5,5)(8,4)
\drawline[AHnb=0](11,4)(11.5,5)(12,4)

\drawline[AHnb=0](9,5.5)(10,5.5)

\drawline[AHnb=0](8,5)(9,5.5)(8,6)
\drawline[AHnb=0](11,5)(10,5.5)(11,6)

\drawline[AHnb=0](5,5)(5.5,4)
\drawline[AHnb=0](6,5)(5.5,4)

\drawline[AHnb=0](5,1)(5.5,0)
\drawline[AHnb=0](6,1)(5.5,0)
\drawline[AHnb=0](9,1)(9.5,0)
\drawline[AHnb=0](10,1)(9.5,0)

\drawline[AHnb=0](0,1)(1,1.5)
\drawline[AHnb=0](0,2)(1,1.5)
\drawline[AHnb=0](4,1)(5,1.5)
\drawline[AHnb=0](4,2)(5,1.5)
\drawline[AHnb=0](8,1)(9,1.5)
\drawline[AHnb=0](8,2)(9,1.5)

\drawline[AHnb=0](1,2)(1.5,3)
\drawline[AHnb=0](2,2)(1.5,3)
\drawline[AHnb=0](5,2)(5.5,3)
\drawline[AHnb=0](6,2)(5.5,3)
\drawline[AHnb=0](9,2)(9.5,3)
\drawline[AHnb=0](10,2)(9.5,3)

\drawline[AHnb=0](1,3)(0,3.5)
\drawline[AHnb=0](1,4)(0,3.5)
\drawline[AHnb=0](5,3)(4,3.5)
\drawline[AHnb=0](5,4)(4,3.5)
\drawline[AHnb=0](9,3)(8,3.5)
\drawline[AHnb=0](9,4)(8,3.5)

\drawline[AHnb=0](3,0)(3.5,1)
\drawline[AHnb=0](4,0)(3.5,1)
\drawline[AHnb=0](7,0)(7.5,1)
\drawline[AHnb=0](8,0)(7.5,1)
\drawline[AHnb=0](11,0)(11.5,1)
\drawline[AHnb=0](12,0)(11.5,1)

\drawline[AHnb=0](2,1.5)(3,1)
\drawline[AHnb=0](2,1.5)(3,2)
\drawline[AHnb=0](6,1.5)(7,1)
\drawline[AHnb=0](6,1.5)(7,2)
\drawline[AHnb=0](10,1.5)(11,1)
\drawline[AHnb=0](10,1.5)(11,2)

\drawline[AHnb=0](3.5,2)(3,3)
\drawline[AHnb=0](3.5,2)(4,3)
\drawline[AHnb=0](7.5,2)(7,3)
\drawline[AHnb=0](7.5,2)(8,3)
\drawline[AHnb=0](11.5,2)(11,3)
\drawline[AHnb=0](11.5,2)(12,3)

\drawline[AHnb=0](3,3.5)(2,4)
\drawline[AHnb=0](3,3.5)(2,3)
\drawline[AHnb=0](7,3.5)(6,4)
\drawline[AHnb=0](7,3.5)(6,3)
\drawline[AHnb=0](11,3.5)(10,4)
\drawline[AHnb=0](11,3.5)(10,3)
\drawline[AHnb=0](7,3.5)(8,3.5)
\drawline[AHnb=0](11,3.5)(12,3.5)
\drawline[AHnb=0](11.5,6)(11.5,5)


\drawline[AHnb=0](1,1.5)(2,1.5)
\drawline[AHnb=0](3.5,1)(3.5,2)
\drawline[AHnb=0](5,1.5)(6,1.5)
\drawline[AHnb=0](7.5,1)(7.5,2)
\drawline[AHnb=0](9,1.5)(10,1.5)
\drawline[AHnb=0](11.5,1)(11.5,2)

\drawline[AHnb=0](1.5,3)(1.5,4)
\drawline[AHnb=0](3,3.5)(4,3.5)
\drawline[AHnb=0](5.5,3)(5.5,4)
\drawline[AHnb=0](7,3.5)(7,3.5)
\drawline[AHnb=0](9.5,3)(9.5,4)
\drawline[AHnb=0](11,3.5)(11,3.5)

\put(-.4,-.5){\scriptsize {\tiny $a_{11}$}}
\put(.4,-.5){\scriptsize {\tiny $a_{12}$}}
\put(1.3,-.5){\scriptsize {\tiny $a_{13}$}}
\put(10.8,-.5){\scriptsize {\tiny $a_{1i}$}}
\put(12,-.5){\scriptsize {\tiny $a_{11}$}}

\put(-.85,1){\scriptsize {\tiny $a_{21}$}}
\put(12.1,1){\scriptsize {\tiny $a_{21}$}}

\put(-.85,2){\scriptsize {\tiny $a_{31}$}}
\put(12.1,2){\scriptsize {\tiny $a_{31}$}}

\put(-.85,5){\scriptsize {\tiny $a_{j1}$}}
\put(12.1,5){\scriptsize {\tiny $a_{j1}$}}

\put(-1,6.5){\scriptsize {\tiny $a_{1(k+1)}$}}
\put(.7,6.5){\scriptsize {\tiny $a_{1(k+2)}$}}
\put(10.9,6.5){\scriptsize {\tiny $a_{1k}$}}
\put(11.7,6.5){\scriptsize {\tiny $a_{1(k+1)}$}}

\put(-1.4,-.1){\scriptsize {\tiny $Q_1$}}
\put(-1.5,1){\scriptsize {\tiny $Q_2$}}
\put(-1.5,2){\scriptsize {\tiny $Q_3$}}
\put(-1.5,5){\scriptsize {\tiny $Q_j$}}
\put(-1.4,5.8){\scriptsize {\tiny $Q_1$}}

\put(3.4,-.5){\scriptsize $\ldots$}
\put(5.4,-.5){\scriptsize $\ldots$}
\put(7.4,-.5){\scriptsize $\ldots$}
\put(9.4,-.5){\scriptsize $\ldots$}

\put(3.4,6.65){\scriptsize $\ldots$}
\put(6,6.65){\scriptsize $\ldots$}
\put(8.4,6.65){\scriptsize $\ldots$}

\put(-1.4,3.5){\scriptsize $\vdots$}

\drawpolygon[fillcolor=black](0,5)(0,0)(11,0)(11.5,1)(11,1)(10,1.5)(10,1)(9,1)(9,1.5)(8,1)(7,1)(6,1.5)(6,1)(5,1)(5,1.5)(4,1)(3,1)(2,1.5)(2,1)(1,1)(1,1.5)(1,2)(11.5,2)(11,3)(11,3.5)(10,3)(9,3)(8,3.5)(8,3)(7,3)(7,3.5)(6,3)(5,3)(4,3.5)(4,3)(3,3)(3,3.5)(2,3)(1,3)(1,4)(11,4)(11.5,5)(11,5)(10,5.5)(10,5)(9,5)(9,5.5)(8,5)(7,5)(6,5.5)(6,5)(5,5)(5,5.5)(4,5)(3,5)(2,5.5)(2,5)(1,5)(1,5.5)(0,5)(-.1,5.1)(1.1,5.6)(1.1,5.1)(1.9,5.1)(1.9,5.6)(3,5.1)(3.9,5.1)(5.1,5.6)(5.1,5.1)(5.9,5.1)(5.9,5.6)(7,5.1)(8,5.1)(9.1,5.6)(9.1,5.1)(9.9,5.1)(9.9,5.6)(11,5.1)(11.6,5.1)(11.1,3.9)(1.1,3.9)(1.1,3.1)(2,3.1)(3.1,3.65)(3.1,3.1)(3.9,3.1)(3.9,3.65)(5,3.1)(6,3.1)(7.1,3.6)(7.1,3.1)(7.9,3.1)(7.9,3.65)(9,3.1)(10,3.1)(11.1,3.6)(11.1,3.1)(11.7,1.9)(1.1,1.9)(1.1,1.1)(1.9,1.1)(1.9,1.65)(3.1,1.1)(3.9,1.1)(5.1,1.6)(5.1,1.1)(5.9,1.1)(5.9,1.6)(7,1.1)(7.9,1.1)(9.1,1.6)(9.1,1.1)(9.9,1.1)(9.9,1.6)(11,1.1)(11.65,1.1)(11.1,-.1)(-.1,-.1)(-.1,5.1)

\put(-1.25,-1.5){\scriptsize {\tiny {\bf Figure 4.3.3:} $M(i,j=4m+2,k)$ of DSEM of type $[3^3.4^2 : 3^2.4.3.4]_1$}} 

\end{picture}

\vspace{5.75cm}

\begin{lem}  \label{l3.3.1}
	A DSEM $M$ of type $[3^3.4^2: 3^2.4.3.4]_1$ admits an $M(i,j,k)$ representation iff the following holds: $(i)$ $i=5m$, $ m  \in \mathbb N $ and $j$ is even, $(ii)$ number of vertices of $M(i,j,k)=6ij/5 \geq 12$, $(iii)$ if $j=4m + 2$, $m\in \mathbb N\cup \{0\}$, then $k \in \{5r + 3: 0 \leq r < i/5\} $, and if $j=4m$, $ m \in \mathbb N$, then $k \in \{5r: 0\leq r < i/5\} $. 
\end{lem}

\noindent{\bf Proof.} A representation $M(i, j, k)$ of a DSEM $M$ of type $[3^3.4^2: 3^2.4.3.4]_1$ has $j$  disjoint horizontal cycles of $C_{1}$ type having length $i$. Let $Q_0, Q_1, \ldots , Q_{j-1}$ be the horizontal cycles of type $C_{1}$. The number of vertices with face-sequence $(3^3,4^2)$ lying between horizontal cycles $Q_{(2s+1)(mod\,j)}$ and $Q_{(2s+2)(mod\,j)}$, for $0 \leq s \leq j-1$, is $2i/5 \cdot j/2$. Thus, the total number of vertices in $M$ is $n = ij + ij/5= 6ij/5$. If $j=1$, then $M(i,1,k)$ has no vertex with face-sequence $(3^3,4^2)$ or $(3^2,4,3,4)$. So $j \geq 2 $. If $j \geq 2 $ and $j$ is not an even integer then we get some vertices in the base horizontal cycle which does not have the face-sequence $(3^3,4^2)$ or $(3^2,4,3,4)$. So, $j$ is even.

If $j$ is even and $i < 5$ then again the representation $M(i,j,k)$  has some vertices which do not have the face-sequence $(3^3,4^2)$ or $(3^2,4,3,4)$. So, $i \geq 5 $. If $i \geq 5 $ and not a multiple of 5, then $2|V(3^3,4^2)| \neq |V(3^2,4,3,4)|$. This is not possible. So $ i =5m$, where $m \in \mathbb{N}$ and $n=6ij/5 \geq 12$.

If $j = 4m+2$, $m  \in \mathbb N \cup \{0\}$ and  $k \in \{ r : 0 \leq r \leq i-1 \}\setminus \{5r+3: 0 \leq r < i/5 \} $ then we get some vertices which do not have the face-sequence $(3^3,4^2)$ or $(3^2,4,3,4)$. So, $k \in  \{5r+3 : 0 \leq r < i/5\}$  for $j = 4m+2$, $m  \in \mathbb N \cup \{0\}$. Similarly if $j=4m$, $m\in \mathbb N$, then $k \in \{5r: 0 \leq r < i/5 \}$. This completes the proof. \hfill $\Box$

%

\subsection{ DSEMs of types $[3^3.4^2: 3^2.4.3.4]_2$}\label{s3.4}

Let $M$ be a DSEM of type $[3^3. 4^2: 3^2.4.3.4]_2$. It is easy to see that, if $M$ exists then $|V_{(3^3,4^2)}| = |V_{(3^2,4,3,4)}|$. We consider following types of paths in $M$ as follows.

A path $P_{1} = P( \ldots, u_{i},v_{j},u_{i+1},v_{j+1}, \ldots)$ in $M$, say of type $Y_{1}$, indicated by thick black paths in Figure 4.4.1. Here the vertices $u_i$'s and $v_j$'s have face-sequences $(3^2,4,3,4)$ and $(3^3,4^2)$ respectively. 

A path $P_{2} = P( \ldots, u_{i-1},u_{i},u_{i+1}, \ldots)$ in $M$, say of type $Y_{2}$, indicated by green paths in Figure 4.4.1. Here the vertices $u_i$'s of the path has face-sequence either $(3^3,4^2)$ or $(3^2,4,3,4)$.

\vspace{.2cm}

To construct $M(i, j, k)$, we cut $M$ along a cycle of type $Y_1$ and then, take second cut along a cycle of type $Y_2$, assume without loss of generality $Y_2$ whose all the vertices have the face-sequence $(3^2,4,3,4)$. Following, a similar argument given in Lemma \ref{l4.1.3}, we obtain Lemma \ref{l3.4.1}.

\begin{picture}(0,0)(-21,30)
\setlength{\unitlength}{6.5mm}

\drawline[AHnb=0](-2,1)(0,0)(2,0)(4,-1)(4.90,-1)
\drawline[AHnb=0](-2.5,0)(-2,0)(0,-1)(2,-1)(4,-2)(4.90,-2)

\drawline[AHnb=0](-1.5,2)(.5,1)(2.5,1)(4.5,0)
\drawline[AHnb=0](-1,2.75)(.5,2)(2.5,2)(4.5,1)

\drawline[AHnb=0](-2.5,-1)(-.5,-2)(1.5,-2)(3.5,-3)(4.5,-3)

\drawline[AHnb=0](-2,1)(-2,0)
\drawline[AHnb=0](-1,.5)(-1,-.5)
\drawline[AHnb=0](0,0)(0,-1)
\drawline[AHnb=0](1,0)(1,-1)
\drawline[AHnb=0](2,0)(2,-1)
\drawline[AHnb=0](3,-.5)(3,-1.5)
\drawline[AHnb=0](4,-1)(4,-2)

\drawline[AHnb=0](0,0)(.5,1)
\drawline[AHnb=0](1,0)(.5,1)

\drawline[AHnb=0](1,0)(1.5,1)
\drawline[AHnb=0](2,0)(1.5,1)
\drawline[AHnb=0](2,0)(2.5,1)
\drawline[AHnb=0](3,-.5)(3.5,.5)
\drawline[AHnb=0](4,-1)(4.5,0)
\drawline[AHnb=0](-1,.5)(-.5,1.5)
\drawline[AHnb=0](-2,1)(-1.5,2)

\drawline[AHnb=0](-2.5,-1)(-2,0)
\drawline[AHnb=0](-1,-.5)(-1.5,-1.5)
\drawline[AHnb=0](0,-1)(-.5,-2)
\drawline[AHnb=0](0,-1)(.5,-2)
\drawline[AHnb=0](1,-1)(.5,-2)
\drawline[AHnb=0](1,-1)(1.5,-2)
\drawline[AHnb=0](2,-1)(1.5,-2)
\drawline[AHnb=0](3,-1.5)(2.5,-2.5)

\drawline[AHnb=0](4,-2)(3.5,-3)
\drawline[AHnb=0](3,-1.5)(2.5,-2.5)

\drawline[AHnb=0](-1.5,2)(-1.5,2.5)
\drawline[AHnb=0](-.5,1.5)(-.5,2.5)
\drawline[AHnb=0](-.5,2.5)(-.25,3)

\drawline[AHnb=0](.5,1)(.5,2)
\drawline[AHnb=0](1.5,1)(1.5,2)
\drawline[AHnb=0](2.5,1)(2.5,2)
\drawline[AHnb=0](3.5,.5)(3.5,1.5)
\drawline[AHnb=0](4.5,0)(4.5,1)

\drawline[AHnb=0](-.5,-2)(-.5,-3)
\drawline[AHnb=0](.5,-2)(.5,-3)
\drawline[AHnb=0](1.5,-2)(1.5,-3)

\drawline[AHnb=0](1.5,-3)(2.5,-2.5)
\drawline[AHnb=0](1.5,-2)(1.5,-3)

\drawline[AHnb=0](-.5,-3)(-2.5,-2)
\drawline[AHnb=0](-2.5,-1)(-2.5,-2)
\drawline[AHnb=0](-1.5,-1.5)(-1.5,-2.5)

\drawline[AHnb=0](-.5,-3)(1.5,-3)

\drawline[AHnb=0](-1.5,2)(-.5,2.5)
\drawline[AHnb=0](-.5,1.5)(.5,2)

\drawline[AHnb=0](2.5,1)(3.5,1.5)
\drawline[AHnb=0](3.5,.5)(4.5,1)

\drawline[AHnb=0](4.5,0)(5,0)
\drawline[AHnb=0](4.5,1)(5,1)

\drawline[AHnb=0](-2,0)(-1,.5)

\drawline[AHnb=0](-1,-.5)(0,0)

\drawline[AHnb=0](2,-1)(3,-.5)
\drawline[AHnb=0](3,-1.5)(4,-1)

\drawline[AHnb=0](-2.5,-2)(-1.5,-1.5)
\drawline[AHnb=0](-1.5,-2.5)(-.5,-2)

\drawline[AHnb=0](-.5,-3)(-.25,-3.25)
\drawline[AHnb=0](.5,-3)(.25,-3.25)

\drawline[AHnb=0](1.5,-3)(1.25,-3.25)
\drawline[AHnb=0](.5,-3)(.75,-3.25)

\drawline[AHnb=0](-2.5,-2)(-3,-2.5)
\drawline[AHnb=0](-1.5,-2.5)(-2,-3)
\drawline[AHnb=0](-.5,-3)(-.75,-3.25)

\drawline[AHnb=0](1.5,-3)(2,-3.25)
\drawline[AHnb=0](2.5,-2.5)(2.5,-3.25)
\drawline[AHnb=0](3.5,-3)(3.5,-3.35)
\drawline[AHnb=0](3.5,-3)(3,-3.25)

\drawline[AHnb=0](-2.5,-2)(-3,-2)
\drawline[AHnb=0](-2.5,-1)(-3,-1)

\drawline[AHnb=0](-2,0)(-2.5,0)
\drawline[AHnb=0](-2,1)(-2.5,1)
\drawline[AHnb=0](-1.5,2)(-2,2)

\drawline[AHnb=0](4,-2)(4.25,-2.5)
\drawline[AHnb=0](4.5,0)(4.75,-.5)

\drawline[AHnb=0](.5,2)(.75,2.5)
\drawline[AHnb=0](1.5,2)(1.25,2.5)

\drawline[AHnb=0](.5,2)(.75,2.5)
\drawline[AHnb=0](1.5,2)(1.25,2.5)

\drawline[AHnb=0](1.5,2)(1.75,2.5)
\drawline[AHnb=0](2.5,2)(2.25,2.5)
\drawline[AHnb=0](2.5,2)(2.75,2.5)

\drawline[AHnb=0](3.5,1.5)(3.75,2)
\drawline[AHnb=0](4.5,1)(4.75,1.5)

\drawline[AHnb=0](-2.5,-1)(-2.75,-.5)
\drawline[AHnb=0](-2,1)(-2.25,1.5)

\drawpolygon[fillcolor=black](-3,-1)(-2.5,-1)(-.5,-2)(1.5,-2)(3.5,-3)(4.5,-3)(4.5,-3.1)(3.5,-3.1)(1.5,-2.1)(-.5,-2.1)(-2.5,-1.1)(-3,-1.1)

\drawpolygon[fillcolor=black](-2.5,0)(-2,0)(0,-1)(2,-1)(4,-2)(5,-2)(5,-2.1)(4,-2.1)(2,-1.1)(0,-1.1)(-2,-.1)(-2.5,-0.1)

\drawpolygon[fillcolor=black](-2.5,1)(-2,1)(0,0)(2,0)(4,-1)(5,-1)(5,-1.1)(4,-1.1)(2,-.1)(0,-.1)(-2,.9)(-2.5,0.9)

\drawpolygon[fillcolor=black](-2,2)(-1.5,2)(0.5,1)(2.5,1)(4.5,0)(5.5,0)(5.5,-.1)(4.5,-.1)(2.5,.9)(0.5,.9)(-1.5,1.9)(-2,1.9)

\drawpolygon[fillcolor=green](-2,-3)(-1.5,-2.5)(-1.5,-1.5)(-1,-.5)(-1,.5)(-.5,1.5)(-.5,2.5)(-.25,3)(-.15,3)(-.4,2.5)(-.4,1.5)(-.9,.5)(-.9,-.5)(-1.4,-1.5)(-1.4,-2.5)(-1.85,-3)

\drawpolygon[fillcolor=green](-.75,-3.25)(-.5,-3)(-.5,-2)(0,-1)(0,0)(0.5,1)(0.5,2)(0.75,2.5)(0.85,2.5)(0.6,2)(0.6,1)(0.1,0)(0.1,-1)(-.4,-2)(-.4,-3)(-.65,-3.25)

\drawpolygon[fillcolor=green](.25,-3.25)(.5,-3)(.5,-2)(1,-1)(1,0)(1.5,1)(1.5,2)(1.75,2.5)(1.85,2.5)(1.6,2)(1.6,1)(1.1,0)(1.1,-1)(.6,-2)(.6,-3)(.35,-3.25)

\drawpolygon[fillcolor=green](1.25,-3.25)(1.5,-3)(1.5,-2)(2,-1)(2,0)(2.5,1)(2.5,2)(2.75,2.5)(2.85,2.5)(2.6,2)(2.6,1)(2.1,0)(2.1,-1)(1.6,-2)(1.6,-3)(1.35,-3.25)

\drawpolygon[fillcolor=green](2.25,-3.75)(2.5,-3.5)(2.5,-2.5)(3,-1.5)(3,-.5)(3.5,.5)(3.5,1.5)(3.75,2)(3.85,2)(3.6,1.5)(3.6,.5)(3.1,-.5)(3.1,-1.5)(2.6,-2.5)(2.6,-3.5)(2.35,-3.75)

\put(1.55,2.75) {\scriptsize {\tiny $Y_2$}}
\put(5.5,-2.1) {\scriptsize {\tiny $Y_1$}}

\put(-2.5,-4.4){\scriptsize {\tiny {\bf Figure 4.4.1:} Paths of types $Y_1$ and $Y_2$}} 

\end{picture}

\begin{picture}(0,0)(-77,46)
\setlength{\unitlength}{6mm}

\drawpolygon(0,0)(10,0)(10,8)(0,8)

\drawline[AHnb=0](1,0)(1,8)
\drawline[AHnb=0](2,0)(2,8)
\drawline[AHnb=0](3,0)(3,8)
\drawline[AHnb=0](4,0)(4,8)
\drawline[AHnb=0](5,0)(5,8)
\drawline[AHnb=0](6,0)(6,8)
\drawline[AHnb=0](7,0)(7,8)
\drawline[AHnb=0](8,0)(8,8)
\drawline[AHnb=0](9,0)(9,8)

\drawline[AHnb=0](0,1)(10,1)
\drawline[AHnb=0](0,2)(10,2)
\drawline[AHnb=0](0,3)(10,3)
\drawline[AHnb=0](0,4)(10,4)
\drawline[AHnb=0](0,5)(10,5)
\drawline[AHnb=0](0,6)(10,6)
\drawline[AHnb=0](0,7)(10,7)

\drawline[AHnb=0](0,1)(1,0)
\drawline[AHnb=0](1,1)(2,0)
\drawline[AHnb=0](4,1)(5,0)
\drawline[AHnb=0](5,1)(6,0)
\drawline[AHnb=0](8,1)(9,0)
\drawline[AHnb=0](9,1)(10,0)

\drawline[AHnb=0](0,3)(1,2)
\drawline[AHnb=0](1,3)(2,2)
\drawline[AHnb=0](4,3)(5,2)
\drawline[AHnb=0](5,3)(6,2)
\drawline[AHnb=0](8,3)(9,2)
\drawline[AHnb=0](9,3)(10,2)

\drawline[AHnb=0](0,5)(1,4)
\drawline[AHnb=0](1,5)(2,4)
\drawline[AHnb=0](4,5)(5,4)
\drawline[AHnb=0](5,5)(6,4)
\drawline[AHnb=0](8,5)(9,4)
\drawline[AHnb=0](9,5)(10,4)

\drawline[AHnb=0](0,7)(1,6)
\drawline[AHnb=0](1,7)(2,6)
\drawline[AHnb=0](4,7)(5,6)
\drawline[AHnb=0](5,7)(6,6)
\drawline[AHnb=0](8,7)(9,6)
\drawline[AHnb=0](9,7)(10,6)


\drawline[AHnb=0](2,1)(3,2)
\drawline[AHnb=0](3,1)(4,2)
\drawline[AHnb=0](6,1)(7,2)
\drawline[AHnb=0](7,1)(8,2)

\drawline[AHnb=0](2,3)(3,4)
\drawline[AHnb=0](3,3)(4,4)
\drawline[AHnb=0](6,3)(7,4)
\drawline[AHnb=0](7,3)(8,4)

\drawline[AHnb=0](2,5)(3,6)
\drawline[AHnb=0](3,5)(4,6)
\drawline[AHnb=0](6,5)(7,6)
\drawline[AHnb=0](7,5)(8,6)

\drawline[AHnb=0](2,7)(3,8)
\drawline[AHnb=0](3,7)(4,8)
\drawline[AHnb=0](6,7)(7,8)
\drawline[AHnb=0](7,7)(8,8)

\drawpolygon[fillcolor=black](.95,0)(1,.05)(.1,.95)(0,.95)
\drawpolygon[fillcolor=black](1.95,0)(2,.05)(1.1,.95)(1,.95)
\drawpolygon[fillcolor=black](4.95,0)(5,.05)(4.1,.95)(4,.95)
\drawpolygon[fillcolor=black](5.95,0)(6,.05)(5.1,.95)(5,.95)
\drawpolygon[fillcolor=black](8.95,0)(9,.05)(8.1,.95)(8,.95)
\drawpolygon[fillcolor=black](9.95,0)(10,.05)(9.1,.95)(9,.95)

\drawpolygon[fillcolor=black](.95,2)(1,2.05)(.1,2.95)(0,2.95)
\drawpolygon[fillcolor=black](1.95,2)(2,2.05)(1.1,2.95)(1,2.95)
\drawpolygon[fillcolor=black](4.95,2)(5,2.05)(4.1,2.95)(4,2.95)
\drawpolygon[fillcolor=black](5.95,2)(6,2.05)(5.1,2.95)(5,2.95)
\drawpolygon[fillcolor=black](8.95,2)(9,2.05)(8.1,2.95)(8,2.95)
\drawpolygon[fillcolor=black](9.95,2)(10,2.05)(9.1,2.95)(9,2.95)

\drawpolygon[fillcolor=black](.95,4)(1,4.05)(.1,4.95)(0,4.95)
\drawpolygon[fillcolor=black](1.95,4)(2,4.05)(1.1,4.95)(1,4.95)
\drawpolygon[fillcolor=black](4.95,4)(5,4.05)(4.1,4.95)(4,4.95)
\drawpolygon[fillcolor=black](5.95,4)(6,4.05)(5.1,4.95)(5,4.95)
\drawpolygon[fillcolor=black](8.95,4)(9,4.05)(8.1,4.95)(8,4.95)
\drawpolygon[fillcolor=black](9.95,4)(10,4.05)(9.1,4.95)(9,4.95)

\drawpolygon[fillcolor=black](.95,6)(1,6.05)(.1,6.95)(0,6.95)
\drawpolygon[fillcolor=black](1.95,6)(2,6.05)(1.1,6.95)(1,6.95)
\drawpolygon[fillcolor=black](4.95,6)(5,6.05)(4.1,6.95)(4,6.95)
\drawpolygon[fillcolor=black](5.95,6)(6,6.05)(5.1,6.95)(5,6.95)
\drawpolygon[fillcolor=black](8.95,6)(9,6.05)(8.1,6.95)(8,6.95)
\drawpolygon[fillcolor=black](9.95,6)(10,6.05)(9.1,6.95)(9,6.95)

\drawpolygon[fillcolor=black](2.08,1)(3,1.93)(2.95,2)(2,1.05)
\drawpolygon[fillcolor=black](3.08,1)(4,1.93)(3.95,2)(3,1.05)
\drawpolygon[fillcolor=black](6.08,1)(7,1.93)(6.95,2)(6,1.05)
\drawpolygon[fillcolor=black](7.08,1)(8,1.93)(7.95,2)(7,1.05)

\drawpolygon[fillcolor=black](2.08,3)(3,3.93)(2.95,4)(2,3.05)
\drawpolygon[fillcolor=black](3.08,3)(4,3.93)(3.95,4)(3,3.05)
\drawpolygon[fillcolor=black](6.08,3)(7,3.93)(6.95,4)(6,3.05)
\drawpolygon[fillcolor=black](7.08,3)(8,3.93)(7.95,4)(7,3.05)

\drawpolygon[fillcolor=black](2.08,5)(3,5.93)(2.95,6)(2,5.05)
\drawpolygon[fillcolor=black](3.08,5)(4,5.93)(3.95,6)(3,5.05)
\drawpolygon[fillcolor=black](6.08,5)(7,5.93)(6.95,6)(6,5.05)
\drawpolygon[fillcolor=black](7.08,5)(8,5.93)(7.95,6)(7,5.05)

\drawpolygon[fillcolor=black](2.08,7)(3,7.93)(2.95,8)(2,7.05)
\drawpolygon[fillcolor=black](3.08,7)(4,7.93)(3.95,8)(3,7.05)
\drawpolygon[fillcolor=black](6.08,7)(7,7.93)(6.95,8)(6,7.05)
\drawpolygon[fillcolor=black](7.08,7)(8,7.93)(7.95,8)(7,7.05)

\put(-.2,-.4){\scriptsize {\tiny $a_{11}$}}
\put(.8,-.4){\scriptsize {\tiny $a_{12}$}}
\put(1.8,-.4){\scriptsize {\tiny $a_{13}$}}
\put(2.8,-.4){\scriptsize {\tiny $a_{14}$}}
\put(7.3,-.4){\scriptsize {\tiny $a_{1(i-1)}$}}
\put(8.8,-.4){\scriptsize {\tiny $a_{1i}$}}
\put(9.8,-.4){\scriptsize {\tiny $a_{1i}$}}

\put(-.6,.9){\scriptsize {\tiny $a_{21}$}}
\put(-.6,1.9){\scriptsize {\tiny $a_{31}$}}
\put(-.6,2.9){\scriptsize {\tiny $a_{41}$}}
\put(-.6,6.9){\scriptsize {\tiny $a_{j1}$}}
\put(-1.1,8.2){\scriptsize {\tiny $a_{1(k+1)}$}}
\put(.5,8.2){\scriptsize {\tiny $a_{1(k+2)}$}}

\put(8.5,8.2){\scriptsize {\tiny $a_{1k}$}}
\put(9.5,8.2){\scriptsize {\tiny $a_{1(k+1)}$}}

\put(10.1,.9){\scriptsize {\tiny $a_{21}$}}
\put(10.1,1.9){\scriptsize {\tiny $a_{31}$}}
\put(10.1,2.9){\scriptsize {\tiny $a_{41}$}}
\put(10.1,6.9){\scriptsize {\tiny $a_{j1}$}}

\put(-1.5,-.15){\scriptsize {\tiny $Q_{1}$}}
\put(-1.5,.85){\scriptsize {\tiny $Q_{2}$}}
\put(-1.5,1.85){\scriptsize {\tiny $Q_{3}$}}
\put(-1.5,2.85){\scriptsize {\tiny $Q_{4}$}}
\put(-1.5,5.85){\scriptsize {\tiny $Q_{j-1}$}}
\put(-1.5,6.85){\scriptsize {\tiny $Q_{j}$}}
\put(-1.5,7.85){\scriptsize {\tiny $Q_{1}$}}

\put(4,-.35){\scriptsize $\ldots$}
\put(6,-.35){\scriptsize $\ldots$}

\put(4,8.25){\scriptsize $\ldots$}
\put(6,8.25){\scriptsize $\ldots$}

\put(-.75,-1.25){\scriptsize {\tiny {\bf Figure 4.4.2:} :  $M(i,j,k)$  of DSEM of type $[3^3.4^2 : 3^2.4.3.4]_2$}} 

\end{picture}

\vspace{5.25cm}

\begin{lem}  \label{l3.4.1}
	A DSEM $M$ of type $[3^3.4^2: 3^2.4.3.4]_2$ admits an $M(i,j,k)$ representation iff the following holds: $(i)$ $ j \geq 2 $ and $j$ even, $(ii)$ if $j = 2$ then $i \geq 8$, and  if $j \geq 4$ then $i \geq 4$, also $i= 4m$, where $m \in \mathbb{N}$,  $(iii)$  $ij \geq 16$, $(iv)$ if $j=2$ then  $k \in \{4r: 0 < r < i/4\} $, and if  $j \geq 4 $ then $k \in \{4r: 0\leq r < i/4\}$.
\end{lem}

%

\subsection{DSEMs of type $[3^6: 3^2.4.3.4]$}\label{s3.5}

Let $M$ be a DSEM of type $[3^6: 3^2.4.3.4]$. We consider the following types of paths in $M$ as follows.

A path $P_{1} = P( \ldots u_{i-1},u_{i},u_{i+1}, \ldots)$ in $M$, say of type $E_{1}$, indicated by thick black paths in Figure 4.5.1. The vertices $u_i$'s have the face-sequence $(3^2,4,3,4)$.

A path $P_{2} = P(\ldots,u_{i},v_{j},v_{j+1},v_{j+2}$, $v_{j+3},u_{i+1}, \ldots)$ in $M$, say of type $E_{2}$, indicated by green paths in Figure 4.5.1. The vertices $u_i$'s and $v_j$'s have the face-sequences $(3^6)$ and $(3^2,4,3,4)$ respectively.

An $M(i,j,k)$ representation of $M$ follows by we first cutting $M$ along an $E_{1}$ type cycle and then along an $E_{2}$ type cycle. Now, proceeding similarly, as in Lemma \ref{l3.3.1}, we get the following lemma. 

\newpage
\begin{picture}(0,0)(-10.5,27)
\setlength{\unitlength}{4mm}

\drawline[AHnb=0](0,0)(1,0)(2,.5)(3,0)(4,0)(5,.5)(6,0)(7,0)(8,.5)(8.5,.25)
\drawline[AHnb=0](0,-1)(1,-1)(2,-1.5)(3,-1)(4,-1)(5,-1.5)(6,-1)(7,-1)(8,-1.5)(8.5,-1.25)

\drawline[AHnb=0](-.5,1.5)(.5,1)(1.5,1.5)(2.5,1.5)(3.5,1)(4.5,1.5)(5.5,1.5)(6.5,1)(7.5,1.5)(8.5,1.5)

\drawline[AHnb=0](-.5,2.5)(.5,3)(1.5,2.5)(2.5,2.5)(3.5,3)(4.5,2.5)(5.5,2.5)(6.5,3)(7.5,2.5)(8.5,2.5)

\drawline[AHnb=0](0,4)(1,4)(2,3.5)(3,4)(4,4)(5,3.5)(6,4)(7,4)(8,3.5)(8.5,3.75)

\drawline[AHnb=0](0,5)(1,5)(2,5.5)(3,5)(4,5)(5,5.5)(6,5)(7,5)(8,5.5)(8.5,5.25)

\drawline[AHnb=0](0,4)(0,5)
\drawline[AHnb=0](1,4)(1,5)
\drawline[AHnb=0](2,3.5)(2,5.5)
\drawline[AHnb=0](3,4)(3,5)
\drawline[AHnb=0](4,4)(4,5)
\drawline[AHnb=0](5,3.5)(5,5.5)
\drawline[AHnb=0](6,4)(6,5)
\drawline[AHnb=0](7,4)(7,5)
\drawline[AHnb=0](8,3.5)(8,5.5)

\drawline[AHnb=0](1,4)(3,5)
\drawline[AHnb=0](1,5)(3,4)

\drawline[AHnb=0](4,4)(6,5)
\drawline[AHnb=0](4,5)(6,4)

\drawline[AHnb=0](7,4)(8,4.5)
\drawline[AHnb=0](7,5)(8,4.5)

\drawline[AHnb=0](8,4.5)(8.5,4.75)
\drawline[AHnb=0](8,4.5)(8.5,4.25)

\drawline[AHnb=0](-.5,1.5)(-.5,2.5)
\drawline[AHnb=0](.5,1)(.5,3)
\drawline[AHnb=0](1.5,1.5)(1.5,2.5)
\drawline[AHnb=0](2.5,1.5)(2.5,2.5)

\drawline[AHnb=0](3.5,1)(3.5,3)
\drawline[AHnb=0](4.5,1.5)(4.5,2.5)
\drawline[AHnb=0](5.5,1.5)(5.5,2.5)
\drawline[AHnb=0](6.5,1)(6.5,3)
\drawline[AHnb=0](7.5,1.5)(7.5,2.5)

\drawline[AHnb=0](8.5,1.5)(8.5,2.5)

\drawline[AHnb=0](-.5,1.5)(1.5,2.5)
\drawline[AHnb=0](-.5,2.5)(1.5,1.5)

\drawline[AHnb=0](2.5,1.5)(4.5,2.5)
\drawline[AHnb=0](2.5,2.5)(4.5,1.5)

\drawline[AHnb=0](5.5,1.5)(7.5,2.5)
\drawline[AHnb=0](5.5,2.5)(7.5,1.5)

\drawline[AHnb=0](0,0)(0,-1)
\drawline[AHnb=0](1,0)(1,-1)
\drawline[AHnb=0](2,.5)(2,-1.5)
\drawline[AHnb=0](3,0)(3,-1)
\drawline[AHnb=0](4,0)(4,-1)
\drawline[AHnb=0](5,.5)(5,-1.5)
\drawline[AHnb=0](6,0)(6,-1)
\drawline[AHnb=0](7,0)(7,-1)
\drawline[AHnb=0](8,.5)(8,-1.5)

\drawline[AHnb=0](1,-1)(3,0)
\drawline[AHnb=0](1,0)(3,-1)

\drawline[AHnb=0](4,-1)(6,0)
\drawline[AHnb=0](4,0)(6,-1)

\drawline[AHnb=0](7,0)(8,-.5)
\drawline[AHnb=0](7,-1)(8,-.5)

\drawline[AHnb=0](8.5,-.25)(8,-.5)
\drawline[AHnb=0](8.5,-.75)(8,-.5)

\drawline[AHnb=0](-.5,-.25)(0,0)
\drawline[AHnb=0](-.5,.25)(0,0)

\drawline[AHnb=0](-.5,-.75)(0,-1)

\drawline[AHnb=0](-.5,-1.25)(0,-1)

\drawline[AHnb=0](0,4)(-.5,4.25)
\drawline[AHnb=0](0,5)(-.5,5.25)

\drawline[AHnb=0](0,4)(-.5,3.75)
\drawline[AHnb=0](0,5)(-.5,4.75)

\drawline[AHnb=0](0,4)(.5,3)
\drawline[AHnb=0](1,4)(.5,3)

\drawline[AHnb=0](3,4)(3.5,3)
\drawline[AHnb=0](4,4)(3.5,3)

\drawline[AHnb=0](6,4)(6.5,3)
\drawline[AHnb=0](7,4)(6.5,3)

\drawline[AHnb=0](0,0)(.5,1)
\drawline[AHnb=0](1,0)(.5,1)

\drawline[AHnb=0](3,0)(3.5,1)
\drawline[AHnb=0](4,0)(3.5,1)

\drawline[AHnb=0](6,0)(6.5,1)
\drawline[AHnb=0](7,0)(6.5,1)

\drawline[AHnb=0](1.5,1.5)(2,.5)
\drawline[AHnb=0](2.5,1.5)(2,.5)

\drawline[AHnb=0](4.5,1.5)(5,.5)
\drawline[AHnb=0](5.5,1.5)(5,.5)

\drawline[AHnb=0](7.5,1.5)(8,.5)
\drawline[AHnb=0](8.5,1.5)(8,.5)

\drawline[AHnb=0](1.5,2.5)(2,3.5)
\drawline[AHnb=0](2.5,2.5)(2,3.5)

\drawline[AHnb=0](4.5,2.5)(5,3.5)
\drawline[AHnb=0](5.5,2.5)(5,3.5)

\drawline[AHnb=0](7.5,2.5)(8,3.5)
\drawline[AHnb=0](8.5,2.5)(8,3.5)

\drawline[AHnb=0](0,-1)(.5,-2)
\drawline[AHnb=0](1,-1)(.5,-2)

\drawline[AHnb=0](3,-1)(3.5,-2)
\drawline[AHnb=0](4,-1)(3.5,-2)

\drawline[AHnb=0](6,-1)(6.5,-2)
\drawline[AHnb=0](7,-1)(6.5,-2)

\drawline[AHnb=0](2,-1.5)(1.5,-2.5)
\drawline[AHnb=0](2,-1.5)(2.5,-2.5)

\drawline[AHnb=0](5,-1.5)(4.5,-2.5)
\drawline[AHnb=0](5,-1.5)(5.5,-2.5)

\drawline[AHnb=0](.5,-2)(1.5,-2.5)(2.5,-2.5)(3.5,-2)(4.5,-2.5)(5.5,-2.5)(6.5,-2)(7.5,-2.5)(8,-1.5)

\drawline[AHnb=0](7.5,-2.5)(8.5,-2.5)(8,-1.5)

\drawline[AHnb=0](.5,-2)(-.5,-2.5)
\drawline[AHnb=0](.5,-2)(.5,-2.75)

\drawline[AHnb=0](1.5,-2.5)(1.5,-2.85)
\drawline[AHnb=0](2.5,-2.5)(2.5,-2.85)

\drawline[AHnb=0](3.5,-2)(3.5,-2.85)
\drawline[AHnb=0](4.5,-2.5)(4.5,-2.85)
\drawline[AHnb=0](5.5,-2.5)(5.5,-2.85)

\drawline[AHnb=0](6.5,-2)(6.5,-2.85)
\drawline[AHnb=0](7.5,-2.5)(7.5,-2.85)
\drawline[AHnb=0](8.5,-2.5)(8.5,-2.85)

\drawline[AHnb=0](7.5,-2.5)(7,-2.75)
\drawline[AHnb=0](4.5,-2.5)(4,-2.75)
\drawline[AHnb=0](5.5,-2.5)(6,-2.75)

\drawline[AHnb=0](1.5,-2.5)(1,-2.75)
\drawline[AHnb=0](2.5,-2.5)(3,-2.75)

\drawline[AHnb=0](-.5,-2.5)(-.5,-2.85)

\drawline[AHnb=0](0,5)(.5,6)
\drawline[AHnb=0](1,5)(.5,6)

\drawline[AHnb=0](3,5)(3.5,6)
\drawline[AHnb=0](4,5)(3.5,6)

\drawline[AHnb=0](6,5)(6.5,6)
\drawline[AHnb=0](7,5)(6.5,6)

\drawline[AHnb=0](2,5.5)(2.25,6)
\drawline[AHnb=0](2,5.5)(1.75,6)

\drawline[AHnb=0](5,5.5)(5.25,6)
\drawline[AHnb=0](5,5.5)(4.75,6)

\drawline[AHnb=0](8,5.5)(8.25,6)
\drawline[AHnb=0](8,5.5)(7.75,6)

\drawline[AHnb=0](.5,6)(0,6.25)
\drawline[AHnb=0](.5,6)(1,6.25)

\drawline[AHnb=0](3.5,6)(3,6.25)
\drawline[AHnb=0](3.5,6)(4,6.25)

\drawline[AHnb=0](6.5,6)(6,6.25)
\drawline[AHnb=0](6.5,6)(7,6.25)

\drawline[AHnb=0](-.5,2.5)(-.75,3)
\drawline[AHnb=0](-.5,-2.5)(-.75,-2)

\drawline[AHnb=0](-.5,1.5)(-.75,1)

\drawline[AHnb=0](8.5,2.5)(8.7,2.5)
\drawline[AHnb=0](8.5,1.5)(8.7,1.5)

\drawpolygon[fillcolor=black](-.5,-1.25)(0,-1)(1,-1)(2,-1.5)(3,-1)(4,-1)(5,-1.5)(6,-1)(7,-1)(8,-1.5)(8.5,-1.25)(8.5,-1.15)(8,-1.4)(7,-.9)(6,-.9)(5,-1.4)(4,-.9)(3,-.9)(2,-1.4)(1,-.9)(0,-.9)(-.5,-1.15)

\drawpolygon[fillcolor=black](-.5,3.65)(0,4)(1,4)(2,3.5)(3,4)(4,4)(5,3.5)(6,4)(7,4)(8,3.5)(8.5,3.75)(8,3.6)(7,4.1)(6,4.1)(5,3.6)(4,4.1)(3,4.1)(2,3.6)(1,4.1)(0,4.1)(-.5,3.75)

\drawpolygon[fillcolor=black](-.5,.25)(0,0)(1,0)(2,.5)(3,0)(4,0)(5,.5)(6,0)(7,0)(8,.5)(8.5,.25)(8.5,.35)(8,.6)(7,.1)(6,.1)(5,.6)(4,.1)(3,.1)(2,.6)(1,.1)(0,.1)(-.5,.35)

\drawpolygon[fillcolor=black](-.5,5.25)(0,5)(1,5)(2,5.5)(3,5)(4,5)(5,5.5)(6,5)(7,5)(8,5.5)(8.5,5.25)(8.5,5.35)(8,5.6)(7,5.1)(6,5.1)(5,5.6)(4,5.1)(3,5.1)(2,5.6)(1,5.1)(0,5.1)(-.5,5.35)

\drawpolygon[fillcolor=black](-.5,1.5)(.5,1)(1.5,1.5)(2.5,1.5)(3.5,1)(4.5,1.5)(5.5,1.5)(6.5,1)(7.5,1.5)(8.6,1.5)(8.6,1.6)(7.5,1.6)(6.5,1.1)(5.5,1.6)(4.5,1.6)(3.5,1.1)(2.5,1.6)(1.5,1.6)(.5,1.1)(-.5,1.6)

\drawpolygon[fillcolor=black](-.5,2.5)(.5,3)(1.5,2.5)(2.5,2.5)(3.5,3)(4.5,2.5)(5.5,2.5)(6.5,3)(7.5,2.5)(8.6,2.5)(8.6,2.6)(7.5,2.6)(6.5,3.1)(5.5,2.6)(4.5,2.6)(3.5,3.1)(2.5,2.6)(1.5,2.6)(.5,3.1)(-.5,2.6)

\drawpolygon[fillcolor=black](-.5,-2.5)(.5,-2)(1.5,-2.5)(2.5,-2.5)(3.5,-2)(4.5,-2.5)(5.5,-2.5)(6.5,-2)(7.5,-2.5)(8.6,-2.5)(8.6,-2.4)(7.5,-2.4)(6.5,-1.9)(5.5,-2.4)(4.5,-2.4)(3.5,-1.9)(2.5,-2.4)(1.5,-2.4)(.5,-1.9)(-.5,-2.4)

\drawpolygon[fillcolor=green](.5,-3)(.5,-2)(1,-1)(1,0)(.5,1)(.5,3)(1,4)(1,5)(.5,6)(.6,6)(1.1,5)(1.1,4)(.6,3)(.6,1)(1.1,0)(1.1,-1)(.6,-2)(.6,-3)

\drawpolygon[fillcolor=green](3.5,-3)(3.5,-2)(4,-1)(4,0)(3.5,1)(3.5,3)(4,4)(4,5)(3.5,6)(3.6,6)(4.1,5)(4.1,4)(3.6,3)(3.6,1)(4.1,0)(4.1,-1)(3.6,-2)(3.6,-3)

\drawpolygon[fillcolor=green](6.5,-3)(6.5,-2)(7,-1)(7,0)(6.5,1)(6.5,3)(7,4)(7,5)(6.5,6)(6.6,6)(7.1,5)(7.1,4)(6.6,3)(6.6,1)(7.1,0)(7.1,-1)(6.6,-2)(6.6,-3)

\drawpolygon[fillcolor=green](.5,-3)(.5,-2)(0,-1)(0,0)(.5,1)(.5,3)(0,4)(0,5)(.5,6)(.4,6)(-.1,5)(-.1,4)(.4,3)(.4,1)(-.1,0)(-.1,-1)(.4,-2)(.4,-3)

\drawpolygon[fillcolor=green](3.5,-3)(3.5,-2)(3,-1)(3,0)(3.5,1)(3.5,3)(3,4)(3,5)(3.5,6)(3.4,6)(2.9,5)(2.9,4)(3.4,3)(3.4,1)(2.9,0)(2.9,-1)(3.4,-2)(3.4,-3)

\drawpolygon[fillcolor=green](6.5,-3)(6.5,-2)(6,-1)(6,0)(6.5,1)(6.5,3)(6,4)(6,5)(6.5,6)(6.4,6)(5.9,5)(5.9,4)(6.4,3)(6.4,1)(5.9,0)(5.9,-1)(6.4,-2)(6.4,-3)

\drawpolygon[fillcolor=green](2.5,-2.85)(2.5,-2.5)(2,-1.5)(2,.5)(2.5,1.5)(2.5,2.5)(2,3.5)(2,5.5)(2.25,6)(2.35,6)(2.1,5.5)(2.1,3.5)(2.6,2.5)(2.6,1.5)(2.1,.5)(2.1,-1.5)(2.6,-2.5)(2.6,-2.85)

\drawpolygon[fillcolor=green](5.5,-2.85)(5.5,-2.5)(5,-1.5)(5,.5)(5.5,1.5)(5.5,2.5)(5,3.5)(5,5.5)(5.25,6)(5.35,6)(5.1,5.5)(5.1,3.5)(5.6,2.5)(5.6,1.5)(5.1,.5)(5.1,-1.5)(5.6,-2.5)(5.6,-2.85)

\drawpolygon[fillcolor=green](1.5,-2.85)(1.5,-2.5)(2,-1.5)(2,.5)(1.5,1.5)(1.5,2.5)(2,3.5)(2,5.5)(1.75,6)(1.65,6)(1.9,5.5)(1.9,3.5)(1.4,2.5)(1.4,1.5)(1.9,.5)(1.9,-1.5)(1.4,-2.5)(1.4,-2.85)

\drawpolygon[fillcolor=green](4.5,-2.85)(4.5,-2.5)(5,-1.5)(5,.5)(4.5,1.5)(4.5,2.5)(5,3.5)(5,5.5)(4.75,6)(4.65,6)(4.9,5.5)(4.9,3.5)(4.4,2.5)(4.4,1.5)(4.9,.5)(4.9,-1.5)(4.4,-2.5)(4.4,-2.85)

\drawpolygon[fillcolor=green](7.5,-2.85)(7.5,-2.5)(8,-1.5)(8,.5)(7.5,1.5)(7.5,2.5)(8,3.5)(8,5.5)(7.75,6)(7.65,6)(7.9,5.5)(7.9,3.5)(7.4,2.5)(7.4,1.5)(7.9,.5)(7.9,-1.5)(7.4,-2.5)(7.4,-2.85)

\drawpolygon[fillcolor=green](8.5,-2.85)(8.5,-2.5)(8,-1.5)(8,.5)(8.5,1.5)(8.5,2.5)(8,3.5)(8,5.5)(8.25,6)(8.35,6)(8.1,5.5)(8.1,3.5)(8.6,2.5)(8.6,1.5)(8.1,.5)(8.1,-1.5)(8.6,-2.5)(8.6,-2.85)

\put(9,.25) {\scriptsize {\tiny $E_1$}}
\put(1.6,6.5) {\scriptsize {\tiny $E_2$}}

\put(-1,-4){\scriptsize {\tiny {\bf Figure 4.5.1:} Paths of types $E_1$ and $E_2$}}

\end{picture}

\begin{picture}(0,0)(-75,33)
\setlength{\unitlength}{5.8mm}

\drawpolygon(0,0)(12,0)(12.5,1)(12.5,2)(12,3)(12,4)(12.5,5)(12.5,6)(.5,6)(.5,5)(0,4)(0,3)(.5,2)(.5,1)

\drawline[AHnb=0](12,3)(11.5,2)
\drawline[AHnb=0](11,3)(10.5,2)

\drawline[AHnb=0](11.5,5)(12,4)
\drawline[AHnb=0](10.5,5)(11,4)
\drawline[AHnb=0](11.5,1)(12,0)
\drawline[AHnb=0](10.5,1)(11,0)

\drawline[AHnb=0](10,3)(10.5,2)
\drawline[AHnb=0](9,3)(9.5,2)

\drawline[AHnb=0](9.5,5)(9,4)
\drawline[AHnb=0](10.5,5)(10,4)
\drawline[AHnb=0](9.5,1)(9,0)
\drawline[AHnb=0](10.5,1)(10,0)


\drawpolygon(1,0)(2,0)(1.5,1)
\drawpolygon(1.5,2)(2,3)(1,3)
\drawpolygon(1,4)(2,4)(1.5,5)

\drawpolygon(2.5,1)(3,0)(3.5,1)
\drawpolygon(2.5,2)(3,3)(3.5,2)
\drawpolygon(2.5,5)(3,4)(3.5,5)

\drawpolygon(5.5,1)(6.5,1)(6,0)
\drawpolygon(7,0)(8,0)(7.5,1)
\drawpolygon(5.5,5)(6.5,5)(6,4)

\drawpolygon(5.5,2)(6.5,2)(6,3)
\drawpolygon(7,3)(8,3)(7.5,2)
\drawpolygon(7,4)(8,4)(7.5,5)

\drawpolygon(4,0)(4.5,1)(5,0)
\drawpolygon(4,3)(4.5,2)(5,3)
\drawpolygon(4,4)(4.5,5)(5,4)

\drawline[AHnb=0](9,0)(8.5,1)
\drawline[AHnb=0](8.5,2)(9,3)

\drawpolygon[fillcolor=black](0,0)(11,0)(10.5,1)(11.5,1)(10.5,1.5)(9.5,1)(8.5,1)(7.5,1.5)(7.5,1)(8.5,1)(5.5,1)(4.5,1.5)(4.5,1)(3.5,1)(1.5,1)(1.5,2)(10.5,2)(11,3)(10,3)(9,3.5)(9,3)(7,3)(6,3.5)(6,3)(4,3)(3,3.5)(3,3)(1,3)(1,4)(11,4)(10.5,5)(11.5,5)(10.5,5.5)(9.5,5)(8.5,5)(7.5,5.5)(7.5,5)(5.5,5)(4.5,5.5)(4.5,5)(2.5,5)(1.5,5.5)(1.5,5)(.5,5)(.5,4.9)(1.6,4.9)(1.6,5.4)(2.5,4.9)(4.6,4.9)(4.6,5.35)(5.5,4.9)(7.6,4.9)(7.6,5.35)(8.4,4.9)(9.5,4.9)(10.5,5.4)(11.2,5.1)(10.35,5.1)(10.85,4.1)(.9,4.1)(.9,2.9)(3.1,2.9)(3.1,3.35)(3.9,2.9)(6.1,2.9)(6.1,3.35)(6.9,2.9)(9.1,2.9)(9.1,3.4)(9.9,2.9)(10.9,2.9)(10.4,2.1)(1.4,2.1)(1.4,.9)(4.6,.9)(4.6,1.35)(5.4,.9)(7.6,.9)(7.6,1.4)(8.4,.9)(9.6,.9)(10.5,1.35)(11,1.1)(10.35,1.1)(10.8,.1)(0,.1)

\drawpolygon[fillcolor=black](12,0)(12.5,1)(12.5,2)(11.5,2)(12,3)(12,4)(12.5,5)(12.6,5)(12.1,4)(12.1,2.9)(11.7,2.1)(12.6,2.1)(12.6,.9)(12.1,0)



\drawline[AHnb=0](.5,1)(12.5,1)
\drawline[AHnb=0](.5,2)(12.5,2)
\drawline[AHnb=0](0,3)(12,3)
\drawline[AHnb=0](0,4)(12,4)
\drawline[AHnb=0](.5,5)(12.5,5)


\drawline[AHnb=0](1.5,1)(1.5,2)
\drawline[AHnb=0](2.5,1)(2.5,2)
\drawline[AHnb=0](3.5,1)(3.5,2)
\drawline[AHnb=0](4.5,1)(4.5,2)
\drawline[AHnb=0](5.5,1)(5.5,2)
\drawline[AHnb=0](6.5,1)(6.5,2)
\drawline[AHnb=0](7.5,1)(7.5,2)
\drawline[AHnb=0](8.5,1)(8.5,2)

\drawline[AHnb=0](9.5,1)(9.5,2)
\drawline[AHnb=0](10.5,1)(10.5,2)
\drawline[AHnb=0](11.5,1)(11.5,2)

\drawline[AHnb=0](9,3)(9,4)
\drawline[AHnb=0](10,3)(10,4)
\drawline[AHnb=0](11,3)(11,4)

\drawline[AHnb=0](8.5,5)(8.5,6)
\drawline[AHnb=0](9.5,5)(9.5,6)
\drawline[AHnb=0](10.5,5)(10.5,6)
\drawline[AHnb=0](11.5,5)(11.5,6)

\drawline[AHnb=0](1,3)(1,4)
\drawline[AHnb=0](2,3)(2,4)
\drawline[AHnb=0](3,3)(3,4)
\drawline[AHnb=0](4,3)(4,4)
\drawline[AHnb=0](5,3)(5,4)
\drawline[AHnb=0](6,3)(6,4)
\drawline[AHnb=0](7,3)(7,4)
\drawline[AHnb=0](8,3)(8,4)

\drawline[AHnb=0](8.5,5)(8.5,6)
\drawline[AHnb=0](9.5,5)(9.5,6)
\drawline[AHnb=0](10.5,5)(10.5,6)
\drawline[AHnb=0](11.5,5)(11.5,6)

\drawline[AHnb=0](1.5,5)(1.5,6)
\drawline[AHnb=0](2.5,5)(2.5,6)
\drawline[AHnb=0](3.5,5)(3.5,6)
\drawline[AHnb=0](4.5,5)(4.5,6)
\drawline[AHnb=0](5.5,5)(5.5,6)
\drawline[AHnb=0](6.5,5)(6.5,6)
\drawline[AHnb=0](7.5,5)(7.5,6)

\drawline[AHnb=0](9,4)(8.5,5)(8.5,6)


\drawline[AHnb=0](.5,1)(2.5,2)
\drawline[AHnb=0](2.5,1)(.5,2)

\drawline[AHnb=0](3.5,1)(5.5,2)
\drawline[AHnb=0](3.5,2)(5.5,1)

\drawline[AHnb=0](6.5,1)(8.5,2)
\drawline[AHnb=0](8.5,1)(6.5,2)

\drawline[AHnb=0](9.5,1)(11.5,2)
\drawline[AHnb=0](11.5,1)(9.5,2)
\drawline[AHnb=0](9.5,5)(11.5,6)
\drawline[AHnb=0](11.5,5)(9.5,6)

\drawline[AHnb=0](8,3)(10,4)
\drawline[AHnb=0](10,3)(8,4)
\drawline[AHnb=0](5,3)(7,4)
\drawline[AHnb=0](7,3)(5,4)

\drawline[AHnb=0](0,3.5)(1,4)
\drawline[AHnb=0](0,3.5)(1,3)

\drawline[AHnb=0](12,3.5)(11,4)
\drawline[AHnb=0](12,3.5)(11,3)

\drawline[AHnb=0](3,3.5)(2,4)
\drawline[AHnb=0](3,3.5)(2,3)

\drawline[AHnb=0](3,3.5)(4,4)
\drawline[AHnb=0](3,3.5)(4,3)

\drawline[AHnb=0](.5,5)(2.5,6)
\drawline[AHnb=0](2.5,5)(.5,6)

\drawline[AHnb=0](3.5,5)(5.5,6)
\drawline[AHnb=0](3.5,6)(5.5,5)

\drawline[AHnb=0](6.5,5)(8.5,6)
\drawline[AHnb=0](8.5,5)(6.5,6)

\drawline[AHnb=0](9,3.5)(8,4)
\drawline[AHnb=0](9,3.5)(8,3)

\put(-.2,-.4){\scriptsize {\tiny $a_{11}$}}
\put(.8,-.4){\scriptsize {\tiny $a_{12}$}}
\put(1.9,-.4){\scriptsize {\tiny $a_{13}$}}
\put(10.8,-.4){\scriptsize {\tiny $a_{1i}$}}
\put(12,-.4){\scriptsize {\tiny $a_{11}$}}

\put(-.3,.95){\scriptsize {\tiny $a_{21}$}}
\put(12.75,1){\scriptsize {\tiny $a_{21}$}}

\put(-.4,2){\scriptsize {\tiny $a_{31}$}}
\put(12.7,2){\scriptsize {\tiny $a_{31}$}}

\put(-.5,5){\scriptsize {\tiny $a_{j1}$}}
\put(12.8,5){\scriptsize {\tiny $a_{j1}$}}

\put(-.5,6.4){\scriptsize {\tiny $a_{1(k+1)}$}}
\put(1,6.4){\scriptsize {\tiny $a_{1(k+2)}$}}
\put(11,6.4){\scriptsize {\tiny $a_{1k}$}}
\put(12,6.4){\scriptsize {\tiny $a_{1(k+1)}$}}

\put(-1.4,0){\scriptsize {\tiny $Q_1$}}
\put(-1.4,1){\scriptsize {\tiny $Q_2$}}
\put(-1.4,5){\scriptsize {\tiny $Q_j$}}
\put(-1.4,6){\scriptsize {\tiny $Q_1$}}

\put(-1.2,3.15){\scriptsize $\vdots$}

\put(4.75,-.5){\scriptsize $\ldots$}
\put(7.25,-.5){\scriptsize $\ldots$}
\put(9.25,-.5){\scriptsize $\ldots$}

\put(4,6.5){\scriptsize $\ldots$}
\put(7,6.5){\scriptsize $\ldots$}
\put(9,6.5){\scriptsize $\ldots$}

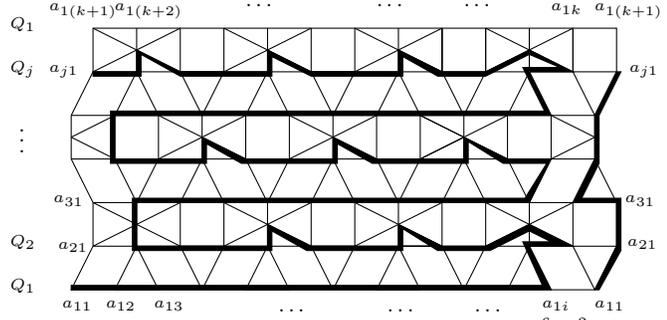
\put(0,-1){\scriptsize {\tiny {\bf Figure 4.5.2:} $M(i,j=4m+2,k)$ of DSEM of type $[3^6 : 3^2.4.3.4]$}} 
\end{picture}

\begin{picture}(0,0)(-40,85)
\setlength{\unitlength}{5.35mm}

\drawpolygon(0,0)(12,0)(12.5,1)(12.5,2)(12,3)(12,4)(12.5,5)(12.5,6)(12,7)(12,8)(0,8)(0,7)(.5,6)(.5,5)(0,4)(0,3)(.5,2)(.5,1)


\drawpolygon(1,0)(2,0)(1.5,1)
\drawpolygon(1.5,2)(2,3)(1,3)
\drawpolygon(1,4)(2,4)(1.5,5)
\drawpolygon(1.5,6)(2,7)(1,7)

\drawpolygon(2.5,1)(3,0)(3.5,1)
\drawpolygon(2.5,2)(3,3)(3.5,2)
\drawpolygon(2.5,5)(3,4)(3.5,5)
\drawpolygon(2.5,6)(3,7)(3.5,6)

\drawpolygon(5.5,1)(6.5,1)(6,0)
\drawpolygon(7,0)(8,0)(7.5,1)
\drawpolygon(5.5,5)(6.5,5)(6,4)
\drawpolygon(5.5,6)(6.5,6)(6,7)

\drawpolygon(5.5,2)(6.5,2)(6,3)
\drawpolygon(7,3)(8,3)(7.5,2)
\drawpolygon(7,4)(8,4)(7.5,5)
\drawpolygon(7,7)(8,7)(7.5,6)

\drawpolygon(4,0)(4.5,1)(5,0)
\drawpolygon(4,3)(4.5,2)(5,3)
\drawpolygon(4,4)(4.5,5)(5,4)
\drawpolygon(4,7)(4.5,6)(5,7)

\drawline[AHnb=0](9,0)(8.5,1)
\drawline[AHnb=0](8.5,2)(9,3)
\drawline[AHnb=0](8.5,6)(9,7)

\drawline[AHnb=0](12,7)(11.5,6)
\drawline[AHnb=0](11,7)(10.5,6)
\drawline[AHnb=0](12,3)(11.5,2)
\drawline[AHnb=0](11,3)(10.5,2)

\drawline[AHnb=0](11.5,5)(12,4)
\drawline[AHnb=0](10.5,5)(11,4)
\drawline[AHnb=0](11.5,1)(12,0)
\drawline[AHnb=0](10.5,1)(11,0)
\drawline[AHnb=0](10,7)(9.5,6)
\drawline[AHnb=0](9,7)(9.5,6)

\drawline[AHnb=0](10,3)(10.5,2)
\drawline[AHnb=0](9,3)(9.5,2)

\drawline[AHnb=0](9.5,5)(9,4)
\drawline[AHnb=0](10.5,5)(10,4)
\drawline[AHnb=0](9.5,1)(9,0)
\drawline[AHnb=0](10.5,1)(10,0)

\drawline[AHnb=0](6,7)(5.5,6)
\drawline[AHnb=0](5,7)(4.5,6)
\drawline[AHnb=0](9,7)(8.5,6)
\drawline[AHnb=0](8,7)(7.5,6)


\put(-.2,-.4){\scriptsize {\tiny $a_{11}$}}
\put(.8,-.4){\scriptsize {\tiny $a_{12}$}}
\put(1.9,-.4){\scriptsize {\tiny $a_{13}$}}
\put(10.8,-.4){\scriptsize {\tiny $a_{1i}$}}
\put(12,-.4){\scriptsize {\tiny $a_{11}$}}

\put(-.3,.95){\scriptsize {\tiny $a_{21}$}}
\put(12.75,1){\scriptsize {\tiny $a_{21}$}}

\put(-.4,2){\scriptsize {\tiny $a_{31}$}}
\put(12.7,2){\scriptsize {\tiny $a_{31}$}}



\put(-.7,7){\scriptsize {\tiny $a_{j1}$}}
\put(12.2,7){\scriptsize {\tiny $a_{j1}$}}

\put(-.95,8.4){\scriptsize {\tiny $a_{1(k+1)}$}}
\put(.7,8.4){\scriptsize {\tiny $a_{1(k+2)}$}}
\put(10.65,8.4){\scriptsize {\tiny $a_{1k}$}}
\put(11.45,8.4){\scriptsize {\tiny $a_{1(k+1)}$}}


\drawline[AHnb=0](.5,1)(12.5,1)
\drawline[AHnb=0](.5,2)(12.5,2)
\drawline[AHnb=0](0,3)(12,3)
\drawline[AHnb=0](0,4)(12,4)
\drawline[AHnb=0](.5,5)(12.5,5)
\drawline[AHnb=0](.5,6)(12.5,6)
\drawline[AHnb=0](0,7)(12,7)

\drawline[AHnb=0](1.5,1)(1.5,2)
\drawline[AHnb=0](2.5,1)(2.5,2)
\drawline[AHnb=0](3.5,1)(3.5,2)
\drawline[AHnb=0](4.5,1)(4.5,2)
\drawline[AHnb=0](5.5,1)(5.5,2)
\drawline[AHnb=0](6.5,1)(6.5,2)
\drawline[AHnb=0](7.5,1)(7.5,2)
\drawline[AHnb=0](8.5,1)(8.5,2)

\drawline[AHnb=0](1,3)(1,4)
\drawline[AHnb=0](2,3)(2,4)
\drawline[AHnb=0](3,3)(3,4)
\drawline[AHnb=0](4,3)(4,4)
\drawline[AHnb=0](5,3)(5,4)
\drawline[AHnb=0](6,3)(6,4)
\drawline[AHnb=0](7,3)(7,4)
\drawline[AHnb=0](8,3)(8,4)

\drawline[AHnb=0](1.5,5)(1.5,6)
\drawline[AHnb=0](2.5,5)(2.5,6)
\drawline[AHnb=0](3.5,5)(3.5,6)
\drawline[AHnb=0](4.5,5)(4.5,6)
\drawline[AHnb=0](5.5,5)(5.5,6)
\drawline[AHnb=0](6.5,5)(6.5,6)
\drawline[AHnb=0](7.5,5)(7.5,6)
\drawline[AHnb=0](9,4)(8.5,5)(8.5,6)

\drawline[AHnb=0](1,7)(1,8)
\drawline[AHnb=0](2,7)(2,8)
\drawline[AHnb=0](3,7)(3,8)
\drawline[AHnb=0](4,7)(4,8)
\drawline[AHnb=0](5,7)(5,8)
\drawline[AHnb=0](6,7)(6,8)
\drawline[AHnb=0](7,7)(7,8)
\drawline[AHnb=0](8,7)(8,8)

\drawline[AHnb=0](9.5,1)(9.5,2)
\drawline[AHnb=0](10.5,1)(10.5,2)
\drawline[AHnb=0](11.5,1)(11.5,2)

\drawline[AHnb=0](9,3)(9,4)
\drawline[AHnb=0](10,3)(10,4)
\drawline[AHnb=0](11,3)(11,4)

\drawline[AHnb=0](9,7)(9,8)
\drawline[AHnb=0](10,7)(10,8)
\drawline[AHnb=0](11,7)(11,8)

\drawline[AHnb=0](8.5,5)(8.5,6)
\drawline[AHnb=0](9.5,5)(9.5,6)
\drawline[AHnb=0](10.5,5)(10.5,6)
\drawline[AHnb=0](11.5,5)(11.5,6)


\drawline[AHnb=0](12,3.5)(11,4)
\drawline[AHnb=0](12,3.5)(11,3)

\drawline[AHnb=0](12,7.5)(11,8)
\drawline[AHnb=0](12,7.5)(11,7)

\drawline[AHnb=0](9.5,1)(11.5,2)
\drawline[AHnb=0](11.5,1)(9.5,2)
\drawline[AHnb=0](9.5,5)(11.5,6)
\drawline[AHnb=0](11.5,5)(9.5,6)

\drawline[AHnb=0](8,7)(10,8)
\drawline[AHnb=0](10,7)(8,8)
\drawline[AHnb=0](8,3)(10,4)
\drawline[AHnb=0](10,3)(8,4)

\drawline[AHnb=0](.5,1)(2.5,2)
\drawline[AHnb=0](2.5,1)(.5,2)

\drawline[AHnb=0](3.5,1)(5.5,2)
\drawline[AHnb=0](3.5,2)(5.5,1)

\drawline[AHnb=0](6.5,1)(8.5,2)
\drawline[AHnb=0](8.5,1)(6.5,2)

\drawline[AHnb=0](0,3.5)(1,4)
\drawline[AHnb=0](0,3.5)(1,3)

\drawline[AHnb=0](5,3)(7,4)
\drawline[AHnb=0](7,3)(5,4)

\drawline[AHnb=0](3,3.5)(2,4)
\drawline[AHnb=0](3,3.5)(2,3)

\drawline[AHnb=0](3,3.5)(4,4)
\drawline[AHnb=0](3,3.5)(4,3)

\drawline[AHnb=0](.5,5)(2.5,6)
\drawline[AHnb=0](2.5,5)(.5,6)

\drawline[AHnb=0](3.5,5)(5.5,6)
\drawline[AHnb=0](3.5,6)(5.5,5)

\drawline[AHnb=0](6.5,5)(8.5,6)
\drawline[AHnb=0](8.5,5)(6.5,6)

\drawline[AHnb=0](9,3.5)(8,4)
\drawline[AHnb=0](9,3.5)(8,3)

\drawline[AHnb=0](0,7.5)(1,8)
\drawline[AHnb=0](0,7.5)(1,7)

\drawline[AHnb=0](3,7.5)(2,8)
\drawline[AHnb=0](3,7.5)(2,7)

\drawline[AHnb=0](3,7.5)(4,8)
\drawline[AHnb=0](3,7.5)(4,7)

\drawline[AHnb=0](5,7)(7,8)
\drawline[AHnb=0](5,8)(7,7)

\drawline[AHnb=0](9,7.5)(8,8)
\drawline[AHnb=0](9,7.5)(8,7)

\drawpolygon[fillcolor=black](0,0)(11,0)(10.5,1)(11.5,1)(10.5,1.5)(9.5,1)(8.5,1)(7.5,1.5)(7.5,1)(8.5,1)(5.5,1)(4.5,1.5)(4.5,1)(3.5,1)(1.5,1)(1.5,2)(10.5,2)(11,3)(10,3)(9,3.5)(9,3)(7,3)(6,3.5)(6,3)(4,3)(3,3.5)(3,3)(1,3)(1,4)(11,4)(10.5,5)(11.5,5)(10.5,5.5)(9.5,5)(8.5,5)(7.5,5.5)(7.5,5)(5.5,5)(4.5,5.5)(4.5,5)(1.5,5)(1.5,6)(10.5,6)(11,7)(10,7)(9,7.5)(9,7)(8,7)(7,7)(6,7.5)(6,7)(4,7)(3,7.5)(3,7)(1,7)(0,7.5)(0,7.4)(.95,6.9)(3.1,6.9)(3.1,7.35)(3.9,6.9)(6.1,6.9)(6.1,7.35)(6.9,6.9)(9.1,6.9)(9.1,7.4)(9.9,6.9)(10.9,6.9)(10.4,6.1)(1.4,6.1)(1.4,4.9)(2.5,4.9)(4.6,4.9)(4.6,5.35)(5.5,4.9)(7.6,4.9)(7.6,5.35)(8.4,4.9)(9.5,4.9)(10.5,5.4)(11.2,5.1)(10.35,5.1)(10.85,4.1)(.9,4.1)(.9,2.9)(3.1,2.9)(3.1,3.35)(3.9,2.9)(6.1,2.9)(6.1,3.35)(6.9,2.9)(9.1,2.9)(9.1,3.4)(9.9,2.9)(10.9,2.9)(10.4,2.1)(1.4,2.1)(1.4,.9)(4.6,.9)(4.6,1.35)(5.4,.9)(7.6,.9)(7.6,1.4)(8.4,.9)(9.6,.9)(10.5,1.35)(11,1.1)(10.35,1.1)(10.8,.1)(0,.1)

\drawpolygon[fillcolor=black](12,0)(12.5,1)(12.5,2)(11.5,2)(12,3)(12,4)(12.5,5)(12.5,6)(11.5,6)(12,7)(12,7.5)(12.1,7.5)(12.1,7)(11.65,6.1)(12.6,6.1)(12.6,5)(12.1,4)(12.1,2.9)(11.7,2.1)(12.6,2.1)(12.6,.9)(12.1,0)

\put(-1.4,0){\scriptsize {\tiny $Q_1$}}
\put(-1.4,1){\scriptsize {\tiny $Q_2$}}
\put(-1.4,2){\scriptsize {\tiny $Q_3$}}
\put(-1.4,6.8){\scriptsize {\tiny $Q_j$}}
\put(-1.4,7.8){\scriptsize {\tiny $Q_1$}}

\put(-1.2,4.5){\scriptsize $\vdots$}

\put(5,-.5){\scriptsize $\ldots$}
\put(6.5,-.5){\scriptsize $\ldots$}
\put(8,-.5){\scriptsize $\ldots$}

\put(3,8.5){\scriptsize $\ldots$}
\put(6,8.5){\scriptsize $\ldots$}
\put(9.25,8.5){\scriptsize $\ldots$}

\put(-1,-1.2){\scriptsize {\tiny {\bf Figure 4.5.3:} $M(i,j=4m+4,k)$  of DSEM of type $[3^6 : 3^2.4.3.4]$}} 
\end{picture}

\vspace{9.35cm}

\begin{lem}  \label{l3.5.1} 
	A DSEM $M$ of type $[3^6: 3^2.4.3.4]$ admits an $M(i,j,k)$ representation iff the following holds: $(i)$ $ j \geq 2 $ even, $(ii)$ if $j=2$ then  $i \geq 9 $ and if $ j \geq 4$ then $i \geq 6 $, also $i=3m$, $m  \in \mathbb N$, $(iii)$ number of vertices of $M(i,j,k) = 7ij/6 \geq 21$, $(iv)$ if $j=2$ then $k \in \{3r + 2: 0 < r < (i-3)/3\},$ if $j=4m+2$ then $k\in\{3r + 2: 0 \leq r < i/3\} $, and if $j=4m$ then $k \in \{3r: 0\leq r < i/3 \} $, $m \in \mathbb N $.  
\end{lem}

%

\subsection{DSEMs of type $[3.4^2.6:3.6.3.6]_1$ and $[3.4^2.6:3.6.3.6]_2$} \label{s3.6}

Let $M^r$ be a DSEM of type $[3.4^2.6:3.6.3.6]_r$, where $r \in \{1,2\}$. In $M^r$, consider following types of paths as follows. 

A path $P_1 = P(\ldots, y_{i-1}, y_{i}, y_{i+1}, \ldots)$ in $M^r$, say of type $D_{1}$, 
indicated by thick black paths in Figure 4.6.1. The vertices $y_i$'s have the face sequence $(3,4^2,6)$. 

A path $P_{2} = P( \ldots, y_{i},y_{i+1},z_{i},y_{i+2},y_{i+3},z_{i+1}, \ldots)$ in $M^r$, say  of type $D_{2}$, indicated by green paths in Figure 4.6.1. The vertices $y_{i}$'s and $z_{i}$'s have the face-sequences $(3,4^2,6)$ and  $(3,6,3,6)$ respectively.	

To construct an $M^r(i, j, k)$ representation of $M^r$, we cut $M^r$ along a cycle of type $D_1$ and then cut along a cycle of type $D_2$. Without loss of generality, let the starting adjacent face to the base horizontal cycle $D_1$ is a 3-gon.


\begin{picture}(0,0)(-30,23.5)
\setlength{\unitlength}{4mm}


\drawline[AHnb=0](6,-.5)(6,0)(7,2)(7,3)(6,5)(6,6)

\drawline[AHnb=0](0,6)(0,5)(1,3)(1,2)(0,0)

\drawline[AHnb=0](2,-.5)(2,0)(3,2)(3,3)(2,5)(2,6)

\drawline[AHnb=0](4,-.5)(4,0)(5,2)(5,3)(4,5)(4,6)


\drawline[AHnb=0](1,-.5)(1,0)(0,2)(0,3)(1,5)(1,6)

\drawline[AHnb=0](3,-.5)(3,0)(2,2)(2,3)(3,5)(3,6)

\drawline[AHnb=0](5,-.5)(5,0)(4,2)(4,3)(5,5)(5,6)

\drawline[AHnb=0](7,-.5)(7,0)(6,2)(6,3)(7,5)(7,6)

\drawline[AHnb=0](1,-.5)(1,0)(0.5,1)(.25,1.5)

\drawline[AHnb=0](0.25,3.5)(0.5,4)(1,5)(1,6)

\drawline[AHnb=0](0,0)(0,-.5)


\drawline[AHnb=0](-.5,0)(7.5,0)

\drawline[AHnb=0](-.5,2)(7.5,2)

\drawline[AHnb=0](-.5,3)(7.5,3)

\drawline[AHnb=0](-.5,5)(7.5,5)

\drawpolygon[fillcolor=black](-.5,0)(7.5,0)(7.5,0.15)(-.5,0.15)
\drawpolygon[fillcolor=black](-.5,2)(7.5,2)(7.5,2.15)(-.5,2.15)
\drawpolygon[fillcolor=black](-.5,3)(7.5,3)(7.5,3.15)(-.5,3.15)
\drawpolygon[fillcolor=black](-.5,5)(7.5,5)(7.5,5.15)(-.5,5.15)

\drawpolygon[fillcolor=green](4,-.5)(4,0)(5,2)(5,3)(4,5)(4,6)(4.15,6)(4.15,5)(5.15,3)(5.15,2)(4.15,0)(4.15,-.5)

\drawpolygon[fillcolor=green](2,-.5)(2,0)(3,2)(3,3)(2,5)(2,6)(2.15,6)(2.15,5)(3.15,3)(3.15,2)(2.15,0)(2.15,-.5)

\drawpolygon[fillcolor=green](0,-.5)(0,0)(1,2)(1,3)(0,5)(0,6)(0.15,6)(0.15,5)(1.15,3)(1.15,2)(0.15,0)(0.15,-.5)

\drawpolygon[fillcolor=green](6,-.5)(6,0)(7,2)(7,3)(6,5)(6,6)(6.15,6)(6.15,5)(7.15,3)(7.15,2)(6.15,0)(6.15,-.5)

\put(7.7,1.9){\scriptsize {\tiny {\bf } $D_1$ }} 

\put(3.7,6.1){\scriptsize {\tiny {\bf }  $D_2$ }}

\put(-1,-1.8){\scriptsize {\tiny {\bf Figure 4.6.1:} Paths of type $D_1,D_2$ in $[3.4^2.6:3.6.3.6]_1$ and $[3.4^2.6:3.6.3.6]_2$}} 

\end{picture}

\begin{picture}(0,0)(-85,18)
\setlength{\unitlength}{4.3mm}

\drawpolygon(0,0)(6,0)(7,2)(7,3)(8,5)(2,5)(1,3)(1,2)

\drawline[AHnb=0](2,-.5)(2,0)(3,2)(3,3)(4,5)(4,5.5)

\drawline[AHnb=0](4,-.5)(4,0)(5,2)(5,3)(6,5)(6,5.5)

\drawline[AHnb=0](0,-.5)(0,0)

\drawline[AHnb=0](6,-.5)(6,0)

\drawline[AHnb=0](2,5)(2,5.5)

\drawline[AHnb=0](8,5)(8,5.5)

\drawline[AHnb=0](1,-.5)(1,0)(0,2)(0,3)(-.5,4)(0,5)(0,5.5)

\drawline[AHnb=0](3,-.5)(3,0)(2,2)(2,3)(1,5)(1,5.5)

\drawline[AHnb=0](5,-.5)(5,0)(4,2)(4,3)(3,5)(3,5.5)

\drawline[AHnb=0](7,-.5)(7,0)(6.5,1)(6,2)(6,3)(5,5)(5,5.5)

\drawline[AHnb=0](8,-.5)(8,0)(8.5,1)(8,2)(8,3)(7.5,4)(7,5)(7,5.5)


\drawline[AHnb=0](-.5,0)(8.5,0)

\drawline[AHnb=0](-.5,2)(8.5,2)

\drawline[AHnb=0](-.5,3)(8.5,3)

\drawline[AHnb=0](-.5,5)(8.5,5)

\drawpolygon[fillcolor=black](-.5,0)(8.5,0)(8.5,0.15)(-.5,0.15)
\drawpolygon[fillcolor=black](-.5,2)(8.5,2)(8.5,2.15)(-.5,2.15)
\drawpolygon[fillcolor=black](-.5,3)(8.5,3)(8.5,3.15)(-.5,3.15)
\drawpolygon[fillcolor=black](-.5,5)(8.5,5)(8.5,5.15)(-.5,5.15)

\drawpolygon[fillcolor=green](0,-.5)(0,0)(1,2)(1,3)(2,5)(2,5.5)(2.15,5.5)(2.15,5)(1.15,3)(1.15,2)(0.15,0)(.15,-.5)

\drawpolygon[fillcolor=green](2,-.5)(2,0)(3,2)(3,3)(4,5)(4,5.5)(4.15,5.5)(4.15,5)(3.15,3)(3.15,2)(2.15,0)(2.15,-.5)

\drawpolygon[fillcolor=green](4,-.5)(4,0)(5,2)(5,3)(6,5)(6,5.5)(6.15,5.5)(6.15,5)(5.15,3)(5.15,2)(4.15,0)(4.15,-.5)

\drawpolygon[fillcolor=green](6,-.5)(6,0)(7,2)(7,3)(8,5)(8,5.5)(8.15,5.5)(8.15,5)(7.15,3)(7.15,2)(6.15,0)(6.15,-.5)

\put(8.5,1.8){\scriptsize {\tiny {\bf }  $D_1$ }}

\put(5.5,5.8){\scriptsize {\tiny {\bf }  $D_2$ }}


\end{picture}

\newpage
\begin{picture}(0,0)(-6,48)
\setlength{\unitlength}{4mm}

\drawpolygon(0,0)(8,0)(9,2)(9,3)(8,5)(8,6)(9,8)(9,9)(8,11)(8,12)(0,12)(0,11)(1,9)(1,8)(0,6)(0,5)(1,3)(1,2)

\drawline[AHnb=0](2,0)(3,2)(3,3)(2,5)(2,6)(3,8)(3,9)(2,11)(2,12)

\drawline[AHnb=0](4,0)(5,2)(5,3)(4,5)(4,6)(5,8)(5,9)(4,11)(4,12)

\drawline[AHnb=0](6,0)(7,2)(7,3)(6,5)(6,6)(7,8)(7,9)(6,11)(6,12)
\drawline[AHnb=0](3,0)(2,2)(2,3)(3,5)(3,6)(2,8)(2,9)(3,11)(3,12)

\drawline[AHnb=0](5,0)(4,2)(4,3)(5,5)(5,6)(4,8)(4,9)(5,11)(5,12)
\drawline[AHnb=0](7,0)(6,2)(6,3)(7,5)(7,6)(6,8)(6,9)(7,11)(7,12)

\drawline[AHnb=0](6.5,1)(6,2)(6,3)(6.5,4)
\drawline[AHnb=0](8.5,1)(8,2)(8,3)(8.5,4)

\drawline[AHnb=0](1,0)(0.5,1)

\drawline[AHnb=0](0.5,4)(1,5)(1,6)(.5,7)

\drawline[AHnb=0](8.5,7)(8,8)(8,9)(8.5,10)

\drawline[AHnb=0](.5,10)(1,11)(1,12)

\drawline[AHnb=0](1,2)(9,2)

\drawline[AHnb=0](1,3)(9,3)

\drawline[AHnb=0](0,5)(8,5)

\drawline[AHnb=0](0,6)(8,6)

\drawline[AHnb=0](1,8)(9,8)

\drawline[AHnb=0](0,11)(8,11)

\drawline[AHnb=0](1,9)(9,9)

\put(-.6,-.6){\scriptsize {\tiny $a_{11}$}}
\put(.7,-.6){\scriptsize {\tiny $a_{12}$}}
\put(1.9,-.6){\scriptsize {\tiny $a_{13}$}}
\put(6.8,-.6){\scriptsize {\tiny $a_{1i}$}}
\put(7.8,-.6){\scriptsize {\tiny $a_{11}$}}

\put(-.8,.9){\scriptsize {\tiny $x_{11}$}}
\put(8.8,.9){\scriptsize {\tiny $x_{11}$}}

\put(-.15,2.15){\scriptsize {\tiny $a_{21}$}}
\put(9.2,2.15){\scriptsize {\tiny $a_{21}$}}

\put(-.2,3){\scriptsize {\tiny $a_{31}$}}
\put(9.2,3){\scriptsize {\tiny $a_{31}$}}

\put(-.65,4){\scriptsize {\tiny $x_{31}$}}
\put(8.9,4){\scriptsize {\tiny $x_{31}$}}

\put(-1.15,5.15){\scriptsize {\tiny $a_{41}$}}
\put(8.25,5.2){\scriptsize {\tiny $a_{41}$}}

\put(-1.5,11){\scriptsize {\tiny $a_{j1}$}}
\put(8.4,11){\scriptsize {\tiny $a_{j1}$}}

\put(-2,12.5){\scriptsize {\tiny $a_{1(k+1)}$}}
\put(.3,12.5){\scriptsize {\tiny $a_{1(k+2)}$}}
\put(6.25,12.5){\scriptsize {\tiny $a_{1k}$}}
\put(7.5,12.5){\scriptsize {\tiny $a_{1(k+1)}$}}

\put(-1.5,-.2){\scriptsize {\tiny $Q_1$}}
\put(-1.5,2.1){\scriptsize {\tiny $Q_2$}}
\put(-1.5,3.1){\scriptsize {\tiny $Q_{3}$}}
\put(-2,5){\scriptsize {\tiny $Q_4$}}
\put(-2.5,10.5){\scriptsize {\tiny $Q_{j}$}}
\put(-2.5,11.65){\scriptsize {\tiny $Q_{1}$}}

\put(3.25,-.65){\scriptsize $\ldots$}
\put(4.65,-.65){\scriptsize $\ldots$}
\put(3.5,12.6){\scriptsize $\ldots$}

\put(-2.35,8){\scriptsize $\vdots$}

\drawpolygon[fillcolor=black](8,0)(8.5,1)(8,2)(9,2)(9,3)(8,3)(8.5,4)(8,5)(8,6)(8.5,7)(8,8)(9,8)(9.1,7.9)(8.15,7.9)(8.65,7)(8.15,6)(8.15,5)(8.65,4)(8.15,3.1)(9.1,3.1)(9.1,1.9)(8.15,1.9)(8.65,1)(8.1,0)

\drawpolygon[fillcolor=black](0,11)(7,11)(6.5,10)(7,9)(5,9)(4.5,10)(4,9)(3,9)(2.5,10)(2,9)(2,8)(1.9,8)(1.9,9.1)(2.5,10.1)(3.1,9.1)(3.9,9.1)(4.5,10.1)(5.1,9.1)(6.85,9.1)(6.35,10)(6.85,10.9)(0,10.9)

\drawpolygon[fillcolor=black](8,11)(8.5,10)(8,9)(9,9)(9,8)(9.1,7.9)(9.1,9.1)(8.15,9.1)(8.6,10)(8.1,11)

\drawpolygon[fillcolor=black](0,0)(2,0)(2.5,1)(3,0)(4,0)(4.5,1)(5,0)(7,0)(6.5,1)(7,2)(2,2)(2,3)(2.5,4)(3,3)(4,3)(4.5,4)(5,3)(7,3)(6.5,4)(7,5)(1,5)(1,6)(2,6)(2.5,7)(3,6)(4,6)(4.5,7)(5,6)(6,6)(7,6)(6.5,7)(7,8)(2,8)(2,8.1)(7.15,8.1)(6.65,7)(7.15,5.9)(4.9,5.9)(4.5,6.8)(4.1,5.9)(2.9,5.9)(2.5,6.85)(2.1,5.9)(1.1,5.9)(1.1,5.1)(7.2,5.1)(6.65,4)(7.2,2.9)(4.9,2.9)(4.5,3.7)(4.1,2.9)(2.9,2.9)(2.5,3.85)(2.15,3)(2.1,2.1)(7.15,2.1)(6.65,1)(7.2,-.1)(4.9,-.1)(4.5,.8)(4.1,-.1)(2.9,-.1)(2.5,.85)(2.1,-.1)(0,-.1)

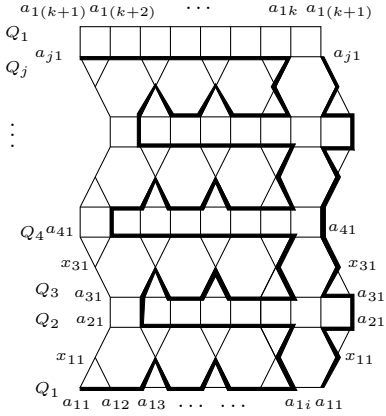
\put(-1,-2){\scriptsize {\tiny {\bf Figure 4.6.2:} $M(i,j=4m+4,k)$ and $M(i,j=4m+2,k)$  of type $[3.4^2.6:3.6.3.6]_1$}}

\end{picture}

\begin{picture}(0,0)(-57,42)
\setlength{\unitlength}{4mm}

\drawpolygon(0,0)(8,0)(9,2)(9,3)(8,5)(8,6)(9,8)(1,8)(0,6)(0,5)(1,3)(1,2)

\drawline[AHnb=0](2,0)(3,2)(3,3)(2,5)(2,6)(3,8)(3,9)

\drawline[AHnb=0](4,0)(5,2)(5,3)(4,5)(4,6)(5,8)(5,9)

\drawline[AHnb=0](6,0)(7,2)(7,3)(6,5)(6,6)(7,8)(7,9)
\drawline[AHnb=0](3,0)(2,2)(2,3)(3,5)(3,6)(2,8)(2,9)
\drawline[AHnb=0](5,0)(4,2)(4,3)(5,5)(5,6)(4,8)(4,9)
\drawline[AHnb=0](7,0)(6,2)(6,3)(7,5)(7,6)(6,8)(6,9)

\drawline[AHnb=0](6.5,1)(6,2)(6,3)(6.5,4)
\drawline[AHnb=0](8.5,1)(8,2)(8,3)(8.5,4)

\drawline[AHnb=0](1,0)(0.5,1)

\drawline[AHnb=0](0.5,4)(1,5)(1,6)(.5,7)

\drawline[AHnb=0](8.5,7)(8,8)(8,9)

\drawline[AHnb=0](1,2)(9,2)

\drawline[AHnb=0](1,3)(9,3)

\drawline[AHnb=0](0,5)(8,5)

\drawline[AHnb=0](0,6)(8,6)

\drawline[AHnb=0](1,8)(9,8)

\drawline[AHnb=0](1,8)(1,9)(9,9)(9,8)

\put(-.4,-.6){\scriptsize {\tiny $a_{11}$}}
\put(.7,-.6){\scriptsize {\tiny $a_{12}$}}
\put(1.9,-.6){\scriptsize {\tiny $a_{13}$}}
\put(6.8,-.6){\scriptsize {\tiny $a_{1i}$}}
\put(7.8,-.6){\scriptsize {\tiny $a_{11}$}}

\put(-.8,.9){\scriptsize {\tiny $x_{11}$}}
\put(8.8,.9){\scriptsize {\tiny $x_{11}$}}

\put(-.15,2.15){\scriptsize {\tiny $a_{21}$}}
\put(9.2,2.15){\scriptsize {\tiny $a_{21}$}}

\put(-.2,3){\scriptsize {\tiny $a_{31}$}}
\put(9.2,3){\scriptsize {\tiny $a_{31}$}}

\put(-.65,4){\scriptsize {\tiny $x_{31}$}}
\put(8.9,4){\scriptsize {\tiny $x_{31}$}}

\put(-1.15,5.15){\scriptsize {\tiny $a_{41}$}}
\put(8.25,5.2){\scriptsize {\tiny $a_{41}$}}

\put(-.5,8){\scriptsize {\tiny $a_{j1}$}}
\put(9.5,8){\scriptsize {\tiny $a_{j1}$}}

\put(-1.5,-.2){\scriptsize {\tiny $Q_1$}}
\put(-1.5,2.1){\scriptsize {\tiny $Q_2$}}
\put(-1.5,3.1){\scriptsize {\tiny $Q_{3}$}}
\put(-2,5){\scriptsize {\tiny $Q_4$}}
\put(-1.5,8.1){\scriptsize {\tiny $Q_{j}$}}
\put(-1.5,9.1){\scriptsize {\tiny $Q_{1}$}}

\put(3.25,-.65){\scriptsize $\ldots$}
\put(4.65,-.65){\scriptsize $\ldots$}
\put(5,9.65){\scriptsize $\ldots$}

\put(-1.35,6.5){\scriptsize $\vdots$}

\put(-.75,9.5){\scriptsize {\tiny $a_{1(k+1)}$}}
\put(1.5,9.5){\scriptsize {\tiny $a_{1(k+2)}$}}

\put(7.65,9.5){\scriptsize {\tiny $a_{1k}$}}
\put(8.85,9.5){\scriptsize {\tiny $a_{1(k+1)}$}}

\drawpolygon[fillcolor=black](8,0)(8.5,1)(8,2)(9,2)(9,3)(8,3)(8.5,4)(8,5)(8,6)(8.5,7)(8,8)(9,8)(9,7.9)(8.15,7.9)(8.65,7)(8.15,6)(8.15,5)(8.65,4)(8.15,3.1)(9.1,3.1)(9.1,1.9)(8.15,1.9)(8.65,1)(8.1,0)

\drawpolygon[fillcolor=black](0,0)(2,0)(2.5,1)(3,0)(4,0)(4.5,1)(5,0)(7,0)(6.5,1)(7,2)(2,2)(2,3)(2.5,4)(3,3)(4,3)(4.5,4)(5,3)(7,3)(6.5,4)(7,5)(1,5)(1,6)(2,6)(2.5,7)(3,6)(4,6)(4.5,7)(5,6)(6,6)(7,6)(6.5,7)(7,8)(1,8)(1,8.1)(7.15,8.1)(6.65,7)(7.15,5.9)(4.9,5.9)(4.5,6.8)(4.1,5.9)(2.9,5.9)(2.5,6.85)(2.1,5.9)(1.1,5.9)(1.1,5.1)(7.2,5.1)(6.65,4)(7.2,2.9)(4.9,2.9)(4.5,3.7)(4.1,2.9)(2.9,2.9)(2.5,3.85)(2.15,3)(2.1,2.1)(7.15,2.1)(6.65,1)(7.2,-.1)(4.9,-.1)(4.5,.8)(4.1,-.1)(2.9,-.1)(2.5,.85)(2.1,-.1)(0,-.1)

\end{picture}

\begin{picture}(0,0)(-107,37.5)
\setlength{\unitlength}{4mm}

\drawpolygon(0,0)(8,0)(9,2)(9,3)(10,5)(10,6)(11,8)(11,9)(3,9)(3,8)(2,6)(2,5)(1,3)(1,2)

\drawline[AHnb=0](2,0)(3,2)(3,3)(4,5)(4,6)(5,8)(5,9)

\drawline[AHnb=0](4,0)(5,2)(5,3)(6,5)(6,6)(7,8)(7,9)

\drawline[AHnb=0](6,0)(7,2)(7,3)(8,5)(8,6)(9,8)(9,9)

\drawline[AHnb=0](3,0)(2,2)(2,3)(1.5,4)

\drawline[AHnb=0](5,0)(4,2)(4,3)(3,5)(3,6)(2.5,7)

\drawline[AHnb=0](7,0)(6,2)(6,3)(5,5)(5,6)(4,8)(4,9)

\drawline[AHnb=0](1,0)(0.5,1)

\drawline[AHnb=0](10.5,7)(10,8)(10,9)

\drawline[AHnb=0](8.5,1)(8,2)(8,3)(7,5)(7,6)(6,8)(6,9)

\drawline[AHnb=0](9.5,4)(9,5)(9,6)(8,8)(8,9)
\drawline[AHnb=0](1,2)(9,2)

\drawline[AHnb=0](1,3)(9,3)

\drawline[AHnb=0](2,5)(10,5)
\drawline[AHnb=0](2,6)(10,6)
\drawline[AHnb=0](3,8)(11,8)

\put(-.6,-.6){\scriptsize {\tiny $a_{11}$}}
\put(.5,-.6){\scriptsize {\tiny $a_{12}$}}
\put(1.7,-.6){\scriptsize {\tiny $a_{13}$}}
\put(6.8,-.6){\scriptsize {\tiny $a_{1i}$}}
\put(7.8,-.6){\scriptsize {\tiny $a_{11}$}}

\put(-.8,.9){\scriptsize {\tiny $x_{11}$}}
\put(8.8,.9){\scriptsize {\tiny $x_{11}$}}

\put(-.15,2.15){\scriptsize {\tiny $a_{21}$}}
\put(9.4,2.15){\scriptsize {\tiny $a_{21}$}}

\put(-.2,3){\scriptsize {\tiny $a_{31}$}}
\put(9.4,3){\scriptsize {\tiny $a_{31}$}}

\put(1.5,8){\scriptsize {\tiny $a_{j1}$}}
\put(11.5,8){\scriptsize {\tiny $a_{j1}$}}

\put(.95,9.5){\scriptsize {\tiny $a_{1(k+1)}$}}
\put(3.35,9.5){\scriptsize {\tiny $a_{1(k+2)}$}}

\put(9.55,9.5){\scriptsize {\tiny $a_{1k}$}}
\put(10.7,9.5){\scriptsize {\tiny $a_{1(k+1)}$}}

\drawpolygon[fillcolor=black](8,0)(8.5,1)(8,2)(9,2)(9,3)(9.5,4)(9,5)(10,5)(10,6)(10.1,6)(10.1,4.9)(9.15,4.9)(9.65,4)(9.15,3)(9.15,1.9)(8.25,1.9)(8.65,1)(8.15,0)

\drawpolygon[fillcolor=black](0,0)(2,0)(2.5,1)(3,0)(4,0)(4.5,1)(5,0)(7,0)(6.5,1)(7,2)(2,2)(2,3)(3,3)(3.5,4)(4,3)(5,3)(5.5,4)(6,3)(8,3)(7.5,4)(8,5)(3,5)(3,6)(3.1,6)(3.1,5.1)(8.2,5.1)(7.65,4)(8.2,2.9)(5.9,2.9)(5.5,3.8)(5.1,2.9)(3.9,2.9)(3.5,3.7)(3.15,2.9)(2.1,2.9)(2.1,2.1)(7.15,2.1)(6.65,1)(7.2,-.1)(4.9,-.1)(4.5,.8)(4.1,-.1)(2.9,-.1)(2.5,.8)(2.1,-.1)(0,-.1)

\drawpolygon[fillcolor=black](3,6)(4,6)(4.5,7)(5,6)(6,6)(6.5,7)(7,6)(9,6)(8.5,7)(9,8)(3,8)(3.1,8.1)(9.15,8.1)(8.65,7)(9.2,5.9)(6.9,5.9)(6.5,6.9)(6.1,5.9)(4.9,5.9)(4.5,6.9)(4.1,5.9)(3,5.9)

\drawpolygon[fillcolor=black](11,8)(10,8)(10.5,7)(10,6)(10.15,6)(10.65,7)(10.15,7.9)(11,7.9)

\put(-2,-.2){\scriptsize {\tiny $Q_1$}}
\put(-1.5,2){\scriptsize {\tiny $Q_2$}}
\put(-1.25,3){\scriptsize {\tiny $Q_{3}$}}
\put(.25,8){\scriptsize {\tiny $Q_{j}$}}
\put(.25,9){\scriptsize {\tiny $Q_{1}$}}

\put(3.8,-.65){\scriptsize $\ldots$}
\put(5.4,-.65){\scriptsize $\ldots$}
\put(7,9.5){\scriptsize $\ldots$}

\put(.1,6.1){\scriptsize $.$}
\put(0,5.8){\scriptsize $.$}
\put(-.1,5.5){\scriptsize $.$}

\put(-1.5,-2.3){\scriptsize {\tiny {\bf Figure 4.6.3:} $M(i,j,k)$ of type $[3.4^2.6 : 3.6.3.6]_2$ }}
\end{picture}

\vspace{5cm}

\begin{lem} \label{l3.6.1}
A DSEM $M^r$ of type $[3.4^2.6:3.6.3.6]_r$, for $r \in \{1,2\}$, admits an $M^r(i,j,k)$-representation iff the following holds: $(i)$ $i \geq 6$ and $i, j$ even, $(ii)$ number of vertices of $M(i,j,k) = 5ij/4 \geq 15 $, $(iii)$  if $j=2$ then $ 4 \leq k \leq i-2 $, and if $j \geq 4$ then $ 0 \leq k \leq i-1$.
\end{lem}

\noindent{\bf Proof.}	Let $M^r$ be a DSEM of type $[3.4^2.6:3.6.3.6]_r$ with $n$ vertices. Its $M^r(i,j,k)$ representation has $j$ disjoint horizontal cycles of $D_{1}$ type, say $Q_0, Q_1,\ldots, Q_{j-1}$, of length $i$. Note that the number of vertices between the horizontal cycles $Q_{2s(mod\,j)}$ and $Q_{(2s+1)(mod\,j)}$ for $0 \leq s \leq j-1$ with the face-sequence $(3,6,3,6)$ is $i/2\cdot j/2$. Therefore, the number of vertices in $M^r$ is $n = ij + ij/4 = 5ij/4$.

If $j=1$, then $M^r(i,1,k)$ has no vertex with face-sequence $(3,4^2,6)$. This is not possible. Thus, $j \geq 2 $. If $j \geq 2 $ and $j$ is not an even integer, then we see that there is no vertex in the base horizontal cycle which has face-sequence $(3,4 ^2,6)$ after identifying the boundaries of $M^r(i, j, k)$. Thus $j$ is even.

For  $i \leq 4$, $M^r$ is not a map. So $i > 4$. If $i > 4 $ and not an even integer, $M^r$ is not a map. So $ i \geq 6$, and $i, j$ are even. Thus $n=5ij/4 \geq 15$.

If $j = 2$ and $k \in \{ r : 0 \leq r \leq i-1 \} \setminus \{0,1,2,3,i-1 \} $, then we get some vertices whose link can not be constructed. So, for $j=2$, we get $4 \leq k \leq i-2$. Similarly, if $j \geq 4$, then $ 0 \leq k \leq i-1$. Thus the proof. \hfill $\Box$


%


\subsection{DSEMs of type $[3^2.6^2:3.6.3.6]$} \label{s3.7}

Let $M$ be a DSEM of type $[3^2.6^2:3.6.3.6]$. We consider the following types of paths in $M$ as follows.

A path $P_{1} = P( \ldots, y_{i-1},y_{i},y_{i+1}, \ldots)$ in $M$, say  of type $F_{1}$ indicated by thick black paths, shown in Figure 4.7.1. The vertices $y_i$'s have the face-sequence $(3^2,6^2)$. 

A path $P_{2} = P( \ldots, y_{i},z_{i},y_{i+1},z_{i+1}$, $\ldots)$ in $M$, say of type $F_{2}$, indicated by green paths, shown Figure 4.7.1. The vertices $y_i$'s and $z_i$'s have the face-sequences $(3^2,6^2)$ and $(3,6,3,6)$ respectively.

\vspace{-.15cm}

\begin{picture}(0,0)(-2,27)
\setlength{\unitlength}{2.75mm}


\drawpolygon(0,0)(8,0)(9,2)(8,4)(9,6)(8,8)(0,8)(1,6)(0,4)(1,2)


\drawline[AHnb=0](.25,-.5)(0,0)

\drawline[AHnb=0](8.25,-.5)(8,0)

\drawline[AHnb=0](.25,8.5)(0,8)

\drawline[AHnb=0](8.25,8.5)(8,8)

\drawline[AHnb=0](2.25,-.5)(2,0)(3,2)(2,4)(3,6)(2,8)(2.25,8.5)


\drawline[AHnb=0](6.25,-.5)(6,0)(7,2)(6,4)(7,6)(6,8)(6.25,8.5)


\drawline[AHnb=0](.75,-.5)(1,0)(0,2)(1,4)(0,6)(1,8)(.75,8.5)

\drawline[AHnb=0](2.75,-.5)(3,0)(2,2)(3,4)(2,6)(3,8)(2.75,8.5)

\drawline[AHnb=0](4.75,-.5)(5,0)(4,2)(5,4)(4,6)(5,8)(4.75,8.5)

\drawline[AHnb=0](6.75,-.5)(7,0)(6,2)(7,4)(6,6)(7,8)(6.75,8.5)

\drawline[AHnb=0](8.75,-.5)(9,0)(8,2)(9,4)(8,6)(9,8)(8.75,8.5)


\drawline[AHnb=0](-.5,0)(9.5,0)
\drawline[AHnb=0](-.5,2)(9.5,2)
\drawline[AHnb=0](-.5,4)(9.5,4)
\drawline[AHnb=0](-.5,6)(9.5,6)
\drawline[AHnb=0](-.5,8)(9.5,8)


\drawpolygon[fillcolor=black](-.5,0)(9.5,0)(9.5,.15)(-.5,.15)

\drawpolygon[fillcolor=black](-.5,2)(9.5,2)(9.5,2.15)(-.5,2.15)

\drawpolygon[fillcolor=black](-.5,4)(9.5,4)(9.5,4.15)(-.5,4.15)

\drawpolygon[fillcolor=black](-.5,6)(9.5,6)(9.5,6.15)(-.5,6.15)

\drawpolygon[fillcolor=black](-.5,8)(9.5,8)(9.5,8.15)(-.5,8.15)


\drawpolygon[fillcolor=green](8.25,-.5)(8,0)(9,2)(8,4)(9,6)(8,8)(8.25,8.5)(8.26,8.5)(8.15,8)(9.15,6)(8.15,4)(9.15,2)(8.15,0)(8.27,-.5)

\drawpolygon[fillcolor=green](6.25,-.5)(6,0)(7,2)(6,4)(7,6)(6,8)(6.25,8.5)(6.26,8.5)(6.15,8)(7.15,6)(6.15,4)(7.15,2)(6.15,0)(6.27,-.5)

\drawpolygon[fillcolor=green](4.25,-.5)(4,0)(5,2)(4,4)(5,6)(4,8)(4.25,8.5)(4.26,8.5)(4.15,8)(5.15,6)(4.15,4)(5.15,2)(4.15,0)(4.27,-.5)

\drawpolygon[fillcolor=green](2.25,-.5)(2,0)(3,2)(2,4)(3,6)(2,8)(2.25,8.5)(2.26,8.5)(2.15,8)(3.15,6)(2.15,4)(3.15,2)(2.15,0)(2.27,-.5)

\drawpolygon[fillcolor=green](.25,-.5)(0,0)(1,2)(0,4)(1,6)(0,8)(.25,8.5)(.26,8.5)(.15,8)(1.15,6)(.15,4)(1.15,2)(.15,0)(.27,-.5)

\put(3.8,8.9){\scriptsize {\tiny {\bf } $F_2$ }}

\put(9.6,1.9){\scriptsize {\tiny {\bf } $F_1$ }}

\put(-2,-2){\scriptsize {\tiny {\bf Figure 4.7.1:} Paths of type $F_1, F_2$ }}

\end{picture}

\begin{picture}(0,0)(0,48)
\setlength{\unitlength}{4.15mm}


\drawpolygon(0,0)(8,0)(9,2)(1,2)


\drawline[AHnb=0](2,0)(3,2)
\drawline[AHnb=0](4,0)(5,2)
\drawline[AHnb=0](6,0)(7,2)
\drawline[AHnb=0](3,0)(2,2)
\drawline[AHnb=0](5,0)(4,2)
\drawline[AHnb=0](7,0)(6,2)

\drawline[AHnb=0](1,0)(.5,1)

\drawline[AHnb=0](8.5,1)(8,2)


\drawline[AHnb=0](1,2)(9,2)


\drawpolygon[fillcolor=black](0,0)(.5,1)(1,0)(2,0)(2.5,1)(3,0)(4,0)(4.5,1)(5,0)(6,0)(6.5,1)(7,0)(8,0)(8,.1)(7.1,.1)(6.5,1.1)(5.9,.1)(5.1,.1)(4.5,1.1)(3.9,.1)(3.1,.1)(2.5,1.1)(1.9,.1)(1.1,.1)(.5,1.1)(-.1,0)

\put(-.4,-.6){\scriptsize {\tiny $a_{11}$}}
\put(.7,-.6){\scriptsize {\tiny $a_{12}$}}
\put(1.7,-.6){\scriptsize {\tiny $a_{13}$}}
\put(6.8,-.6){\scriptsize {\tiny $a_{1i}$}}
\put(7.8,-.6){\scriptsize {\tiny $a_{11}$}}

\put(-.7,1){\scriptsize {\tiny $x_{11}$}}
\put(8.8,.9){\scriptsize {\tiny $x_{11}$}}

\put(-.75,2.5){\scriptsize {\tiny $a_{1(k+1)}$}}
\put(1.55,2.5){\scriptsize {\tiny $a_{1(k+2)}$}}
\put(7.6,2.5){\scriptsize {\tiny $a_{1k}$}}
\put(8.8,2.5){\scriptsize {\tiny $a_{1(k+1)}$}}

\put(4,-.5){\scriptsize $\ldots$}
\put(4.25,2.5){\scriptsize $\ldots$}
\put(6,2.5){\scriptsize $\ldots$}

\put(.2,-1.7){\scriptsize {\tiny {\bf Figure 4.7.2:} $M(i,1,k)$ }}

\end{picture}

\begin{picture}(0,0)(-57,44)
\setlength{\unitlength}{4mm}


\drawpolygon(0,0)(8,0)(9,2)(8,4)(9,6)(8,8)(0,8)(1,6)(0,4)(1,2)

\drawpolygon(0,8)(8,8)(9,10)(8,12)(0,12)(1,10)


\drawline[AHnb=0](2,0)(3,2)(2,4)(3,6)(2,8)

\drawline[AHnb=0](4,0)(5,2)(4,4)(5,6)(4,8)

\drawline[AHnb=0](6,0)(7,2)(6,4)(7,6)(6,8)

\drawline[AHnb=0](3,0)(2,2)(3,4)(2,6)(3,8)

\drawline[AHnb=0](5,0)(4,2)(5,4)(4,6)(5,8)

\drawline[AHnb=0](7,0)(6,2)(7,4)(6,6)(7,8)

\drawline[AHnb=0](2,8)(3,10)(2,12)
\drawline[AHnb=0](4,8)(5,10)(4,12)
\drawline[AHnb=0](6,8)(7,10)(6,12)

\drawline[AHnb=0](3,8)(2,10)(3,12)
\drawline[AHnb=0](5,8)(4,10)(5,12)
\drawline[AHnb=0](7,8)(6,10)(7,12)
\drawline[AHnb=0](8.5,9)(8,10)(8.5,11)

\drawline[AHnb=0](1,8)(.5,9)
\drawline[AHnb=0](.5,11)(1,12)

\drawline[AHnb=0](1,0)(.5,1)
\drawline[AHnb=0](1,4)(.5,5)
\drawline[AHnb=0](1,4)(.5,3)
\drawline[AHnb=0](1,8)(.5,7)

\drawline[AHnb=0](8.5,1)(8,2)
\drawline[AHnb=0](8.5,3)(8,2)

\drawline[AHnb=0](8.5,5)(8,6)
\drawline[AHnb=0](8.5,7)(8,6)



\drawline[AHnb=0](1,2)(9,2)
\drawline[AHnb=0](0,4)(8,4)
\drawline[AHnb=0](1,6)(9,6)
\drawline[AHnb=0](0,8)(8,8)
\drawline[AHnb=0](1,12)(8,12)
\drawline[AHnb=0](1,10)(9,10)


\drawpolygon[fillcolor=black](0,0)(2,0)(2.5,1)(3,0)(4,0)(4.5,1)(5,0)(7,0)(6.5,1)(7,2)(6.5,3)(6,2)(5,2)(4.5,3)(4,2)(2,2)(2.5,3)(2,4)(2.5,5)(3,4)(4,4)(4.5,5)(5,4)(7,4)(6.5,5)(7,6)(6.5,7)(6,6)(5,6)(4.5,7)(4,6)(2,6)(2.5,7)(2,8)(1.9,8)(2.4,7)(1.9,5.9)(4.1,5.9)(4.5,6.9)(4.9,5.9)(6.1,5.9)(6.5,6.9)(6.9,6)(6.4,5)(6.8,4.1)(5.1,4.1)(4.5,5.1)(3.9,4.1)(3.1,4.1)(2.5,5.1)(1.85,4)(2.4,3)(1.85,1.9)(4.1,1.9)(4.5,2.9)(4.9,1.9)(6.1,1.9)(6.5,2.9)(6.9,2)(6.4,1)(6.85,.1)(5.1,.1)(4.5,1.1)(3.9,.1)(3.1,.1)(2.5,1.1)(1.9,.1)(0,.1)

\drawpolygon[fillcolor=black](2,8)(2.5,9)(3,8)(4,8)(4.5,9)(5,8)(7,8)(6.5,9)(7,10)(6.5,11)(6,10)(5,10)(4.5,11)(4,10)(3,10)(2.5,11)(2,10)(1,10)(.5,11)(.4,11)(.9,9.9)(2.1,9.9)(2.5,10.9)(2.9,9.9)(4.1,9.9)(4.5,10.9)(4.9,9.9)(6.1,9.9)(6.5,10.9)(6.9,10)(6.4,9)(6.8,8.1)(5.1,8.1)(4.5,9.1)(3.9,8.1)(3.1,8.1)(2.5,9.1)(1.85,8)

\drawpolygon[fillcolor=black](8,0)(8.5,1)(8,2)(9,2)(8.5,3)(8.45,2.9)(8.85,2.1)(7.85,2.1)(8.4,1)(7.9,0)

\drawpolygon[fillcolor=black](.5,3)(0,4)(1,4)(.5,5)(.55,5.1)(1.15,3.9)(.15,3.9)(.55,3.1)

\drawpolygon[fillcolor=black](8.5,5)(8,6)(9,6)(8.5,7)(8.55,7.1)(9.15,5.9)(8.15,5.9)(8.6,5)

\drawpolygon[fillcolor=black](.5,7)(0,8)(1,8)(0.5,9)(.6,9.1)(1.15,7.9)(.15,7.9)(.55,7.2)

\drawpolygon[fillcolor=black](8.5,9)(8,10)(8.5,11)(8.55,10.9)(8.1,10)(8.55,9.1)

\put(-.4,-.6){\scriptsize {\tiny $a_{11}$}}
\put(.7,-.6){\scriptsize {\tiny $a_{12}$}}
\put(1.7,-.6){\scriptsize {\tiny $a_{13}$}}
\put(6.8,-.6){\scriptsize {\tiny $a_{1i}$}}
\put(7.8,-.6){\scriptsize {\tiny $a_{11}$}}

\put(-.8,.9){\scriptsize {\tiny $x_{11}$}}
\put(8.8,.9){\scriptsize {\tiny $x_{11}$}}

\put(-.35,2){\scriptsize {\tiny $a_{21}$}}
\put(9.2,2){\scriptsize {\tiny $a_{21}$}}

\put(-.6,3){\scriptsize {\tiny $x_{21}$}}
\put(9,3){\scriptsize {\tiny $x_{21}$}}

\put(-.9,4.2){\scriptsize {\tiny $a_{31}$}}
\put(8.4,4.2){\scriptsize {\tiny $a_{31}$}}

\put(-1.05,11){\scriptsize {\tiny $x_{j1}$}}
\put(8.75,11){\scriptsize {\tiny $x_{j1}$}}

\put(-.5,10){\scriptsize {\tiny $a_{j1}$}}
\put(9.25,10){\scriptsize {\tiny $a_{j1}$}}

\put(-1.95,12.5){\scriptsize {\tiny $a_{1(k+1)}$}}
\put(.45,12.5){\scriptsize {\tiny $a_{1(k+2)}$}}
\put(6.25,12.5){\scriptsize {\tiny $a_{1k}$}}
\put(7.5,12.5){\scriptsize {\tiny $a_{1(k+1)}$}}

\put(-2,-.2){\scriptsize {\tiny $Q_1$}}
\put(-2.2,1.8){\scriptsize {\tiny $Q_2$}}
\put(-2.2,3.8){\scriptsize {\tiny $Q_3$}}
\put(-2.2,10){\scriptsize {\tiny $Q_{j}$}}
\put(-2.4,11.8){\scriptsize {\tiny $Q_{1}$}}

\put(4.25,-.5){\scriptsize $\ldots$}
\put(3,12.5){\scriptsize $\ldots$}
\put(4.75,12.5){\scriptsize $\ldots$}
\put(-2.25,6.25){\scriptsize $\vdots$}

\put(.2,-1.5){\scriptsize {\tiny {\bf Figure 4.7.3:} $M(i,j=2m,k)$ }}

\end{picture}

\begin{picture}(0,0)(-112,37)
\setlength{\unitlength}{4.3mm}


\drawpolygon(0,0)(8,0)(9,2)(8,4)(9,6)(8,8)(9,10)(1,10)(0,8)(1,6)(0,4)(1,2)


\drawline[AHnb=0](2,0)(3,2)(2,4)(3,6)(2,8)(3,10)

\drawline[AHnb=0](4,0)(5,2)(4,4)(5,6)(4,8)(5,10)

\drawline[AHnb=0](6,0)(7,2)(6,4)(7,6)(6,8)(7,10)

\drawline[AHnb=0](3,0)(2,2)(3,4)(2,6)(3,8)(2,10)

\drawline[AHnb=0](5,0)(4,2)(5,4)(4,6)(5,8)(4,10)

\drawline[AHnb=0](7,0)(6,2)(7,4)(6,6)(7,8)(6,10)

\drawline[AHnb=0](1,0)(.5,1)
\drawline[AHnb=0](1,4)(.5,5)
\drawline[AHnb=0](1,4)(.5,3)
\drawline[AHnb=0](1,8)(.5,7)
\drawline[AHnb=0](1,8)(.5,9)

\drawline[AHnb=0](8.5,1)(8,2)
\drawline[AHnb=0](8.5,3)(8,2)

\drawline[AHnb=0](8.5,5)(8,6)
\drawline[AHnb=0](8.5,7)(8,6)

\drawline[AHnb=0](8.5,9)(8,10)


\drawline[AHnb=0](1,2)(9,2)
\drawline[AHnb=0](0,4)(8,4)
\drawline[AHnb=0](1,6)(9,6)
\drawline[AHnb=0](0,8)(8,8)


\drawpolygon[fillcolor=black](0,0)(2,0)(2.5,1)(3,0)(4,0)(4.5,1)(5,0)(7,0)(6.5,1)(7,2)(6.5,3)(6,2)(5,2)(4.5,3)(4,2)(2,2)(2.5,3)(2,4)(2.5,5)(3,4)(4,4)(4.5,5)(5,4)(7,4)(6.5,5)(7,6)(6.5,7)(6,6)(5,6)(4.5,7)(4,6)(2,6)(2.5,7)(2,8)(1.9,8)(2.4,7)(1.9,5.9)(4.1,5.9)(4.5,6.9)(4.9,5.9)(6.1,5.9)(6.5,6.9)(6.9,6)(6.4,5)(6.8,4.1)(5.1,4.1)(4.5,5.1)(3.9,4.1)(3.1,4.1)(2.5,5.1)(1.85,4)(2.4,3)(1.85,1.9)(4.1,1.9)(4.5,2.9)(4.9,1.9)(6.1,1.9)(6.5,2.9)(6.9,2)(6.4,1)(6.85,.1)(5.1,.1)(4.5,1.1)(3.9,.1)(3.1,.1)(2.5,1.1)(1.9,.1)(0,.1)

\drawpolygon[fillcolor=black](2,8)(2.5,9)(3,8)(4,8)(4.5,9)(5,8)(6,8)(6.5,9)(7,8)(8,8)(8.5,9)(8.4,9)(7.9,8.1)(7.1,8.1)(6.5,9.1)(5.9,8.1)(5.1,8.1)(4.5,9.1)(3.9,8.1)(3.1,8.1)(2.5,9.1)(1.85,8)

\drawpolygon[fillcolor=black](8,0)(8.5,1)(8,2)(9,2)(8.5,3)(8.45,2.9)(8.85,2.1)(7.9,2.1)(8.4,1)(7.9,0)

\drawpolygon[fillcolor=black](.5,3)(0,4)(1,4)(.5,5)(.55,5.1)(1.15,3.9)(.15,3.9)(.55,3.1)

\drawpolygon[fillcolor=black](8.5,5)(8,6)(9,6)(8.5,7)(8.55,7.1)(9.15,5.9)(8.15,5.9)(8.6,5)

\drawpolygon[fillcolor=black](.5,7)(1,8)(.5,9)(.55,9.1)(1.1,8)(.55,6.9)

\put(-.4,-.6){\scriptsize {\tiny $a_{11}$}}
\put(.7,-.6){\scriptsize {\tiny $a_{12}$}}
\put(1.7,-.6){\scriptsize {\tiny $a_{13}$}}
\put(6.8,-.6){\scriptsize {\tiny $a_{1i}$}}
\put(7.8,-.6){\scriptsize {\tiny $a_{11}$}}

\put(-.8,.9){\scriptsize {\tiny $x_{11}$}}
\put(8.8,.9){\scriptsize {\tiny $x_{11}$}}

\put(-.35,2){\scriptsize {\tiny $a_{21}$}}
\put(9.2,2){\scriptsize {\tiny $a_{21}$}}

\put(-.6,3){\scriptsize {\tiny $x_{21}$}}
\put(9,3){\scriptsize {\tiny $x_{21}$}}

\put(-.9,4.2){\scriptsize {\tiny $a_{31}$}}
\put(8.4,4.2){\scriptsize {\tiny $a_{31}$}}

\put(-1,8){\scriptsize {\tiny $a_{j1}$}}
\put(8.5,8){\scriptsize {\tiny $a_{j1}$}}

\put(-.5,9){\scriptsize {\tiny $x_{j1}$}}
\put(9,9){\scriptsize {\tiny $x_{j1}$}}

\put(-.7,10.5){\scriptsize {\tiny $a_{1(k+1)}$}}
\put(1.5,10.5){\scriptsize {\tiny $a_{1(k+2)}$}}
\put(7.5,10.5){\scriptsize {\tiny $a_{1k}$}}
\put(8.6,10.5){\scriptsize {\tiny $a_{1(k+1)}$}}

\put(-2,-.2){\scriptsize {\tiny $Q_1$}}
\put(-2,2){\scriptsize {\tiny $Q_2$}}
\put(-2,8){\scriptsize {\tiny $Q_j$}}
\put(-2,10){\scriptsize {\tiny $Q_{1}$}}

\put(4,-.5){\scriptsize $\ldots$}
\put(4,10.5){\scriptsize $\ldots$}
\put(5.75,10.5){\scriptsize $\ldots$}
\put(-2,5.25){\scriptsize $\vdots$}

\put(.2,-1.75){\scriptsize {\tiny {\bf Figure 4.7.4:} $M(i,j=2m+1,k)$ }}

\end{picture}

\vspace{5.2cm}

A construction of an $M(i,j,k)$ representation for $M$ follows by first cutting $M$ along a cycle of type $F_1$ and then cutting it along a cycle of type $F_2$. This gives the following result.

\begin{lem}\label{l3.7.1} 
	A DSEM $M$ of type $[3^2.6^2:3.6.3.6]$ admits an $M(i,j,k)$-representation iff the following holds: $(i)$ $j \geq 1$ and $i$ even,  $(ii)$  number of vertices of $M(i,j,k) = 3ij/2 \geq 15 $, $(iii)$ $i \geq 10$ if $j=1$, and $i \geq 6$ if $j \geq 2$, $(iv)$ if $j=1$ then $ k \in \{2r+5: 0 \leq r < (i-8)/2\} \setminus \{(i-8)/2+5\}$, if $j=2$ then $ k \in \{2r: 0 < r < i/2\}$, if $j=2m+1, \text{ where } m \in \mathbb{N}$ then $ k \in \{2r+1: 0 \leq r < i/2\}$, and if $j=2m+2, \text{ where } m \in \mathbb{N}$ then $ k \in \{2r: 0 \leq r < i/2\}$.
\end{lem}

\noindent{\bf Proof.} Let $M$ be a DSEM of type $[3^2.6^2:3.6.3.6]$ with $n$ vertices. Then $M(i,j,k)$ of $M$ has $j$ disjoint horizontal cycles of type $F_{1}$, say $Q_0, Q_1,\ldots, Q_{j-1}$, of length $i$. Note that the number of adjacent vertices with face-sequence $(3,6,3,6)$ which are lying on one side of horizontal cycles $Q_s$ and not belonging to any horizontal cycles $Q_s$ for $0 \leq s \leq j-1$ is $i/2 \cdot j$. So, $n = ij + ij/2 = 3ij/2$. 

If $j = 1$ and $i$ is not an even integer, then some vertex in the base horizontal cycle do not have face-sequence $(3^2,6^2)$, which is not possible. So, if $j = 1$, then $i$ is an even integer. Similarly, if $j > 1 $, then $i$ is an even integer.

If $j = 1$ and $ i < 10 $, then we get some vertex in $M(i,1,k)$ 
whose link can not be constructed. So, for $j = 1$, we get $i \geq 10$. Similarly, as above, we get that $ i \geq 6$ if $j \geq 2$. Thus $n = 3ij/2 \geq 15$.

If $j =1$ and $ k \in \{r:0 \leq r \leq i-1\} \setminus (\{2r+5: 0 \leq r < (i-8)/2\} \setminus \{(i-8)/2+5\})$, then we get some vertices whose link can not be constructed. So, for $j = 1$, we get $ k \in \{2r+5: 0 \leq r < (i-8)/2\} \setminus \{(i-8)/2+5\}$.  If  $j=2$ and $ k \in \{r:0 \leq r \leq i-1\} \setminus \{2r: 0 < r < i/2\}$, then  some vertices in the lower horizontal cycle do not follow the face-sequence $(3^2,6^2)$ after identifying the boundaries of $M(i,j,k)$. Which is a contradiction. So, for $j = 2$, we get $ k \in \{2r: 0 < r < i/2\}$ . Similarly as above, we see that if $j=2m+1, \text{ where } m \in \mathbb{N}$, then $ k \in \{2r+1: 0 \leq r < i/2\}$, and if $j=2m+2, \text{ where } m \in \mathbb{N}$, then $ k \in \{2r: 0 \leq r < i/2\}$. This completes the proof. \hfill $\Box$

%

\subsection {DSEMs of type $[3^6:3^2.6^2]$} \label{s3.8}

Let $M$ be a DSEM of type $[3^6:3^2.6^2]$. In $M$, we consider a path as follows.

A path $P_{1} = P( \ldots, y_{i},z_{i},z_{i+1},y_{i+1}, \ldots)$ in $M$, say of type $H_{1}$, indicated by thick black or green paths, shown in Figure 4.8.1. Here the vertices $y_i$'s and $z_i$'s have the face-sequences $(3^6)$ and $(3^2,6^2)$ respectively.

\begin{picture}(0,0)(-2,35)
\setlength{\unitlength}{4.5mm}

\drawpolygon(3,0)(6,0)(9,6)(6,6)

\drawline[AHnb=0](1.75,3.5)(3.75,7.5)


\drawline[AHnb=0](5.25,-1.5)(9.75,7.5)

\drawline[AHnb=0](8.25,-1.5)(10.25,2.5)


\drawline[AHnb=0](1.75,3)(10.5,3)

\drawline[AHnb=0](1.5,0)(10.5,0)

\drawline[AHnb=0](1.5,6)(10.5,6)

\drawline[AHnb=0](3.75,-1.5)(1.75,2.5)

\drawline[AHnb=0](6.75,-1.5)(2.25,7.5)

\drawline[AHnb=0](9.75,-1.5)(5.25,7.5)

\drawline[AHnb=0](10.25,3.5)(8.25,7.5)


\drawline[AHnb=0](2,0)(2.5,1)(3.5,1)(4,0)(3.5,-1)(2.5,-1)(2,0)

\drawline[AHnb=0](5,0)(5.5,1)(6.5,1)(7,0)(6.5,-1)(5.5,-1)(5,0)

\drawline[AHnb=0](8,0)(8.5,1)(9.5,1)(10,0)(9.5,-1)(8.5,-1)(8,0)

\drawline[AHnb=0](1.5,2)(2,2)(2.5,3)(2,4)(1.5,4)

\drawline[AHnb=0](3.5,3)(4,2)(5,2)(5.5,3)(5,4)(4,4)(3.5,3)

\drawline[AHnb=0](7,2)(8,2)(8.5,3)(8,4)(7,4)(6.5,3)(7,2)

\drawline[AHnb=0](10.5,2)(10,2)(9.5,3)(10,4)(10.5,4)

\drawline[AHnb=0](5,6)(5.5,5)(6.5,5)(7,6)(6.5,7)(5.5,7)(5,6)

\drawline[AHnb=0](8,6)(8.5,5)(9.5,5)(10,6)(9.5,7)(8.5,7)(8,6)

\drawline[AHnb=0](2.5,5)(3.5,5)(4,6)(3.5,7)(2.5,7)(2,6)(2.5,5)

\drawline[AHnb=0](5.25,-1.5)(9.75,7.5)

\drawpolygon[fillcolor=black](1.5,6)(10.5,6)(10.5,6.1)(1.5,6.1)

\drawpolygon[fillcolor=black](1.5,3)(10.5,3)(10.5,3.1)(1.5,3.1)

\drawpolygon[fillcolor=black](1.5,0)(10.5,0)(10.5,.15)(1.5,.15)

\drawpolygon[fillcolor=green](1.5,3)(3.75,7.5)(3.85,7.5)(1.6,3)

\drawpolygon[fillcolor=green](2.25,-1.5)(3,0)(6,6)(6.75,7.5)(6.85,7.5)(6.1,6)(3.1,0)(2.35,-1.5)

\drawpolygon[fillcolor=green](5.25,-1.5)(6,0)(9,6)(9.75,7.5)(9.85,7.5)(9.1,6)(6.1,0)(5.35,-1.5)

\drawpolygon[fillcolor=green](8.25,-1.5)(9,0)(10.5,3)(10.6,3)(9.1,0)(8.35,-1.5)


\put(2,-2.75){\scriptsize {\tiny {\bf Figure 4.8.1:}  Paths of type $H_1$ }}

\end{picture}

\begin{picture}(0,0)(-70,36)
\setlength{\unitlength}{4.5mm}

\drawpolygon(0,0)(9,0)(13.5,9)(4.5,9)


\drawline[AHnb=0](3,0)(7.5,9)

\drawline[AHnb=0](6,0)(10.5,9)


\drawline[AHnb=0](1.5,3)(10.5,3)
\drawline[AHnb=0](3,6)(12,6)
\drawline[AHnb=0](4.5,9)(13.5,9)

\drawline[AHnb=0](3,0)(1.5,3)

\drawline[AHnb=0](6,0)(3,6)

\drawline[AHnb=0](9,0)(4.5,9)

\drawline[AHnb=0](10.5,3)(7.5,9)


\drawline[AHnb=0](.5,1)(1,0)

\drawline[AHnb=0](2.5,1)(2,0)
\drawline[AHnb=0](2.5,1)(3.5,1)
\drawline[AHnb=0](3.5,1)(4,0)

\drawline[AHnb=0](5,0)(5.5,1)
\drawline[AHnb=0](5.5,1)(6.5,1)
\drawline[AHnb=0](6.5,1)(7,0)

\drawline[AHnb=0](8,0)(8.5,1)
\drawline[AHnb=0](8.5,1)(9.5,1)

\drawline[AHnb=0](1,2)(2,2)
\drawline[AHnb=0](2,2)(2.5,3)
\drawline[AHnb=0](2.5,3)(2,4)

\drawline[AHnb=0](3.5,3)(4,2)
\drawline[AHnb=0](4,2)(5,2)
\drawline[AHnb=0](5,2)(5.5,3)
\drawline[AHnb=0](5.5,3)(5,4)
\drawline[AHnb=0](5,4)(4,4)
\drawline[AHnb=0](4,4)(3.5,3)

\drawline[AHnb=0](6.5,3)(7,4)
\drawline[AHnb=0](7,2)(8,2)
\drawline[AHnb=0](8,2)(8.5,3)
\drawline[AHnb=0](8.5,3)(8,4)
\drawline[AHnb=0](8,4)(7,4)
\drawline[AHnb=0](7,2)(6.5,3)

\drawline[AHnb=0](9.5,3)(10,4)
\drawline[AHnb=0](10,4)(11,4)
\drawline[AHnb=0](9.5,3)(10,2)

\drawline[AHnb=0](2.5,5)(3.5,5)
\drawline[AHnb=0](3.5,5)(4,6)

\drawline[AHnb=0](3.5,7)(4,6)

\drawline[AHnb=0](5,6)(5.5,7)
\drawline[AHnb=0](5.5,7)(6.5,7)
\drawline[AHnb=0](6.5,7)(7,6)

\drawline[AHnb=0](5,6)(5.5,5)
\drawline[AHnb=0](5.5,5)(6.5,5)
\drawline[AHnb=0](6.5,5)(7,6)

\drawline[AHnb=0](8,6)(8.5,7)
\drawline[AHnb=0](8.5,7)(9.5,7)
\drawline[AHnb=0](9.5,7)(10,6)

\drawline[AHnb=0](8,6)(8.5,5)
\drawline[AHnb=0](8.5,5)(9.5,5)
\drawline[AHnb=0](9.5,5)(10,6)

\drawline[AHnb=0](11,6)(11.5,7)
\drawline[AHnb=0](11.5,7)(12.5,7)

\drawline[AHnb=0](11,6)(11.5,5)
\drawline[AHnb=0](10,4)(11,4)
\drawline[AHnb=0](9.5,3)(10,2)

\drawline[AHnb=0](6.5,9)(7,8)(8,8)(8.5,9)
\drawline[AHnb=0](9.5,9)(10,8)(11,8)(11.5,9)

\drawline[AHnb=0](10.5,9)(12,6)
\drawline[AHnb=0](12.5,9)(13,8)

\drawline[AHnb=0](5.5,9)(5,8)
\drawline[AHnb=0](4,8)(5,8)

\put(-.2,-.6){\scriptsize {\tiny $a_{11}$}}
\put(.7,-.6){\scriptsize {\tiny $a_{12}$}}
\put(1.7,-.6){\scriptsize {\tiny $a_{13}$}}
\put(7.8,-.6){\scriptsize {\tiny $a_{1i}$}}
\put(8.8,-.6){\scriptsize {\tiny $a_{11}$}}

\put(-.5,.9){\scriptsize {\tiny $x_{11}$}}
\put(-.3,2){\scriptsize {\tiny $x_{12}$}}
\put(9.8,.9){\scriptsize {\tiny $x_{11}$}}
\put(10.3,2){\scriptsize {\tiny $x_{12}$}}

\put(.5,3){\scriptsize {\tiny $a_{21}$}}
\put(10.8,3){\scriptsize {\tiny $a_{21}$}}

\put(2,6){\scriptsize {\tiny $a_{j1}$}}
\put(12.5,6){\scriptsize {\tiny $a_{j1}$}}

\put(2.4,7){\scriptsize {\tiny $x_{j1}$}}
\put(3.2,8.3){\scriptsize {\tiny $x_{j2}$}}
\put(13,6.9){\scriptsize {\tiny $x_{j1}$}}
\put(13.4,8){\scriptsize {\tiny $x_{j2}$}}

\put(2.75,9.4){\scriptsize {\tiny $a_{1(k+1)}$}}
\put(4.85,9.4){\scriptsize {\tiny $a_{1(k+2)}$}}

\put(12.2,9.4){\scriptsize {\tiny $a_{1k}$}}
\put(13.2,9.4){\scriptsize {\tiny $a_{1(k+1)}$}}

\drawpolygon[fillcolor=black](0,0)(8,0)(8.5,1)(8,2)(7,2)(6.5,1)(5.5,1)(5,2)(4,2)(3.5,1)(2.5,1)(2,2)(2.5,3)(9.5,3)(10,4)(9.5,5)(8.5,5)(8,4)(7,4)(6.5,5)(5.5,5)(5,4)(4,4)(3.5,5)(4,6)(11,6)(11.5,7)(11,8)(10,8)(9.5,7)(8.5,7)(8,8)(7,8)(6.5,7)(5.5,7)(5,8)(4,8)(0,0)(-.1,0)(3.9,8.1)(5.1,8.1)(5.6,7.1)(6.4,7.1)(6.9,8.1)(8.1,8.1)(8.6,7.1)(9.4,7.1)(9.9,8.1)(11.1,8.1)(11.6,7)(11.1,5.9)(4.1,5.9)(3.65,5)(4.1,4.1)(4.9,4.1)(5.4,5.1)(6.6,5.1)(7.1,4.1)(7.9,4.1)(8.4,5.1)(9.6,5.1)(10.1,4)(9.6,2.9)(2.6,2.9)(2.1,2)(2.6,1.1)(3.4,1.1)(3.9,2.1)(5.1,2.1)(5.6,1.1)(6.4,1.1)(6.9,2.1)(8.1,2.1)(8.6,1)(8.1,-.1)(-.1,-.1)

\put(-3,-.1){\scriptsize {\tiny $Q_1$}}
\put(-.5,6){\scriptsize {\tiny $Q_{j}$}}
\put(.5,9){\scriptsize {\tiny $Q_{1}$}}

\put(-.7,5){\scriptsize $.$}
\put(-.8,4.8){\scriptsize $.$}
\put(-.9,4.6){\scriptsize $.$}

\put(-1.6,3){\scriptsize $.$}
\put(-1.7,2.8){\scriptsize $.$}
\put(-1.8,2.6){\scriptsize $.$}

\put(4,-.5){\scriptsize $\ldots$}
\put(5.5,-.5){\scriptsize $\ldots$}
\put(8.5,9.5){\scriptsize $\ldots$}

\put(2.5,-1.5){\scriptsize {\tiny {\bf Figure 4.8.2:} $M(i,j,k)$ }}

\end{picture}

\vspace{4.5cm}

Now, construct $M(i,j,k)$ representation of $M$ by first cutting $M$ along a black colored cycle of type $H_1$ through a vertex $v$ with face-sequence $(3^6)$ and then cutting it along a green colored cycle of type $H_1$. Note that both these cycles are non-homologous. Then we have the following lemma.

\begin{lem}\label{l3.8.1} 
	
	A DSEM $M$ of type $[3^6:3^2.6^2]$ admits an $M(i,j,k)$-representation iff the following holds: $(i)$ $i \geq 9$ and $i=3m$, $m \in \mathbb{N} $ for $j=1$, $(ii)$ $i \geq 6$ and $i=3m$, $m \in \mathbb{N} $ for $j > 1$, $(iii)$ number of vertices of $M(i,j,k) = 7ij/3 $, $(iv)$ if $j=1$ then $ k \in \{3r: 1 < r < i/3\}$, and if $j > 1$ then $ k \in \{3r: 0 \leq r < i/3\}$.
\end{lem}

\noindent{\bf Proof.} Let $M$ be a DSEM of type $[3^6:3^2.6^2]$ with $n$ vertices. Clearly $M(i,j,k)$ of $M$ has $j$ disjoint horizontal cycles of $H_{1}$ type, say $Q_0, Q_1,\ldots, Q_{j-1}$, of length $i$. Note that the number of vertices with face-sequence $(3^2,6^2)$ which are above the horizontal cycles $Q_s$ and not belonging to any horizontal cycles $Q_s$ for $0 \leq s \leq j-1$ is $4i/3 \cdot j$. In $M$, therefore  $n = ij + 4ij/3 = 7ij/3$. 

If $j = 1$ and $ i < 9 $, then the link of some vertices in $M(i,1,k)$ can not be completed. This is not possible. So, for $j = 1$, we get $i \geq 9$. Also, if $i \geq 9$ and $i$ is not a multiple of $3$, then after identifying the boundaries of $M(i,j,k)$, there is no vertex in the horizontal base cycle having face-sequences $(3^6)$ and $(3^2,6^2)$. Which is not possible. So, $i \geq 9$ and $i=3m$, $m \in \mathbb{N} $ for $j=1$. Similarly, as above,  for $j > 1$ we get that $ i \geq 6$ and $i = 3m, \text{ where } m \in \mathbb{N} $. Thus $ n = 7ij/3 \geq 21$.

If $j =1$ and $ k \in \{r:0 \leq r \leq i-1\} \setminus \{3r: 1 < r < i/3\}$, then we get some vertices whose link can not be completed or does not follow the face-sequences $(3^6)$ and $(3^2,6^2)$. So, for $j = 1$, we get $ k \in \{3r: 1 < r < i/3\}$. Similarly, as above, we see that if $j > 1$, then $ k \in \{3r: 0 \leq r < i/3\}$. Thus the proof. \hfill $\Box$

%

\subsection{DSEMs of type $[3^4.6:3^2.6^2]$} \label{s3.9}

Let $M$ be a DSEM of type $[3^4.6:3^2.6^2]$. Consider the following types of paths as follows.

A path $P_{1} = P( \ldots, y_{i},z_{i},y_{i+1},z_{i+1}, \ldots)$ in $M$, say of type $I_{1}$, indicated by thick black paths, shown in Figure 4.9.1. Here, the vertices $y_i$'s and $z_i$'s have the face-sequences $(3^2,6^2)$ and $(3^4,6)$ respectively.

 A path $P_{2} = P( \ldots,y_{i},y_{i+1},z_{i},z_{i+1},y_{i+2},y_{i+3}, \ldots)$ in $M$, say of type $I_{2}$, indicated by green paths, shown in Figure 4.9.1. Here the vertices $y_i$'s and $z_i$'s have the face-sequences $(3^2,6^2)$ and $(3^4,6)$ respectively.

\begin{picture}(0,0)(5,45)
\setlength{\unitlength}{5.5mm}


\drawline[AHnb=0](2.25,1.75)(3,2)(3,1)
\drawline[AHnb=0](3,1.5)(5,2)(5,1)
\drawline[AHnb=0](5,1.5)(7,2)(7,1)
\drawline[AHnb=0](7,1.5)(9,2)(9,1)
\drawline[AHnb=0](9,1.5)(11,2)(11,1)

\drawline[AHnb=0](2.25,3)(4,3.5)(4,1.75)
\drawline[AHnb=0](4,3)(6,3.5)(6,1.75)
\drawline[AHnb=0](6,3)(8,3.5)(8,1.75)
\drawline[AHnb=0](8,3)(10,3.5)(10,1.75)
\drawline[AHnb=0](10,3)(11.75,3.5)

\drawline[AHnb=0](2.25,4.75)(3,4.5)(5,5)(5,3.25)
\drawline[AHnb=0](5,4.5)(7,5)(7,3.25)
\drawline[AHnb=0](7,4.5)(9,5)(9,3.25)
\drawline[AHnb=0](9,4.5)(11,5)(11,3.25)
\drawline[AHnb=0](11,4.5)(11.75,4.75)

\drawline[AHnb=0](2.25,6.5)(4,6)(6,6.5)(6,4.75)
\drawline[AHnb=0](6,6)(8,6.5)(8,4.75)
\drawline[AHnb=0](8,6)(10,6.5)(10,4.75)
\drawline[AHnb=0](10,6)(11.75,6.5)


\drawline[AHnb=0](2.25,1.75)(3,1.5)
\drawline[AHnb=0](3,2)(5,1.5)
\drawline[AHnb=0](5,2)(7,1.5)
\drawline[AHnb=0](7,2)(9,1.5)
\drawline[AHnb=0](9,2)(11,1.5)

\drawline[AHnb=0](2.25,3.5)(4,3)
\drawline[AHnb=0](4,3.5)(6,3)
\drawline[AHnb=0](6,3.5)(8,3)
\drawline[AHnb=0](8,3.5)(10,3)
\drawline[AHnb=0](10,3.5)(11.75,3)

\drawline[AHnb=0](2.25,4.75)(3,5)(5,4.5)
\drawline[AHnb=0](5,5)(7,4.5)
\drawline[AHnb=0](7,5)(9,4.5)
\drawline[AHnb=0](9,5)(11,4.5)
\drawline[AHnb=0](11,5)(11.75,4.75)

\drawline[AHnb=0](4,6.5)(2.25,6)
\drawline[AHnb=0](4,6.5)(6,6)
\drawline[AHnb=0](6,6.5)(8,6)
\drawline[AHnb=0](8,6.5)(10,6)
\drawline[AHnb=0](10,6.5)(11.75,6)


\drawline[AHnb=0](3,3.25)(3,4)
\drawline[AHnb=0](4,4.75)(4,5.5)
\drawline[AHnb=0](6,4.75)(6,5.5)
\drawline[AHnb=0](8,4.75)(8,5.5)
\drawline[AHnb=0](10,4.75)(10,5.5)

\drawline[AHnb=0](3,6.25)(3,7)
\drawline[AHnb=0](5,6.25)(5,7)
\drawline[AHnb=0](7,6.25)(7,7)
\drawline[AHnb=0](9,6.25)(9,7)
\drawline[AHnb=0](11,6.25)(11,7)


\drawline[AHnb=0](2.25,2.35)(3,2)(4,2.5)(5,2)(6,2.5)(7,2)(8,2.5)(9,2)(10,2.5)(11,2)(11.75,2.25)

\drawline[AHnb=0](2.25,3.75)(3,4)(4,3.5)(5,4)(6,3.5)(7,4)(8,3.5)(9,4)(10,3.5)(11,4)(11.75,3.75)

\drawline[AHnb=0](2.25,5.5)(3,5)(4,5.5)(5,5)(6,5.5)(7,5)(8,5.5)(9,5)(10,5.5)(11,5)(11.75,5.25)

\drawline[AHnb=0](4,5.5)(4,6.5)(5,7)(6,6.5)(7,7)(8,6.5)(9,7)(10,6.5)(11,7)(11.75,6.75)

\drawline[AHnb=0](10,5.5)(9,5)(8,5.5)(7,5)(6,5.5)(5,5)(4,5.5)(3,5)(3,4)
\drawline[AHnb=0](4,6.5)(3,7)(2.25,6.5)

\drawpolygon[fillcolor=black](2.25,3.75)(3,4)(4,3.5)(5,4)(6,3.5)(7,4)(8,3.5)(9,4)(10,3.5)(11,4)(11.75,3.75)(11.75,3.85)(11,4.1)(10,3.6)(9,4.1)(8,3.6)(7,4.1)(6,3.6)(5,4.1)(4,3.6)(3,4.1)(2.25,3.85)

\drawpolygon[fillcolor=black](2.25,6.5)(3,7)(4,6.5)(5,7)(6,6.5)(7,7)(8,6.5)(9,7)(10,6.5)(11,7)(11.75,6.75)(11.75,6.85)(11,7.1)(10,6.6)(9,7.1)(8,6.6)(7,7.1)(6,6.6)(5,7.1)(4,6.6)(3,7.1)(2.25,6.65)

\drawpolygon[fillcolor=black](2.25,5.5)(3,5)(4,5.5)(5,5)(6,5.5)(7,5)(8,5.5)(9,5)(10,5.5)(11,5)(12,5.5)(12,5.6)(11,5.1)(10,5.6)(9,5.1)(8,5.6)(7,5.1)(6,5.6)(5,5.1)(4,5.6)(3,5.1)(2.25,5.6)

\drawpolygon[fillcolor=black](2.25,2.35)(3,2)(4,2.5)(5,2)(6,2.5)(7,2)(8,2.5)(9,2)(10,2.5)(11,2)(11.75,2.25)(11.75,2.35)(11,2.1)(10,2.6)(9,2.1)(8,2.6)(7,2.1)(6,2.6)(5,2.1)(4,2.6)(3,2.1)(2.25,2.45)

\drawpolygon[fillcolor=black](2.25,6)(3,6.25)(4,6)(5,6.25)(6,6)(7,6.25)(8,6)(9,6.25)(10,6)(11,6.25)(11.75,6)(11.75,6.15)(11,6.35)(10,6.1)(9,6.35)(8,6.1)(7,6.35)(6,6.1)(5,6.35)(4,6.1)(3,6.35)(2.25,6.1)

\drawpolygon[fillcolor=black](2.25,3)(3,3.25)(4,3)(5,3.25)(6,3)(7,3.25)(8,3)(9,3.25)(10,3)(11,3.25)(11.75,3)(11.75,3.15)(11,3.35)(10,3.1)(9,3.35)(8,3.1)(7,3.35)(6,3.1)(5,3.35)(4,3.1)(3,3.35)(2.25,3.1)

\drawpolygon[fillcolor=black](2.25,4.75)(3,4.5)(4,4.75)(5,4.5)(6,4.75)(7,4.5)(8,4.75)(9,4.5)(10,4.75)(11,4.5)(12,4.75)(12,4.85)(11,4.6)(10,4.85)(9,4.6)(8,4.85)(7,4.6)(6,4.85)(5,4.6)(4,4.85)(3,4.6)(2.25,4.85)

\drawpolygon[fillcolor=black](2.25,1.75)(3,1.5)(4,1.75)(5,1.5)(6,1.75)(7,1.5)(8,1.75)(9,1.5)(10,1.75)(11,1.5)(12,1.75)(12,1.85)(11,1.6)(10,1.85)(9,1.6)(8,1.85)(7,1.6)(6,1.85)(5,1.6)(4,1.85)(3,1.6)(2.25,1.85)

\drawpolygon[fillcolor=green](2.25,3.75)(3,4)(3,5)(4,5.5)(4,6.5)(5,7)(5,7.5)(5.1,7.5)(5.1,6.95)(4.1,6.45)(4.1,5.45)(3.1,4.95)(3.1,3.95)(2.26,3.65)

\drawpolygon[fillcolor=green](3,1)(3,2)(4,2.5)(4,3.5)(5,4)(5,5)(6,5.5)(6,6.5)(7,7)(7,7.5)(7.1,7.5)(7.1,6.95)(6.1,6.45)(6.1,5.45)(5.1,4.95)(5.1,3.95)(4.1,3.45)(4.1,2.45)(3.1,1.95)(3.1,1)

\drawpolygon[fillcolor=green](5,1)(5,2)(6,2.5)(6,3.5)(7,4)(7,5)(8,5.5)(8,6.5)(9,7)(9,7.5)(9.1,7.5)(9.1,6.95)(8.1,6.45)(8.1,5.45)(7.1,4.95)(7.1,3.95)(6.1,3.45)(6.1,2.45)(5.1,1.95)(5.1,1)

\drawpolygon[fillcolor=green](7,1)(7,2)(8,2.5)(8,3.5)(9,4)(9,5)(10,5.5)(10,6.5)(11,7)(11,7.5)(11.1,7.5)(11.1,6.95)(10.1,6.45)(10.1,5.45)(9.1,4.95)(9.1,3.95)(8.1,3.45)(8.1,2.45)(7.1,1.95)(7.1,1)

\drawpolygon[fillcolor=green](9,1)(9,2)(10,2.5)(10,3.5)(11,4)(11,5)(12,5.5)(12.1,5.45)(11.1,4.95)(11.1,3.95)(10.1,3.45)(10.1,2.45)(9.1,1.95)(9.1,1)

\put(8.7,8){\scriptsize {\tiny  $I_2$ }}

\put(12,3.8){\scriptsize {\tiny $I_1$ }}

\put(3.5,0){\scriptsize {\tiny {\bf Figure 4.9.1:} Paths of type $I_1, I_2$ }}

\end{picture}

\begin{picture}(0,0)(-76,8)
\setlength{\unitlength}{4.75mm}


\drawline[AHnb=0](2,2.5)(1,2)(1,1)(2,.5)(3,1)(4,.5)(5,1)(6,.5)(7,1)(8,.5)(9,1)(10,.5)(11,1)(11,2)

\drawline[AHnb=0](1,1.5)(3,2)(3,1)
\drawline[AHnb=0](3,1.5)(5,2)(5,1)
\drawline[AHnb=0](5,1.5)(7,2)(7,1)
\drawline[AHnb=0](7,1.5)(9,2)(9,1)
\drawline[AHnb=0](9,1.5)(11,2)(11,1)

\drawline[AHnb=0](1,2)(3,1.5)
\drawline[AHnb=0](3,2)(5,1.5)
\drawline[AHnb=0](5,2)(7,1.5)
\drawline[AHnb=0](7,2)(9,1.5)
\drawline[AHnb=0](9,2)(11,1.5)

\drawline[AHnb=0](2,2.5)(3,2)(4,2.5)(5,2)(6,2.5)(7,2)(8,2.5)(9,2)(10,2.5)(11,2)


\drawpolygon[fillcolor=black](1,1.5)(2,1.75)(2,2.5)(3,2)(3,1.5)(4,1.75)(4,2.5)(5,2)(5,1.5)(6,1.75)(6,2.5)(7,2)(7,1.5)(8,1.75)(8,2.5)(9,2)(9,1.5)(10,1.75)(10,2.5)(11,2)(11,1.5)(11.1,1.5)(11.1,2.1)(9.9,2.6)(9.9,1.8)(9.1,1.6)(9.1,2.1)(7.9,2.6)(7.9,1.8)(7.1,1.6)(7.1,2.1)(5.9,2.6)(5.9,1.8)(5.1,1.6)(5.1,2.1)(3.9,2.6)(3.9,1.8)(3.1,1.6)(3.1,2.1)(1.9,2.6)(1.9,1.85)(1,1.6)

\put(.5,.6){\scriptsize {\tiny $a_{11}$}}
\put(1.7,.2){\scriptsize {\tiny $a_{12}$}}
\put(2.7,.5){\scriptsize {\tiny $a_{13}$}}
\put(9.7,0){\scriptsize {\tiny $a_{1i}$}}
\put(10.7,.5){\scriptsize {\tiny $a_{11}$}}

\put(0,1.4){\scriptsize {\tiny $a_{21}$}}
\put(11.2,1.2){\scriptsize {\tiny $a_{21}$}}

\put(-.2,2.45){\scriptsize {\tiny $a_{1(k+1)}$}}
\put(1.6,2.85){\scriptsize {\tiny $a_{1(k+2)}$}}
\put(10,2.7){\scriptsize {\tiny $a_{1k}$}}
\put(11.2,2){\scriptsize {\tiny $a_{1(k+1)}$}}

\put(4,0){\scriptsize $\ldots$}
\put(6.5,0){\scriptsize $\ldots$}
\put(4.5,2.9){\scriptsize $\ldots$}
\put(7.5,2.9){\scriptsize $\ldots$}

\put(3,-.9){\scriptsize {\tiny {\bf Figure 4.9.2:} $M(i,2,k)$ }}

\end{picture}

\begin{picture}(0,0)(-74,34)
\setlength{\unitlength}{4.75mm}


\drawline[AHnb=0](3,4)(2,3.5)(2,2.5)(1,2)(1,1)(2,.5)(3,1)(4,.5)(5,1)(6,.5)(7,1)(8,.5)(9,1)(10,.5)(11,1)(11,2)(12,2.5)(12,3.5)

\drawline[AHnb=0](1,1.5)(3,2)(3,1)
\drawline[AHnb=0](3,1.5)(5,2)(5,1)
\drawline[AHnb=0](5,1.5)(7,2)(7,1)
\drawline[AHnb=0](7,1.5)(9,2)(9,1)
\drawline[AHnb=0](9,1.5)(11,2)(11,1)

\drawline[AHnb=0](2,3)(4,3.5)(4,1.75)
\drawline[AHnb=0](4,3)(6,3.5)(6,1.75)
\drawline[AHnb=0](6,3)(8,3.5)(8,1.75)
\drawline[AHnb=0](8,3)(10,3.5)(10,1.75)
\drawline[AHnb=0](10,3)(12,3.5)


\drawline[AHnb=0](2,2.5)(3,2)(4,2.5)(5,2)(6,2.5)(7,2)(8,2.5)(9,2)(10,2.5)(11,2)

\drawline[AHnb=0](3,4)(4,3.5)(5,4)(6,3.5)(7,4)(8,3.5)(9,4)(10,3.5)(11,4)(12,3.5)

\drawline[AHnb=0](1,2)(3,1.5)
\drawline[AHnb=0](3,2)(5,1.5)
\drawline[AHnb=0](5,2)(7,1.5)
\drawline[AHnb=0](7,2)(9,1.5)
\drawline[AHnb=0](9,2)(11,1.5)

\drawline[AHnb=0](2,3.5)(4,3)
\drawline[AHnb=0](4,3.5)(6,3)
\drawline[AHnb=0](6,3.5)(8,3)
\drawline[AHnb=0](8,3.5)(10,3)
\drawline[AHnb=0](10,3.5)(12,3)

\drawline[AHnb=0](2,1.75)(2,2.5)

\drawpolygon[fillcolor=black](1,1.5)(3,2)(3,1.5)(4,1.75)(4,2.5)(5,2)(5,1.5)(6,1.75)(6,2.5)(7,2)(7,1.5)(8,1.75)(8,2.5)(9,2)(9,1.5)(10,1.75)(10,3)(11,3.25)(11,4)(10,3.5)(9,3.25)(9,4)(8,3.5)(8,3)(7,3.25)(7,4)(6,3.5)(6,3)(5,3.25)(5,4)(4,3.5)(4,3)(3,3.25)(3,4)(2,3.5)(2,2.5)(1,2)(1,1.5)(1.1,1.5)(1.1,1.9)(2.1,2.4)(2.1,3.4)(2.9,3.8)(2.9,3.2)(4.1,2.9)(4.1,3.4)(4.9,3.8)(4.9,3.2)(6.1,2.9)(6.1,3.45)(6.9,3.85)(6.9,3.2)(8.1,2.9)(8.1,3.4)(8.9,3.85)(8.9,3.1)(10.1,3.4)(10.9,3.85)(10.9,3.3)(9.9,3.1)(9.9,1.8)(9.1,1.6)(9.1,2.1)(7.9,2.6)(7.9,1.8)(7.1,1.6)(7.1,2.1)(5.9,2.6)(5.9,1.8)(5.1,1.6)(5.1,2.1)(3.9,2.6)(3.9,1.8)(3.1,1.6)(3.1,2.1)(1.1,1.6)

\put(.5,.6){\scriptsize {\tiny $a_{11}$}}
\put(1.7,.2){\scriptsize {\tiny $a_{12}$}}
\put(2.7,.5){\scriptsize {\tiny $a_{13}$}}
\put(9.7,0){\scriptsize {\tiny $a_{1i}$}}
\put(10.7,.5){\scriptsize {\tiny $a_{11}$}}

\put(0,1.2){\scriptsize {\tiny $a_{21}$}}
\put(11.2,1.2){\scriptsize {\tiny $a_{21}$}}

\put(0,2){\scriptsize {\tiny $a_{31}$}}
\put(2.2,2.45){\scriptsize {\tiny $a_{32}$}}
\put(11.1,1.8){\scriptsize {\tiny $a_{31}$}}
\put(12.1,2.2){\scriptsize {\tiny $a_{32}$}}

\put(1,2.8){\scriptsize {\tiny $a_{41}$}}
\put(12.1,2.7){\scriptsize {\tiny $a_{41}$}}

\put(0.25,3.9){\scriptsize {\tiny $a_{1(k+1)}$}}
\put(2.5,4.3){\scriptsize {\tiny $a_{1(k+2)}$}}
\put(11,4.1){\scriptsize {\tiny $a_{1k}$}}
\put(12.1,3.45){\scriptsize {\tiny $a_{1(k+1)}$}}

\put(4.5,0){\scriptsize $\ldots$}
\put(6.5,0){\scriptsize $\ldots$}

\put(5.5,4.5){\scriptsize $\ldots$}
\put(8.5,4.5){\scriptsize $\ldots$}

\put(4,-1){\scriptsize {\tiny {\bf Figure 4.9.3:} $M(i,4,k)$ }}

\end{picture}


\begin{picture}(0,0)(-2.5,80)
\setlength{\unitlength}{4.5mm}


\drawpolygon(1,1)(2,.5)(3,1)(4,.5)(5,1)(6,.5)(7,1)(8,.5)(9,1)(10,.5)(11,1)(11,2)(12,2.5)(12,3.5)(13,4)(13,5)(12,5.5)(11,5)(10,5.5)(9,5)(8,5.5)(7,5)(6,5.5)(5,5)(4,5.5)(3,5)(3,4)(2,3.5)(2,2.5)(1,2)

\drawpolygon(5,7)(6,6.5)(7,7)(8,6.5)(9,7)(10,6.5)(11,7)(12,6.5)(13,7)(14,6.5)(15,7)(15,8)(14,8.5)(13,8)(12,8.5)(11,8)(10,8.5)(9,8)(8,8.5)(7,8)(6,8.5)(5,8)

\drawline[AHnb=0](7,7)(7,8)(5,7.5)
\drawline[AHnb=0](9,7)(9,8)(7,7.5)
\drawline[AHnb=0](11,7)(11,8)(9,7.5)
\drawline[AHnb=0](13,7)(13,8)(11,7.5)

\drawline[AHnb=0](7,7.5)(5,8)
\drawline[AHnb=0](9,7.5)(7,8)
\drawline[AHnb=0](11,7.5)(9,8)
\drawline[AHnb=0](13,7.5)(11,8)
\drawline[AHnb=0](15,7.5)(13,8)
\drawline[AHnb=0](15,8)(13,7.5)

\drawline[AHnb=0](8,7.75)(8,8.5)
\drawline[AHnb=0](10,7.75)(10,8.5)
\drawline[AHnb=0](12,7.75)(12,8.5)
\drawline[AHnb=0](14,7.75)(14,8.5)
\drawline[AHnb=0](6,7.75)(6,8.5)


\drawline[AHnb=0](1,1.5)(3,2)(3,1)
\drawline[AHnb=0](3,1.5)(5,2)(5,1)
\drawline[AHnb=0](5,1.5)(7,2)(7,1)
\drawline[AHnb=0](7,1.5)(9,2)(9,1)
\drawline[AHnb=0](9,1.5)(11,2)(11,1)

\drawline[AHnb=0](2,3)(4,3.5)(4,1.75)
\drawline[AHnb=0](4,3)(6,3.5)(6,1.75)
\drawline[AHnb=0](6,3)(8,3.5)(8,1.75)
\drawline[AHnb=0](8,3)(10,3.5)(10,1.75)
\drawline[AHnb=0](10,3)(12,3.5)

\drawline[AHnb=0](3,4.5)(5,5)(5,3.25)
\drawline[AHnb=0](5,4.5)(7,5)(7,3.25)
\drawline[AHnb=0](7,4.5)(9,5)(9,3.25)
\drawline[AHnb=0](9,4.5)(11,5)(11,3.25)
\drawline[AHnb=0](11,4.5)(13,5)

\drawline[AHnb=0](4,6)(6,6.5)(6,4.75)
\drawline[AHnb=0](6,6)(8,6.5)(8,4.75)
\drawline[AHnb=0](8,6)(10,6.5)(10,4.75)
\drawline[AHnb=0](10,6)(12,6.5)(12,4.75)
\drawline[AHnb=0](12,6)(14,6.5)


\drawline[AHnb=0](1,2)(3,1.5)
\drawline[AHnb=0](3,2)(5,1.5)
\drawline[AHnb=0](5,2)(7,1.5)
\drawline[AHnb=0](7,2)(9,1.5)
\drawline[AHnb=0](9,2)(11,1.5)

\drawline[AHnb=0](2,3.5)(4,3)
\drawline[AHnb=0](4,3.5)(6,3)
\drawline[AHnb=0](6,3.5)(8,3)
\drawline[AHnb=0](8,3.5)(10,3)
\drawline[AHnb=0](10,3.5)(12,3)

\drawline[AHnb=0](3,5)(5,4.5)
\drawline[AHnb=0](5,5)(7,4.5)
\drawline[AHnb=0](7,5)(9,4.5)
\drawline[AHnb=0](9,5)(11,4.5)
\drawline[AHnb=0](11,5)(13,4.5)

\drawline[AHnb=0](4,6.5)(6,6)
\drawline[AHnb=0](6,6.5)(8,6)
\drawline[AHnb=0](8,6.5)(10,6)
\drawline[AHnb=0](10,6.5)(12,6)
\drawline[AHnb=0](12,6.5)(14,6)

\drawline[AHnb=0](2,1.75)(2,2.5)
\drawline[AHnb=0](3,3.25)(3,4)
\drawline[AHnb=0](4,4.75)(4,5.5)
\drawline[AHnb=0](6,4.75)(6,5.5)
\drawline[AHnb=0](8,4.75)(8,5.5)
\drawline[AHnb=0](10,4.75)(10,5.5)
\drawline[AHnb=0](12,4.75)(12,5.5)

\drawline[AHnb=0](5,6.25)(5,7)
\drawline[AHnb=0](7,6.25)(7,7)
\drawline[AHnb=0](9,6.25)(9,7)
\drawline[AHnb=0](11,6.25)(11,7)
\drawline[AHnb=0](13,6.25)(13,7)



\drawline[AHnb=0](2,2.5)(3,2)(4,2.5)(5,2)(6,2.5)(7,2)(8,2.5)(9,2)(10,2.5)(11,2)

\drawline[AHnb=0](3,4)(4,3.5)(5,4)(6,3.5)(7,4)(8,3.5)(9,4)(10,3.5)(11,4)(12,3.5)

\drawline[AHnb=0](4,5.5)(4,6.5)(5,7)(6,6.5)(7,7)(8,6.5)(9,7)(10,6.5)(11,7)(12,6.5)(13,7)(14,6.5)(14,5.5)(13,5)


\put(.5,.6){\scriptsize {\tiny $a_{11}$}}
\put(1.7,.2){\scriptsize {\tiny $a_{12}$}}
\put(2.7,.5){\scriptsize {\tiny $a_{13}$}}
\put(9.7,0){\scriptsize {\tiny $a_{1i}$}}
\put(10.7,.5){\scriptsize {\tiny $a_{11}$}}

\put(0,1.2){\scriptsize {\tiny $a_{21}$}}
\put(11.2,1.2){\scriptsize {\tiny $a_{21}$}}

\put(0,2){\scriptsize {\tiny $a_{31}$}}
\put(2.2,2.45){\scriptsize {\tiny $a_{32}$}}
\put(11.1,1.8){\scriptsize {\tiny $a_{31}$}}
\put(12.1,2.2){\scriptsize {\tiny $a_{32}$}}

\put(1,2.8){\scriptsize {\tiny $a_{41}$}}
\put(12.1,2.7){\scriptsize {\tiny $a_{41}$}}

\put(4,7.5){\scriptsize {\tiny $a_{j1}$}}
\put(15.2,7.5){\scriptsize {\tiny $a_{j1}$}}

\put(-1.75,.5){\scriptsize {\tiny $Q_1$}}
\put(-1.25,1.25){\scriptsize {\tiny $Q_2$}}
\put(-.95,2){\scriptsize {\tiny $Q_{3}$}}
\put(-.45,2.75){\scriptsize {\tiny $Q_{4}$}}
\put(1.25,7.25){\scriptsize {\tiny $Q_{j}$}}
\put(1.75,8.15){\scriptsize {\tiny $Q_{1}$}}

\put(3.25,8.65){\scriptsize {\tiny $a_{1(k+1)}$}}
\put(5.35,9.15){\scriptsize {\tiny $a_{1(k+2)}$}}
\put(13.8,9.25){\scriptsize {\tiny $a_{1k}$}}
\put(14.75,8.75){\scriptsize {\tiny $a_{1(k+1)}$}}

\put(4.5,0){\scriptsize $\ldots$}
\put(7,0){\scriptsize $\ldots$}
\put(8,9.25){\scriptsize $\ldots$}
\put(11,9.25){\scriptsize $\ldots$}

\put(.7,5.6){\scriptsize $.$}
\put(.6,5.4){\scriptsize $.$}
\put(0.5,5.2){\scriptsize $.$}

\drawpolygon[fillcolor=black](1,1.5)(3,2)(3,1.5)(4,1.75)(4,2.5)(5,2)(5,1.5)(6,1.75)(6,2.5)(7,2)(7,1.5)(8,1.75)(8,2.5)(9,2)(9,1.5)(10,1.75)(10,3)(11,3.25)(11,4)(10,3.5)(9,3.25)(9,4)(8,3.5)(8,3)(7,3.25)(7,4)(6,3.5)(6,3)(5,3.25)(4,3)(4,3.5)(5,4)(5,5)(6,4.75)(6,5.5)(7,5)(7,4.5)(8,4.75)(8,5.5)(9,5)(9,4.5)(10,4.75)(10,5.5)(11,5)(11,4.5)(12,4.75)(12,6)(13,6.25)(13,7)(12,6.5)(11,6.25)(11,7)(10,6.5)(10,6)(9,6.25)(9,7)(8,6.5)(8,6)(7,6.25)(6,6)(6,6.5)(7,7)(7,6.9)(6.1,6.5)(6.1,6.1)(7,6.3)(7.9,6.1)(7.9,6.6)(9.1,7.1)(9.1,6.26)(9.9,6.1)(9.9,6.6)(11.1,7.1)(11.1,6.35)(12,6.6)(13.1,7.1)(13.1,6.24)(12.1,5.9)(12.1,4.65)(10.9,4.4)(10.9,4.9)(10.1,5.3)(10.1,4.65)(8.9,4.4)(8.9,4.9)(8.1,5.3)(8.1,4.65)(6.9,4.4)(6.9,4.9)(6.1,5.3)(6.1,4.65)(5.1,4.9)(5.1,3.9)(4.1,3.4)(4.1,3.1)(5,3.4)(5.9,3.1)(5.9,3.6)(7.1,4.1)(7.1,3.28)(7.9,3.1)(7.9,3.6)(9.1,4.1)(9.1,3.35)(9.9,3.6)(11.1,4.1)(11.1,3.1)(10.1,2.9)(10.1,1.6)(8.9,1.4)(8.9,2)(8.1,2.3)(8.1,1.6)(6.9,1.4)(6.9,2)(6.1,2.3)(6.1,1.6)(4.9,1.4)(4.9,2)(4.1,2.3)(4.1,1.6)(2.9,1.4)(2.9,1.9)(1,1.4)

\drawpolygon[fillcolor=black](5,7)(4,6.5)(5,6.25)(4,6)(4,5.5)(3,5)(4,4.75)(3,4.5)(3,4)(2,3.5)(3,3.25)(2,3)(2,2.5)(1,2)(1,1.4)(.9,1.4)(.9,2.1)(1.9,2.6)(1.9,3.1)(2.7,3.25)(1.8,3.48)(2.9,4.1)(2.9,4.6)(3.7,4.75)(2.8,4.98)(3.9,5.6)(3.9,6.1)(4.7,6.25)(3.8,6.5)(5,7.1)

\drawpolygon[fillcolor=black](15,8)(14,8.5)(14,7.75)(13,7.5)(13,8)(12,8.5)(12,7.75)(11,7.5)(11,8)(10,8.5)(10,7.75)(9,7.5)(9,8)(8,8.5)(8,7.75)(7,8)(7,7)(6.9,7)(6.9,8.1)(7.9,7.85)(7.9,8.6)(9.1,8.1)(9.1,7.6)(9.9,7.76)(9.9,8.6)(11.1,8.1)(11.1,7.6)(11.9,7.85)(11.91,8.6)(13.1,8.1)(13.1,7.6)(13.9,7.85)(13.9,8.6)(15,8.1)

\drawpolygon[fillcolor=black](5,7)(5,7.5)(6,7.75)(6,8.5)(5,8)(5,8.1)(6.1,8.65)(6.1,7.7)(5.1,7.4)(5.1,7)

\put(2.5,-1){\scriptsize {\tiny {\bf Figure 4.9.4:} $M(i,j=4m+2,k)$ }}
\end{picture}

\begin{picture}(0,0)(-72,74.75)
\setlength{\unitlength}{4.65mm}


\drawpolygon(1,1)(2,.5)(3,1)(4,.5)(5,1)(6,.5)(7,1)(8,.5)(9,1)(10,.5)(11,1)(11,2)(12,2.5)(12,3.5)(13,4)(13,5)(12,5.5)(11,5)(10,5.5)(9,5)(8,5.5)(7,5)(6,5.5)(5,5)(4,5.5)(3,5)(3,4)(2,3.5)(2,2.5)(1,2)

\drawline[AHnb=0](1,1.5)(3,2)(3,1)
\drawline[AHnb=0](3,1.5)(5,2)(5,1)
\drawline[AHnb=0](5,1.5)(7,2)(7,1)
\drawline[AHnb=0](7,1.5)(9,2)(9,1)
\drawline[AHnb=0](9,1.5)(11,2)(11,1)

\drawline[AHnb=0](2,3)(4,3.5)(4,1.75)
\drawline[AHnb=0](4,3)(6,3.5)(6,1.75)
\drawline[AHnb=0](6,3)(8,3.5)(8,1.75)
\drawline[AHnb=0](8,3)(10,3.5)(10,1.75)
\drawline[AHnb=0](10,3)(12,3.5)

\drawline[AHnb=0](3,4.5)(5,5)(5,3.25)
\drawline[AHnb=0](5,4.5)(7,5)(7,3.25)
\drawline[AHnb=0](7,4.5)(9,5)(9,3.25)
\drawline[AHnb=0](9,4.5)(11,5)(11,3.25)
\drawline[AHnb=0](11,4.5)(13,5)

\drawline[AHnb=0](4,6)(6,6.5)(6,4.75)
\drawline[AHnb=0](6,6)(8,6.5)(8,4.75)
\drawline[AHnb=0](8,6)(10,6.5)(10,4.75)
\drawline[AHnb=0](10,6)(12,6.5)(12,4.75)
\drawline[AHnb=0](12,6)(14,6.5)


\drawline[AHnb=0](1,2)(3,1.5)
\drawline[AHnb=0](3,2)(5,1.5)
\drawline[AHnb=0](5,2)(7,1.5)
\drawline[AHnb=0](7,2)(9,1.5)
\drawline[AHnb=0](9,2)(11,1.5)

\drawline[AHnb=0](2,3.5)(4,3)
\drawline[AHnb=0](4,3.5)(6,3)
\drawline[AHnb=0](6,3.5)(8,3)
\drawline[AHnb=0](8,3.5)(10,3)
\drawline[AHnb=0](10,3.5)(12,3)

\drawline[AHnb=0](3,5)(5,4.5)
\drawline[AHnb=0](5,5)(7,4.5)
\drawline[AHnb=0](7,5)(9,4.5)
\drawline[AHnb=0](9,5)(11,4.5)
\drawline[AHnb=0](11,5)(13,4.5)

\drawline[AHnb=0](4,6.5)(6,6)
\drawline[AHnb=0](6,6.5)(8,6)
\drawline[AHnb=0](8,6.5)(10,6)
\drawline[AHnb=0](10,6.5)(12,6)
\drawline[AHnb=0](12,6.5)(14,6)

\drawline[AHnb=0](2,1.75)(2,2.5)
\drawline[AHnb=0](3,3.25)(3,4)
\drawline[AHnb=0](4,4.75)(4,5.5)
\drawline[AHnb=0](6,4.75)(6,5.5)
\drawline[AHnb=0](8,4.75)(8,5.5)
\drawline[AHnb=0](10,4.75)(10,5.5)
\drawline[AHnb=0](12,4.75)(12,5.5)

\drawline[AHnb=0](5,6.25)(5,7)
\drawline[AHnb=0](7,6.25)(7,7)
\drawline[AHnb=0](9,6.25)(9,7)
\drawline[AHnb=0](11,6.25)(11,7)
\drawline[AHnb=0](13,6.25)(13,7)



\drawline[AHnb=0](2,2.5)(3,2)(4,2.5)(5,2)(6,2.5)(7,2)(8,2.5)(9,2)(10,2.5)(11,2)

\drawline[AHnb=0](3,4)(4,3.5)(5,4)(6,3.5)(7,4)(8,3.5)(9,4)(10,3.5)(11,4)(12,3.5)

\drawline[AHnb=0](4,5.5)(4,6.5)(5,7)(6,6.5)(7,7)(8,6.5)(9,7)(10,6.5)(11,7)(12,6.5)(13,7)(14,6.5)(14,5.5)(13,5)


\drawpolygon[fillcolor=black](1,1.5)(3,2)(3,1.5)(4,1.75)(4,2.5)(5,2)(5,1.5)(6,1.75)(6,2.5)(7,2)(7,1.5)(8,1.75)(8,2.5)(9,2)(9,1.5)(10,1.75)(10,3)(11,3.25)(11,4)(10,3.5)(9,3.25)(9,4)(8,3.5)(8,3)(7,3.25)(7,4)(6,3.5)(6,3)(5,3.25)(4,3)(4,3.5)(5,4)(5,5)(6,4.75)(6,5.5)(7,5)(7,4.5)(8,4.75)(8,5.5)(9,5)(9,4.5)(10,4.75)(10,5.5)(11,5)(11,4.5)(12,4.75)(12,6)(13,6.25)(13,7)(12,6.5)(11,6.25)(11,7)(10,6.5)(10,6)(9,6.25)(9,7)(8,6.5)(8,6)(7,6.25)(7,7)(6,6.5)(6,6)(5,6.25)(5,7)(4,6.5)(4,5.5)(3.9,5.5)(3.9,6.6)(5.1,7.1)(5.1,6.35)(5.9,6.1)(5.9,6.6)(7.1,7.1)(7.1,6.35)(7.9,6.1)(7.9,6.6)(9.1,7.1)(9.1,6.26)(9.9,6.1)(9.9,6.6)(11.1,7.1)(11.1,6.35)(12,6.6)(13.1,7.1)(13.1,6.24)(12.1,5.9)(12.1,4.65)(10.9,4.4)(10.9,4.9)(10.1,5.3)(10.1,4.65)(8.9,4.4)(8.9,4.9)(8.1,5.3)(8.1,4.65)(6.9,4.4)(6.9,4.9)(6.1,5.3)(6.1,4.65)(5.1,4.9)(5.1,3.9)(4.1,3.4)(4.1,3.1)(5,3.4)(5.9,3.1)(5.9,3.6)(7.1,4.1)(7.1,3.28)(7.9,3.1)(7.9,3.6)(9.1,4.1)(9.1,3.35)(9.9,3.6)(11.1,4.1)(11.1,3.1)(10.1,2.9)(10.1,1.6)(8.9,1.4)(8.9,2)(8.1,2.3)(8.1,1.6)(6.9,1.4)(6.9,2)(6.1,2.3)(6.1,1.6)(4.9,1.4)(4.9,2)(4.1,2.3)(4.1,1.6)(2.9,1.4)(2.9,1.9)(1,1.4)

\drawpolygon[fillcolor=black](4,5.5)(3,5)(4,4.75)(3,4.5)(3,4)(2,3.5)(3,3.25)(2,3)(2,2.5)(1,2)(1,1.4)(.9,1.4)(.9,2.1)(1.9,2.6)(1.9,3.1)(2.75,3.25)(1.8,3.47)(2.9,4.1)(2.9,4.6)(3.75,4.75)(2.85,4.98)(4,5.6)

\put(.5,.6){\scriptsize {\tiny $a_{11}$}}
\put(1.7,.2){\scriptsize {\tiny $a_{12}$}}
\put(2.7,.5){\scriptsize {\tiny $a_{13}$}}
\put(9.7,0){\scriptsize {\tiny $a_{1i}$}}
\put(10.7,.5){\scriptsize {\tiny $a_{11}$}}

\put(0,1.2){\scriptsize {\tiny $a_{21}$}}
\put(11.2,1.2){\scriptsize {\tiny $a_{21}$}}

\put(0,2){\scriptsize {\tiny $a_{31}$}}
\put(2.2,2.45){\scriptsize {\tiny $a_{32}$}}
\put(11.1,1.8){\scriptsize {\tiny $a_{31}$}}
\put(12.1,2.2){\scriptsize {\tiny $a_{32}$}}

\put(1,2.8){\scriptsize {\tiny $a_{41}$}}
\put(12.1,2.7){\scriptsize {\tiny $a_{41}$}}

\put(3,5.9){\scriptsize {\tiny $a_{j1}$}}
\put(14.1,5.9){\scriptsize {\tiny $a_{j1}$}}

\put(2.5,7){\scriptsize {\tiny $a_{1(k+1)}$}}
\put(4.25,7.5){\scriptsize {\tiny $a_{1(k+2)}$}}
\put(12.8,7.5){\scriptsize {\tiny $a_{1k}$}}
\put(13.75,7){\scriptsize {\tiny $a_{1(k+1)}$}}

\put(-1.75,.5){\scriptsize {\tiny $Q_1$}}
\put(-1.25,1.25){\scriptsize {\tiny $Q_2$}}
\put(-.95,2){\scriptsize {\tiny $Q_{3}$}}
\put(-.5,2.75){\scriptsize {\tiny $Q_{4}$}}
\put(1.5,6){\scriptsize {\tiny $Q_{j}$}}
\put(1.75,6.5){\scriptsize {\tiny $Q_{1}$}}

\put(4,0){\scriptsize $\ldots$}
\put(7,0){\scriptsize $\ldots$}
\put(7.5,7.5){\scriptsize $\ldots$}
\put(10.5,7.5){\scriptsize $\ldots$}

\put(1,4.8){\scriptsize $.$}
\put(.9,4.6){\scriptsize $.$}
\put(0.8,4.4){\scriptsize $.$}

\put(2.5,-1){\scriptsize {\tiny {\bf Figure 4.9.5:} $M(i,j=4m+4,k)$ }}

\end{picture}

\vspace{8.25cm}

As in previous sections, we construct $M(i,j,k)$ representation of $M$ by first cutting $M$ along a cycle of type $I_1$ and then cutting it along a cycle of type $I_2$ where the beginning adjacent face to the cycle $I_1$ is a hexagon. This gives the following lemma.

\begin{lem}\label{l3.9.1}
	
	A DSEM $M$ of type $[3^4.6:3^2.6^2]$ admits an $M(i,j,k)$-representation iff the following holds: $(i)$ $i \geq 6$ and $i, j$ even, $(ii)$ number of vertices in $M(i,j,k) = ij \geq 12 $, $(iii)$ if $j = 2$ then $k \in \{2r+3: 0 \leq r \leq (i-6)/2\}$, and if $j \geq 4$ then $k \in \{2r+1: 0 \leq r \leq (i-2)/2\}$.
	
\end{lem}

\noindent{\bf Proof.} Let $M$ be a DSEM of type $[3^4.6:3^2.6^2]$ with $n$ vertices. Then its $M(i,j,k)$ has $j$  disjoint horizontal cycles  of $I_{1}$ type having length $i$. Thus the number of vertices in $M$ is $n = ij$.

If $j=1$, then $M(i,1,k)$ has no vertex with face-sequence $(3^4,6)$. It is not possible. Therefore, $j \geq 2 $. If $j \geq 2 $ and $j$ is not an even integer, then no vertex in the base horizontal cycle follows the face-sequence $(3 ^ 4,6)$ after identifying the boundaries of $M(i, j , k)$. Therefore, $j$ is even.

If $j$ is an even integer and $ i < 6 $, then we get some vertex in $M(i,1,k)$ whose link can not be constructed. So, $i \geq 6$. Also, if $i \geq 6$ and $i$ is not an even integer, then similarly as above, there is no vertex in the horizontal base cycle having face-sequences $(3^4,6)$ and $(3^2,6^2)$. Which is not possible. So, $i \geq 6$ and $i, j$ are even integers. Thus, $ n = ij \geq 12$.

If $j = 2$ and $ k \in \{r:0 \leq r \leq i-1\} \setminus (\{2r+3: 0 \leq r \leq (i-6)/2\}$, then some vertices in the lower horizontal cycle do not follow the face-sequences $(3^4,6)$ and $(3^2,6^2)$. Which is a contradiction. So, for $j = 2$, we get $ k \in \{2r+3: 0 \leq r \leq (i-6)/2\}$. Similarly, as above, we see that if $j \geq 4$, then $ k \in \{2r+1: 0 \leq r \leq (i-2)/2\}$. This completes the proof. \hfill $\Box$

%
%
%

%
%
%
%

\subsection{DSEMs of type $[3^2.4.3.4:3.4.6.4]$} \label{s3.10}

Let $M$ be a DSEM of type $[3^2.4.3.4:3.4.6.4]$. We consider a path in $M$ as follows. A path $P_{1} = P( \ldots, y_{i},y_{i+1},z_{i},z_{i+1}, \ldots)$ in $M$, say of type $J_{1}$, indicated by thick black or green paths, shown in Figure 4.10.1. The inner vertices $y_i$'s and $z_i$'s have the face-sequences $(3,4,6,4)$ and $(3^2,4,3,4)$ respectively.

\begin{picture}(0,0)(8,42)
\setlength{\unitlength}{4.5mm}


\drawline[AHnb=0](2,4.5)(3,5)(4,5)(4,4)(5,3.5)(6,3.5)(7,4)(8,4)(9,3.5)(10,3.5)(11,4)(12,4)(13,3.5)(14,3.5)(14.75,3.75)

\drawline[AHnb=0](4.5,6)(4,5)(5,4.5)(6,4.5)(7,5)(8,5)(9,4.5)(10,4.5)(11,5)(12,5)(13,4.5)(14,4.5)(14.75,4.75)

\drawline[AHnb=0](2,7.5)(3,7)(4,7)(5,7.5)(6,7.5)(7,7)(8,7)(9,7.5)(10,7.5)(11,7)(12,7)(13,7.5)(14,7.5)(14,8.5)(14.5,8.25)

\drawline[AHnb=0](2,8.5)(3,8)(4,8)(5,8.5)(6,8.5)(7,8)(8,8)(9,8.5)(10,8.5)(11,8)(12,8)(13,8.5)(14,8.5)(14,7.5)(14.5,7.25)

\drawline[AHnb=0](14,7.5)(14.5,7.75)
\drawline[AHnb=0](3,8)(3,7)
\drawline[AHnb=0](4,8)(4,7)
\drawline[AHnb=0](5,8.5)(5,7.5)(5.5,6.5)

\drawline[AHnb=0](10.5,2.5)(11.5,2)(11.5,3)(12.5,2.5)(11.5,2)

\drawline[AHnb=0](12.5,6)(13.5,5.5)(13.5,6.5)(14.5,6)(13.5,5.5)(13,4.5)(14,4.5)(13.5,5.5)

\drawline[AHnb=0](2.5,2.5)(2,3.5)(3,4)(3.5,3)

\drawline[AHnb=0](2,5.75)(2.5,6)(2,6.25)

\drawline[AHnb=0](5,3.5)(5,4.5)(5.5,5.5)(6,4.5)(6,3.5)
\drawline[AHnb=0](5.5,6.5)(6,7.5)(6,8.5)

\drawline[AHnb=0](9,3.5)(9,4.5)(9.5,5.5)(10,4.5)(10,3.5)
\drawline[AHnb=0](9,8.5)(9,7.5)(9.5,6.5)(10,7.5)(10,8.5)
\drawline[AHnb=0](7,5)(7,4)(7.5,3)(8,4)(8,5)

\drawline[AHnb=0](8.5,6)(9.5,5.5)(10.5,6)(9.5,6.5)(8.5,6)
\drawline[AHnb=0](4.5,6)(5.5,5.5)(6.5,6)(5.5,6.5)(4.5,6)


\drawline[AHnb=0](3,4)(3,5)(2.5,6)(3,7)(4,7)(4.5,6)(4,5)(4,4)(3,4)(3,5)(4,5)

\drawline[AHnb=0](7,4)(7,5)(6.5,6)(7,7)(8,7)(8.5,6)(8,5)(8,4)

\drawline[AHnb=0](11,4)(11,5)(10.5,6)(11,7)(12,7)(12.5,6)(12,5)(12,4)

\drawline[AHnb=0](5,1.5)(4.5,2.5)(5,3.5)(6,3.5)(6.5,2.5)(6,1.5)(5,1.5)

\drawline[AHnb=0](9,1.5)(8.5,2.5)(9,3.5)(10,3.5)(10.5,2.5)(10,1.5)(9,1.5)

\drawline[AHnb=0](13,1.5)(12.5,2.5)(13,3.5)(14,3.5)(14.5,2.5)(14,1.5)(13,1.5)

\drawline[AHnb=0](3.5,2)(3.5,3)(2.5,2.5)(2,1.5)

\drawline[AHnb=0](2.5,2.5)(3.5,2)(4.5,2.5)(3.5,3)(4,4)

\drawline[AHnb=0](6.5,2.5)(7.5,2)(8.5,2.5)(7.5,3)(6.5,2.5)

\drawline[AHnb=0](7.5,2)(7.5,3)
\drawline[AHnb=0](7,8)(7,7)
\drawline[AHnb=0](8,8)(8,7)
\drawline[AHnb=0](11,8)(11,7)
\drawline[AHnb=0](12,8)(12,7)
\drawline[AHnb=0](13,8.5)(13,7.5)(13.5,6.5)
\drawline[AHnb=0](12.5,6)(13.5,6.5)(14,7.5)
\drawline[AHnb=0](10.5,2.5)(11.5,3)(12,4)(11,4)(11.5,3)

\drawline[AHnb=0](13,3.5)(13,4.5)
\drawline[AHnb=0](14,3.5)(14,4.5)

\drawline[AHnb=0](12,4)(13,4.5)
\drawline[AHnb=0](8,4)(9,4.5)
\drawline[AHnb=0](4,4)(5,4.5)

\drawline[AHnb=0](12,8)(13,7.5)
\drawline[AHnb=0](8,8)(9,7.5)
\drawline[AHnb=0](4,8)(5,7.5)

\drawline[AHnb=0](14.75,4)(14,4.5)
\drawline[AHnb=0](11,4)(10,4.5)
\drawline[AHnb=0](7,4)(6,4.5)
\drawline[AHnb=0](3,4)(2,4.5)

\drawline[AHnb=0](11,8)(10,7.5)
\drawline[AHnb=0](7,8)(6,7.5)
\drawline[AHnb=0](3,8)(2,7.5)

\drawline[AHnb=0](9.5,5.5)(9.5,6.5)
\drawline[AHnb=0](5.5,5.5)(5.5,6.5)

\drawline[AHnb=0](3.5,2)(3.25,1.5)
\drawline[AHnb=0](3.5,2)(3.75,1.5)

\drawline[AHnb=0](7.5,2)(7.25,1.5)
\drawline[AHnb=0](7.5,2)(7.75,1.5)

\drawline[AHnb=0](11.5,2)(11.25,1.5)
\drawline[AHnb=0](11.5,2)(11.75,1.5)


\drawline[AHnb=0](4.5,1.25)(5,1.5)(5,1)
\drawline[AHnb=0](6.5,1.25)(6,1.5)(6,1)

\drawline[AHnb=0](8.5,1.25)(9,1.5)(9,1)
\drawline[AHnb=0](10.5,1.25)(10,1.5)(10,1)

\drawline[AHnb=0](12.5,1.25)(13,1.5)(13,1)
\drawline[AHnb=0](14.5,1.25)(14,1.5)(14,1)


\drawline[AHnb=0](3,8)(3.25,8.5)
\drawline[AHnb=0](7,8)(7.25,8.5)
\drawline[AHnb=0](11,8)(11.25,8.5)

\drawline[AHnb=0](2,8.5)(2.25,9)
\drawline[AHnb=0](6,8.5)(6.25,9)
\drawline[AHnb=0](10,8.5)(10.25,9)
\drawline[AHnb=0](14,8.5)(14.25,9)

\drawline[AHnb=0](4,8)(3.75,8.5)
\drawline[AHnb=0](8,8)(7.75,8.5)
\drawline[AHnb=0](12,8)(11.75,8.5)

\drawline[AHnb=0](5,8.5)(4.75,9)
\drawline[AHnb=0](9,8.5)(8.75,9)
\drawline[AHnb=0](13,8.5)(12.75,9)

\drawpolygon[fillcolor=black](2,4.5)(3,5)(4,5)(5,4.5)(6,4.5)(7,5)(8,5)(9,4.5)(10,4.5)(11,5)(12,5)(13,4.5)(14,4.5)(14.75,4.75)(14.75,4.85)(14,4.6)(13,4.6)(12,5.1)(11,5.1)(10,4.6)(9,4.6)(8,5.1)(7,5.1)(6,4.6)(5,4.6)(4,5.1)(3,5.1)(2,4.6)

\drawpolygon[fillcolor=black](2,3.5)(3,4)(4,4)(5,3.5)(6,3.5)(7,4)(8,4)(9,3.5)(10,3.5)(11,4)(12,4)(13,3.5)(14,3.5)(14.75,3.75)(14.75,3.85)(14,3.6)(13,3.6)(12,4.1)(11,4.1)(10,3.6)(9,3.6)(8,4.1)(7,4.1)(6,3.6)(5,3.6)(4,4.1)(3,4.1)(2,3.6)

\drawpolygon[fillcolor=black](2,7.5)(3,7)(4,7)(5,7.5)(6,7.5)(7,7)(8,7)(9,7.5)(10,7.5)(11,7)(12,7)(13,7.5)(14,7.5)(14.75,7.75)(14.75,7.85)(14,7.6)(13,7.6)(12,7.1)(11,7.1)(10,7.6)(9,7.6)(8,7.1)(7,7.1)(6,7.6)(5,7.6)(4,7.1)(3,7.1)(2,7.6)

\drawpolygon[fillcolor=black](2,8.5)(3,8)(4,8)(5,8.5)(6,8.5)(7,8)(8,8)(9,8.5)(10,8.5)(11,8)(12,8)(13,8.5)(14,8.5)(14.75,8.75)(14.75,8.85)(14,8.6)(13,8.6)(12,8.1)(11,8.1)(10,8.6)(9,8.6)(8,8.1)(7,8.1)(6,8.6)(5,8.6)(4,8.1)(3,8.1)(2,8.6)

\drawpolygon[fillcolor=green](14.25,9)(14,8.5)(14,7.5)(13.5,6.5)(12.5,6)(12,5)(12,4)(11.5,3)(10.5,2.5)(10,1.5)(10,1)(10.1,1)(10.1,1.4)(10.6,2.4)(11.6,2.9)(12.1,3.9)(12.1,4.9)(12.6,5.9)(13.6,6.4)(14.1,7.5)(14.1,8.5)(14.35,9)

\drawpolygon[fillcolor=green](10.25,9)(10,8.5)(10,7.5)(9.5,6.5)(8.5,6)(8,5)(8,4)(7.5,3)(6.5,2.5)(6,1.5)(6,1)(6.1,1)(6.1,1.4)(6.6,2.4)(7.6,2.9)(8.1,3.9)(8.1,4.9)(8.6,5.9)(9.6,6.4)(10.1,7.5)(10.1,8.5)(10.35,9)

\drawpolygon[fillcolor=green](6.25,9)(6,8.5)(6,7.5)(5.5,6.5)(4.5,6)(4,5)(4,4)(3.5,3)(2.5,2.5)(2,1.5)(2,1)(2.1,1)(2.1,1.4)(2.6,2.4)(3.6,2.9)(4.1,3.9)(4.1,4.9)(4.6,5.9)(5.6,6.4)(6.1,7.5)(6.1,8.5)(6.35,9)

\drawpolygon[fillcolor=green](11.25,8.5)(11,8)(11,7)(10.5,6)(9.5,5.5)(9,4.5)(9,3.5)(8.5,2.5)(7.5,2)(7,1)(7,.5)(7.1,.5)(7.1,.9)(7.6,1.9)(8.6,2.4)(9.1,3.4)(9.1,4.4)(9.6,5.4)(10.6,5.9)(11.1,7)(11.1,8)(11.35,8.5)

\drawpolygon[fillcolor=green](7.25,8.5)(7,8)(7,7)(6.5,6)(5.5,5.5)(5,4.5)(5,3.5)(4.5,2.5)(3.5,2)(3,1)(3,.5)(3.1,.5)(3.1,.9)(3.6,1.9)(4.6,2.4)(5.1,3.4)(5.1,4.4)(5.6,5.4)(6.6,5.9)(7.1,7)(7.1,8)(7.35,8.5)

\drawpolygon[fillcolor=green](15,8)(15,7)(14.5,6)(13.5,5.5)(13,4.5)(13,3.5)(12.5,2.5)(11.5,2)(11.6,1.9)(12.6,2.4)(13.1,3.4)(13.1,4.4)(13.6,5.4)(14.6,5.9)(15.1,7)(15.1,8)


\put(4.3,-.5){\scriptsize {\tiny {\bf Figure 4.10.1:} Path of type $J_1$ }} 

\end{picture}

\begin{picture}(0,0)(-61.5,35)
\setlength{\unitlength}{5mm}


\drawline[AHnb=0](4,5)(4,4)(5,3.5)(6,3.5)(7,4)(8,4)(9,3.5)(10,3.5)(11,4)(12,4)(13,3.5)(14,3.5)(15,4)(16,4)(16,5)(15,5)(14,4.5)(13,4.5)(12,5)

\drawline[AHnb=0](6,8.5)(6,7.5)(7,7)(8,7)(9,7.5)(10,7.5)(11,7)(12,7)(13,7.5)(14,7.5)(15,7)(16,7)(17,7.5)(18,7.5)(18,8.5)(17,8.5)(16,8)(15,8)(14,8.5)

\drawline[AHnb=0](6,8.5)(6,7.5)(7,7)(8,7)(9,7.5)(10,7.5)(11,7)(12,7)(13,7.5)(14,7.5)(15,7)(16,7)(17,7.5)(18,7.5)(18,8.5)(17,8.5)(16,8)(15,8)(14,8.5)

\drawline[AHnb=0](4.5,6)(4,5)(5,4.5)(6,4.5)(7,5)(8,5)(9,4.5)(10,4.5)(11,5)(12,5)

\drawline[AHnb=0](6,7.5)(7,7)(8,7)(9,7.5)(10,7.5)(11,7)(12,7)(13,7.5)(14,7.5)(14,8.5)

\drawline[AHnb=0](6,8.5)(7,8)(8,8)(9,8.5)(10,8.5)(11,8)(12,8)(13,8.5)(14,8.5)

\drawline[AHnb=0](5,3.5)(5,4.5)(5.5,5.5)(6,4.5)(6,3.5)
\drawline[AHnb=0](5.5,6.5)(6,7.5)(6,8.5)

\drawline[AHnb=0](9,3.5)(9,4.5)(9.5,5.5)(10,4.5)(10,3.5)
\drawline[AHnb=0](9,8.5)(9,7.5)(9.5,6.5)(10,7.5)(10,8.5)
\drawline[AHnb=0](7,5)(7,4)(7.5,3)(8,4)(8,5)
\drawline[AHnb=0](13,3.5)(13,4.5)(13.5,5.5)(14,4.5)(14,3.5)
\drawline[AHnb=0](13,8.5)(13,7.5)(13.5,6.5)(14,7.5)(14,8.5)
\drawline[AHnb=0](11,5)(11,4)(11.5,3)(12,4)(12,5)

\drawline[AHnb=0](8.5,6)(9.5,5.5)(10.5,6)(9.5,6.5)(8.5,6)
\drawline[AHnb=0](4.5,6)(5.5,5.5)(6.5,6)(5.5,6.5)(4.5,6)
\drawline[AHnb=0](12.5,6)(13.5,5.5)(14.5,6)(13.5,6.5)(12.5,6)
\drawline[AHnb=0](10.5,2.5)(11.5,2)(12.5,2.5)(11.5,3)(10.5,2.5)

\drawline[AHnb=0](17,7.5)(17,8.5)
\drawline[AHnb=0](16,7)(16,8)(17,7.5)
\drawline[AHnb=0](15,7)(15,8)(14,7.5)

\drawline[AHnb=0](11.5,2)(12,1)(11,1)(11.5,2)(11.5,3)

\drawline[AHnb=0](12,1)(13,1.5)
\drawline[AHnb=0](11,1)(10,1.5)

\drawline[AHnb=0](13.5,5.5)(13.5,6.5)
\drawline[AHnb=0](7,4)(7,5)(6.5,6)(7,7)(8,7)(8.5,6)(8,5)(8,4)

\drawline[AHnb=0](11,4)(11,5)(10.5,6)(11,7)(12,7)(12.5,6)(12,5)(12,4)

\drawline[AHnb=0](15,4)(15,5)(14.5,6)(15,7)(16,7)(16.5,6)(16,5)(16,4)

\drawline[AHnb=0](5,1.5)(4.5,2.5)(5,3.5)(6,3.5)(6.5,2.5)(6,1.5)(5,1.5)
\drawline[AHnb=0](9,1.5)(8.5,2.5)(9,3.5)(10,3.5)(10.5,2.5)(10,1.5)(9,1.5)
\drawline[AHnb=0](13,1.5)(12.5,2.5)(13,3.5)(14,3.5)(14.5,2.5)(14,1.5)(13,1.5)

\drawline[AHnb=0](3.5,2)(3.5,3)(2.5,2.5)(2,1.5)(3,1)(3.5,2)(4,1)(5,1.5)(6,1.5)(7,1)(7.5,2)(8,1)(9,1.5)

\drawline[AHnb=0](8,4)(9,4.5)
\drawline[AHnb=0](4,4)(5,4.5)

\drawline[AHnb=0](12,8)(13,7.5)
\drawline[AHnb=0](8,8)(9,7.5)

\drawline[AHnb=0](11,4)(10,4.5)
\drawline[AHnb=0](7,4)(6,4.5)

\drawline[AHnb=0](11,8)(10,7.5)
\drawline[AHnb=0](7,8)(6,7.5)

\drawline[AHnb=0](9.5,5.5)(9.5,6.5)
\drawline[AHnb=0](5.5,5.5)(5.5,6.5)

\drawline[AHnb=0](2.5,2.5)(3.5,2)(4.5,2.5)(3.5,3)(4,4)

\drawline[AHnb=0](6.5,2.5)(7.5,2)(8.5,2.5)(7.5,3)(6.5,2.5)

\drawline[AHnb=0](3,1)(4,1)
\drawline[AHnb=0](7,1)(8,1)
\drawline[AHnb=0](7.5,2)(7.5,3)
\drawline[AHnb=0](9.5,5.5)(9.5,6.5)
\drawline[AHnb=0](7,8)(7,7)
\drawline[AHnb=0](8,8)(8,7)
\drawline[AHnb=0](11,8)(11,7)
\drawline[AHnb=0](12,8)(12,7)
\drawline[AHnb=0](13,8.5)(13,7.5)(13.5,6.5)
\drawline[AHnb=0](12.5,6)(13.5,6.5)(14,7.5)
\drawline[AHnb=0](10.5,2.5)(11.5,3)(12,4)(11,4)(11.5,3)

\drawline[AHnb=0](15.5,3)(15,4)
\drawline[AHnb=0](17.5,6.5)(17,7.5)

\drawpolygon[fillcolor=black](17,7.5)(18,7.5)(17.5,6.5)(16.5,6)(16,5)(16,4)(15.5,3)(14.5,2.5)(14,1.5)(14.1,1.4)(14.6,2.4)(15.6,2.9)(16.1,3.9)(16.1,4.9)(16.6,5.9)(17.6,6.4)(18.15,7.6)(17,7.6)

\drawpolygon[fillcolor=black](14,1.5)(13,1.5)(12,1)(11,1)(11.5,2)(10.5,2.5)(10,1.5)(9,1.5)(8,1)(7,1)(7.5,2)(6.5,2.5)(6,1.5)(5,1.5)(4,1)(3,1)(3.5,2)(4.5,2.5)(5,3.5)(6,3.5)(7,4)(8,4)(7.5,3)(8.5,2.5)(9,3.5)(10,3.5)(11,4)(12,4)(11.5,3)(12.5,2.5)(13,3.5)(14,3.5)(15,4)(15,5)(14,4.5)(13,4.5)(13.5,5.5)(12.5,6)(12,5)(11,5)(10,4.5)(9,4.5)(9.5,5.5)(8.5,6)(8,5)(7,5)(6,4.5)(5,4.5)(5.5,5.5)(6.5,6)(7,7)(8,7)(9,7.5)(10,7.5)(9.5,6.5)(10.5,6)(11,7)(12,7)(13,7.5)(14,7.5)(13.5,6.5)(14.5,6)(15,7)(16,7)(17,7.5)(16.9,7.6)(16,7.1)(14.9,7.1)(14.4,6.15)(13.65,6.6)(14.15,7.6)(12.9,7.6)(11.9,7.1)(10.9,7.1)(10.4,6.15)(9.65,6.6)(10.2,7.6)(8.9,7.6)(7.9,7.1)(6.9,7.1)(6.4,6.1)(5.4,5.6)(4.8,4.4)(6.1,4.4)(7.1,4.9)(8.1,4.9)(8.6,5.85)(9.3,5.4)(8.85,4.4)(10.1,4.4)(11,4.9)(12.1,4.9)(12.6,5.85)(13.35,5.4)(12.8,4.4)(14,4.4)(14.9,4.9)(14.9,4.1)(14,3.6)(12.9,3.6)(12.45,2.6)(11.65,3.1)(12.15,4.1)(11,4.1)(9.9,3.6)(8.9,3.6)(8.4,2.65)(7.7,3.1)(8.1,4.1)(6.9,4.1)(5.9,3.6)(4.9,3.6)(4.4,2.6)(3.4,2.1)(2.9,.9)(4.1,.9)(5.1,1.4)(6.1,1.4)(6.6,2.3)(7.3,1.9)(6.8,.9)(8.1,.9)(9.1,1.4)(10.1,1.4)(10.6,2.35)(11.35,1.9)(10.8,.9)(12.1,.9)(13.1,1.4)(14,1.4)

\put(1.5,.8){\scriptsize {\tiny $a_{11}$}}
\put(2.5,.4){\scriptsize {\tiny $a_{12}$}}
\put(3.6,.4){\scriptsize {\tiny $a_{13}$}}
\put(4.6,.8){\scriptsize {\tiny $a_{14}$}}

\put(12.6,.8){\scriptsize {\tiny $a_{1i}$}}
\put(13.8,.8){\scriptsize {\tiny $a_{11}$}}

\put(1.55,2.6){\scriptsize {\tiny $x_{11}$}}
\put(14.9,2){\scriptsize {\tiny $x_{11}$}}

\put(2.35,3.1){\scriptsize {\tiny $x_{12}$}}
\put(15.8,2.7){\scriptsize {\tiny $x_{12}$}}

\put(2.8,3.9){\scriptsize {\tiny $a_{21}$}}
\put(16.2,3.7){\scriptsize {\tiny $a_{21}$}}

\put(4.8,7.4){\scriptsize {\tiny $a_{j1}$}}
\put(18.25,7.3){\scriptsize {\tiny $a_{j1}$}}

\put(4.25,9){\scriptsize {\tiny $a_{1(k+1)}$}}
\put(6.5,8.65){\scriptsize {\tiny $a_{1(k+2)}$}}
\put(16.4,9.2){\scriptsize {\tiny $a_{1k}$}}
\put(17.6,9.2){\scriptsize {\tiny $a_{1(k+1)}$}}

\put(.5,1.5){\scriptsize {\tiny $Q_1$}}
\put(1.5,4.25){\scriptsize {\tiny $Q_2$}}
\put(3,7.65){\scriptsize {\tiny $Q_{j}$}}
\put(3.5,8.75){\scriptsize {\tiny $Q_{1}$}}

\put(7.5,.5){\scriptsize $\ldots$}
\put(9.5,.5){\scriptsize $\ldots$}
\put(10,9.05){\scriptsize $\ldots$}
\put(14,9.05){\scriptsize $\ldots$}

\put(2.7,6.2){\scriptsize $.$}
\put(2.6,6){\scriptsize $.$}
\put(2.5,5.8){\scriptsize $.$}

\put(4.75,-1){\scriptsize {\tiny {\bf Figure 4.10.2:} $M(i,j,k)$ }} 

\end{picture}

\vspace{4.3cm}

Now, construct $M(i,j,k)$ representation of $M$ by first cutting $M$ along a black colored cycle of type $J_1$ and then cutting it along the green colored cycle of type $J_1$, where without loss of generality, let the beginning adjacent face to the base horizontal cycle is a quadrangle.

\begin{lem}\label{l3.10.1} A DSEM $M$ of type $[3^2.4.3.4:3.4.6.4]$ admits an $M(i,j,k)$-representation iff the following holds: $(i)$ $j$ even and $i=4m$, $m \in \mathbb{N} $, $(ii)$ number of vertices of $M(i,j,k) = 3ij/2 \geq 12 $, $(iii)$ $ k \in \{4r: 0 \leq r < i/4\}$.
	
\end{lem}

\noindent{\bf Proof.} Let $M$ be a DSEM of type $[3^2.4.3.4:3.4.6.4]$ having $n$ vertices. Then its $M(i,j,k)$  has $j$ disjoint horizontal cycles of $J_{1}$ type, say $Q_0, Q_1,\ldots, Q_{j-1}$, of length $i$. Let $Q_{0}=C(w_{0,0},w_{0,1}, \linebreak \ldots,w_{0,i-1}),Q_{1}=C(w_{1,0},w_{1,1},\ldots,w_{1,i-1}),\ldots,Q_{j-1}=C(w_{j-1,0},w_{j-1,1},\ldots,w_{j-1,i-1})$ be the list of horizontal cycles. Observe that the vertices having face-sequence $(3^2,4,3,4)$ lying between horizontal cycles $Q_{2s(mod\,j)}$ and $Q_{(2s+1)(mod\,j)}$ for $0 \leq s \leq j-1$ is $4i/4\cdot j/2$. Therefore, $n = ij + ij/2 = 3ij/2$. 

If $j=1$, then  no vertex in the base horizontal cycle follow the face-sequence $(3^2,4,3,4)$ after identifying the boundaries of $M(i,1,k)$. Therefore, $j \geq 2$. If $j \geq 2$ and $j$ is not an even integer, then as above no vertex in the base horizontal cycle follows the face-sequence $(3^2,4,3,4)$. So $j$ is an even integer. If $j$ is even and $i \neq 4m$, where $m \in \mathbb{N}$, then the ${\rm lk}(w_{0,0})$ is not of type $(3,4,6,4)$. Which is not possible. So, $j$ is an even integer and $i = 4m$, where $m \in \mathbb{N}$. Thus, $3ij/ \geq 12$.

If $j$ is an even integer and $ k \in \{r:0 \leq r \leq i-1\} \setminus (k \in \{4r: 0 \leq r < i/4\})$, then some vertex in $M(i,j,k)$ do not follow the face-sequences $(3^2,4,3,4)$ and $(3,4,6,4)$. So, $ k \in \{4r: 0 \leq r < i/4\}$. This completes the proof. \hfill $\Box$

%
%

\subsection{DSEMs of type $[3^6:3^4.6]_1$} \label{s3.11}

Let $M$ be a DSEM of type $[3^6:3^4.6]_1$. Then for the existence of the map $M$, we see that $|V_{(3^6)}|=|V_{(3^4,6)}|$. Now, consider a fixed type path in $M$ as follows. A path $P_{1} = P( \ldots,z_{i},y_{i},y_{i+1},y_{i+2},z_{i+1}, \linebreak \ldots)$ in $M$, say of type $K_{1}$, indicated by thick black or green paths, shown in Figure 4.11.1. The vertices $z_i$'s and $y_i$'s have the face-sequences $(3^6)$ and $(3^4,6)$ respectively or vertices $z_i's$ and $y_i's$ have the face-sequence $(3^6)$.

\vspace{1cm}
\begin{picture}(0,0)(-12,44)
\setlength{\unitlength}{3.5mm}


\drawline[AHnb=0](-1.25,.5)(-3,3.75)
\drawline[AHnb=0](-.5,1)(-1,2)(-3,6)
\drawline[AHnb=0](.5,1)(0,2)(-1,4)
\drawline[AHnb=0](-2,6)(-3,8)
\drawline[AHnb=0](.5,3)(-3,10)
\drawline[AHnb=0](1.75,.5)(1.5,1)(2.5,1)(-2,10)
\drawline[AHnb=0](3.5,1)(3,2)(0,8)
\drawline[AHnb=0](4.5,1)(4,2)(2.5,5)
\drawline[AHnb=0](1.5,7)(0,10)
\drawline[AHnb=0](4,4)(1,10)
\drawline[AHnb=0](5,4)(2,10)
\drawline[AHnb=0](5,6)(3.5,9)
\drawline[AHnb=0](5,8)(4,10)


\drawline[AHnb=0](2.5,1)(4.5,1)
\drawline[AHnb=0](-1.5,1)(-.5,1)(1.5,1)(2,2)(5,2)
\drawline[AHnb=0](-3,2)(-1,2)(0,2)(.5,3)(1.5,3)(2.5,3)(3.5,3)(4,4)(5,4)
\drawline[AHnb=0](-3,3)(-.5,3)(0,4)(1,4)(2,4)(3,4)(3.5,5)(5,5)
\drawline[AHnb=0](-3,4)(-1,4)(-.5,5)(.5,5)(1.5,5)(2.5,5)(3,6)(5,6)
\drawline[AHnb=0](-3,5)(-2.5,5)(-2,6)(0,6)(1,6)(1.5,7)(5,7)
\drawline[AHnb=0](-3,6)(-2.5,7)(-.5,7)(.5,7)(1,8)(5,8)
\drawline[AHnb=0](-3,8)(-1,8)(0,8)(.5,9)(5,9)
\drawline[AHnb=0](-3,9)(-1.5,9)(-1,10)(5,10)
\drawline[AHnb=0](-3,10)(-1,10)

\drawline[AHnb=0](2.25,.5)(2.5,1)(3.5,3)
\drawline[AHnb=0](1.25,.5)(2,2)(3,4)(4,4)(4.5,5)
\drawline[AHnb=0](1.5,3)(2.5,5)(3.5,5)(5,8)
\drawline[AHnb=0](-1.75,.5)(-1,2)(-.5,3)(.5,3)(1.5,5)
\drawline[AHnb=0](-2,2)(-1.5,3)(-1,4)(0,4)(1,6)

\drawline[AHnb=0](-.5,5)(.5,7)(1.5,7)(2.5,9)
\drawline[AHnb=0](-1,6)(0,8)(1,8)(2,10)
\drawline[AHnb=0](-3,2.25)(-2,4)

\drawline[AHnb=0](-3.5,5)(-3,6)(-2,6)(-1,8)
\drawline[AHnb=0](-3,7)(-2.5,7)(-1.5,9)
\drawline[AHnb=0](-3,8)(-2,10)

\drawline[AHnb=0](3.5,1)(4,2)
\drawline[AHnb=0](4.5,1)(5,2)
\drawline[AHnb=0](3,6)(4.5,9)
\drawline[AHnb=0](2.5,7)(3.5,9)

\drawline[AHnb=0](-.5,1)(0,2)
\drawline[AHnb=0](4.5,5)(5,6)
\drawline[AHnb=0](.5,9)(1,10)

\drawline[AHnb=0](-.75,.5)(-.5,1)(-.25,.5)
\drawline[AHnb=0](.25,.5)(.5,1)(.75,.5)
\drawline[AHnb=0](2.25,.5)(2.5,1)(2.75,.5)
\drawline[AHnb=0](3.25,.5)(3.5,1)(3.75,.5)
\drawline[AHnb=0](4.25,.5)(4.5,1)(4.75,.5)


\drawpolygon[fillcolor=black](-3,3)(-.5,3)(0,4)(.5,4)(2,4)(3,4)(3.5,5)(5,5)(5,5.1)(3.4,5.1)(2.9,4.1)(2,4.1)(.4,4.1)(-.1,4.1)(-.6,3.1)(-3,3.1)

\drawpolygon[fillcolor=black](-3,4)(-1,4)(-.5,5)(.5,5)(1.5,5)(2.5,5)(3,6)(5,6)(5,6.1)(3,6.1)(2.4,5.1)(1.4,5.1)(.5,5.1)(-.6,5.1)(-1.1,4.1)(-3,4.1)

\drawpolygon[fillcolor=black](-3,5)(-2.5,5)(-2,6)(1,6)(1.5,7)(2,7)(3.5,7)(4.5,7)(5,8)(5.5,8)(5.5,8.1)(4.9,8.1)(4.4,7.1)(3.5,7.1)(1.9,7.1)(1.4,7.1)(.9,6.1)(-2.15,6.1)(-2.6,5.1)(-3,5.1)

\drawpolygon[fillcolor=black](-3,6)(-2.5,7)(.5,7)(1,8)(1.5,8)(3,8)(4,8)(4.5,9)(5,9)(5,9.1)(4.4,9.1)(3.9,8.1)(3,8.1)(1.4,8.1)(.9,8.1)(.4,7.1)(-2.55,7.1)(-3.1,6)

\drawpolygon[fillcolor=black](-3,8)(0,8)(.5,9)(1.5,9)(2.5,9)(3.5,9)(4,10)(5,10)(5,10.1)(4,10.1)(3.4,9.1)(2.4,9.1)(1.5,9.1)(.4,9.1)(-.1,8.1)(-3,8.1)

\drawpolygon[fillcolor=black](-3,2)(0,2)(.5,3)(1,3)(2.,3)(3.5,3)(4,4)(5.5,4)(5.5,4.1)(3.9,4.1)(3.4,3.1)(2.5,3.1)(.9,3.1)(.4,3.1)(-.1,2.1)(-3,2.1)

\drawpolygon[fillcolor=black](-3,2)(0,2)(.5,3)(1,3)(2.,3)(3.5,3)(4,4)(5.5,4)(5.5,4.1)(3.9,4.1)(3.4,3.1)(2.5,3.1)(.9,3.1)(.4,3.1)(-.1,2.1)(-3,2.1)

\drawpolygon[fillcolor=green](3,0)(1.4,3.1)(.4,3.1)(-1,6)(-2,6)(-3.5,9)(-3.35,9)(-1.9,6.15)(-.9,6.15)(0.5,3.3)(1.5,3.3)(3.1,.15)

\drawpolygon[fillcolor=green](4.75,.5)(4.5,1)(3.5,1)(1.9,4.1)(.9,4.1)(-.5,7)(-1.5,7)(-3,10)(-2.85,10)(-1.4,7.15)(-.4,7.15)(1,4.3)(2,4.3)(3.6,1.15)(4.6,1.15)(4.9,.5)

\drawpolygon[fillcolor=green](5.25,1.5)(5,2)(4,2)(2.4,5.1)(1.4,5.1)(0,8)(-1,8)(-2,10)(-1.85,10)(-.9,8.15)(.1,8.15)(1.5,5.3)(2.5,5.3)(4.1,2.15)(5.1,2.15)(5.4,1.5)

\drawpolygon[fillcolor=green](1.75,.5)(1.5,1)(.5,1)(-1.1,4.1)(-2.1,4.1)(-3.5,7)(-3.4,7.15)(-2,4.3)(-1,4.3)(.7,1.15)(1.6,1.15)(1.9,.5)

\drawpolygon[fillcolor=green](-.5,1)(-1.5,3)(-2.5,3)(-3,4)(-3.15,4)(-2.65,2.85)(-1.65,2.85)(-.65,1)


\drawpolygon[fillcolor=green](5.5,4.1)(3.85,4.1)(2.5,7)(1.5,7)(0,10)(.15,10)(1.6,7.15)(2.6,7.15)(4,4.25)(5.5,4.25)

\drawpolygon[fillcolor=green](4.85,6)(3.5,9)(2.5,9)(2,10)(2.15,10)(2.6,9.15)(3.6,9.15)(5,6)

\drawpolygon[fillcolor=green](5,5.1)(4.35,5.1)(3,8)(2,8)(1,10)(1.15,10)(2.1,8.15)(3.1,8.15)(4.5,5.25)(5,5.25)



\put(-3.5,-1){\scriptsize {\tiny {\bf Figure 4.11.1:} Paths of type $K_1$ }}

\end{picture}

\begin{picture}(0,0)(-64,30)
\setlength{\unitlength}{4mm}

\drawline[AHnb=0](-5,6)(-6.5,9)
\drawline[AHnb=0](-2.5,3)(-5.5,9)
\drawline[AHnb=0](0,0)(-.5,1)(-1,2)(-3.5,7)
\drawline[AHnb=0](1,0)(.5,1)(0,2)(-1,4)
\drawline[AHnb=0](-2,6)(-4,10)
\drawline[AHnb=0](.5,3)(-3,10)
\drawline[AHnb=0](2,0)(1.5,1)(2.5,1)(-2,10)
\drawline[AHnb=0](3.5,1)(3,2)(0,8)
\drawline[AHnb=0](4.5,1)(4,2)(2.5,5)
\drawline[AHnb=0](1.5,7)(-.5,11)
\drawline[AHnb=0](4,4)(.5,11)
\drawline[AHnb=0](5.5,1)(5,2)(6,2)(5.5,3)(3,8)
\drawline[AHnb=0](7,2)(5.5,5)


\drawline[AHnb=0](0,0)(2,0)(2.5,1)(5.5,1)(6,2)(7,2)
\drawline[AHnb=0](-.5,1)(1.5,1)(2,2)(5,2)(5.5,3)(6.5,3)
\drawline[AHnb=0](-1,2)(0,2)(.5,3)(1.5,3)(2.5,3)(3.5,3)(4,4)(6,4)
\drawline[AHnb=0](-2.5,3)(-.5,3)(0,4)(1,4)(2,4)(3,4)(3.5,5)(5.5,5)
\drawline[AHnb=0](-3,4)(-1,4)(-.5,5)(.5,5)(1.5,5)(2.5,5)(3,6)(4,6)
\drawline[AHnb=0](-3.5,5)(-2.5,5)(-2,6)(0,6)(1,6)(1.5,7)(3.5,7)
\drawline[AHnb=0](-5,6)(-3,6)(-2.5,7)(-.5,7)(.5,7)(1,8)(3,8)
\drawline[AHnb=0](-5.5,7)(-3.5,7)(-3,8)(-1,8)(0,8)(.5,9)(1.5,9)
\drawline[AHnb=0](-6,8)(-5,8)(-4.5,9)(-1.5,9)(-1,10)(1,10)
\drawline[AHnb=0](-6.5,9)(-5.5,9)(-5,10)(-2,10)(-1.5,11)(.5,11)


\drawline[AHnb=0](2.5,1)(3.5,3)
\drawline[AHnb=0](2,2)(3,4)(4,4)(4.5,5)
\drawline[AHnb=0](1.5,3)(2.5,5)(3.5,5)(4,6)
\drawline[AHnb=0](-1,2)(-.5,3)(.5,3)(1.5,5)
\drawline[AHnb=0](-1.5,3)(-1,4)(0,4)(1,6)
\drawline[AHnb=0](-.5,5)(.5,7)(1.5,7)(2,8)
\drawline[AHnb=0](-1,6)(0,8)(1,8)(1.5,9)
\drawline[AHnb=0](-3.5,5)(-3,6)(-2,6)(-1,8)
\drawline[AHnb=0](-4,6)(-3.5,7)(-2.5,7)(-1.5,9)
\drawline[AHnb=0](-3,8)(-2,10)

\drawline[AHnb=0](0,0)(.5,1)
\drawline[AHnb=0](1,0)(1.5,1)

\drawline[AHnb=0](3.5,1)(4,2)
\drawline[AHnb=0](4.5,1)(5,2)
\drawline[AHnb=0](6,2)(6.5,3)
\drawline[AHnb=0](5.5,3)(6,4)
\drawline[AHnb=0](5,4)(5.5,5)
\drawline[AHnb=0](3,6)(3.5,7)
\drawline[AHnb=0](2.5,7)(3,8)

\drawline[AHnb=0](.5,9)(1,10)
\drawline[AHnb=0](0,10)(.5,11)
\drawline[AHnb=0](-1,10)(-.5,11)
\drawline[AHnb=0](-4.5,9)(-4,10)
\drawline[AHnb=0](-6,8)(-5.5,9)
\drawline[AHnb=0](-5.5,7)(-5,8)
\drawline[AHnb=0](-5,6)(-4.5,7)
\drawline[AHnb=0](-3,4)(-2.5,5)
\drawline[AHnb=0](-2.5,3)(-2,4)
\drawline[AHnb=0](-.5,1)(0,2)
\drawline[AHnb=0](-5.5,9)(-4.5,9)(-5,10)
\drawline[AHnb=0](-2.5,10)(-1,10)(-1.5,11)
\drawline[AHnb=0](-3.5,9)(-3,10)


\drawpolygon[fillcolor=black](0,0)(2,0)(2.5,1)(5.5,1)(6,2)(5.5,3)(5,2)(2,2)(1.5,1)(0.5,1)(0,2)(.5,3)(1.5,3)(2.5,3)(3.5,3)(4,4)(3.5,5)(3,4)(2,4)(1,4)(0,4)(-.5,3)(-1,4)(-.5,5)(.5,5)(1.5,5)(2.5,5)(3,6)(2.5,7)(1.5,7)(1,6)(0,6)(-2,6)(-2.5,5)(-3,6)(-2.5,7)(-.5,7)(.5,7)(1,8)(.5,9)(0,8)(-1,8)(-3,8)(-3.5,7)(-4.5,7)(-5,8)(-4.5,9)(-1.5,9)(-1,10)(1,10)(1,9.9)(-.9,9.9)(-1.4,8.9)(-4.4,8.9)(-4.9,8)(-4.4,7.1)(-3.6,7.1)(-3.1,8.1)(-.1,8.1)(.5,9.1)(1.1,8.1)(.6,6.9)(-2.4,6.9)(-2.9,6)(-2.5,5.1)(-2.1,6.1)(.9,6.1)(1.4,7.1)(2.6,7.1)(3.1,6)(2.6,4.9)(-.4,4.9)(-.9,4)(-.5,3.1)(-.1,4.1)(2.9,4.1)(3.5,5.1)(4.1,4.1)(3.6,2.9)(.6,2.9)(.15,1.9)(.6,1.1)(1.4,1.1)(1.9,2.1)(4.9,2.1)(5.5,3.1)(6.15,1.9)(5.6,.9)(2.6,.9)(2.1,-.1)(0,-.1)

\drawpolygon[fillcolor=black](0,0)(-1,2)(-.85,2)(.15,0)

\drawpolygon[fillcolor=black](-1.5,3)(-2.5,3)(-2,4)(-3,4)(-4,6)(-5,6)(-6,8)(-5.9,8)(-4.9,6.1)(-3.9,6.1)(-2.9,4.1)(-1.9,4.1)(-2.3,3.1)(-1.4,3.1)

\drawpolygon[fillcolor=black](5.5,5)(5,4)(6,4)(6,4.1)(5.15,4.1)(5.6,5)

\put(-.4,-.5){\scriptsize {\tiny $a_{11}$}}
\put(.7,-.5){\scriptsize {\tiny $a_{12}$}}
\put(1.7,-.5){\scriptsize {\tiny $a_{13}$}}
\put(6.1,1.4){\scriptsize {\tiny $a_{1i}$}}
\put(7.1,1.4){\scriptsize {\tiny $a_{11}$}}

\put(-1.5,.7){\scriptsize {\tiny $a_{21}$}}
\put(6.9,2.8){\scriptsize {\tiny $a_{21}$}}

\put(-2,1.6){\scriptsize {\tiny $a_{31}$}}
\put(6.4,3.9){\scriptsize {\tiny $a_{31}$}}

\put(-3.6,2.9){\scriptsize {\tiny $a_{41}$}}
\put(-2.45,2.5){\scriptsize {\tiny $a_{42}$}}
\put(4.4,5.25){\scriptsize {\tiny $a_{41}$}}
\put(5.9,4.9){\scriptsize {\tiny $a_{42}$}}

\put(-4,3.8){\scriptsize {\tiny $a_{51}$}}
\put(4.3,6.1){\scriptsize {\tiny $a_{51}$}}

\put(-7.1,7.6){\scriptsize {\tiny $a_{j1}$}}
\put(1.4,9.8){\scriptsize {\tiny $a_{j1}$}}

\put(-8,9.4){\scriptsize {\tiny $a_{1(k+1)}$}}
\put(-7.25,10.4){\scriptsize {\tiny $a_{1(k+2)}$}}
\put(-1,11.5){\scriptsize {\tiny $a_{1k}$}}
\put(0,11.5){\scriptsize {\tiny $a_{1(k+1)}$}}

\put(-2.5,-.35){\scriptsize {\tiny $Q_1$}}
\put(-3,.75){\scriptsize {\tiny $Q_2$}}
\put(-8.5,8.55){\scriptsize {\tiny $Q_1$}}
\put(-8.25,7.5){\scriptsize {\tiny $Q_{j}$}}

\put(4.5,0.5){\scriptsize $\ldots$}
\put(-4,10.5){\scriptsize $\ldots$}

\put(-5.7,4.9){\scriptsize $.$}
\put(-5.6,4.7){\scriptsize $.$}
\put(-5.5,4.5){\scriptsize $.$}

\put(-2,-1.75){\scriptsize {\tiny {\bf Figure 4.11.2:} $M(i,j=6m+3,k)$ }}

\end{picture}

\begin{picture}(0,0)(-125,35)
\setlength{\unitlength}{4mm}

\drawline[AHnb=0](-5,6)(-6.5,9)
\drawline[AHnb=0](-2.5,3)(-5.5,9)
\drawline[AHnb=0](0,0)(-.5,1)(-1,2)(-3.5,7)
\drawline[AHnb=0](1,0)(.5,1)(0,2)(-1,4)
\drawline[AHnb=0](-2,6)(-4,10)
\drawline[AHnb=0](.5,3)(-3,10)
\drawline[AHnb=0](2,0)(1.5,1)(2.5,1)(-2,10)
\drawline[AHnb=0](3.5,1)(3,2)(0,8)
\drawline[AHnb=0](4.5,1)(4,2)(2.5,5)
\drawline[AHnb=0](1.5,7)(-.5,11)
\drawline[AHnb=0](4,4)(.5,11)
\drawline[AHnb=0](5.5,1)(5,2)(6,2)(5.5,3)(3,8)
\drawline[AHnb=0](7,2)(5.5,5)

\drawline[AHnb=0](-8,10)(-8.5,11)(-9,12)
\drawline[AHnb=0](-5.5,9)(-6,10)(-6.5,9)(-7,10)(-7.5,11)(-8,12)

\drawline[AHnb=0](-1,12)(-1.5,13)(-2,14)

\drawline[AHnb=0](-7.5,9)(-7,10)

\drawline[AHnb=0](-.5,11)(-1,12)(-1.5,11)(-2,12)(-2.5,13)(-3,14)(-3.5,13)
\drawline[AHnb=0](-4.5,13)(-3.5,13)(-4,14)
\drawline[AHnb=0](-6,10)(-5,10)(-5.5,11)

\drawline[AHnb=0](-3.5,11)(-4,12)(-4.5,13)(-5,12)

\drawline[AHnb=0](-3.5,9)(-4,10)(-4.5,11)(-5,10)

\drawline[AHnb=0](-2.5,9)(-3,10)(-3.5,11)(-4,10)

\drawline[AHnb=0](-1.5,9)(-2,10)(-2.5,11)(-3,10)

\drawline[AHnb=0](-4.5,11)(-5,12)(-5.5,13)(-6,12)
\drawline[AHnb=0](-5,12)(-5.5,11)(-6,12)(-6.5,13)

\drawline[AHnb=0](-2.5,11)(-1.5,11)

\drawline[AHnb=0](-6.5,13)(-7,12)(-7.5,13)

\drawline[AHnb=0](-8.5,11)(-8,12)(-6.5,12)
\drawline[AHnb=0](-7.5,11)(-8,10)
\drawline[AHnb=0](-4,12)(-4.5,11)
\drawline[AHnb=0](-1.5,13)(-2,12)
\drawline[AHnb=0](-2,14)(-2.5,13)


\drawline[AHnb=0](0,0)(2,0)(2.5,1)(5.5,1)(6,2)(7,2)
\drawline[AHnb=0](-.5,1)(1.5,1)(2,2)(5,2)(5.5,3)(6.5,3)
\drawline[AHnb=0](-1,2)(0,2)(.5,3)(1.5,3)(2.5,3)(3.5,3)(4,4)(6,4)
\drawline[AHnb=0](-2.5,3)(-.5,3)(0,4)(1,4)(2,4)(3,4)(3.5,5)(5.5,5)
\drawline[AHnb=0](-3,4)(-1,4)(-.5,5)(.5,5)(1.5,5)(2.5,5)(3,6)(4,6)
\drawline[AHnb=0](-3.5,5)(-2.5,5)(-2,6)(0,6)(1,6)(1.5,7)(3.5,7)
\drawline[AHnb=0](-5,6)(-3,6)(-2.5,7)(-.5,7)(.5,7)(1,8)(3,8)
\drawline[AHnb=0](-5.5,7)(-3.5,7)(-3,8)(-1,8)(0,8)(.5,9)(1.5,9)
\drawline[AHnb=0](-6,8)(-5,8)(-4.5,9)(-1.5,9)(-1,10)(1,10)
\drawline[AHnb=0](-8,10)(-7.5,9)(-5.5,9)(-5,10)(-2,10)(-1.5,11)(.5,11)

\drawline[AHnb=0](-8,10)(-6,10)(-5.5,11)(-2.5,11)(-2,12)(-1,12)
\drawline[AHnb=0](-8.5,11)(-7.5,11)(-7,12)(-6,12)(-5,12)(-4,12)(-3.5,13)(-1.5,13)
\drawline[AHnb=0](-9,12)(-8,12)(-7.5,13)(-6.5,13)(-5.5,13)(-4.5,13)(-4,14)(-2,14)


\drawline[AHnb=0](2.5,1)(3.5,3)
\drawline[AHnb=0](2,2)(3,4)(4,4)(4.5,5)
\drawline[AHnb=0](1.5,3)(2.5,5)(3.5,5)(4,6)
\drawline[AHnb=0](-1,2)(-.5,3)(.5,3)(1.5,5)
\drawline[AHnb=0](-1.5,3)(-1,4)(0,4)(1,6)
\drawline[AHnb=0](-.5,5)(.5,7)(1.5,7)(2,8)
\drawline[AHnb=0](-1,6)(0,8)(1,8)(1.5,9)
\drawline[AHnb=0](-3.5,5)(-3,6)(-2,6)(-1,8)
\drawline[AHnb=0](-4,6)(-3.5,7)(-2.5,7)(-1.5,9)
\drawline[AHnb=0](-3,8)(-2,10)

\drawline[AHnb=0](0,0)(.5,1)
\drawline[AHnb=0](1,0)(1.5,1)

\drawline[AHnb=0](3.5,1)(4,2)
\drawline[AHnb=0](4.5,1)(5,2)
\drawline[AHnb=0](6,2)(6.5,3)
\drawline[AHnb=0](5.5,3)(6,4)
\drawline[AHnb=0](5,4)(5.5,5)
\drawline[AHnb=0](3,6)(3.5,7)
\drawline[AHnb=0](2.5,7)(3,8)

\drawline[AHnb=0](.5,9)(1,10)
\drawline[AHnb=0](0,10)(.5,11)
\drawline[AHnb=0](-1,10)(-.5,11)
\drawline[AHnb=0](-4.5,9)(-4,10)
\drawline[AHnb=0](-6,8)(-5.5,9)
\drawline[AHnb=0](-5.5,7)(-5,8)
\drawline[AHnb=0](-5,6)(-4.5,7)
\drawline[AHnb=0](-3,4)(-2.5,5)
\drawline[AHnb=0](-2.5,3)(-2,4)
\drawline[AHnb=0](-.5,1)(0,2)
\drawline[AHnb=0](-5.5,9)(-4.5,9)(-5,10)
\drawline[AHnb=0](-2.5,10)(-1,10)(-1.5,11)
\drawline[AHnb=0](-3.5,9)(-3,10)


\drawpolygon[fillcolor=black](0,0)(2,0)(2.5,1)(5.5,1)(6,2)(5.5,3)(5,2)(2,2)(1.5,1)(0.5,1)(0,2)(.5,3)(1.5,3)(2.5,3)(3.5,3)(4,4)(3.5,5)(3,4)(2,4)(1,4)(0,4)(-.5,3)(-1,4)(-.5,5)(.5,5)(1.5,5)(2.5,5)(3,6)(2.5,7)(1.5,7)(1,6)(0,6)(-2,6)(-2.5,5)(-3,6)(-2.5,7)(-.5,7)(.5,7)(1,8)(.5,9)(0,8)(-1,8)(-3,8)(-3.5,7)(-4.5,7)(-5,8)(-4.5,9)(-1.5,9)(-1,10)(-1.5,11)(-2,10)(-5,10)(-5.5,9)(-6,10)(-5.5,11)(-2.5,11)(-2,12)
(-2.5,13)(-3.5,13)(-4,12)(-7,12)(-7.5,11)(-8.5,11)(-8.5,11.1)(-7.6,11.1)(-7.1,12.1)(-4.1,12.1)(-3.6,13.1)(-2.4,13.1)(-1.9,12)(-2.4,10.9)(-5.4,10.9)(-5.9,10)(-5.5,9.15)(-5.1,10.1)(-2.1,10.1)(-1.5,11.1)(-.85,10)(-1.4,8.9)(-4.4,8.9)(-4.9,8)(-4.4,7.1)(-3.6,7.1)(-3.1,8.1)(-.1,8.1)(.5,9.1)(1.1,8.1)(.6,6.9)(-2.4,6.9)(-2.9,6)(-2.5,5.1)(-2.1,6.1)(.9,6.1)(1.4,7.1)(2.6,7.1)(3.1,6)(2.6,4.9)(-.4,4.9)(-.9,4)(-.5,3.1)(-.1,4.1)(2.9,4.1)(3.5,5.1)(4.1,4.1)(3.6,2.9)(.6,2.9)(.15,1.9)(.6,1.1)(1.4,1.1)(1.9,2.1)(4.9,2.1)(5.5,3.1)(6.15,1.9)(5.6,.9)(2.6,.9)(2.1,-.1)(0,-.1)

\drawpolygon[fillcolor=black](0,0)(-1,2)(-.85,2)(.15,0)

\drawpolygon[fillcolor=black](-1.5,3)(-2.5,3)(-2,4)(-3,4)(-4,6)(-5,6)(-6,8)(-5.9,8)(-4.9,6.1)(-3.9,6.1)(-2.9,4.1)(-1.9,4.1)(-2.3,3.1)(-1.4,3.1)

\drawpolygon[fillcolor=black](5.5,5)(5,4)(6,4)(6,4.1)(5.15,4.1)(5.6,5)

\drawpolygon[fillcolor=black](.5,11)(0,10)(1,10)(1,10.1)(.15,10.1)(.6,11)

\drawpolygon[fillcolor=black](-6.5,9)(-7.5,9)(-7,10)(-8,10)(-8.5,11)(-8.4,11)(-7.9,10.1)(-6.85,10.1)(-7.35,9.1)(-6.5,9.1)

\put(-.4,-.5){\scriptsize {\tiny $a_{11}$}}
\put(.7,-.5){\scriptsize {\tiny $a_{12}$}}
\put(1.7,-.5){\scriptsize {\tiny $a_{13}$}}
\put(6.1,1.4){\scriptsize {\tiny $a_{1i}$}}
\put(7.1,1.4){\scriptsize {\tiny $a_{11}$}}

\put(-1.5,.7){\scriptsize {\tiny $a_{21}$}}
\put(6.9,2.8){\scriptsize {\tiny $a_{21}$}}

\put(-2,1.6){\scriptsize {\tiny $a_{31}$}}
\put(6.4,3.9){\scriptsize {\tiny $a_{31}$}}

\put(-3.6,2.9){\scriptsize {\tiny $a_{41}$}}
\put(-2.45,2.5){\scriptsize {\tiny $a_{42}$}}
\put(4.4,5.25){\scriptsize {\tiny $a_{41}$}}
\put(5.9,4.9){\scriptsize {\tiny $a_{42}$}}

\put(-4,3.8){\scriptsize {\tiny $a_{51}$}}
\put(4.3,6.1){\scriptsize {\tiny $a_{51}$}}

\put(-4.4,4.6){\scriptsize {\tiny $a_{61}$}}
\put(3.9,7){\scriptsize {\tiny $a_{61}$}}

\put(-6,5.8){\scriptsize {\tiny $a_{71}$}}
\put(-5,5.5){\scriptsize {\tiny $a_{72}$}}
\put(2.1,8.25){\scriptsize {\tiny $a_{71}$}}
\put(3.25,8){\scriptsize {\tiny $a_{72}$}}

\put(-6.5,6.6){\scriptsize {\tiny $a_{81}$}}
\put(1.8,9){\scriptsize {\tiny $a_{81}$}}

\put(-7.1,7.6){\scriptsize {\tiny $a_{91}$}}

\put(1.4,9.8){\scriptsize {\tiny $a_{91}$}}

\put(-9.75,10.85){\scriptsize {\tiny $a_{j1}$}}
\put(-1.25,12.85){\scriptsize {\tiny $a_{j1}$}}

\put(-10.5,12.4){\scriptsize {\tiny $a_{1(k+1)}$}}
\put(-9.95,13.4){\scriptsize {\tiny $a_{1(k+2)}$}}
\put(-3.5,14.5){\scriptsize {\tiny $a_{1k}$}}
\put(-2.5,14.5){\scriptsize {\tiny $a_{1(k+1)}$}}

\put(-2.5,-.35){\scriptsize {\tiny $Q_1$}}
\put(-3.5,.75){\scriptsize {\tiny $Q_2$}}
\put(-11.5,11.65){\scriptsize {\tiny $Q_1$}}
\put(-11,10.85){\scriptsize {\tiny $Q_{j}$}}

\put(4.5,0.5){\scriptsize $\ldots$}
\put(-6.5,13.5){\scriptsize $\ldots$}

\put(-8.7,8.9){\scriptsize $.$}
\put(-8.6,8.7){\scriptsize $.$}
\put(-8.5,8.5){\scriptsize $.$}

\put(-2,-1.75){\scriptsize {\tiny {\bf Figure 4.11.3:} $M(i,j=6m+6,k)$ }}

\end{picture}

\vspace{4.5cm}

As in previous sections, $M$ has an $M(i, j, k)$ representation. The admissible relations among $i, j, k$ of $M(i, j, k)$ are given below.

\begin{lem}\label{l3.11.1} A DSEM $M$ of type $[3^6:3^4.6]_1$ admits an $M(i,j,k)$-representation iff the following holds: $(i)$ $j = 3m$ and $i=4m$, where $m \in \mathbb{N} $, $(ii)$ number of vertices of $M(i,j,k) = ij \geq 12 $, $(iii)$ $ k \in \{4r+1: 0 \leq r < i/4\}$.
	
\end{lem}

\noindent{\bf Proof.} Let $M$ be a DSEM of type $[3^6:3^4.6]_1$ having $n$ vertices. Its $M(i,j,k)$  has $j$  disjoint horizontal cycles of $K_{1}$ type, say $Q_0, Q_1,\ldots, Q_{j-1}$, of length $i$. Since the vertices of $M$ are in these cycles, the vertices in $M$ is $n = ij$. If $j = 1,2$ then $M$ is not a map. So $j \geq 3$. If $j \geq 3$ and $j \neq 3m,\ m \in \mathbb{N}$, then  no vertex in the base horizontal cycle follows the face-sequence $(3^6)$ after identifying the boundaries of $M(i,j,k)$. So $j = 3m, \text{ where } m \in \mathbb{N}$. Clearly if $i \leq 3$, $M$ is not a map. So $i \geq 4$. If $j=3m$ and $i \neq 4m$, where $m \in \mathbb{N}$, then either $|V_{(3^6)}| \neq |V_{(3^4,6)}|$ or $M$ is not a map. So, $i = 4m$, where $m \in \mathbb{N}$. Thus $n = ij \geq 12$.

If $ k \in \{r:0 \leq r \leq i-1\} \setminus (\{4r+1: 0 \leq r < i/4\})$, then the link of some vertex in $M(i,j,k)$ can not be completed. So, $ k \in \{4r+1: 0 \leq r < i/4\}$. This completes the proof. \hfill $\Box$

%
%

\subsection{DSEMs of type $[3^6:3^4.6]_2$} \label{s3.12}

Let $M$ be a DSEM of type $[3^6:3^4.6]_2$. We consider a fixed type path as follows. 

A path $P_{1} = P( \ldots, y_{i},z_{i},z_{i+1}$, $y_{i+1}, \ldots)$ in $M$, say of type $K_{1}'$, indicated by thick black or green paths, shown in Figure 4.12.1. The vertices $y_i's$ and $z_i's$ have the face-sequences $(3^6)$ and $(3^4,6)$ respectively.

Now, construct $M(i,j,k)$ representation of $M$ by first cutting $M$ along a black colored cycle of type $K'_1$ through a vertex with face-sequence $(3^6)$ and then cutting it along the green colored cycle of type $K'_1$.

\begin{picture}(0,0)(-15,35)
\setlength{\unitlength}{4.5mm}


\drawline[AHnb=0](-.2,0)(8.5,0)
\drawline[AHnb=0](6,0)(8.35,4.75)
\drawline[AHnb=0](3,6)(8,6)
\drawline[AHnb=0](0,0)(3,6)


\drawline[AHnb=0](2,0)(5,6)
\drawline[AHnb=0](3,0)(6,6)
\drawline[AHnb=0](5,0)(8,6)


\drawline[AHnb=0](.25,2)(8.5,2)
\drawline[AHnb=0](.25,3)(8.5,3)
\drawline[AHnb=0](.25,5)(8.5,5)
\drawline[AHnb=0](.25,6)(8.5,6)


\drawline[AHnb=0](3,0)(.5,5)

\drawline[AHnb=0](6,0)(3,6)

\drawline[AHnb=0](8.5,1)(6,6)


\drawline[AHnb=0](.5,1)(1,0)

\drawline[AHnb=0](2,0)(2.5,1)(3.5,1)(4,0)

\drawline[AHnb=0](5,0)(5.5,1)(6.5,1)(7,0)

\drawline[AHnb=0](7,0)(8,0)(8.5,.75)

\drawline[AHnb=0](1,2)(2,2)(2.5,3)(2,4)(1,4)

\drawline[AHnb=0](3.5,3)(4,2)(5,2)(5.5,3)(5,4)(4,4)(3.5,3)

\drawline[AHnb=0](7,4)(8,4)(8.4,3.25)

\drawline[AHnb=0](8,2)(8.4,2.75)

\drawline[AHnb=0](2.5,5)(3.5,5)(4,6)

\drawline[AHnb=0](5,6)(5.5,5)(6.5,5)(7,6)

\drawline[AHnb=0](8,6)(8.4,5.25)

\drawline[AHnb=0](6.5,3)(7,2)

\drawline[AHnb=0](2.5,3)(3,2)
\drawline[AHnb=0](5.5,3)(6,2)

\drawline[AHnb=0](4,6)(4.5,5)
\drawline[AHnb=0](7,6)(7.5,5)


\drawline[AHnb=0](3,2)(3.5,1)
\drawline[AHnb=0](6,2)(6.5,1)

\drawline[AHnb=0](4.5,5)(5,4)
\drawline[AHnb=0](7.5,5)(8,4)
\drawline[AHnb=0](.5,5)(1,6)(2,4)

\drawline[AHnb=0](1,2)(.5,3)(2,6)(2.5,5)

\drawline[AHnb=0](0.75,6.5)(1,6)
\drawline[AHnb=0](2,6)(2.25,6.5)

\drawline[AHnb=0](-0.25,-.5)(0,0)(.25,-.5)
\drawline[AHnb=0](0.75,-.5)(1,0)(1.25,-.5)
\drawline[AHnb=0](1.75,-.5)(2,0)(2.25,-.5)
\drawline[AHnb=0](2.75,-.5)(3,0)(3.25,-.5)
\drawline[AHnb=0](3.75,-.5)(4,0)(4.25,-.5)
\drawline[AHnb=0](4.75,-.5)(5,0)(5.25,-.5)
\drawline[AHnb=0](5.75,-.5)(6,0)(6.25,-.5)
\drawline[AHnb=0](6.75,-.5)(7,0)(7.25,-.5)
\drawline[AHnb=0](7.75,-.5)(8,0)(8.25,-.5)

\drawline[AHnb=0](2.75,6.5)(3,6)(3.25,6.5)
\drawline[AHnb=0](5.75,6.5)(6,6)(6.25,6.5)

\drawline[AHnb=0](3.75,6.5)(4,6)
\drawline[AHnb=0](5,6)(5.25,6.5)

\drawline[AHnb=0](6.75,6.5)(7,6)
\drawline[AHnb=0](8,6)(8.25,6.5)

\drawpolygon[fillcolor=black](.25,5)(8.5,5)(8.5,5.1)(.25,5.1)
\drawpolygon[fillcolor=black](.25,3)(8.5,3)(8.5,3.1)(.25,3.1)
\drawpolygon[fillcolor=black](.25,2)(8.5,2)(8.5,2.1)(.25,2.1)
\drawpolygon[fillcolor=black](-.25,0)(8.5,0)(8.5,.1)(-.25,.1)

\drawpolygon[fillcolor=green](-.25,-.5)(0,0)(3,6)(3.25,6.5)(3.35,6.5)(3.1,6)(.1,0)(-.15,-.5)

\drawpolygon[fillcolor=green](1.75,-.5)(2,0)(5,6)(5.25,6.5)(5.35,6.5)(5.1,6)(2.1,0)(1.85,-.5)

\drawpolygon[fillcolor=green](2.75,-.5)(3,0)(6,6)(6.25,6.5)(6.35,6.5)(6.1,6)(3.1,0)(2.85,-.5)

\drawpolygon[fillcolor=green](4.75,-.5)(5,0)(8,6)(8.25,6.5)(8.35,6.5)(8.1,6)(5.1,0)(4.85,-.5)

\drawpolygon[fillcolor=green](5.75,-.5)(6,0)(9,6)(9.25,6.5)(9.35,6.5)(9.1,6)(6.1,0)(5.85,-.5)

\put(.5,-1.5){\scriptsize {\tiny {\bf Figure 4.12.1:} Paths of type $K_1'$ }}

\end{picture}

\begin{picture}(0,0)(-70,33)
\setlength{\unitlength}{4mm}

\drawpolygon(0,0)(9,0)(13.5,9)(4.5,9)

\drawline[AHnb=0](2,0)(5,6)
\drawline[AHnb=0](3,0)(6,6)
\drawline[AHnb=0](5,0)(8,6)
\drawline[AHnb=0](6,0)(9,6)
\drawline[AHnb=0](8,0)(11,6)

\drawline[AHnb=0](5,6)(6,8)
\drawline[AHnb=0](6,6)(7,8)
\drawline[AHnb=0](8,6)(9,8)
\drawline[AHnb=0](9,6)(10,8)
\drawline[AHnb=0](11,6)(12,8)


\drawline[AHnb=0](1.5,2)(10,2)
\drawline[AHnb=0](1.5,3)(10.5,3)
\drawline[AHnb=0](2.5,5)(11.5,5)
\drawline[AHnb=0](3,6)(12,6)
\drawline[AHnb=0](4,8)(13,8)


\drawline[AHnb=0](3,0)(1.5,3)
\drawline[AHnb=0](6,0)(3,6)
\drawline[AHnb=0](9,0)(6,6)
\drawline[AHnb=0](10.5,3)(9,6)

\drawline[AHnb=0](.5,1)(1,0)

\drawline[AHnb=0](2.5,1)(2,0)
\drawline[AHnb=0](2.5,1)(3.5,1)
\drawline[AHnb=0](3.5,1)(4,0)

\drawline[AHnb=0](5,0)(5.5,1)
\drawline[AHnb=0](5.5,1)(6.5,1)
\drawline[AHnb=0](6.5,1)(7,0)

\drawline[AHnb=0](8,0)(8.5,1)
\drawline[AHnb=0](8.5,1)(9.5,1)

\drawline[AHnb=0](1,2)(2,2)
\drawline[AHnb=0](2,2)(2.5,3)
\drawline[AHnb=0](2.5,3)(2,4)

\drawline[AHnb=0](3.5,3)(4,2)
\drawline[AHnb=0](4,2)(5,2)
\drawline[AHnb=0](5,2)(5.5,3)
\drawline[AHnb=0](5.5,3)(5,4)
\drawline[AHnb=0](5,4)(4,4)
\drawline[AHnb=0](4,4)(3.5,3)

\drawline[AHnb=0](6.5,3)(7,4)
\drawline[AHnb=0](7,2)(8,2)
\drawline[AHnb=0](8,2)(8.5,3)
\drawline[AHnb=0](8.5,3)(8,4)
\drawline[AHnb=0](8,4)(7,4)
\drawline[AHnb=0](7,2)(6.5,3)

\drawline[AHnb=0](9.5,3)(10,4)
\drawline[AHnb=0](10,4)(11,4)
\drawline[AHnb=0](9.5,3)(10,2)

\drawline[AHnb=0](2.5,5)(3.5,5)
\drawline[AHnb=0](3.5,5)(4,6)

\drawline[AHnb=0](5,6)(5.5,5)
\drawline[AHnb=0](5.5,5)(6.5,5)
\drawline[AHnb=0](6.5,5)(7,6)

\drawline[AHnb=0](8,6)(8.5,5)
\drawline[AHnb=0](8.5,5)(9.5,5)
\drawline[AHnb=0](9.5,5)(10,6)

\drawline[AHnb=0](11,6)(11.5,5)
\drawline[AHnb=0](10,4)(11,4)
\drawline[AHnb=0](9.5,3)(10,2)


\drawline[AHnb=0](2.5,3)(3,2)
\drawline[AHnb=0](5.5,3)(6,2)
\drawline[AHnb=0](8.5,3)(9,2)

\drawline[AHnb=0](4,6)(4.5,5)
\drawline[AHnb=0](7,6)(7.5,5)
\drawline[AHnb=0](10,6)(10.5,5)


\drawline[AHnb=0](3,2)(3.5,1)
\drawline[AHnb=0](6,2)(6.5,1)
\drawline[AHnb=0](9,2)(9.5,1)

\drawline[AHnb=0](4.5,5)(5,4)
\drawline[AHnb=0](7.5,5)(8,4)
\drawline[AHnb=0](10.5,5)(11,4)

\drawline[AHnb=0](6,8)(6.5,7)(5.5,7)(5,8)(5.5,9)(6,8)(6.5,9)(7,8)
\drawline[AHnb=0](9,8)(9.5,7)(8.5,7)(8,8)(8.5,9)(9,8)(9.5,9)(10,8)
\drawline[AHnb=0](12,8)(12.5,7)(11.5,7)(11,8)(11.5,9)(12,8)(12.5,9)(13,8)

\drawline[AHnb=0](12,8)(12.5,7)
\drawline[AHnb=0](5,8)(4.5,9)

\drawline[AHnb=0](8,8)(7.5,9)(7,8)
\drawline[AHnb=0](11,8)(10.5,9)(10,8)

\drawline[AHnb=0](3.5,7)(4,6)
\drawline[AHnb=0](5.5,7)(6,6)
\drawline[AHnb=0](6.5,7)(7,6)
\drawline[AHnb=0](8.5,7)(9,6)
\drawline[AHnb=0](9.5,7)(10,6)
\drawline[AHnb=0](11.5,7)(12,6)


\put(-.4,-.5){\scriptsize {\tiny $a_{11}$}}
\put(.7,-.5){\scriptsize {\tiny $a_{12}$}}
\put(1.7,-.5){\scriptsize {\tiny $a_{13}$}}
\put(7.85,-.5){\scriptsize {\tiny $a_{1i}$}}
\put(8.8,-.5){\scriptsize {\tiny $a_{11}$}}

\put(-.6,.9){\scriptsize {\tiny $x_{11}$}}
\put(9.8,.9){\scriptsize {\tiny $x_{11}$}}

\put(-.1,2.2){\scriptsize {\tiny $a_{21}$}}
\put(10.4,2){\scriptsize {\tiny $a_{21}$}}

\put(2.55,8){\scriptsize {\tiny $a_{j1}$}}
\put(13.25,8){\scriptsize {\tiny $a_{j1}$}}

\put(3.55,9.5){\scriptsize {\tiny $a_{1(k+1)}$}}
\put(13,9.5){\scriptsize {\tiny $a_{1(k+1)}$}}

\drawpolygon[fillcolor=black](4,8)(0,0)(1,0)(2,0)(2.5,1)(3,0)(3.5,1)(4,0)(5,0)(5.5,1)(6,0)(6.5,1)(7,0)(8,0)(9,2)(2,2)(2.5,3)(3.5,3)(4,4)(4.5,3)(5,4)(5.5,3)(6.5,3)(7,4)(7.5,3)(8,4)(8.5,3)(9.5,3)(10.5,5)(3.5,5)(4,6)(5,6)(5.5,7)(6,6)(6.5,7)(7,6)(8,6)(8.5,7)(9,6)(9.5,7)(10,6)(11,6)(12,8)(3.9,8)(3.9,8.1)(12.15,8.1)(11.1,5.9)(9.9,5.9)(9.5,6.85)(9,5.85)(8.5,6.85)(8.1,5.9)(6.9,5.9)(6.5,6.85)(6,5.85)(5.5,6.85)(5.1,5.9)(4.1,5.9)(3.6,5.1)(10.65,5.1)(9.56,2.9)(8.4,2.9)(8,3.8)(7.5,2.8)(7,3.8)(6.6,2.9)(5.4,2.9)(5,3.8)(4.5,2.8)(4,3.8)(3.6,2.9)(2.6,2.9)(2.15,2.1)(9.15,2.1)(8.1,-.1)(6.9,-.1)(6.5,.8)(6,-.25)(5.5,.8)(5.1,-.1)(3.9,-.1)(3.5,.8)(3,-.25)(2.5,.8)(2.1,-.1)(-.15,-.1)(3.9,8.1)

\put(-2.25,-.2){\scriptsize {\tiny $Q_1$}}
\put(-1.4,2){\scriptsize {\tiny $Q_2$}}
\put(1,8){\scriptsize {\tiny $Q_{j}$}}
\put(1.5,9){\scriptsize {\tiny $Q_{1}$}}

\put(4,-.5){\scriptsize $\ldots$}
\put(5.5,-.5){\scriptsize $\ldots$} 
\put(7.25,9.5){\scriptsize $\ldots$} 
\put(10,9.5){\scriptsize $\ldots$} 

\put(.3,5.8){\scriptsize $.$}
\put(0.2,5.6){\scriptsize $.$}
\put(.1,5.4){\scriptsize $.$}

\put(1.25,-1.5){\scriptsize {\tiny {\bf Figure 4.12.2:} $M(i,j,k)$ }}

\end{picture}

\vspace{4.25cm}

We get admissible relations among $i, j, k$ of $M(i, j, k)$
in the next lemma and its proof follows from the similar arguments as in Lemma \ref{l3.10.1}.

\begin{lem}\label{l3.12.1} A DSEM $M$ of type $[3^6:3^4.6]_2$ admits an $M(i,j,k)$-representation iff the following holds: $(i)$ $j$ even and $i=3m$, $m \in \mathbb{N} $, $(ii)$ $i \geq 9$  for $j=2$ and $i \geq 6$ for $j \geq 4$, $(iii)$ number of vertices of $M(i,j,k) = 4ij/3 \geq 24 $, $(iv)$ if $j=2$ then $ 4 \leq k \leq i-1$, and if $j \geq 4$ then $ 0 \leq k \leq i-1$.
\end{lem}

%
%

\subsection{DSEMs of type $[3^6:3^2.4.12]$} \label{s3.13}

Let $M$ be a DSEM of type $[3^6:3^2.4.12]$.  Consider a fixed type  path $P_{1} = P( \ldots,y_{i},z_{i},z_{i+1},z_{i+2},z_{i+3}, \linebreak y_{i+1}, \ldots)$ in $M$, say of type $L_{1}$, indicated by thick black or green paths, shown in Figure 4.13.1. The vertices $y_i$'s and $z_i$'s have the face-sequences $(3^6)$ and $(3^2,4,12)$ respectively.

\smallskip

An $M(i,j,k)$ representation of $M$ follows by first cutting $M$ along a black colored cycle of type $L_1$ through a vertex with face-sequence $(3^6)$ and then cutting it along a green colored cycle of type $L_1$.

\vspace{-.75cm}

\begin{picture}(0,0)(16,43)
\setlength{\unitlength}{3.35mm}

\drawline[AHnb=0](16,4)(15,4.5)(14,4.5)(13,4)(11,4)(10,4.5)(9,4.5)(8,4)(6,4)(5,4.5)(4.5,4.5)

\drawline[AHnb=0](4.25,.5)(5,.5)(6,1)(8,1)(9,.5)(10,.5)(11,1)(12,1)(13,1)(14,.5)(15,.5)(16,1)(16.5,2)(16.5,3)(16,4)(17.5,4)

\drawline[AHnb=0](16.5,9)(15.5,8.5)(13.5,8.5)(12.5,9)(11.5,9)(10.5,8.5)(8.5,8.5)(7.5,9)(6.5,9)(5.5,8.5)(5,7.5)(5,6.5)(4.5,5.5)(5.5,5.5)(6.5,5)(7.5,5)(8.5,5.5)(10.5,5.5)(11.5,5)(12.5,5)(13.5,5.5)(14.5,5.5)(15.5,5.5)(16.5,5)


\drawline[AHnb=0](6,1)(6.5,2)(6.5,3)(6,4)(8,4)(7.5,3)(7.5,2)(8,1)

\drawline[AHnb=0](11,1)(11.5,2)(11.5,3)(11,4)(13,4)(12.5,3)(12.5,2)(13,1)

\drawline[AHnb=0](8.5,5.5)(9,6.5)(9,7.5)(8.5,8.5)(10.5,8.5)(10,7.5)(10,6.5)(10.5,5.5)

\drawline[AHnb=0](13.5,5.5)(14,6.5)(14,7.5)(13.5,8.5)(15.5,8.5)(15,7.5)(15,6.5)(15.5,5.5)

\drawline[AHnb=0](4.5,7.5)(5,7.5)
\drawline[AHnb=0](4.5,6.5)(5,6.5)

\drawline[AHnb=0](4.5,5.5)(5,4.5)

\drawline[AHnb=0](4.5,8.5)(5.5,8.5)

\drawline[AHnb=0](4.5,8.25)(5,7.5)

\drawline[AHnb=0](4.5,8.75)(5,9.5)(4.35,9.5)
\drawline[AHnb=0](6.5,2)(7.5,2)(7,1)(6.5,2)
\drawline[AHnb=0](11.5,2)(12.5,2)(12,1)(11.5,2)

\drawline[AHnb=0](16.5,2)(17,1)(17.5,1)

\drawline[AHnb=0](17.5,1.75)(16.75,.5)
\drawline[AHnb=0](16.25,.5)(16,1)(17,1)
\drawline[AHnb=0](16.5,2)(17.5,2)

\drawline[AHnb=0](9,6.5)(10,6.5)(9.5,5.5)(9,6.5)
\drawline[AHnb=0](14,6.5)(15,6.5)(14.5,5.5)(14,6.5)

\drawline[AHnb=0](6.5,3)(7.5,3)(7,4)(6.5,3)
\drawline[AHnb=0](11.5,3)(12.5,3)(12,4)(11.5,3)
\drawline[AHnb=0](16.5,3)(17.5,3)(17,4)(16.5,3)

\drawline[AHnb=0](9,7.5)(10,7.5)(9.5,8.5)(9,7.5)
\drawline[AHnb=0](14,7.5)(15,7.5)(14.5,8.5)(14,7.5)

\drawline[AHnb=0](6.5,5)(7.5,5)(7,4)(6.5,5)
\drawline[AHnb=0](11.5,5)(12.5,5)(12,4)(11.5,5)
\drawline[AHnb=0](16.5,5)(17.5,5)(17,4)(16.5,5)

\drawline[AHnb=0](9,9.5)(10,9.5)(9.5,8.5)(9,9.5)
\drawline[AHnb=0](14,9.5)(15,9.5)(14.5,8.5)(14,9.5)

\drawline[AHnb=0](9,4.5)(10,4.5)(9.5,5.5)(9,4.5)
\drawline[AHnb=0](14,4.5)(15,4.5)(14.5,5.5)(14,4.5)
\drawline[AHnb=0](6.5,9)(7.5,9)(7,10)(6.5,9)
\drawline[AHnb=0](11.5,9)(12.5,9)(12,10)(11.5,9)
\drawline[AHnb=0](16.5,9)(17.5,9)(17,10)(16.5,9)


\drawline[AHnb=0](7,10)(8,10)
\drawline[AHnb=0](11,10)(13,10)
\drawline[AHnb=0](16,10)(18,10)
\drawline[AHnb=0](5,6.5)(5.5,5.5)


\drawline[AHnb=0](5,4.5)(5.5,5.5)
\drawline[AHnb=0](6,4)(6.5,5)
\drawline[AHnb=0](10,4.5)(10.5,5.5)
\drawline[AHnb=0](11,4)(11.5,5)
\drawline[AHnb=0](15,4.5)(15.5,5.5)
\drawline[AHnb=0](16,4)(16.5,5)

\drawline[AHnb=0](8,4)(7.5,5)
\drawline[AHnb=0](9,4.5)(8.5,5.5)
\drawline[AHnb=0](13,4)(12.5,5)
\drawline[AHnb=0](14,4.5)(13.5,5.5)

\drawline[AHnb=0](7.5,9)(8,10)(9,9.5)(8.5,8.5)
\drawline[AHnb=0](12.5,9)(13,10)(14,9.5)(13.5,8.5)
\drawline[AHnb=0](17.5,9)(18,10)

\drawline[AHnb=0](5.5,8.5)(5,9.5)(6,10)(6.5,9)
\drawline[AHnb=0](10.5,8.5)(10,9.5)(11,10)(11.5,9)
\drawline[AHnb=0](15.5,8.5)(15,9.5)(16,10)(16.5,9)

\drawline[AHnb=0](6.75,.5)(7,1)(7.25,.5)
\drawline[AHnb=0](11.75,.5)(12,1)(12.25,.5)
\drawline[AHnb=0](16.75,.5)(17,1)(17.25,.5)

\drawline[AHnb=0](4.75,0)(5,.5)(5.25,0)
\drawline[AHnb=0](8.75,0)(9,.5)(9.25,0)
\drawline[AHnb=0](9.75,0)(10,.5)(10.25,0)

\drawline[AHnb=0](13.75,0)(14,.5)(14.25,0)
\drawline[AHnb=0](14.75,0)(15,.5)(15.25,0)

\drawline[AHnb=0](8,1)(7.75,.5)
\drawline[AHnb=0](9,0.5)(8.75,0)
\drawline[AHnb=0](13,1)(12.75,.5)
\drawline[AHnb=0](14,.5)(13.75,0)

\drawline[AHnb=0](6,1)(6.25,.5)
\drawline[AHnb=0](11,1)(11.25,.5)
\drawline[AHnb=0](10,.5)(10.25,0)
\drawline[AHnb=0](15,.5)(15.25,0)

\drawline[AHnb=0](6,10)(7,10)

\drawline[AHnb=0](6.75,10.5)(7,10)(7.25,10.5)
\drawline[AHnb=0](11.75,10.5)(12,10)(12.25,10.5)
\drawline[AHnb=0](16.75,10.5)(17,10)(17.25,10.5)

\drawline[AHnb=0](8,10)(7.75,10.5)
\drawline[AHnb=0](13,10)(12.75,10.5)
\drawline[AHnb=0](18,10)(17.75,10.5)

\drawline[AHnb=0](6,10)(6.25,10.5)
\drawline[AHnb=0](11,10)(11.25,10.5)
\drawline[AHnb=0](16,10)(16.25,10.5)

\drawpolygon[fillcolor=green](6.75,0.5)(7.5,2)(7.5,3)(8,4)(9,4.5)(9.5,5.5)(10,6.5)(10,7.5)(10.5,8.5)(11.5,9)(12,10)(12.25,10.5)(12.4,10.5)(12.15,9.9)(11.65,8.9)(10.65,8.4)(10.15,7.5)(10.15,6.5)(9.7,5.5)(9.15,4.4)(8.15,3.9)(7.65,3)(7.65,2)(6.9,0.5)

\drawpolygon[fillcolor=green](11.75,0.5)(12.5,2)(12.5,3)(13,4)(14,4.5)(14.5,5.5)(15,6.5)(15,7.5)(15.5,8.5)(16.5,9)(17,10)(17.25,10.5)(17.4,10.5)(17.15,9.9)(16.65,8.9)(15.65,8.4)(15.15,7.5)(15.15,6.5)(14.7,5.5)(14.15,4.4)(13.15,3.9)(12.65,3)(12.65,2)(11.9,0.5)

\drawpolygon[fillcolor=black](4.5,5.5)(5.5,5.5)(6.5,5)(7.5,5)(8.5,5.5)(10.5,5.5)(11.5,5)(12.5,5)(13.5,5.5)(14.5,5.5)(15.5,5.5)(16.5,5)(17.5,5)(17.5,5.1)(16.5,5.1)(15.5,5.6)(14.5,5.6)(13.5,5.6)(12.5,5.1)(11.5,5.1)(10.5,5.6)(8.5,5.6)(7.5,5.1)(6.5,5.1)(5.5,5.6)(4.5,5.6)

\drawpolygon[fillcolor=black](4.5,8.5)(5.5,8.5)(6.5,9)(7.5,9)(8.5,8.5)(10.5,8.5)(11.5,9)(12.5,9)(13.5,8.5)(14.5,8.5)(15.5,8.5)(16.5,9)(17.5,9)(17.5,9.1)(16.5,9.1)(15.5,8.6)(14.5,8.6)(13.5,8.6)(12.5,9.1)(11.5,9.1)(10.5,8.6)(8.5,8.6)(7.5,9.1)(6.5,9.1)(5.5,8.6)(4.5,8.6)

\drawpolygon[fillcolor=green](6,1)(6.5,2)(6.5,3)(7.5,5)(8.5,5.5)(9,6.5)(9,7.5)(10,9.5)(11,10)(11.25,10.5)(11.35,10.5)(11.15,9.9)(10.15,9.4)(9.15,7.5)(9.15,6.5)(8.65,5.45)(7.65,4.9)(6.65,3)(6.65,2)(6.15,1)
\drawpolygon[fillcolor=green](11,1)(11.5,2)(11.5,3)(12.5,5)(13.5,5.5)(14,6.5)(14,7.5)(15,9.5)(16,10)(16.25,10.5)(16.35,10.5)(16.15,9.9)(15.15,9.4)(14.15,7.5)(14.15,6.5)(13.65,5.45)(12.65,4.9)(11.65,3)(11.65,2)(11.15,1)

\drawpolygon[fillcolor=black](4,4.5)(5,4.5)(6,4)(7,4)(8,4)(9,4.5)(10,4.5)(11,4)(12,4)(13,4)(14,4.5)(15,4.5)(16,4)(17.5,4)(17.5,4.1)(16,4.1)(15,4.6)(14,4.6)(13,4.1)(12,4.1)(11,4.1)(10,4.6)(9,4.6)(8,4.1)(7,4.1)(6,4.1)(5,4.6)(4,4.6)

\drawpolygon[fillcolor=black](4,9.5)(5,9.5)(6,10)(7,10)(8,10)(9,9.5)(10,9.5)(11,10)(12,10)(13,10)(14,9.5)(15,9.5)(16,10)(18,10)(18,10.1)(16,10.1)(15,9.6)(14,9.6)(13,10.1)(12,10.1)(11,10.1)(10,9.6)(9,9.6)(8,10.1)(7,10.1)(6,10.1)(5,9.6)(4,9.6)

\drawpolygon[fillcolor=black](4,0.5)(5,.5)(6,1)(7,1)(8,1)(9,.5)(10,.5)(11,1)(12,1)(13,1)(14,.5)(15,.5)(16,1)(18,1)(18,1.1)(16,1.1)(15,.6)(14,.6)(13,1.1)(12,1.1)(11,1.1)(10,.6)(9,.6)(8,1.1)(7,1.1)(6,1.1)(5,.6)(4,.6)

\put(6.2,-1.5){\scriptsize {\tiny {\bf Figure 4.13.1:} Paths of type $L_1$ }} 

\end{picture}

\begin{picture}(0,0)(-50.5,41)
\setlength{\unitlength}{4mm}

\drawline[AHnb=0](16,4)(15,4.5)(14,4.5)(13,4)(11,4)(10,4.5)(9,4.5)(8,4)(6,4)(5,4.5)(4,4.5)(3,4)(2.5,3)(2.5,2)(2,1)(3,1)(4,.5)(5,.5)(6,1)(8,1)(9,.5)(10,.5)(11,1)(12,1)(13,1)(14,.5)(15,.5)(16,1)(17,1)(17.5,2)(17.5,3)(18,4)(16,4)(16.5,3)(16.5,2)(16,1)

\drawline[AHnb=0](18.5,8.5)(17.5,9)(16.5,9)(15.5,8.5)(13.5,8.5)(12.5,9)(11.5,9)(10.5,8.5)(8.5,8.5)(7.5,9)(6.5,9)(5.5,8.5)(5,7.5)(5,6.5)(4.5,5.5)(5.5,5.5)(6.5,5)(7.5,5)(8.5,5.5)(10.5,5.5)(11.5,5)(12.5,5)(13.5,5.5)(14.5,5.5)(15.5,5.5)(16.5,5)(17.5,5)(18.5,5.5)(19.5,5.5)(20,6.5)(20,7.5)(20.5,8.5)(18.5,8.5)(19,7.5)(19,6.5)(18.5,5.5)

\drawline[AHnb=0](6,1)(6.5,2)(6.5,3)(6,4)(8,4)(7.5,3)(7.5,2)(8,1)

\drawline[AHnb=0](11,1)(11.5,2)(11.5,3)(11,4)(13,4)(12.5,3)(12.5,2)(13,1)

\drawline[AHnb=0](8.5,5.5)(9,6.5)(9,7.5)(8.5,8.5)(10.5,8.5)(10,7.5)(10,6.5)(10.5,5.5)

\drawline[AHnb=0](13.5,5.5)(14,6.5)(14,7.5)(13.5,8.5)(15.5,8.5)(15,7.5)(15,6.5)(15.5,5.5)


\drawline[AHnb=0](6.5,2)(7.5,2)(7,1)(6.5,2)
\drawline[AHnb=0](11.5,2)(12.5,2)(12,1)(11.5,2)
\drawline[AHnb=0](16.5,2)(17.5,2)(17,1)(16.5,2)

\drawline[AHnb=0](9,6.5)(10,6.5)(9.5,5.5)(9,6.5)
\drawline[AHnb=0](14,6.5)(15,6.5)(14.5,5.5)(14,6.5)
\drawline[AHnb=0](19,6.5)(20,6.5)(19.5,5.5)(19,6.5)

\drawline[AHnb=0](6.5,3)(7.5,3)(7,4)(6.5,3)
\drawline[AHnb=0](11.5,3)(12.5,3)(12,4)(11.5,3)
\drawline[AHnb=0](16.5,3)(17.5,3)(17,4)(16.5,3)

\drawline[AHnb=0](9,7.5)(10,7.5)(9.5,8.5)(9,7.5)
\drawline[AHnb=0](14,7.5)(15,7.5)(14.5,8.5)(14,7.5)
\drawline[AHnb=0](19,7.5)(20,7.5)(19.5,8.5)(19,7.5)

\drawline[AHnb=0](6.5,5)(7.5,5)(7,4)(6.5,5)
\drawline[AHnb=0](11.5,5)(12.5,5)(12,4)(11.5,5)
\drawline[AHnb=0](16.5,5)(17.5,5)(17,4)(16.5,5)

\drawline[AHnb=0](9,9.5)(10,9.5)(9.5,8.5)(9,9.5)
\drawline[AHnb=0](14,9.5)(15,9.5)(14.5,8.5)(14,9.5)
\drawline[AHnb=0](19,9.5)(20,9.5)(19.5,8.5)(19,9.5)
\drawline[AHnb=0](21,10)(22,10)(21.5,9)(21,10)

\drawline[AHnb=0](4,4.5)(5,4.5)(4.5,5.5)(4,4.5)
\drawline[AHnb=0](9,4.5)(10,4.5)(9.5,5.5)(9,4.5)
\drawline[AHnb=0](14,4.5)(15,4.5)(14.5,5.5)(14,4.5)
\drawline[AHnb=0](6.5,9)(7.5,9)(7,10)(6.5,9)
\drawline[AHnb=0](11.5,9)(12.5,9)(12,10)(11.5,9)
\drawline[AHnb=0](16.5,9)(17.5,9)(17,10)(16.5,9)

\drawline[AHnb=0](2.5,2)(3,1)
\drawline[AHnb=0](18,4)(19,4.5)(19.5,5.5)
\drawline[AHnb=0](7,10)(8,10)
\drawline[AHnb=0](11,10)(13,10)
\drawline[AHnb=0](16,10)(18,10)

\drawline[AHnb=0](18,4)(19,4.5)(19.5,5.5)
\drawline[AHnb=0](5,6.5)(5.5,5.5)


\drawline[AHnb=0](5,4.5)(5.5,5.5)
\drawline[AHnb=0](6,4)(6.5,5)
\drawline[AHnb=0](10,4.5)(10.5,5.5)
\drawline[AHnb=0](11,4)(11.5,5)
\drawline[AHnb=0](15,4.5)(15.5,5.5)
\drawline[AHnb=0](16,4)(16.5,5)

\drawline[AHnb=0](8,4)(7.5,5)
\drawline[AHnb=0](9,4.5)(8.5,5.5)
\drawline[AHnb=0](13,4)(12.5,5)
\drawline[AHnb=0](14,4.5)(13.5,5.5)
\drawline[AHnb=0](18,4)(17.5,5)
\drawline[AHnb=0](19,4.5)(18.5,5.5)

\drawline[AHnb=0](7.5,9)(8,10)(9,9.5)(8.5,8.5)
\drawline[AHnb=0](12.5,9)(13,10)(14,9.5)(13.5,8.5)
\drawline[AHnb=0](17.5,9)(18,10)(19,9.5)(18.5,8.5)
\drawline[AHnb=0](10.5,8.5)(10,9.5)(11,10)(11.5,9)
\drawline[AHnb=0](15.5,8.5)(15,9.5)(16,10)(16.5,9)
\drawline[AHnb=0](20.5,8.5)(20,9.5)(21,10)(21.5,9)(20.5,8.5)

\drawpolygon[fillcolor=black](2,1)(3,1)(4,.5)(5,.5)(6,1)(6.5,2)(7,1)(7.5,2)(8,1)(9,.5)(10,.5)(11,1)(11.5,2)(12,1)(12.5,2)(13,1)(14,.5)(15,.5)(16,1)(16.5,2)(16.5,3)(17,4)(16,4)(15,4.5)(14,4.5)(13,4)(12.5,3)(12,4)(11.5,3)(11,4)(10,4.5)(9,4.5)(8,4)(7.5,3)(7,4)(6.5,3)(6,4)(5,4.5)(5.5,5.5)(6.5,5)(7.5,5)(8.5,5.5)(9,6.5)(9.5,5.5)(10,6.5)(10.5,5.5)(11.5,5)(12.5,5)(13.5,5.5)(14,6.5)(14.5,5.5)(15,6.5)(15.5,5.5)(16.5,5)(17.5,5)(18.5,5.5)(19,6.5)(19,7.5)(19.5,8.5)(18.5,8.5)(17.5,9)(16.5,9)(15.5,8.5)(15,7.5)(14.5,8.5)(14,7.5)(13.5,8.5)(12.5,9)(11.5,9)(10.5,8.5)(10,7.5)(9.5,8.5)(9,7.5)(8.5,8.5)(7.5,9)(6.5,9)(6.5,9.1)(7.5,9.1)(8.65,8.6)(9,7.8)(9.5,8.75)(10,7.8)(10.35,8.6)(11.4,9.15)(12.6,9.15)(13.6,8.65)(14,7.75)(14.5,8.75)(15,7.75)(15.4,8.65)(16.5,9.15)(17.5,9.15)(18.6,8.6)(19.75,8.6)(19.15,7.5)(19.15,6.4)(18.6,5.4)(17.6,4.85)(16.5,4.85)(15.4,5.4)(15,6.15)(14.5,5.3)(14,6.15)(13.6,5.4)(12.5,4.85)(11.5,4.85)(10.4,5.35)(10,6.15)(9.5,5.3)(9,6.15)(8.6,5.4)(7.6,4.85)(6.5,4.85)(5.55,5.3)(5.2,4.6)(6.1,4.15)(6.5,3.25)(7,4.25)(7.5,3.25)(7.9,4.1)(9,4.65)(10.1,4.65)(11.1,4.1)(11.5,3.25)(12,4.2)(12.5,3.25)(12.9,4.1)(14,4.65)(15,4.65)(16.1,4.15)(17.2,4.15)(16.65,3)(16.65,2)(16.1,.9)(15.1,0.35)(14,.35)(12.95,.8)(12.5,1.6)(12,.8)(11.52,1.65)(11.1,.85)(10,.35)(9,.35)(7.9,.9)(7.5,1.65)(7,.85)(6.5,1.65)(6.1,.85)(5.1,.35)(4,.35)(3,.85)(1.9,.85)

\drawpolygon[fillcolor=black](6.5,9)(5.5,8.5)(5,7.5)(5,6.5)(4,4.5)(3,4)(2.5,3)(2.5,2)(2,1)(1.9,.9)(2.35,2)(2.35,3)(2.85,4.1)(3.9,4.6)(4.85,6.5)(4.85,7.5)(5.35,8.6)(6.5,9.1)

\put(1.7,.3){\scriptsize {\tiny $a_{11}$}}
\put(2.7,.3){\scriptsize {\tiny $a_{12}$}}
\put(3.8,-.1){\scriptsize {\tiny $a_{13}$}}
\put(4.9,-.1){\scriptsize {\tiny $a_{14}$}}
\put(5.8,.3){\scriptsize {\tiny $a_{15}$}}
\put(6.8,.3){\scriptsize {\tiny $a_{16}$}}
\put(15.9,.3){\scriptsize {\tiny $a_{1i}$}}
\put(16.8,.3){\scriptsize {\tiny $a_{11}$}}

\put(1.8,4){\scriptsize {\tiny $a_{21}$}}
\put(3,4.8){\scriptsize {\tiny $a_{22}$}}
\put(18.15,3.65){\scriptsize {\tiny $a_{21}$}}
\put(19.2,4.4){\scriptsize {\tiny $a_{22}$}}

\put(3.25,5.4){\scriptsize {\tiny $a_{31}$}}
\put(19.8,5.4){\scriptsize {\tiny $a_{31}$}}

\put(1.35,1.8){\scriptsize {\tiny $x_{11}$}}
\put(17.8,1.8){\scriptsize {\tiny $x_{11}$}}

\put(1.35,2.8){\scriptsize {\tiny $x_{12}$}}
\put(17.8,2.8){\scriptsize {\tiny $x_{12}$}}

\put(4.35,8.5){\scriptsize {\tiny $a_{j1}$}}
\put(5.5,9.25){\scriptsize {\tiny $a_{j2}$}}
\put(20.65,8.15){\scriptsize {\tiny $a_{j1}$}}
\put(21.65,8.5){\scriptsize {\tiny $a_{j2}$}}

\put(5,10.5){\scriptsize {\tiny $a_{1(k+1)}$}}
\put(7.45,10.5){\scriptsize {\tiny $a_{1(k+2)}$}}
\put(20.5,10.5){\scriptsize {\tiny $a_{1k}$}}
\put(21.5,10.5){\scriptsize {\tiny $a_{1(k+1)}$}}

\put(-.5,.9){\scriptsize {\tiny $Q_1$}}
\put(.5,4){\scriptsize {\tiny $Q_2$}}
\put(1.2,5.5){\scriptsize {\tiny $Q_{3}$}}
\put(2.4,8.5){\scriptsize {\tiny $Q_{j}$}}
\put(3,10){\scriptsize {\tiny $Q_{1}$}}

\put(8.8,0){\scriptsize $\ldots$}
\put(11.2,0){\scriptsize $\ldots$} 

\put(10.5,10.5){\scriptsize $\ldots$} 
\put(14.2,10.5){\scriptsize $\ldots$} 
\put(17.5,10.5){\scriptsize $\ldots$}

\put(2,7.2){\scriptsize $.$}
\put(1.9,7){\scriptsize $.$}
\put(1.8,6.8){\scriptsize $.$}

\put(7.7,-1){\scriptsize {\tiny {\bf Figure 4.13.2:} $M(i,j,k)$ }} 

\end{picture}

\vspace{4.25cm}

\begin{lem}\label{l3.13.1} A DSEM $M$ of type $[3^6:3^2.4.12]$ admits an $M(i,j,k)$-representation iff the following holds: $(i)$ $j$ even and $i=5m$, $m \in \mathbb{N} $, $(ii)$ $i \geq 15$ for $j=2$, $(ii)$ $i \geq 10$  for $j \geq 4$, $(iii)$ number of vertices of $M(i,j,k) = 7ij/5\geq 42 $, $(iv)$ if $j=2$ then $ k \in \{5r: 1 < r < i/5\}$, and if $j \geq 4$ then $ k \in \{5r: 0 \leq r < i/5\}$.
\end{lem}

\noindent{\bf Proof.} Let $M$ be a DSEM of type $[3^6:3^2.4.12]$ having $n$ vertices. An $M(i,j,k)$  has $j$  disjoint horizontal cycles of $L_{1}$ type, say $Q_0, Q_1,\ldots, Q_{j-1}$, of length $i$. Let $Q_{0}=C(w_{0,0},w_{0,1}, \linebreak \ldots,w_{0,i-1}),Q_{1}=C(w_{1,0},w_{1,1},\ldots,w_{1,i-1}),\ldots,Q_{j-1}=C(w_{j-1,0},w_{j-1,1},\ldots,w_{j-1,i-1})$ be the list of horizontal cycles. Note that the vertices with face-sequence $(3^2,4,12)$ lying between cycles $Q_{2s(mod\,j)}$ and $Q_{(2s+1)(mod\,j)}$ for $0 \leq s \leq j-1$ is $4i/5\cdot j/2$. Therefore, $n = ij + 2ij/5 = 7ij/5$.

If $j=1$, then  no vertex in the base horizontal cycle follows the face-sequence $(3^2,4,12)$ after identifying the boundaries of $M(i,1,k)$. Therefore, $j \geq 2$. If $j \geq 2$ and $j$ is not an even integer, then similarly as above no vertex in the base horizontal cycle follows the face-sequence $(3^2,4,12)$. So $j$ is an even integer. If $j$ is even and $i \neq 5m$, where $m \in \mathbb{N}$, then the ${\rm lk}(w_{0,0})$ is not of type $(3^6)$. Which is not possible. So, $i = 5m$, where $m \in \mathbb{N}$.

If $j=2$ and $ i < 15 $, then we get some vertex in $M(i,2,k)$ whose link can not be constructed. This is not possible. So, for $j=2$, we get $i \geq 15$. Similarly, $i \geq 10$ for $ j \geq 4$. Thus $n = 7ij/5 \geq 42$.

If $j = 2$ and $ k \in \{r:0 \leq r \leq i-1\} \setminus (\{5r: 1 < r < i/5\})$, then we get some vertex in $M(i,2,k)$ whose link can not be constructed. So, if $j=2$, then $ k \in \{5r: 1 < r < i/5\}$. Similarly, as above, we see that if $j \geq 4$, then $ k \in \{5r: 0 \leq r < i/5\}$. This completes the proof. \hfill $\Box$  

%
%
%

\subsection{DSEMs of type $[3.4.3.12:3.12^2]$} \label{s3.14}

Let $M$ be a DSEM of type $[3.4.3.12:3.12^2]$. In $M$, consider a fixed type path $P_{1} = P( \ldots, y_{i},z_{i},z_{i+1}$, $y_{i+1}, \ldots)$, say of type $N_{1}$, indicated by thick black or green paths, shown in Figure 4.14.1. The vertices $y_i$'s and $z_i$'s have the face-sequences $(3,4,3,12)$ and $(3,12^2)$ respectively.

\begin{picture}(0,0)(-6,33)
\setlength{\unitlength}{2.5mm}

\drawline[AHnb=0](-.5,.25)(0,0)(1,0)(2,.5)(3,.5)(4,0)(5,0)(6,.5)(7,.5)(8,0)(9,0)(10,.5)(11,.5)(12,0)(13,0)(14,.5)(14.5,.5)

\drawline[AHnb=0](-.5,4.25)(0,4)(1,4)(2,4.5)(3,4.5)(4,4)(5,4)(6,4.5)(7,4.5)(8,4)(9,4)(10,4.5)(11,4.5)(12,4)(13,4)(14,4.5)(14.5,4.5)

\drawline[AHnb=0](-.5,8.25)(0,8)(1,8)(2,8.5)(3,8.5)(4,8)(5,8)(6,8.5)(7,8.5)(8,8)(9,8)(10,8.5)(11,8.5)(12,8)(13,8)(14,8.5)(14.5,8.5)

\drawline[AHnb=0](-.5,12.25)(0,12)(1,12)(2,12.5)(3,12.5)(4,12)(5,12)(6,12.5)(7,12.5)(8,12)(9,12)(10,12.5)(11,12.5)(12,12)(13,12)(14,12.5)(14.5,12.5)

\drawline[AHnb=0](0,0)(0,1)(.5,2)(.5,3)(0,4)

\drawline[AHnb=0](4,0)(4,1)(4.5,2)(4.5,3)(4,4)

\drawline[AHnb=0](8,0)(8,1)(8.5,2)(8.5,3)(8,4)

\drawline[AHnb=0](12,0)(12,1)(12.5,2)(12.5,3)(12,4)

\drawline[AHnb=0](0,4)(0,5)(.5,6)(.5,7)(0,8)

\drawline[AHnb=0](4,4)(4,5)(4.5,6)(4.5,7)(4,8)

\drawline[AHnb=0](8,4)(8,5)(8.5,6)(8.5,7)(8,8)

\drawline[AHnb=0](12,4)(12,5)(12.5,6)(12.5,7)(12,8)

\drawline[AHnb=0](0,8)(0,9)(.5,10)(.5,11)(0,12)

\drawline[AHnb=0](4,8)(4,9)(4.5,10)(4.5,11)(4,12)

\drawline[AHnb=0](8,8)(8,9)(8.5,10)(8.5,11)(8,12)

\drawline[AHnb=0](12,8)(12,9)(12.5,10)(12.5,11)(12,12)


\drawline[AHnb=0](1,0)(1,1)(.5,2)

\drawline[AHnb=0](5,0)(5,1)(4.5,2)

\drawline[AHnb=0](9,0)(9,1)(8.5,2)

\drawline[AHnb=0](13,0)(13,1)(12.5,2)

\drawline[AHnb=0](1,4)(1,5)(.5,6)

\drawline[AHnb=0](5,4)(5,5)(4.5,6)

\drawline[AHnb=0](9,4)(9,5)(8.5,6)

\drawline[AHnb=0](13,4)(13,5)(12.5,6)

\drawline[AHnb=0](1,8)(1,9)(.5,10)

\drawline[AHnb=0](5,8)(5,9)(4.5,10)

\drawline[AHnb=0](9,8)(9,9)(8.5,10)

\drawline[AHnb=0](13,8)(13,9)(12.5,10)


\drawline[AHnb=0](0,1)(1,1)(2,.5)

\drawline[AHnb=0](3,.5)(4,1)(5,1)(6,.5)

\drawline[AHnb=0](7,.5)(8,1)(9,1)(10,.5)

\drawline[AHnb=0](11,.5)(12,1)(13,1)(14,.5)

\drawline[AHnb=0](0,5)(1,5)(2,4.5)

\drawline[AHnb=0](3,4.5)(4,5)(5,5)(6,4.5)

\drawline[AHnb=0](7,4.5)(8,5)(9,5)(10,4.5)

\drawline[AHnb=0](11,4.5)(12,5)(13,5)(14,4.5)

\drawline[AHnb=0](0,9)(1,9)(2,8.5)

\drawline[AHnb=0](3,8.5)(4,9)(5,9)(6,8.5)

\drawline[AHnb=0](7,8.5)(8,9)(9,9)(10,8.5)

\drawline[AHnb=0](11,8.5)(12,9)(13,9)(14,8.5)


\drawline[AHnb=0](.5,3)(1,4)

\drawline[AHnb=0](4.5,3)(5,4)

\drawline[AHnb=0](8.5,3)(9,4)

\drawline[AHnb=0](12.5,3)(13,4)

\drawline[AHnb=0](.5,7)(1,8)

\drawline[AHnb=0](4.5,7)(5,8)

\drawline[AHnb=0](8.5,7)(9,8)

\drawline[AHnb=0](12.5,7)(13,8)

\drawline[AHnb=0](.5,11)(1,12)

\drawline[AHnb=0](4.5,11)(5,12)

\drawline[AHnb=0](8.5,11)(9,12)

\drawline[AHnb=0](12.5,11)(13,12)


\drawline[AHnb=0](0,1)(-.5,0.75)
\drawline[AHnb=0](0,5)(-.5,4.75)
\drawline[AHnb=0](0,9)(-.5,8.75)


\drawline[AHnb=0](0,0)(.25,-0.5)
\drawline[AHnb=0](1,0)(.75,-0.5)

\drawline[AHnb=0](4,0)(4.25,-0.5)
\drawline[AHnb=0](5,0)(4.75,-0.5)

\drawline[AHnb=0](8,0)(8.25,-0.5)
\drawline[AHnb=0](9,0)(8.75,-0.5)

\drawline[AHnb=0](12,0)(12.25,-0.5)
\drawline[AHnb=0](13,0)(12.75,-0.5)


\drawline[AHnb=0](0,12)(0,12.5)
\drawline[AHnb=0](1,12)(1,12.5)

\drawline[AHnb=0](4,12)(4,12.5)
\drawline[AHnb=0](5,12)(5,12.5)

\drawline[AHnb=0](8,12)(8,12.5)
\drawline[AHnb=0](9,12)(9,12.5)

\drawline[AHnb=0](12,12)(12,12.5)
\drawline[AHnb=0](13,12)(13,12.5)


\drawline[AHnb=0](2,12.5)(1.5,12.75)
\drawline[AHnb=0](6,12.5)(5.5,12.75)
\drawline[AHnb=0](10,12.5)(9.5,12.75)
\drawline[AHnb=0](14,12.5)(13.5,12.75)

\drawline[AHnb=0](3,12.5)(3.5,12.75)
\drawline[AHnb=0](7,12.5)(7.5,12.75)
\drawline[AHnb=0](11,12.5)(11.5,12.75)

\drawpolygon[fillcolor=black](.75,-.5)(1,0)(2,.5)(3,.5)(4,1)(4.5,2)(4.5,3)(5,4)(6,4.5)(7,4.5)(8,5)(8.5,6)(8.5,7)(9,8)(10,8.5)(11,8.5)(12,9)(12.5,10)(12.5,11)(13,12)(14,12.5)(14.5,12.5)(14.5,12.4)(14,12.4)(13.1,12)(12.6,11)(12.6,10)(12.1,9)(11.1,8.4)(10.1,8.4)(9.1,8)(8.6,7)(8.6,6)(8.1,5)(7.1,4.4)(6.1,4.4)(5.1,4)(4.6,3)(4.6,2)(4.1,1)(3.1,.4)(2.1,.4)(1.1,0)(.85,-.5)

\drawpolygon[fillcolor=black](.75,3.5)(1,4)(2,4.5)(3,4.5)(4,5)(4.5,6)(4.5,7)(5,8)(6,8.5)(7,8.5)(8,9)(8.5,10)(8.5,11)(9,12)(10,12.5)(11,12.5)(11.1,12.4)(10.1,12.4)(9.1,12)(8.6,11)(8.6,10)(8.1,9)(7.1,8.4)(6.1,8.4)(5.1,8)(4.6,7)(4.6,6)(4.1,5)(3.1,4.4)(2.1,4.4)(1.1,4)(.85,3.5)

\drawpolygon[fillcolor=black](.75,7.5)(1,8)(2,8.5)(3,8.5)(4,9)(4.5,10)(4.5,11)(5,12)(6,12.5)(7,12.5)(7.1,12.4)(6.1,12.4)(5.1,12)(4.6,11)(4.6,10)(4.1,9)(3.1,8.4)(2.1,8.4)(1.1,8)(.85,7.5)

\drawpolygon[fillcolor=black](8.75,-.5)(9,0)(10,.5)(11,.5)(12,1)(12.5,2)(12.5,3)(13,4)(14,4.5)(15,4.5)(15.1,4.4)(14.1,4.4)(13.1,4)(12.6,3)(12.6,2)(12.1,1)(11.1,.4)(10.1,.4)(9.1,0)(8.85,-.5)

\drawpolygon[fillcolor=black](4.75,-.5)(5,0)(6,.5)(7,.5)(8,1)(8.5,2)(8.5,3)(9,4)(10,4.5)(11,4.5)(12,5)(12.5,6)(12.5,7)(13,8)(14,8.5)(14.45,8.5)(14.45,8.4)(14,8.4)(13,7.9)(12.6,7)(12.6,6)(12,4.9)(11.1,4.4)(10.1,4.4)(9.1,4)(8.6,3)(8.6,2)(8.1,1)(7.1,.4)(6.1,.4)(5.1,0)(4.85,-.5)

\drawpolygon[fillcolor=green](8.25,-.5)(8,0)(7,.5)(6,.5)(5,1)(4.5,2)(4.5,3)(4,4)(3,4.5)(2,4.5)(1,5)(.5,6)(.5,7)(0,8)(-1,8.5)(-.9,8.65)(.15,8.15)(.65,7)(.65,6)(1.15,5.15)(2,4.7)(3,4.7)(4.1,4.25)(4.65,3)(4.65,2)(5.15,1.15)(6,.65)(7.1,.65)(8.15,0.1)(8.35,-.5)

\drawpolygon[fillcolor=green](12.25,-.5)(12,0)(11,.5)(10,.5)(9,1)(8.5,2)(8.5,3)(8,4)(7,4.5)(6,4.5)(5,5)(4.5,6)(4.5,7)(4,8)(3,8.5)(2,8.5)(1,9)(.5,10)(.5,11)(0,12)(-.5,12.25)(-.5,12.45)(.15,12.1)(.65,11)(.65,10)(1.15,9.1)(2,8.65)(3.1,8.65)(4.15,8.15)(4.65,7)(4.65,6)(5.15,5.15)(6,4.7)(7,4.7)(8.1,4.25)(8.65,3)(8.65,2)(9.15,1.15)(10.1,.65)(11.1,.65)(12.15,0.1)(12.35,-.5)

\drawpolygon[fillcolor=green](15,.5)(14,.5)(13,1)(12.5,2)(12.5,3)(12,4)(11,4.5)(10,4.5)(9,5)(8.5,6)(8.5,7)(8,8)(7,8.5)(6,8.5)(5,9)(4.5,10)(4.5,11)(4,12)(3.5,12.25)(3.5,12.45)(4.15,12.1)(4.65,11)(4.65,10)(5.15,9.1)(6,8.65)(7.1,8.65)(8.15,8.15)(8.65,7)(8.65,6)(9.15,5.15)(10,4.7)(11,4.7)(12.1,4.25)(12.65,3)(12.65,2)(13.15,1.15)(14.1,.65)(15.1,.65)

\drawpolygon[fillcolor=green](15,4.5)(14,4.5)(13,5)(12.5,6)(12.5,7)(12,8)(11,8.5)(10,8.5)(9,9)(8.5,10)(8.5,11)(8,12)(7.5,12.25)(7.5,12.45)(8.15,12.1)(8.65,11)(8.65,10)(9.15,9.1)(10,8.65)(11.1,8.65)(12.15,8.15)(12.65,7)(12.65,6)(13.15,5.15)(14,4.7)(15,4.7)

\put(.5,-2.25){\scriptsize {\tiny {\bf Figure 4.14.1:} Paths of type $N_1$ }}

\end{picture}

\begin{picture}(0,0)(-60,22)
\setlength{\unitlength}{3.6mm}


\drawpolygon(0,0)(1,.5)(2,0)(1.5,1)(2,2)(1,1.5)(0,2)(.5,1)
\drawpolygon(.5,1)(1,.5)(1.5,1)(1,1.5)

\drawpolygon(6,0)(7,.5)(8,0)(7.5,1)(8,2)(7,1.5)(6,2)(6.5,1)
\drawpolygon(6.5,1)(7,.5)(7.5,1)(7,1.5)

\drawpolygon(12,0)(13,.5)(14,0)(13.5,1)(14,2)(13,1.5)(12,2)(12.5,1)
\drawpolygon(12.5,1)(13,.5)(13.5,1)(13,1.5)

\drawpolygon(18,0)(19,.5)(20,0)(19.5,1)(20,2)(19,1.5)(18,2)(18.5,1)
\drawpolygon(18.5,1)(19,.5)(19.5,1)(19,1.5)


\drawline[AHnb=0](2,0)(3,-1)(4,-1.5)(5,-1)(6,0)

\drawline[AHnb=0](8,0)(9,-1)(10,-1.5)(11,-1)(12,0)

\drawline[AHnb=0](14,0)(15,-1)(16,-1.5)(17,-1)(18,0)

\drawline[AHnb=0](20,0)(21,-1)(22,-1.5)(23,-1)(24,0)(24.5,1)(24,2)

\drawline[AHnb=0](3,3)(4,3.5)(5,3)
\drawline[AHnb=0](9,3)(10,3.5)(11,3)
\drawline[AHnb=0](15,3)(16,3.5)(17,3)
\drawline[AHnb=0](21,3)(22,3.5)(23,3)


\drawline[AHnb=0](2,2)(3,3)
\drawline[AHnb=0](5,3)(6,2)

\drawline[AHnb=0](8,2)(9,3)
\drawline[AHnb=0](11,3)(12,2)

\drawline[AHnb=0](14,2)(15,3)
\drawline[AHnb=0](17,3)(18,2)

\drawline[AHnb=0](20,2)(21,3)
\drawline[AHnb=0](23,3)(24,2)

\drawpolygon[fillcolor=black](0,2)(.5,1)(1,1.5)(1.5,1)(2,2)(3,3)(4,3.5)(5,3)(6,2)(6.5,1)(7,1.5)(7.5,1)(8,2)(9,3)(10,3.5)(11,3)(12,2)(12.5,1)(13,1.5)(13.5,1)(14,2)(15,3)(16,3.5)(17,3)(18,2)(18.5,1)(19,1.5)(19.5,1)(20,2)(21,3)(22,3.5)(23,3)(24,2)(24.2,2)(23.2,3)(22,3.7)(21,3.15)(19.9,2.1)(19.5,1.2)(19,1.65)(18.5,1.2)(18.1,2.1)(17.1,3.1)(16,3.7)(14.9,3.1)(13.9,2.1)(13.5,1.2)(13,1.65)(12.5,1.2)(12.1,2.1)(11.1,3.1)(10,3.7)(8.9,3.1)(7.9,2.1)(7.5,1.2)(7,1.65)(6.5,1.2)(6.1,2.1)(5.1,3.1)(4,3.7)(2.9,3.1)(1.9,2.1)(1.5,1.2)(1,1.7)(.55,1.2)(.15,2)

\put(-2.5,2.3){\scriptsize {\tiny $a_{1(k+1)}$}}
\put(-.25,2.7){\scriptsize {\tiny $a_{1(k+2)}$}}
\put(2.2,1.8){\scriptsize {\tiny $a_{1(k+3)}$}}
\put(20.5,4){\scriptsize {\tiny $a_{1(k-1)}$}}
\put(23.1,3.3){\scriptsize {\tiny $a_{1k}$}}
\put(24,2.3){\scriptsize {\tiny $a_{1(k+1)}$}}

\put(-.3,-.5){\scriptsize {\tiny $a_{11}$}}
\put(.6,-.2){\scriptsize {\tiny $a_{12}$}}
\put(1.35,-.8){\scriptsize {\tiny $a_{13}$}}
\put(2.4,-1.55){\scriptsize {\tiny $a_{14}$}}
\put(3.5,-2){\scriptsize {\tiny $a_{15}$}}
\put(22.85,-1.5){\scriptsize {\tiny $a_{1i}$}}
\put(23.9,-.5){\scriptsize {\tiny $a_{11}$}}

\put(-.8,.7){\scriptsize {\tiny $x_{11}$}}
\put(24.7,.7){\scriptsize {\tiny $x_{11}$}}

\put(6.5,-2){\scriptsize $\ldots$}
\put(12.5,-2){\scriptsize $\ldots$}
\put(18.5,-2){\scriptsize $\ldots$}

\put(6.5,4){\scriptsize $\ldots$}
\put(12.5,4){\scriptsize $\ldots$}
\put(18.5,4){\scriptsize $\ldots$}

\put(8,-3.4){\scriptsize {\tiny {\bf Figure 4.14.2:} $M(i,1,k)$ }}

\end{picture}

\begin{picture}(0,0)(-8.5,70)
\setlength{\unitlength}{2.3mm}


\drawpolygon(0,0)(1,.5)(2,0)(1.5,1)(2,2)(1,1.5)(0,2)(.5,1)
\drawpolygon(.5,1)(1,.5)(1.5,1)(1,1.5)

\drawpolygon(6,0)(7,.5)(8,0)(7.5,1)(8,2)(7,1.5)(6,2)(6.5,1)
\drawpolygon(6.5,1)(7,.5)(7.5,1)(7,1.5)

\drawpolygon(12,0)(13,.5)(14,0)(13.5,1)(14,2)(13,1.5)(12,2)(12.5,1)
\drawpolygon(12.5,1)(13,.5)(13.5,1)(13,1.5)

\drawpolygon(18,0)(19,.5)(20,0)(19.5,1)(20,2)(19,1.5)(18,2)(18.5,1)
\drawpolygon(18.5,1)(19,.5)(19.5,1)(19,1.5)

\drawpolygon(3,3)(4,3.5)(5,3)(4.5,4)(5,5)(4,4.5)(3,5)(3.5,4)
\drawpolygon(3.5,4)(4,3.5)(4.5,4)(4,4.5)

\drawpolygon(9,3)(10,3.5)(11,3)(10.5,4)(11,5)(10,4.5)(9,5)(9.5,4)
\drawpolygon(9.5,4)(10,3.5)(10.5,4)(10,4.5)

\drawpolygon(15,3)(16,3.5)(17,3)(16.5,4)(17,5)(16,4.5)(15,5)(15.5,4)
\drawpolygon(15.5,4)(16,3.5)(16.5,4)(16,4.5)

\drawpolygon(21,3)(22,3.5)(23,3)(22.5,4)(23,5)(22,4.5)(21,5)(21.5,4)
\drawpolygon(21.5,4)(22,3.5)(22.5,4)(22,4.5)

\drawpolygon(0,6)(1,6.5)(2,6)(1.5,7)(2,8)(1,7.5)(0,8)(.5,7)
\drawpolygon(.5,7)(1,6.5)(1.5,7)(1,7.5)

\drawpolygon(6,6)(7,6.5)(8,6)(7.5,7)(8,8)(7,7.5)(6,8)(6.5,7)
\drawpolygon(6.5,7)(7,6.5)(7.5,7)(7,7.5)

\drawpolygon(12,6)(13,6.5)(14,6)(13.5,7)(14,8)(13,7.5)(12,8)(12.5,7)
\drawpolygon(12.5,7)(13,6.5)(13.5,7)(13,7.5)

\drawpolygon(18,6)(19,6.5)(20,6)(19.5,7)(20,8)(19,7.5)(18,8)(18.5,7)
\drawpolygon(18.5,7)(19,6.5)(19.5,7)(19,7.5)

\drawpolygon(3,9)(4,9.5)(5,9)(4.5,10)(5,11)(4,10.5)(3,11)(3.5,10)
\drawpolygon(3.5,10)(4,9.5)(4.5,10)(4,10.5)

\drawpolygon(9,9)(10,9.5)(11,9)(10.5,10)(11,11)(10,10.5)(9,11)(9.5,10)
\drawpolygon(9.5,10)(10,9.5)(10.5,10)(10,10.5)

\drawpolygon(15,9)(16,9.5)(17,9)(16.5,10)(17,11)(16,10.5)(15,11)(15.5,10)
\drawpolygon(15.5,10)(16,9.5)(16.5,10)(16,10.5)

\drawpolygon(21,9)(22,9.5)(23,9)(22.5,10)(23,11)(22,10.5)(21,11)(21.5,10)
\drawpolygon(21.5,10)(22,9.5)(22.5,10)(22,10.5)


\drawline[AHnb=0](24,6)(24.5,7)(24,8)

\drawline[AHnb=0](2,0)(3,-1)(4,-1.5)(5,-1)(6,0)

\drawline[AHnb=0](8,0)(9,-1)(10,-1.5)(11,-1)(12,0)

\drawline[AHnb=0](14,0)(15,-1)(16,-1.5)(17,-1)(18,0)

\drawline[AHnb=0](20,0)(21,-1)(22,-1.5)(23,-1)(24,0)(24.5,1)(24,2)


\drawline[AHnb=0](2,2)(3,3)
\drawline[AHnb=0](5,3)(6,2)

\drawline[AHnb=0](8,2)(9,3)
\drawline[AHnb=0](11,3)(12,2)

\drawline[AHnb=0](14,2)(15,3)
\drawline[AHnb=0](17,3)(18,2)

\drawline[AHnb=0](20,2)(21,3)
\drawline[AHnb=0](23,3)(24,2)

\drawline[AHnb=0](2,6)(3,5)
\drawline[AHnb=0](5,5)(6,6)

\drawline[AHnb=0](8,6)(9,5)
\drawline[AHnb=0](11,5)(12,6)

\drawline[AHnb=0](14,6)(15,5)
\drawline[AHnb=0](17,5)(18,6)

\drawline[AHnb=0](20,6)(21,5)
\drawline[AHnb=0](23,5)(24,6)

\drawline[AHnb=0](2,8)(3,9)
\drawline[AHnb=0](5,9)(6,8)

\drawline[AHnb=0](8,8)(9,9)
\drawline[AHnb=0](11,9)(12,8)

\drawline[AHnb=0](14,8)(15,9)
\drawline[AHnb=0](17,9)(18,8)

\drawline[AHnb=0](20,8)(21,9)
\drawline[AHnb=0](23,9)(24,8)

\drawline[AHnb=0](0,2)(-1,3)(-1.5,4)(-1,5)(0,6)
\drawline[AHnb=0](0,8)(-1,9)(-1.5,10)(-1,11)(0,12)(1,12.5)(2,12)(3,11)

\drawline[AHnb=0](5,11)(6,12)(7,12.5)(8,12)(9,11)
\drawline[AHnb=0](11,11)(12,12)(13,12.5)(14,12)(15,11)
\drawline[AHnb=0](17,11)(18,12)(19,12.5)(20,12)(21,11)

\drawpolygon[fillcolor=black](0,2)(.5,1)(1,1.5)(1.5,1)(2,2)(3,3)(3.5,4)(4,3.5)(4.5,4)(5,3)(6,2)(6.5,1)(7,1.5)(7.5,1)(8,2)(9,3)(9.5,4)(10,3.5)(10.5,4)(11,3)(12,2)(12.5,1)(13,1.5)(13.5,1)(14,2)(15,3)(15.5,4)(16,3.5)(16.5,4)(17,3)(18,2)(18.5,1)(19,1.5)(19.5,1)(20,2)(21,3)(22,3.5)(21.5,4)(22,4.5)(21,5)(20,6)(19,6.5)(18,6)(17,5)(16,4.5)(15,5)(14,6)(13,6.5)(12,6)(11,5)(10,4.5)(9,5)(8,6)(7,6.5)(6,6)(5,5)(4,4.5)(3,5)(2,6)(1,6.5)(1.5,7)(1,7.5)(2,8)(3,9)(3.5,10)(4,9.5)(4.5,10)(5,9)(6,8)(6.5,7)(7,7.5)(7.5,7)(8,8)(9,9)(9.5,10)(10,9.5)(10.5,10)(11,9)(12,8)(12.5,7)(13,7.5)(13.5,7)(14,8)(15,9)(15.5,10)(16,9.5)(16.5,10)(17,9)(18,8)(18.5,7)(19,7.5)(19.5,7)(20,8)(21,9)(22,9.5)(21.5,10)(22,10.5)(21,11)(20,12)(19,12.5)(18,12)(17,11)(16,10.5)(15,11)(14,12)(13,12.5)(12,12)(11,11)(10,10.5)(9,11)(8,12)(7,12.5)(6,12)(5,11)(4,10.5)(3,11)(2,12)(1,12.5)(0,12)(-1,11)(-1.5,10)(-1,9)(0,8)(.5,7)(0,6)(-1,5)(-1.5,4)(-1,3)(0,2)(.15,2)(-.8,3)(-1.3,4)(-.8,4.9)(.2,5.9)(.7,6.9)(.2,8)(-.8,9)(-1.3,9.9)(-.8,11)(.2,12)(1,12.35)(1.9,11.9)(3,10.8)(4,10.35)(5.1,10.85)(6.1,11.85)(7,12.35)(7.9,11.9)(9,10.8)(10,10.35)(11.1,10.85)(12.1,11.85)(13,12.35)(13.9,11.9)(15,10.8)(16,10.35)(17.1,10.85)(18.1,11.85)(19,12.35)(19.9,11.9)(21,10.8)(21.8,10.45)(21.3,10)(21.8,9.5)(21,9.15)(19.9,8.1)(19.5,7.2)(19,7.65)(18.5,7.2)(18.1,8.1)(17.1,9.1)(16.6,10.2)(16,9.7)(15.5,10.2)(14.9,9.1)(13.9,8.1)(13.5,7.2)(13,7.65)(12.5,7.2)(12.1,8.1)(11.1,9.1)(10.6,10.2)(10,9.7)(9.5,10.2)(8.9,9.1)(7.9,8.1)(7.5,7.2)(7,7.65)(6.5,7.2)(6.1,8.1)(5.1,9.1)(4.6,10.2)(4,9.7)(3.5,10.2)(2.9,9.1)(1.9,8.1)(.8,7.55)(1.3,7)(.8,6.5)(1.8,5.9)(3,4.8)(4,4.35)(5.1,4.85)(6.1,5.85)(7.1,6.3)(8,5.85)(9,4.8)(10,4.35)(11.1,4.85)(12.1,5.85)(13,6.3)(14,5.85)(15,4.8)(16,4.35)(17.1,4.85)(18.1,5.85)(19,6.35)(20,5.85)(21,4.8)(21.8,4.45)(21.3,4)(21.8,3.5)(21,3.15)(19.9,2.1)(19.5,1.2)(19,1.65)(18.5,1.2)(18.1,2.1)(17.1,3.1)(16.6,4.2)(16,3.7)(15.5,4.2)(14.9,3.1)(13.9,2.1)(13.5,1.2)(13,1.65)(12.5,1.2)(12.1,2.1)(11.1,3.1)(10.6,4.2)(10,3.7)(9.5,4.2)(8.9,3.1)(7.9,2.1)(7.5,1.2)(7,1.65)(6.5,1.2)(6.1,2.1)(5.1,3.1)(4.6,4.2)(4,3.7)(3.5,4.2)(2.9,3.1)(1.9,2.1)(1.5,1.2)(1,1.7)(.55,1.2)(.15,2)

\put(-1,-.5){\scriptsize {\tiny $a_{11}$}}
\put(.4,-.2){\scriptsize {\tiny $a_{12}$}}
\put(2.1,0){\scriptsize {\tiny $a_{13}$}}
\put(1.7,-1.6){\scriptsize {\tiny $a_{14}$}}
\put(3.5,-2.1){\scriptsize {\tiny $a_{15}$}}
\put(22.85,-1.5){\scriptsize {\tiny $a_{1i}$}}
\put(23.9,-.5){\scriptsize {\tiny $a_{11}$}}

\put(-1.4,.9){\scriptsize {\tiny $x_{11}$}}
\put(24.65,.7){\scriptsize {\tiny $x_{11}$}}

\put(-2.5,2.5){\scriptsize {\tiny $a_{21}$}}
\put(-1.75,1.7){\scriptsize {\tiny $a_{22}$}}
\put(21.5,2.7){\scriptsize {\tiny $a_{2i}$}}
\put(23.15,2.9){\scriptsize {\tiny $a_{21}$}}
\put(24,2.1){\scriptsize {\tiny $a_{22}$}}

\put(-2.5,8.4){\scriptsize {\tiny $a_{j1}$}}
\put(-1.65,7.5){\scriptsize {\tiny $a_{j2}$}}

\put(23,9){\scriptsize {\tiny $a_{j1}$}}
\put(24,8){\scriptsize {\tiny $a_{j2}$}}


\put(-3.15,9.8){\scriptsize {\tiny $x_{j1}$}}
\put(23,9.8){\scriptsize {\tiny $x_{j1}$}}

\put(-4.65,11.4){\scriptsize {\tiny $a_{1(k+1)}$}}
\put(-3.5,12.5){\scriptsize {\tiny $a_{1(k+2)}$}}
\put(21.2,11.1){\scriptsize {\tiny $a_{1k}$}}
\put(23.1,11){\scriptsize {\tiny $a_{1(k+1)}$}}

\put(-4,-.2){\scriptsize {\tiny $Q_1$}}

\put(-4,2.5){\scriptsize {\tiny $Q_2$}}


\put(-3.8,5.2){\scriptsize $\vdots$}


\put(-4.2,8.5){\scriptsize {\tiny $Q_{j}$}}

\put(6,-2){\scriptsize $\ldots$}

\put(11.5,-2){\scriptsize $\ldots$}
\put(17.5,-2){\scriptsize $\ldots$}

\put(5,-3.5){\scriptsize {\tiny {\bf Figure 4.14.3:} $M(i,j=2m,k)$ }}


\end{picture}

\begin{picture}(0,0)(-83,65.5)
\setlength{\unitlength}{2.25mm}


\drawpolygon(0,0)(1,.5)(2,0)(1.5,1)(2,2)(1,1.5)(0,2)(.5,1)
\drawpolygon(.5,1)(1,.5)(1.5,1)(1,1.5)

\drawpolygon(6,0)(7,.5)(8,0)(7.5,1)(8,2)(7,1.5)(6,2)(6.5,1)
\drawpolygon(6.5,1)(7,.5)(7.5,1)(7,1.5)

\drawpolygon(12,0)(13,.5)(14,0)(13.5,1)(14,2)(13,1.5)(12,2)(12.5,1)
\drawpolygon(12.5,1)(13,.5)(13.5,1)(13,1.5)

\drawpolygon(18,0)(19,.5)(20,0)(19.5,1)(20,2)(19,1.5)(18,2)(18.5,1)
\drawpolygon(18.5,1)(19,.5)(19.5,1)(19,1.5)

\drawpolygon(3,3)(4,3.5)(5,3)(4.5,4)(5,5)(4,4.5)(3,5)(3.5,4)
\drawpolygon(3.5,4)(4,3.5)(4.5,4)(4,4.5)

\drawpolygon(9,3)(10,3.5)(11,3)(10.5,4)(11,5)(10,4.5)(9,5)(9.5,4)
\drawpolygon(9.5,4)(10,3.5)(10.5,4)(10,4.5)

\drawpolygon(15,3)(16,3.5)(17,3)(16.5,4)(17,5)(16,4.5)(15,5)(15.5,4)
\drawpolygon(15.5,4)(16,3.5)(16.5,4)(16,4.5)

\drawpolygon(21,3)(22,3.5)(23,3)(22.5,4)(23,5)(22,4.5)(21,5)(21.5,4)
\drawpolygon(21.5,4)(22,3.5)(22.5,4)(22,4.5)

\drawpolygon(0,6)(1,6.5)(2,6)(1.5,7)(2,8)(1,7.5)(0,8)(.5,7)
\drawpolygon(.5,7)(1,6.5)(1.5,7)(1,7.5)

\drawpolygon(6,6)(7,6.5)(8,6)(7.5,7)(8,8)(7,7.5)(6,8)(6.5,7)
\drawpolygon(6.5,7)(7,6.5)(7.5,7)(7,7.5)

\drawpolygon(12,6)(13,6.5)(14,6)(13.5,7)(14,8)(13,7.5)(12,8)(12.5,7)
\drawpolygon(12.5,7)(13,6.5)(13.5,7)(13,7.5)

\drawpolygon(18,6)(19,6.5)(20,6)(19.5,7)(20,8)(19,7.5)(18,8)(18.5,7)
\drawpolygon(18.5,7)(19,6.5)(19.5,7)(19,7.5)

\drawpolygon(3,9)(4,9.5)(5,9)(4.5,10)(5,11)(4,10.5)(3,11)(3.5,10)
\drawpolygon(3.5,10)(4,9.5)(4.5,10)(4,10.5)

\drawpolygon(9,9)(10,9.5)(11,9)(10.5,10)(11,11)(10,10.5)(9,11)(9.5,10)
\drawpolygon(9.5,10)(10,9.5)(10.5,10)(10,10.5)

\drawpolygon(15,9)(16,9.5)(17,9)(16.5,10)(17,11)(16,10.5)(15,11)(15.5,10)
\drawpolygon(15.5,10)(16,9.5)(16.5,10)(16,10.5)

\drawpolygon(21,9)(22,9.5)(23,9)(22.5,10)(23,11)(22,10.5)(21,11)(21.5,10)
\drawpolygon(21.5,10)(22,9.5)(22.5,10)(22,10.5)

\drawpolygon(0,12)(1,12.5)(2,12)(1.5,13)(2,14)(1,13.5)(0,14)(.5,13)
\drawpolygon(.5,13)(1,12.5)(1.5,13)(1,13.5)

\drawpolygon(6,12)(7,12.5)(8,12)(7.5,13)(8,14)(7,13.5)(6,14)(6.5,13)
\drawpolygon(6.5,13)(7,12.5)(7.5,13)(7,13.5)

\drawpolygon(12,12)(13,12.5)(14,12)(13.5,13)(14,14)(13,13.5)(12,14)(12.5,13)
\drawpolygon(12.5,13)(13,12.5)(13.5,13)(13,13.5)

\drawpolygon(18,12)(19,12.5)(20,12)(19.5,13)(20,14)(19,13.5)(18,14)(18.5,13)
\drawpolygon(18.5,13)(19,12.5)(19.5,13)(19,13.5)


\drawline[AHnb=0](24,6)(24.5,7)(24,8)

\drawline[AHnb=0](2,0)(3,-1)(4,-1.5)(5,-1)(6,0)

\drawline[AHnb=0](8,0)(9,-1)(10,-1.5)(11,-1)(12,0)

\drawline[AHnb=0](14,0)(15,-1)(16,-1.5)(17,-1)(18,0)

\drawline[AHnb=0](20,0)(21,-1)(22,-1.5)(23,-1)(24,0)(24.5,1)(24,2)

\drawline[AHnb=0](2,14)(3,15)(4,15.5)(5,15)(6,14)

\drawline[AHnb=0](8,14)(9,15)(10,15.5)(11,15)(12,14)

\drawline[AHnb=0](14,14)(15,15)(16,15.5)(17,15)(18,14)

\drawline[AHnb=0](20,14)(21,15)(22,15.5)(23,15)(24,14)(24.5,13)(24,12)(23,11)


\drawline[AHnb=0](2,2)(3,3)
\drawline[AHnb=0](5,3)(6,2)

\drawline[AHnb=0](8,2)(9,3)
\drawline[AHnb=0](11,3)(12,2)

\drawline[AHnb=0](14,2)(15,3)
\drawline[AHnb=0](17,3)(18,2)

\drawline[AHnb=0](20,2)(21,3)
\drawline[AHnb=0](23,3)(24,2)

\drawline[AHnb=0](2,6)(3,5)
\drawline[AHnb=0](5,5)(6,6)

\drawline[AHnb=0](8,6)(9,5)
\drawline[AHnb=0](11,5)(12,6)

\drawline[AHnb=0](14,6)(15,5)
\drawline[AHnb=0](17,5)(18,6)

\drawline[AHnb=0](20,6)(21,5)
\drawline[AHnb=0](23,5)(24,6)

\drawline[AHnb=0](2,8)(3,9)
\drawline[AHnb=0](5,9)(6,8)

\drawline[AHnb=0](8,8)(9,9)
\drawline[AHnb=0](11,9)(12,8)

\drawline[AHnb=0](14,8)(15,9)
\drawline[AHnb=0](17,9)(18,8)

\drawline[AHnb=0](20,8)(21,9)
\drawline[AHnb=0](23,9)(24,8)

\drawline[AHnb=0](0,2)(-1,3)(-1.5,4)(-1,5)(0,6)
\drawline[AHnb=0](0,8)(-1,9)(-1.5,10)(-1,11)(0,12)(1,12.5)(2,12)(3,11)

\drawline[AHnb=0](5,11)(6,12)(7,12.5)(8,12)(9,11)
\drawline[AHnb=0](11,11)(12,12)(13,12.5)(14,12)(15,11)
\drawline[AHnb=0](17,11)(18,12)(19,12.5)(20,12)(21,11)

\drawpolygon[fillcolor=black](2,14)(3,15)(4,15.5)(5,15)(6,14)(7,13.5)(8,14)(9,15)(10,15.5)(11,15)(12,14)(13,13.5)(14,14)(15,15)(16,15.5)(17,15)(18,14)(19,13.5)(20,14)(21,15)(22,15.5)(23,15)(24,14)(24.5,13)(24,12)(23,11)(23.15,11)(24.15,12)(24.65,13)(24.15,14)(23.15,15)(22,15.65)(21,15.15)(20,14.15)(19,13.65)(18,14.15)(17,15.15)(16,15.65)(15,15.15)(14,14.15)(13,13.65)(12,14.15)(11,15.15)(10,15.65)(9,15.15)(8,14.15)(7,13.65)(6,14.15)(5,15.15)(4,15.65)(3,15.15)(2,14.15)

\drawpolygon[fillcolor=black](.5,1)(1,1.5)(1.5,1)(2,2)(3,3)(3.5,4)(4,3.5)(4.5,4)(5,3)(6,2)(6.5,1)(7,1.5)(7.5,1)(8,2)(9,3)(9.5,4)(10,3.5)(10.5,4)(11,3)(12,2)(12.5,1)(13,1.5)(13.5,1)(14,2)(15,3)(15.5,4)(16,3.5)(16.5,4)(17,3)(18,2)(18.5,1)(19,1.5)(19.5,1)(20,2)(21,3)(22,3.5)(21.5,4)(22,4.5)(21,5)(20,6)(19,6.5)(18,6)(17,5)(16,4.5)(15,5)(14,6)(13,6.5)(12,6)(11,5)(10,4.5)(9,5)(8,6)(7,6.5)(6,6)(5,5)(4,4.5)(3,5)(2,6)(1,6.5)(1.5,7)(1,7.5)(2,8)(3,9)(3.5,10)(4,9.5)(4.5,10)(5,9)(6,8)(6.5,7)(7,7.5)(7.5,7)(8,8)(9,9)(9.5,10)(10,9.5)(10.5,10)(11,9)(12,8)(12.5,7)(13,7.5)(13.5,7)(14,8)(15,9)(15.5,10)(16,9.5)(16.5,10)(17,9)(18,8)(18.5,7)(19,7.5)(19.5,7)(20,8)(21,9)(22,9.5)(21.5,10)(22,10.5)(21,11)(20,12)(19.5,13)(19,12.5)(18.5,13)(18,12)(17,11)(16,10.5)(15,11)(14,12)(13.5,13)(13,12.5)(12.5,13)(12,12)(11,11)(10,10.5)(9,11)(8,12)(7.5,13)(7,12.5)(6.5,13)(6,12)(5,11)(4,10.5)(3,11)(2,12)(1,12.5)(1.5,13)(1,13.5)(2,14)(2,14.15)(.9,13.5)(1.35,13)(.85,12.5)(1.9,11.9)(3,10.8)(4,10.35)(5.1,10.85)(6.1,11.85)(6.55,12.75)(7,12.35)(7.5,12.75)(7.9,11.9)(9,10.8)(10,10.35)(11.1,10.85)(12.1,11.85)(12.55,12.75)(13,12.35)(13.5,12.75)(13.9,11.9)(15,10.8)(16,10.35)(17.1,10.85)(18.1,11.85)(18.55,12.75)(19,12.35)(19.5,12.75)(19.9,11.9)(21,10.8)(21.8,10.45)(21.3,10)(21.8,9.5)(21,9.15)(19.9,8.1)(19.5,7.2)(19,7.65)(18.5,7.2)(18.1,8.1)(17.1,9.1)(16.6,10.2)(16,9.7)(15.5,10.2)(14.9,9.1)(13.9,8.1)(13.5,7.2)(13,7.65)(12.5,7.2)(12.1,8.1)(11.1,9.1)(10.6,10.2)(10,9.7)(9.5,10.2)(8.9,9.1)(7.9,8.1)(7.5,7.2)(7,7.65)(6.5,7.2)(6.1,8.1)(5.1,9.1)(4.6,10.2)(4,9.7)(3.5,10.2)(2.9,9.1)(1.9,8.1)(.8,7.55)(1.3,7)(.8,6.5)(1.8,5.9)(3,4.8)(4,4.35)(5.1,4.85)(6.1,5.85)(7.1,6.3)(8,5.85)(9,4.8)(10,4.35)(11.1,4.85)(12.1,5.85)(13,6.3)(14,5.85)(15,4.8)(16,4.35)(17.1,4.85)(18.1,5.85)(19,6.35)(20,5.85)(21,4.8)(21.8,4.45)(21.3,4)(21.8,3.5)(21,3.15)(19.9,2.1)(19.5,1.2)(19,1.65)(18.5,1.2)(18.1,2.1)(17.1,3.1)(16.6,4.2)(16,3.7)(15.5,4.2)(14.9,3.1)(13.9,2.1)(13.5,1.2)(13,1.65)(12.5,1.2)(12.1,2.1)(11.1,3.1)(10.6,4.2)(10,3.7)(9.5,4.2)(8.9,3.1)(7.9,2.1)(7.5,1.2)(7,1.65)(6.5,1.2)(6.1,2.1)(5.1,3.1)(4.6,4.2)(4,3.7)(3.5,4.2)(2.9,3.1)(1.9,2.1)(1.5,1.2)(1,1.7)(.5,1.2)

\drawpolygon[fillcolor=black](24,14)(24.5,13)(24,12)(23,11)(22.5,10)(23,9)(24,8)(24.5,7)(24,6)(23,5)(22.5,4)(23,3)(24,2)(24.5,1)(24.65,1)(24.15,2)(23.15,3)(22.65,4)(23.15,5)(24.15,6)(24.65,7)(24.15,8)(23.15,9)(22.65,10)(23.15,11)(24.15,12)(24.65,13)(24.15,14)

\put(-1,-.5){\scriptsize {\tiny $a_{11}$}}
\put(.4,-.2){\scriptsize {\tiny $a_{12}$}}
\put(2.1,0){\scriptsize {\tiny $a_{13}$}}
\put(1.6,-1.6){\scriptsize {\tiny $a_{14}$}}
\put(3.5,-2){\scriptsize {\tiny $a_{15}$}}
\put(23,-1.5){\scriptsize {\tiny $a_{1i}$}}
\put(24,-.5){\scriptsize {\tiny $a_{11}$}}

\put(-1.4,.9){\scriptsize {\tiny $x_{11}$}}
\put(24.85,.7){\scriptsize {\tiny $x_{11}$}}

\put(-2.5,2.5){\scriptsize {\tiny $a_{21}$}}
\put(-1.75,1.7){\scriptsize {\tiny $a_{22}$}}
\put(21.5,2.7){\scriptsize {\tiny $a_{2i}$}}
\put(23.35,2.9){\scriptsize {\tiny $a_{21}$}}
\put(24.15,2.1){\scriptsize {\tiny $a_{22}$}}


\put(-1.25,12.95){\scriptsize {\tiny $x_{j1}$}}
\put(24.7,12.85){\scriptsize {\tiny $x_{j1}$}}

\put(-2.65,11.3){\scriptsize {\tiny $a_{j1}$}}
\put(-1.85,12){\scriptsize {\tiny $a_{j2}$}}
\put(23.4,10.7){\scriptsize {\tiny $a_{j1}$}}
\put(24.35,11.7){\scriptsize {\tiny $a_{j2}$}}

\put(-4.5,-.2){\scriptsize {\tiny $Q_1$}}
\put(-4.5,2.65){\scriptsize {\tiny $Q_2$}}
\put(-4.5,11.25){\scriptsize {\tiny $Q_{j}$}}
\put(-4.5,13.85){\scriptsize {\tiny $Q_{1}$}}

\put(6.25,-2){\scriptsize $\ldots$}
\put(12,-2){\scriptsize $\ldots$}
\put(18,-2){\scriptsize $\ldots$}

\put(6.25,16){\scriptsize $\ldots$}
\put(12,16){\scriptsize $\ldots$}
\put(18,16){\scriptsize $\ldots$}


\put(-3,14.5){\scriptsize {\tiny $a_{1(k+1)}$}}
\put(-.5,15.75){\scriptsize {\tiny $a_{1(k+2)}$}}
\put(22.75,15.2){\scriptsize {\tiny $a_{1k}$}}
\put(24,14.3){\scriptsize {\tiny $a_{1(k+1)}$}}

\put(-3.8,6.2){\scriptsize $\vdots$}

\put(5.2,-3.75){\scriptsize {\tiny {\bf Figure 4.14.4:} $M(i,j=2m+1,k)$ }}


\end{picture}

\vspace{7.75cm}

Now construct an $M(i,j,k)$ representation of $M$, as in previous subsections, by cutting $M$ through a vertex $v$ with the face-sequence $(3, 12^2)$ one by one along the cycles of type $N_1$. Assume that the beginning adjacent face to the base horizontal black cycle is 3-gon.

\begin{lem}\label{l3.14.1} A DSEM $M$ of type $[3.4.3.12:3.12^2]$ admits an $M(i,j,k)$-representation iff the following holds: $(i)$ $j \geq 1$ and $i=6m$, $m \in \mathbb{N}$, $(ii)$ $i \geq 24$ for $j=1$, and $i \geq 18$ for $j=2$,  $(iii)$ $i \geq 12$  for $j > 2$, $(iv)$ number of vertices of $M(i,j,k) = 8ij/6 \geq 32 $, $(v)$ if $j=1$ then $ k \in \{6r+3: 1 \leq r < (i-6)/6\}$, if $j = 2$ then $ k \in \{6r+5: 1 \leq r < (i-12)/6\}$, if $j = 2m+1$ then $ k \in \{6r+3: 0 \leq r \leq (i-6)/6\}$, and if $j = 2m+2$ then $ k \in \{6r+5: 0 \leq r \leq (i-6)/6\}$.
\end{lem}

\noindent{\bf Proof.} Let $M$ be a DSEM of type $[3.4.3.12:3.12^2]$ having $n$ vertices. An $M(i,j,k)$ of $M$ has $j$ disjoint horizontal cycles of $N_{1}$ type, say $Q_0, Q_1,\ldots, Q_{j-1}$, of length $i$. Let $Q_{0}=C(w_{0,0},w_{0,1}, \linebreak \ldots,w_{0,i-1}),Q_{1}=C(w_{1,0},w_{1,1},\ldots,w_{1,i-1}),\ldots,Q_{j-1}=C(w_{j-1,0},w_{j-1,1},\ldots,w_{j-1,i-1})$ be the list of horizontal cycles. Observe that the vertices having face-sequence $(3,4,3,12)$ which are lie above the cycles $Q_s$ and not belonging to any of these cycles $Q_s$ for $0 \leq s \leq j-1$ is $2i/6 \cdot j$. So, the total number of vertices in $M$ is $n = ij + 2ij/6 = 8ij/6$.

If $j = 1$ and $i \neq 6m$, where $m \in \mathbb{N}$, then the face-sequence of $w_{0,0}$ is not of the type $(3,12^2)$. This is not possible. So, if $j = 1$, then $i=6m$, where $m \in \mathbb{N}$. Similarly, if $j > 1 $, then $i=6m$.

If $j = 1$ and $ i < 24 $, then link of some vertices can not be completed. So, for $j = 1$, we get $i \geq 24$. Similarly, as above, we get that $ i \geq 18$ if $j = 2$ and $i \geq 12$  for $j > 2$. Thus, $ n = 8ij/6 \geq 32$.

If $j =1$ and $ k \in \{r:0 \leq r \leq i-1\} \setminus ( k \in \{6r+3: 1 \leq r < (i-6)/6\})$ then some vertices do not have the given face-sequences. So, $ k \in \{6r+3: 1 \leq r < (i-6)/6\}$ for $j = 1$. Proceeding similarly, we see if $j=2$ then $ k \in \{6r+5: 1 \leq r < (i-12)/6\}$, if $j = 2m+1$ then $ k \in \{6r+3: 0 \leq r \leq (i-6)/6\}$, and if $j = 2m+2$ then $ k \in \{6r+5: 0 \leq r \leq (i-6)/6\}$. Thus the proof. \hfill $\Box$ 

%
%
%


\subsection{DSEMs of type $[3.4.6.4:4.6.12]$} \label{s3.15}

Let $M$ be a DSEM  of type $[3.4.6.4:4.6.12]$. In $M$, consider a fixed type path $P_{1} = P( \ldots, y_{i},y_{i+1}$, $y_{i+2},y_{i+3},z_{i},z_{i+1}, \ldots)$ in $M$, say of type $O_{1}$, indicated by thick black or green paths, shown in Figure 4.15.1. The vertices $y_i$'s and $z_i$'s have the face-sequences $(4,6,12)$ and $(3,4,6,4)$ respectively.

\begin{picture}(0,0)(-17,47)
\setlength{\unitlength}{2.75mm}

\drawline[AHnb=0](-1.5,2.75)(-1,3)(0,3)(1,2.5)(2,3)(3,3)(4,2.5)(5,3)(6,3)(7,2.5)(8,3)(9,3)(10,2.5)(11,3)(11,4)(10,4.5)(9,4)(8,4)(7,4.5)(6,4)(5,4)(4,4.5)(3,4)(2,4)(1,4.5)(0,4)(-1,4)(-1.5,4.25)

\drawline[AHnb=0](-1.5,7.75)(-1,8)(0,8)(1,7.5)(2,8)(3,8)(4,7.5)(5,8)(6,8)(7,7.5)(8,8)(9,8)(10,7.5)(11,8)(11,9)(10,9.5)(9,9)(8,9)(7,9.5)(6,9)(5,9)(4,9.5)(3,9)(2,9)(1,9.5)(0,9)(-1,9)(-1.5,9.25)

\drawline[AHnb=0](-1.5,12.75)(-1,13)(0,13)(1,12.5)(2,13)(3,13)(4,12.5)(5,13)(6,13)(7,12.5)(8,13)(9,13)(10,12.5)(11,13)(11,14)


\drawline[AHnb=0](2,3)(2,4)(2.5,5)(3,4)(3,3)
\drawline[AHnb=0](8,3)(8,4)(8.5,5)(9,4)(9,3)

\drawline[AHnb=0](-1,8)(-1,9)(-.5,10)(0,9)(0,8)
\drawline[AHnb=0](5,8)(5,9)(5.5,10)(6,9)(6,8)
\drawline[AHnb=0](11,8)(11,9)(11.5,10)(12,9.25)

\drawline[AHnb=0](2,13)(2,14)(2.25,14.5)
\drawline[AHnb=0](3,13)(3,14)(2.75,14.5)

\drawline[AHnb=0](8,13)(8,14)(8.25,14.5)
\drawline[AHnb=0](9,13)(9,14)(8.75,14.5)

\drawline[AHnb=0](-1.25,14.25)(-1,14)(0,14)(.25,14.25)
\drawline[AHnb=0](1.75,14.25)(2,14)(3,14)(3.25,14.25)
\drawline[AHnb=0](4.75,14.25)(5,14)(6,14)(6.25,14.25)
\drawline[AHnb=0](7.75,14.25)(8,14)(9,14)(9.25,14.25)
\drawline[AHnb=0](10.75,14.25)(11,14)(12,14)

\drawline[AHnb=0](11,13)(12,13)
\drawline[AHnb=0](10.5,11.5)(11.5,12)(12,11.75)
\drawline[AHnb=0](10.5,10.5)(11.5,10)(12,10.25)
\drawline[AHnb=0](11,9)(12,9)
\drawline[AHnb=0](11,8)(12,8)
\drawline[AHnb=0](11,3)(12,3)
\drawline[AHnb=0](11,4)(12,4)

\drawline[AHnb=0](-1,1.75)(-.5,2)(0,1.75)
\drawline[AHnb=0](5,1.75)(5.5,2)(6,1.75)
\drawline[AHnb=0](11,1.75)(11.5,2)(12,1.75)

\drawline[AHnb=0](1,2.5)(.75,2)
\drawline[AHnb=0](7,2.5)(6.75,2)

\drawline[AHnb=0](4,2.5)(4.25,2)
\drawline[AHnb=0](10,2.5)(10.25,2)


\drawline[AHnb=0](-1,4)(-1,3)(-.5,2)(0,3)(0,4)
\drawline[AHnb=0](5,4)(5,3)(5.5,2)(6,3)(6,4)
\drawline[AHnb=0](11,4)(11,3)(11.5,2)(12,2.75)

\drawline[AHnb=0](2,9)(2,8)(2.5,7)(3,8)(3,9)
\drawline[AHnb=0](8,9)(8,8)(8.5,7)(9,8)(9,9)

\drawline[AHnb=0](-1,14)(-1,13)(-.5,12)(0,13)(0,14)
\drawline[AHnb=0](5,14)(5,13)(5.5,12)(6,13)(6,14)
\drawline[AHnb=0](11,14)(11,13)(11.5,12)(12,13)
\drawline[AHnb=0](1.75,14.25)(2,14)(3,14)(3.25,14.25)


\drawline[AHnb=0](1,4.5)(1.5,5.5)(1.5,6.5)(1,7.5)
\drawline[AHnb=0](4,4.5)(3.5,5.5)(3.5,6.5)(4,7.5)
\drawline[AHnb=0](-2,9.5)(-1.5,10.5)(-1.5,11.5)(-2,12.5)
\drawline[AHnb=0](1,9.5)(.5,10.5)(.5,11.5)(1,12.5)

\drawline[AHnb=0](7,4.5)(7.5,5.5)(7.5,6.5)(7,7.5)
\drawline[AHnb=0](10,4.5)(9.5,5.5)(9.5,6.5)(10,7.5)
\drawline[AHnb=0](4,9.5)(4.5,10.5)(4.5,11.5)(4,12.5)
\drawline[AHnb=0](7,9.5)(6.5,10.5)(6.5,11.5)(7,12.5)
\drawline[AHnb=0](10,9.5)(10.5,10.5)(10.5,11.5)(10,12.5)


\drawline[AHnb=0](1.5,5.5)(2.5,5)(3.5,5.5)
\drawline[AHnb=0](7.5,5.5)(8.5,5)(9.5,5.5)
\drawline[AHnb=0](-1.5,10.5)(-.5,10)(.5,10.5)
\drawline[AHnb=0](4.5,10.5)(5.5,10)(6.5,10.5)

\drawline[AHnb=0](1.5,6.5)(2.5,7)(3.5,6.5)
\drawline[AHnb=0](7.5,6.5)(8.5,7)(9.5,6.5)
\drawline[AHnb=0](-1.5,11.5)(-.5,12)(.5,11.5)
\drawline[AHnb=0](4.5,11.5)(5.5,12)(6.5,11.5)


\drawpolygon[fillcolor=green](-3,8)(-3,9)(-2,9.5)(-1.5,10.5)(-1.5,11.5)(-.5,12)(0,13)(0,14)(.25,14.25)(.35,14.25)(.15,14)(.15,13)(-.35,11.92)(-1.35,11.4)(-1.35,10.4)(-1.85,9.45)(-2.85,8.9)(-2.85,8)

\drawpolygon[fillcolor=green](-.5,2)(0,3)(0,4)(1,4.5)(1.5,5.5)(1.5,6.5)(2.5,7)(3,8)(3,9)(4,9.5)(4.5,10.5)(4.5,11.5)(5.5,12)(6,13)(6,14)(6.25,14.25)(6.35,14.25)(6.15,14)(6.15,13)(5.65,11.92)(4.65,11.4)(4.65,10.4)(4.15,9.45)(3.15,8.9)(3.15,8)(2.65,6.9)(1.65,6.45)(1.65,5.45)(1.15,4.45)(0.15,3.95)(0.15,3)(-.4,2)

\drawpolygon[fillcolor=green](5.5,2)(6,3)(6,4)(7,4.5)(7.5,5.5)(7.5,6.5)(8.5,7)(9,8)(9,9)(10,9.5)(10.5,10.5)(10.5,11.5)(11.5,12)(12,13)(12,14)(12.25,14.25)(12.35,14.25)(12.15,14)(12.15,13)(11.65,11.92)(10.65,11.4)(10.65,10.4)(10.15,9.45)(9.15,8.9)(9.15,8)(8.65,6.9)(7.65,6.45)(7.65,5.45)(7.15,4.45)(6.15,3.95)(6.15,3)(5.6,2)

\drawpolygon[fillcolor=green](.75,2)(1,2.5)(2,3)(2,4)(2.5,5)(3.5,5.5)(3.5,6.5)(4,7.5)(5,8)(5,9)(5.5,10)(6.5,10.5)(6.5,11.5)(7,12.5)(8,13)(8,14)(8.25,14.5)(8.35,14.5)(8.1,14)(8.15,12.9)(7.15,12.4)(6.65,11.5)(6.65,10.4)(5.65,9.9)(5.15,9)(5.15,7.9)(4.15,7.4)(3.65,6.5)(3.65,5.4)(2.65,4.9)(2.15,4)(2.15,2.9)(1.1,2.4)(.9,2)

\drawpolygon[fillcolor=green](-2,7.5)(-1,8)(-1,9)(-.5,10)(.5,10.5)(.5,11.5)(1,12.5)(2,13)(2,14)(2.25,14.5)(2.35,14.5)(2.1,14)(2.15,12.9)(1.15,12.4)(.65,11.5)(.65,10.4)(-.35,9.9)(-.85,9)(-.85,7.9)(-1.85,7.4)

\drawpolygon[fillcolor=green](6.75,2)(7,2.5)(8,3)(8,4)(8.5,5)(9.5,5.5)(9.5,6.5)(10,7.5)(11,8)(11,9)(11.5,10)(12.5,10.5)(12.5,11.5)(12.65,11.5)(12.65,10.4)(11.65,9.9)(11.15,9)(11.15,7.9)(10.15,7.4)(9.65,6.5)(9.65,5.4)(8.65,4.9)(8.15,4)(8.15,2.9)(7.1,2.4)(6.9,2)

\drawpolygon[fillcolor=black](-1.5,7.75)(-1,8)(0,8)(1,7.5)(2,8)(3,8)(4,7.5)(5,8)(6,8)(7,7.5)(8,8)(9,8)(10,7.5)(11,8)(12,8)(12,8.1)(11,8.1)(10,7.6)(9,8.1)(8,8.1)(7,7.6)(6,8.1)(5,8.1)(4,7.6)(3,8.1)(2,8.1)(1,7.6)(0,8.1)(-1,8.1)(-1.5,7.85)

\drawpolygon[fillcolor=black](-1.5,12.75)(-1,13)(0,13)(1,12.5)(2,13)(3,13)(4,12.5)(5,13)(6,13)(7,12.5)(8,13)(9,13)(10,12.5)(11,13)(12,13)(12,13.1)(11,13.1)(10,12.6)(9,13.1)(8,13.1)(7,12.6)(6,13.1)(5,13.1)(4,12.6)(3,13.1)(2,13.1)(1,12.6)(0,13.1)(-1,13.1)(-1.5,12.85)

\drawpolygon[fillcolor=black](-1.5,2.75)(-1,3)(0,3)(1,2.5)(2,3)(3,3)(4,2.5)(5,3)(6,3)(7,2.5)(8,3)(9,3)(10,2.5)(11,3)(12,3)(12,3.1)(11,3.1)(10,2.6)(9,3.1)(8,3.1)(7,2.6)(6,3.1)(5,3.1)(4,2.6)(3,3.1)(2,3.1)(1,2.6)(0,3.1)(-1,3.1)(-1.5,2.85)

\drawpolygon[fillcolor=black](-1.5,9.25)(-1,9)(0,9)(1,9.5)(2,9)(3,9)(4,9.5)(5,9)(6,9)(7,9.5)(8,9)(9,9)(10,9.5)(11,9)(12,9)(12,9.1)(11,9.1)(10,9.6)(9,9.1)(8,9.1)(7,9.6)(6,9.1)(5,9.1)(4,9.6)(3,9.1)(2,9.1)(1,9.6)(0,9.1)(-1,9.1)(-1.5,9.35)

\drawpolygon[fillcolor=black](-1.5,4.25)(-1,4)(0,4)(1,4.5)(2,4)(3,4)(4,4.5)(5,4)(6,4)(7,4.5)(8,4)(9,4)(10,4.5)(11,4)(12,4)(12,4.1)(11,4.1)(10,4.6)(9,4.1)(8,4.1)(7,4.6)(6,4.1)(5,4.1)(4,4.6)(3,4.1)(2,4.1)(1,4.6)(0,4.1)(-1,4.1)(-1.5,4.35)

\put(-2,0){\scriptsize {\tiny {\bf Figure 4.15.1:} Paths of type $O_1$ }} 

\end{picture}

\begin{picture}(0,0)(-64,49)
\setlength{\unitlength}{3.5mm}

\drawline[AHnb=0](12,4)(11,4)(10,4.5)(9,4)(8,4)(7,4.5)(6,4)(5,4)(4,4.5)(3,4)(2,4)(1,4.5)

\drawline[AHnb=0](3,8)(4,7.5)(5,8)(6,8)(7,7.5)(8,8)(9,8)(10,7.5)(11,8)(11,9)(10,9.5)(9,9)(8,9)(7,9.5)(6,9)(5,9)(4,9.5)(3,9)

\drawline[AHnb=0](6,13)(7,12.5)(8,13)(9,13)(10,12.5)(11,13)(11,14)

\drawline[AHnb=0](13,4.5)(14,4)(15,4)(16,4.5)(17,4)(18,4)(19,4.5)(19.5,5.5)(19.5,6.5)(20.5,7)(21,8)(21,9)(22,9.5)(22.5,10.5)(22.5,11.5)(23.5,12)(24,13)(24,14)(23,14)(22,14.5)(21,14)(20,14)(19,14.5)(18,14)(18,13)(19,12.5)(20,13)(21,13)(22,12.5)(23,13)(24,13)

\drawline[AHnb=0](16,4.5)(15.5,5.5)(15.5,6.5)(16,7.5)(17,8)(17,9)(17.5,10)(18.5,10.5)(18.5,11.5)(19,12.5)

\drawline[AHnb=0](17.5,10)(16.5,10.5)

\drawline[AHnb=0](19,9.5)(18.5,10.5)
\drawline[AHnb=0](18,8)(18,9)(17.5,10)
\drawline[AHnb=0](21,13)(21,14)
\drawline[AHnb=0](20,13)(20,14)

\drawline[AHnb=0](15,8)(16,7.5)(17,8)(18,8)(19,7.5)(20,8)(21,8)(21,9)(20,9)(19,9.5)(18,9)(17,9)(16,9.5)

\drawline[AHnb=0](15.5,5.5)(14.5,5)(13.5,5.5)

\drawline[AHnb=0](15.5,6.5)(14.5,7)
\drawline[AHnb=0](18.5,11.5)(17.5,12)

\drawline[AHnb=0](2,4)(2.5,5)(3,4)
\drawline[AHnb=0](8,4)(8.5,5)(9,4)
\drawline[AHnb=0](14,4)(14.5,5)(15,4)

\drawline[AHnb=0](5,8)(5,9)(5.5,10)(6,9)(6,8)
\drawline[AHnb=0](11,8)(11,9)(11.5,10)(12,9)(12,8)

\drawline[AHnb=0](9,13)(9,14)

\drawline[AHnb=0](14,9)(14,8)(14.5,7)(13.5,6.5)(13,7.5)
\drawline[AHnb=0](17,14)(17,13)(17.5,12)(16.5,11.5)(16,12.5)

\drawline[AHnb=0](20,9)(20,8)(20.5,7)(19.5,6.5)(19,7.5)
\drawline[AHnb=0](23,14)(23,13)(23.5,12)(22.5,11.5)(22,12.5)

\drawline[AHnb=0](11.5,12)(11,13)
\drawline[AHnb=0](15,14)(15,13)
\drawline[AHnb=0](12,8)(12,9)

\drawline[AHnb=0](11,13)(12,13)
\drawline[AHnb=0](10.5,11.5)(11.5,12)(12,11.75)
\drawline[AHnb=0](10.5,10.5)(11.5,10)(12,10.25)
\drawline[AHnb=0](11,9)(12,9)
\drawline[AHnb=0](11,8)(12,8)

\drawline[AHnb=0](11,4)(12,4)


\drawline[AHnb=0](8,9)(8,8)(8.5,7)(9,8)(9,9)


\drawline[AHnb=0](4,4.5)(3.5,5.5)(3.5,6.5)(4,7.5)

\drawline[AHnb=0](7,4.5)(7.5,5.5)(7.5,6.5)(7,7.5)
\drawline[AHnb=0](10,4.5)(9.5,5.5)(9.5,6.5)(10,7.5)
\drawline[AHnb=0](7,9.5)(6.5,10.5)(6.5,11.5)(7,12.5)
\drawline[AHnb=0](10,9.5)(10.5,10.5)(10.5,11.5)(10,12.5)

\drawline[AHnb=0](1.5,5.5)(2.5,5)(3.5,5.5)
\drawline[AHnb=0](7.5,5.5)(8.5,5)(9.5,5.5)
\drawline[AHnb=0](4.5,10.5)(5.5,10)(6.5,10.5)

\drawline[AHnb=0](1.5,6.5)(2.5,7)(3.5,6.5)
\drawline[AHnb=0](7.5,6.5)(8.5,7)(9.5,6.5)
\drawline[AHnb=0](4.5,11.5)(5.5,12)(6.5,11.5)


\drawline[AHnb=0](1,4.5)(1.5,5.5)(1.5,6.5)(2.5,7)(3,8)(3,9)(4,9.5)(4.5,10.5)(4.5,11.5)(5.5,12)(6,13)(6,14)(7,14.5)(8,14)(9,14)(10,14.5)(11,14)(12,14)(13,14.5)(14,14)(15,14)(16,14.5)(17,14)(18,14)(18,13)(17,13)(16,12.5)(15,13)(14,13)(13,12.5)(12,13)(11.5,12)(12.5,11.5)(12.5,10.5)(11.5,10)(12,9)(13,9.5)(12.5,10.5)

\drawline[AHnb=0](7,4.5)(7.5,5.5)(7.5,6.5)(8.5,7)(9,8)(9,9)(10,9.5)(10.5,10.5)(10.5,11.5)(11.5,12)(12,13)(12,14)

\drawline[AHnb=0](12,4)(13,4.5)(13.5,5.5)(13.5,6.5)(14.5,7)(15,8)(15,9)(16,9.5)(16.5,10.5)(16.5,11.5)(17.5,12)(18,13)(18,14)

\drawline[AHnb=0](2,4)(2.5,5)(3.5,5.5)(3.5,6.5)(4,7.5)(5,8)(5,9)(5.5,10)(6.5,10.5)(6.5,11.5)(7,12.5)(8,13)(8,14)

\drawline[AHnb=0](8,4)(8.5,5)(9.5,5.5)(9.5,6.5)(10,7.5)(11,8)(11,9)(11.5,10)(12.5,10.5)(12.5,11.5)(13,12.5)(14,13)(14,14)

\drawline[AHnb=0](3,9)(4,9.5)(5,9)(6,9)(7,9.5)(8,9)(9,9)(10,9.5)(11,9)(12,9)(13,9.5)(14,9)(15,9)

\drawline[AHnb=0](3,8)(4,7.5)(5,8)(6,8)(7,7.5)(8,8)(9,8)(10,7.5)(11,8)(12,8)(13,7.5)(14,8)(15,8)

\put(.25,3.75){\scriptsize {\tiny $a_{11}$}}
\put(1.35,3.35){\scriptsize {\tiny $a_{12}$}}
\put(2.65,3.35){\scriptsize {\tiny $a_{13}$}}
\put(3.75,3.65){\scriptsize {\tiny $a_{14}$}}
\put(18,3.35){\scriptsize {\tiny $a_{1i}$}}
\put(19.25,4){\scriptsize {\tiny $a_{11}$}}

\put(.25,5.5){\scriptsize {\tiny $x_{11}$}}
\put(.25,6.35){\scriptsize {\tiny $x_{12}$}}
\put(1.15,7.1){\scriptsize {\tiny $x_{13}$}}
\put(19.85,5){\scriptsize {\tiny $x_{11}$}}
\put(19.85,6){\scriptsize {\tiny $x_{12}$}}
\put(20.85,6.75){\scriptsize {\tiny $x_{13}$}}

\put(1.8,7.8){\scriptsize {\tiny $a_{21}$}}
\put(21.2,7.5){\scriptsize {\tiny $a_{21}$}}

\put(1.85,8.8){\scriptsize {\tiny $a_{31}$}}
\put(2.85,9.7){\scriptsize {\tiny $a_{32}$}}
\put(21.3,8.5){\scriptsize {\tiny $a_{31}$}}
\put(22.15,9){\scriptsize {\tiny $a_{32}$}}

\put(4.85,13.15){\scriptsize {\tiny $a_{j1}$}}
\put(24.2,13){\scriptsize {\tiny $a_{j1}$}}

\put(3.65,14.5){\scriptsize {\tiny $a_{1(k+1)}$}}
\put(6,15.15){\scriptsize {\tiny $a_{1(k+2)}$}}
\put(22.35,14.5){\scriptsize {\tiny $a_{1k}$}}
\put(23.65,14.35){\scriptsize {\tiny $a_{1(k+1)}$}}

\put(-1,4){\scriptsize {\tiny $Q_1$}}
\put(0,7.8){\scriptsize {\tiny $Q_{2}$}}
\put(0.25,8.8){\scriptsize {\tiny $Q_{3}$}}
\put(2,13){\scriptsize {\tiny $Q_{j}$}}
\put(2.25,14){\scriptsize {\tiny $Q_{1}$}}

\put(10.8,15.1){\scriptsize $\ldots$}
\put(14.8,15.1){\scriptsize $\ldots$}
\put(18.8,15.1){\scriptsize $\ldots$}

\put(6.8,3){\scriptsize $\ldots$}
\put(10.8,3){\scriptsize $\ldots$}
\put(14.8,3){\scriptsize $\ldots$}

\put(1.7,11.2){\scriptsize $.$}
\put(1.6,11){\scriptsize $.$}
\put(1.5,10.8){\scriptsize $.$}

\drawpolygon[fillcolor=black](24,13)(23.5,12)(22.5,11.5)(22.5,10.5)(22,9.5)(21,9)(21,8)(20.5,7)(19.5,6.5)(19.5,5.5)(19,4.5)(19.05,4.4)(19.6,5.5)(19.6,6.45)(20.6,6.95)(21.1,8)(21.1,8.95)(22.1,9.45)(22.6,10.5)(22.6,11.45)(23.6,11.95)(24.1,13)

\drawpolygon[fillcolor=black](19,4.5)(18,4)(17,4)(16,4.5)(15,4)(14,4)(13,4.5)(12,4)(11,4)(10,4.5)(9,4)(8,4)(7,4.5)(6,4)(5,4)(4,4.5)(3,4)(2,4)(2.5,5)(3.5,5.5)(3.5,6.5)(4,7.5)(5,8)(6,8)(7,7.5)(8,8)(9,8)(8.5,7)(7.5,6.5)(7.5,5.5)(8.5,5)(9.5,5.5)(9.5,6.5)(10,7.5)(11,8)(12,8)(13,7.5)(14,8)(15,8)(14.5,7)(13.5,6.5)(13.5,5.5)(14.5,5)(15.5,5.5)(15.5,6.5)(16,7.5)(17,8)(18,8)(19,7.5)(20,8)(20,9)(19,9.5)(18,9)(17,9)(16,9.5)(15,9)(14,9)(13,9.5)(12,9)(11,9)(10,9.5)(9,9)(8,9)(7,9.5)(6,9)(5,9)(5.5,10)(6.5,10.5)(6.5,11.5)(7,12.5)(8,13)(9,13)(10,12.5)(11,13)(12,13)(11.5,12)(10.5,11.5)(10.5,10.5)(11.5,10)(12.5,10.5)(12.5,11.5)(13,12.5)(14,13)(15,13)(16,12.5)(17,13)(18,13)(17.5,12)(16.5,11.5)(16.5,10.5)(17.5,10)(18.5,10.5)(18.5,11.5)(19,12.5)(20,13)(21,13)(22,12.5)(23,13)(24.1,13)(24.1,13.1)(22.9,13.1)(22,12.6)(21,13.1)(20,13.1)(18.9,12.55)(18.4,11.5)(18.4,10.6)(17.5,10.1)(16.6,10.6)(16.6,11.4)(17.6,11.9)(18.1,13.1)(17,13.1)(16,12.6)(15,13.1)(14,13.1)(12.9,12.55)(12.4,11.5)(12.4,10.6)(11.5,10.1)(10.6,10.6)(10.6,11.4)(11.6,11.9)(12.1,13.1)(11,13.1)(10,12.6)(9,13.1)(8,13.1)(6.9,12.55)(6.4,11.6)(6.4,10.6)(5.4,10.1)(4.9,8.9)(6.1,8.9)(7,9.4)(7.9,8.9)(9.1,8.9)(10,9.4)(10.9,8.9)(12.1,8.9)(13,9.4)(13.9,8.9)(15,8.9)(16,9.4)(17,8.9)(18,8.9)(19,9.4)(19.9,8.9)(19.9,8.1)(19,7.6)(18,8.1)(17,8.1)(15.9,7.6)(15.4,6.6)(15.4,5.6)(14.6,5.1)(13.6,5.6)(13.6,6.4)(14.6,6.9)(15.1,8.1)(13.9,8.1)(13,7.6)(12,8.1)(11,8.1)(9.9,7.55)(9.4,6.5)(9.4,5.6)(8.5,5.1)(7.6,5.6)(7.6,6.4)(8.6,6.9)(9.15,8.1)(7.9,8.1)(7,7.6)(6,8.1)(5,8.1)(3.9,7.55)(3.4,6.5)(3.4,5.6)(2.4,5.1)(1.9,3.9)(3.1,3.9)(4,4.4)(5,3.9)(6,3.9)(7,4.4)(8,3.9)(9,3.9)(10,4.4)(11,3.9)(12,3.9)(13.05,4.4)(14,3.9)(15,3.9)(16,4.4)(17,3.9)(18,3.9)(19,4.4)

\put(6,1.6){\scriptsize {\tiny {\bf Figure 4.15.2:} $M(i,j,k)$}} 

\end{picture}

\vspace{4.65cm}

We construct $M(i,j,k)$ representation of $M$ by first cutting $M$ along a black coloured cycle of type $O_1$ through a vertex $v$ having face-sequence $(4,6, 12)$ and then cutting it along a green coloured cycle of type $O_1$ where the beginning adjacent face to black coloured cycle is a 4-gon. The admissible relations among $i, j, k$ of $M(i, j, k)$ are given in the next lemma, and its proof follows from the similar arguments as in Lemma \ref{l3.13.1}.

\begin{lem}\label{l3.15.1} 
	A DSEM $M$ of type $[3.4.6.4:4.6.12]$ admits an $M(i,j,k)$-representation iff the following holds: $(i)$ $j$ even,  $(ii)$ $i \geq 12$ and $i=6m$, $m \in \mathbb{N} $, $(iii)$ number of vertices of $M(i,j,k) = 3ij/2\geq 36 $, $(iv)$  if $j=2$ then $ k \in \{3r+2: 1 < r < (i-3)/3\}$, and if $j \geq 4$ then $ k \in \{3r+2: 0 \leq r < i/3\}$.
\end{lem}

%

\subsection{DSEMs of type $[3.4^2.6:3.4.6.4]$} \label{s3.16}

Let $M$ be a DSEM  of type $[3.4^2.6:3.4.6.4]$. In $M$, consider a fixed type path $P_{1} = P( \ldots, y_{i},y_{i+1},y_{i+2}$, $z_{i},z_{i+1},y_{i+3}, \ldots)$, say of type $R_{1}$, indicated by thick black or green paths, shown in Figure 4.16.1. The vertices $y_i$'s and $z_i$'s have the face-sequences $(3,4,6,4)$ and $(3,4^2,6)$ respectively.

\begin{picture}(0,0)(-28.5,37)
\setlength{\unitlength}{3.5mm}

\drawline[AHnb=0](-5,5)(-4,5)(-3,4.5)(-2,5)(-1,5)(-1,6)(-2,6)(-3,6.5)(-4,6)(-5,6)

\drawline[AHnb=0](0,5)(1,5)(2,4.5)(3,5)(4,5)(4,6)(3,6)(2,6.5)(1,6)(0,6)(0,5)

\drawline[AHnb=0](-4.5,1.5)(-4.5,.5)(-2.5,.5)(-1.5,.5)(-.5,0)(.5,.5)(1.5,.5)(3.5,.5)(3.5,1.5)(1.5,1.5)(.5,1.5)(-.5,2)(-1.5,1.5)(-2.5,1.5)(-4.5,1.5)

\drawline[AHnb=0](3.5,.5)(4.5,.25)

\drawline[AHnb=0](4,5)(4.25,5)
\drawline[AHnb=0](4,6)(4.25,6)


\drawline[AHnb=0](-1.5,1.5)(-2.5,3.5)
\drawline[AHnb=0](3.5,1.5)(2.5,3.5)

\drawline[AHnb=0](-4.5,1.5)(-3.5,3.5)
\drawline[AHnb=0](0.5,1.5)(1.5,3.5)

\drawline[AHnb=0](-2,5)(-.5,2)(1,5)(-2,5)

\drawline[AHnb=0](-2,6)(-.5,9)(1,6)(-2,6)

\drawline[AHnb=0](-4.5,7)(-5.5,9)
\drawline[AHnb=0](-3.5,7.5)(-4.5,9.5)
\drawline[AHnb=0](-2.5,7.5)(-1.5,9.5)
\drawline[AHnb=0](1.5,7.5)(.5,9.5)
\drawline[AHnb=0](2.5,7.5)(3.5,9.5)
\drawline[AHnb=0](3.5,7)(4.5,9)

\drawline[AHnb=0](3,5)(4.35,2.25)

\drawline[AHnb=0](-2.25,8.75)(-2,8.5)(-1,8)(0,8)(1,8.5)(1.25,8.75)

\drawline[AHnb=0](2.75,8.75)(3,8.5)(4,8)(4.25,8)

\drawline[AHnb=0](-5,8)(-4,8.5)(-3.75,8.75)


\drawline[AHnb=0](-4.5,.5)(-4,-.5)(-3.5,.5)(-3.5,1.5)(-4,2.5)(-5,3)

\drawline[AHnb=0](.5,1.5)(.5,.5)(1,-.5)(1.5,.5)(1.5,1.5)(1,2.5)(0,3)(-1,3)(-2,2.5)(-2.5,1.5)(-2.5,.5)(-2,-.5)(-1.5,.5)(-1.5,1.5)

\drawline[AHnb=0](4.25,3)(4,3)(3,2.5)(2.5,1.5)(2.5,.5)(3,-.5)(3.5,.5)(3.5,1.5)(4.5,2)

\drawline[AHnb=0](0,6)(.5,7)(1,6)
\drawline[AHnb=0](3,6)(3.5,7)(4,6)

\drawline[AHnb=0](-5,6.25)(-4.5,7)(-4,6)
\drawline[AHnb=0](-2,6)(-1.5,7)(-1,6)


\drawpolygon(1.5,7.5)(2.5,7.5)(2,6.5)

\drawpolygon(-3.5,7.5)(-2.5,7.5)(-3,6.5)

\drawpolygon(1.5,3.5)(2.5,3.5)(2,4.5)
\drawpolygon(-3.5,3.5)(-2.5,3.5)(-3,4.5)

\drawline[AHnb=0](0,5)(.5,4)(1,5)(1,6)
\drawline[AHnb=0](3,6)(3,5)(3.5,4)(4,5)

\drawline[AHnb=0](-5,4.75)(-4.5,4)(-4,5)(-4,6)
\drawline[AHnb=0](-2,6)(-2,5)(-1.5,4)(-1,5)


\drawline[AHnb=0](.5,4)(1.5,3.5)
\drawline[AHnb=0](2.5,3.5)(3.5,4)

\drawline[AHnb=0](-4.5,4)(-3.5,3.5)
\drawline[AHnb=0](-2.5,3.5)(-1.5,4)

\drawline[AHnb=0](.5,7)(1.5,7.5)
\drawline[AHnb=0](2.5,7.5)(3.5,7)

\drawline[AHnb=0](-4.5,7)(-3.5,7.5)
\drawline[AHnb=0](-2.5,7.5)(-1.5,7)


\drawline[AHnb=0](-.75,-.5)(-.5,0)(-.25,-.5)

\drawpolygon[fillcolor=black](-5,5)(-4,5)(-3,4.5)(-2,5)(1,5)(2,4.5)(3,5)(4.5,5)(4.5,5.1)(3,5.1)(2,4.6)(1,5.1)(-2,5.1)(-3,4.6)(-4,5.1)(-5,5.1)

\drawpolygon[fillcolor=black](-5,6)(-4,6)(-3,6.5)(-2,6)(1,6)(2,6.5)(3,6)(4.5,6)(4.5,6.1)(3,6.1)(2,6.6)(1,6.1)(-2,6.1)(-3,6.6)(-4,6.1)(-5,6.1)

\drawpolygon[fillcolor=black](-5.5,2)(-4.5,1.5)(-1.5,1.5)(-.5,2)(0.5,1.5)(3.5,1.5)(4.5,2)(4.5,2.1)(3.5,1.6)(.5,1.6)(-.5,2.1)(-1.5,1.6)(-4.5,1.6)(-5.5,2.1)

\drawpolygon[fillcolor=black](-5.5,0)(-4.5,.5)(-1.5,.5)(-.5,0)(0.5,.5)(3.5,.5)(4.5,0.25)(4.5,.35)(3.5,.6)(.5,.6)(-.5,.1)(-1.5,.6)(-4.5,.6)(-5.5,.1)

\drawpolygon[fillcolor=green](-2,-.5)(-1.5,.5)(-1.5,1.5)(-.5,2)(0,3)(1,5)(1,6)(2,6.5)(3.5,9.5)(3.65,9.5)(2.15,6.4)(1.15,5.9)(1.15,5)(.15,3)(-.4,1.9)(-1.35,1.4)(-1.35,.5)(-1.85,-.5)

\drawpolygon[fillcolor=green](-6.5,1.5)(-5.5,2)(-5,3)(-4,5)(-4,6)(-3,6.5)(-1.5,9.5)(-1.35,9.5)(-2.85,6.4)(-3.85,5.9)(-3.85,5)(-4.85,3)(-5.4,1.9)(-6.35,1.4)

\drawpolygon[fillcolor=green](-.75,-.5)(-.5,0)(.5,.5)(.5,1.5)(1,2.5)(2,4.5)(3,5)(3,6)(4.5,9)(4.65,9)(3.15,6)(3.15,4.95)(2.15,4.45)(1.15,2.5)(.65,1.5)(.65,.4)(-.4,-.1)(-.65,-.5)

\drawpolygon[fillcolor=green](-5.5,0)(-4.5,.5)(-4.5,1.5)(-4,2.5)(-3,4.5)(-2,5)(-2,6)(-.5,9)(-.35,9)(-1.85,6)(-1.85,4.95)(-2.85,4.45)(-3.85,2.5)(-4.35,1.5)(-4.35,.4)(-5.4,-.1)

\put(-5.5,-1.65){\scriptsize {\tiny {\bf Figure 4.16.1:} Paths of type $R_1$ }} 

\end{picture}

\begin{picture}(0,0)(-70,36)
\setlength{\unitlength}{3.5mm}

\drawpolygon(0,0)(3,0)(4,1)(5,0)(8,0)(9,1)(10,0)(13,0)(14,1)(15,0)(17,3)(16,4)(13,4)(12,3)(11,4)(8,4)(7,3)(6,4)(3,4)

\drawpolygon(3,5)(6,5)(7,6)(8,5)(11,5)(12,6)(13,5)(15,8)(14,9)(11,9)(10,8)(9,9)(6,9)(5,8)

\drawpolygon(3,4)(6,4)(7,3)(8,4)(11,4)(12,3)(13,4)(13,5)(12,6)(11,5)(8,5)(7,6)(6,5)(3,5)

\drawpolygon(6,9)(9,9)(10,8)(11,9)(14,9)(15,8)(16,9)(16,10)(15,11)(14,10)(11,10)(10,11)(9,10)(6,10)

\drawline[AHnb=0](13,5)(16,5)(17,6)(18,5)(18,4)(17,3)

\drawline[AHnb=0](16,10)(19,10)(20,11)(21,10)(21,9)(20,8)(19,9)(16,9)

\drawline[AHnb=0](11,0)(10.75,1.25)

\drawline[AHnb=0](9,1)(11,4)
\drawline[AHnb=0](9,1)(8,4)
\drawline[AHnb=0](14,1)(16,4)
\drawline[AHnb=0](14,1)(13,4)

\drawline[AHnb=0](8,0)(7,3)
\drawline[AHnb=0](5,0)(7,3)
\drawline[AHnb=0](13,0)(12,3)
\drawline[AHnb=0](10,0)(12,3)

\drawline[AHnb=0](4,1)(6,4)
\drawline[AHnb=0](4,1)(3,4)
\drawline[AHnb=0](14,1)(16,4)
\drawline[AHnb=0](14,1)(13,4)
\drawline[AHnb=0](9,1)(11,4)
\drawline[AHnb=0](9,1)(8,4)

\drawline[AHnb=0](3,0)(2,3)

\drawline[AHnb=0](12,6)(14,9)
\drawline[AHnb=0](12,6)(11,9)
\drawline[AHnb=0](17,6)(19,9)
\drawline[AHnb=0](17,6)(16,9)

\drawline[AHnb=0](11,5)(10,8)
\drawline[AHnb=0](8,5)(10,8)
\drawline[AHnb=0](16,5)(15,8)
\drawline[AHnb=0](13,5)(15,8)

\drawline[AHnb=0](16,5)(15,8)
\drawline[AHnb=0](13,5)(15,8)
\drawline[AHnb=0](18,5)(20,8)

\drawline[AHnb=0](7,6)(9,9)
\drawline[AHnb=0](7,6)(6,9)
\drawline[AHnb=0](17,6)(19,9)
\drawline[AHnb=0](17,6)(16,9)

\drawline[AHnb=0](12,6)(14,9)
\drawline[AHnb=0](12,6)(11,9)

\drawline[AHnb=0](6,5)(5,8)
\drawline[AHnb=0](11,5)(10,8)
\drawline[AHnb=0](4,4)(4,5)
\drawline[AHnb=0](5,4)(5,5)

\drawline[AHnb=0](6,4)(6,5)
\drawline[AHnb=0](8,4)(8,5)

\drawline[AHnb=0](9,4)(9,5)
\drawline[AHnb=0](10,4)(10,5)
\drawline[AHnb=0](11,4)(11,5)

\drawline[AHnb=0](7,9)(7,10)
\drawline[AHnb=0](8,9)(8,10)

\drawline[AHnb=0](9,9)(9,10)
\drawline[AHnb=0](11,9)(11,10)

\drawline[AHnb=0](12,9)(12,10)
\drawline[AHnb=0](13,9)(13,10)
\drawline[AHnb=0](14,9)(14,10)

\drawline[AHnb=0](2.66,1.15)(3.66,2.15)
\drawline[AHnb=0](2.33,2.15)(3.33,3.15)

\drawline[AHnb=0](6.4,2.15)(5.4,3.15)
\drawline[AHnb=0](5.8,1.15)(4.8,2.15)

\drawline[AHnb=0](7.66,1.15)(8.66,2.15)
\drawline[AHnb=0](7.33,2.15)(8.33,3.15)

\drawline[AHnb=0](11.4,2.15)(10.4,3.15)
\drawline[AHnb=0](10.8,1.15)(9.8,2.15)
\drawline[AHnb=0](5.66,6.15)(6.66,7.15)
\drawline[AHnb=0](5.33,7.15)(6.33,8.15)

\drawline[AHnb=0](9.4,7.15)(8.4,8.15)
\drawline[AHnb=0](8.8,6.15)(7.8,7.15)

\drawline[AHnb=0](10.66,6.15)(11.66,7.15)
\drawline[AHnb=0](10.33,7.15)(11.33,8.15)

\drawline[AHnb=0](14.4,7.15)(13.4,8.15)
\drawline[AHnb=0](13.8,6.15)(12.8,7.15)

\drawline[AHnb=0](.75,1.15)(1,0)
\drawline[AHnb=0](2.66,1.15)(2,0)
\drawline[AHnb=0](7.66,1.15)(7,0)
\drawline[AHnb=0](5.8,1.15)(6,0)

\drawline[AHnb=0](3.75,6.15)(4,5)
\drawline[AHnb=0](5.66,6.15)(5,5)
\drawline[AHnb=0](10.66,6.15)(10,5)
\drawline[AHnb=0](8.8,6.15)(9,5)

\drawline[AHnb=0](1.4,2.15)(2.3,2.15)
\drawline[AHnb=0](3.66,2.15)(4.8,2.15)
\drawline[AHnb=0](6.4,2.15)(7.3,2.15)
\drawline[AHnb=0](8.66,2.15)(9.8,2.15)

\drawline[AHnb=0](9.4,7.15)(10.3,7.15)
\drawline[AHnb=0](6.66,7.15)(7.8,7.15)
\drawline[AHnb=0](4.4,7.15)(5.3,7.15)
\drawline[AHnb=0](11.66,7.15)(12.8,7.15)

\drawline[AHnb=0](12,0)(12.6,1.25)(13.66,2.15)(14.8,2.15)(15.8,1.25)
\drawline[AHnb=0](15,5)(15.6,6.25)(16.66,7.15)(17.8,7.15)(18.8,6.25)

\drawline[AHnb=0](11.45,2.15)(12.25,2.15)(13.35,3)(14.2,4)(14.2,5)(13.8,6.2)
\drawline[AHnb=0](14.45,7.15)(15.25,7.15)(16.35,8)(17.2,9)(17.2,10)

\drawline[AHnb=0](19.5,7.25)(18.5,8.2)(18.2,9)(18.2,10)

\drawline[AHnb=0](16.5,2.25)(15.5,3.2)(15,4)(15,5)

\drawline[AHnb=0](16,4)(16,5)

\drawline[AHnb=0](19,9)(19,10)

\drawline[AHnb=0](3.33,3.15)(4,4)
\drawline[AHnb=0](5.4,3.15)(5,4)
\drawline[AHnb=0](8.33,3.15)(9,4)
\drawline[AHnb=0](10.4,3.15)(10,4)

\drawline[AHnb=0](6.33,8.15)(7,9)
\drawline[AHnb=0](8.4,8.15)(8,9)
\drawline[AHnb=0](11.33,8.15)(12,9)
\drawline[AHnb=0](13.4,8.15)(13,9)

\put(-.7,-.6){\scriptsize {\tiny $a_{11}$}}
\put(.5,-.6){\scriptsize {\tiny $a_{12}$}}
\put(1.6,-.6){\scriptsize {\tiny $a_{13}$}}
\put(2.7,-.6){\scriptsize {\tiny $a_{14}$}}
\put(3.55,.2){\scriptsize {\tiny $a_{15}$}}
\put(13.65,.2){\scriptsize {\tiny $a_{1i}$}}
\put(14.8,-.5){\scriptsize {\tiny $a_{11}$}}

\put(-.5,1.1){\scriptsize {\tiny $x_{11}$}}
\put(15.95,.9){\scriptsize {\tiny $x_{11}$}}

\put(.2,2.2){\scriptsize {\tiny $x_{12}$}}
\put(16.75,2){\scriptsize {\tiny $x_{12}$}}

\put(.85,3){\scriptsize {\tiny $a_{21}$}}
\put(1.85,4){\scriptsize {\tiny $a_{22}$}}

\put(17.45,2.85){\scriptsize {\tiny $a_{21}$}}
\put(18.2,4){\scriptsize {\tiny $a_{22}$}}

\put(1.85,4.75){\scriptsize {\tiny $a_{31}$}}
\put(18.2,4.8){\scriptsize {\tiny $a_{31}$}}

\put(4,8.2){\scriptsize {\tiny $a_{j1}$}}
\put(4.8,9.2){\scriptsize {\tiny$a_{j2}$}}
\put(20.15,7.75){\scriptsize {\tiny $a_{j1}$}}
\put(21.2,9){\scriptsize {\tiny $a_{j2}$}}

\put(3.85,10.5){\scriptsize {\tiny $a_{1(k+1)}$}}
\put(6.45,10.5){\scriptsize {\tiny $a_{1(k+2)}$}}
\put(19.6,11.45){\scriptsize {\tiny $a_{1k}$}}
\put(20.9,10.5){\scriptsize {\tiny $a_{1(k+1)}$}}

\put(-1.8,0){\scriptsize {\tiny $Q_1$}}
\put(-.5,3){\scriptsize {\tiny $Q_2$}}
\put(2.2,8.2){\scriptsize {\tiny $Q_{j}$}}
\put(3.2,10){\scriptsize {\tiny $Q_{1}$}}

\put(5.8,-.65){\scriptsize $\ldots$}
\put(9.2,-.65){\scriptsize $\ldots$}

\put(12,11){\scriptsize $\ldots$}
\put(17.2,11){\scriptsize $\ldots$}

\put(1.6,6.2){\scriptsize $.$}
\put(1.5,6){\scriptsize $.$}
\put(1.4,5.8){\scriptsize $.$}

\drawpolygon[fillcolor=black](0,0)(3,0)(2.66,1.15)(3.66,2.15)(4.8,2.15)(4,1)(5,0)(5.8,1.15)(6,0)(8,0)(7.66,1.15)(8.66,2.15)(9.8,2.15)(9,1)(10,0)(10.8,1.15)(11,0)(13,0)(12.66,1.15)(13.66,2)(14,1)(14.8,2.15)(15.4,3.15)(16,4)(14.2,4)(13.3,3.15)(13,4)(12,3)(12.25,2.15)(11.4,2.15)(10.4,3.15)(11,4)(9,4)(8.3,3.15)(8,4)(7,3)(7.25,2.15)(6.4,2.15)(5.4,3.15)(6,4)(4,4)(4,5)(6,5)(5.66,6.15)(6.66,7.15)(7.8,7.15)(7,6)(8,5)(8.8,6.15)(9,5)(11,5)(10.66,6.15)(11.66,7.15)(12.8,7.15)(12,6)(13,5)(13.8,6.15)(14,5)(16,5)(15.66,6.15)(16.66,7.15)(17,6)(17.8,7.15)(18.4,8.15)(19,9)(17,9)(16.3,8.15)(16,9)(15,8)(15.3,7.15)(14.4,7.15)(13.45,8.15)(14,9)(12,9)(11.3,8.15)(11,9)(10,8)(10.3,7.15)(9.4,7.15)(8.45,8.15)(9,9)(7,9)(6.9,8.85)(8.7,8.85)(8.2,8.1)(9.35,6.95)(10.5,7)(10.2,8)(10.95,8.75)(11.25,7.85)(12.1,8.8)(13.7,8.85)(13.2,8.1)(14.35,6.95)(15.5,7)(15.2,8)(15.95,8.75)(16.25,7.85)(17.1,8.8)(18.7,8.85)(17.1,6.35)(16.7,7.3)(15.45,6.15)(15.8,5.15)(14.1,5.15)(13.85,6.4)(13,5.25)(12.25,6)(13.1,7.3)(11.6,7.3)(10.5,6.25)(10.8,5.15)(9.15,5.15)(8.85,6.4)(8,5.25)(7.25,6)(8.1,7.3)(6.6,7.3)(5.5,6.25)(5.8,5.15)(3.85,5.15)(3.85,3.85)(5.7,3.85)(5.15,3.15)(6.3,2)(7.5,2)(7.2,3)(7.92,3.7)(8.25,2.9)(9.1,3.85)(10.7,3.85)(10.15,3.15)(11.3,2)(12.5,2)(12.2,3)(12.92,3.7)(13.25,2.9)(14.2,3.85)(15.7,3.85)(14.05,1.25)(13.75,2.2)(12.5,1.2)(12.8,.15)(11.1,.15)(10.9,1.45)(10,0.25)(9.2,1)(10.05,2.3)(8.55,2.3)(7.5,1.2)(7.8,.15)(6.15,.15)(5.8,1.3)(4.95,.25)(4.2,1)(5.05,2.3)(3.55,2.3)(2.5,1.2)(2.8,.15)(0.1,.15)

\drawpolygon[fillcolor=black](7,9)(6,9)(6.33,8.15)(5.33,7.15)(5,8)(3,5)(3,4)(3.33,3.15)(2.33,2.15)(2,3)(0,0)(0.1,0)(2,2.65)(2.3,1.9)(3.5,3.1)(3.15,4.1)(3.15,5)(4.95,7.65)(5.3,6.9)(6.5,8.15)(6.15,8.9)(7,8.9)

\put(4,-1.75){\scriptsize {\tiny {\bf Figure 4.16.2:} $M(i,j,k)$ }} 

\end{picture}

\vspace{4.75cm}

We construct $M(i,j,k)$ representation of $M$ by first cutting $M$ along a black colored cycle of type $R_1$ through a vertex $v$ with the face-sequence $(3,4, 6,4)$ and then cutting it along a green colored cycle of type $R_1$, assume that the beginning adjacent face to base horizontal cycle is a 3-gon. The admissible relations among $i, j, k$ are given by the next lemma, and its proof follows from the similar arguments as in Lemma \ref{l3.10.1}.

\begin{lem}\label{l3.16.1} 
	A DSEM $M$ of type $[3.4^2.6:3.4.6.4]$ admits an $M(i,j,k)$-representation iff the following holds: $(i)$ $j$ even,  $(ii)$ $i \geq 5$ and $i=5m$, $m \in \mathbb{N} $, $(iii)$ number of vertices of $M(i,j,k) = 9ij/5\geq 18 $, $(iv)$  $ k \in \{5r: 0 \leq r < i/5\}$.
\end{lem}

%

\subsection{DSEMs of type $[3^3.4^2:3.4.6.4]$} \label{s3.17}

Let $M$ be a DSEM of type $[3^3.4^2:3.4.6.4]$. In $M$, consider a path $P_{1} = P( \ldots, y_{i},z_{i},z_{i+1},z_{i+2},y_{i+1}$, $\ldots)$, say of type $S_{1}$, indicated by thick black or green paths shown in Figure 4.17.1. The vertices $y_i$'s and $z_i$'s have the face-sequences $(3^3,4^2)$ and $(3,4,6,4)$ respectively or $y_i$'s and $z_i$'s have the face-sequence $(3^3,4^2)$.

\begin{picture}(0,0)(0,25)
\setlength{\unitlength}{2.5mm}

\drawline[AHnb=0](0,0)(0,1)(-.5,2)(0,3)(0,5)(-.5,6)(0,7)(0,9)

\drawline[AHnb=0](7,0)(7,1)(6.5,2)(7,3)(7,5)(6.5,6)(7,7)(7,9)

\drawline[AHnb=0](1,9)(1,7)(1.5,6)(1,5)(1,3)(1.5,2)(1,1)(1,0)

\drawline[AHnb=0](7,1)(8,1)(8,0)

\drawline[AHnb=0](7,9)(8,9)(8,7)(7,7)

\drawline[AHnb=0](7,5)(8,5)(8,3)(7,3)

\drawline[AHnb=0](2,.5)(2.5,1.5)

\drawline[AHnb=0](6,.5)(5.5,1.5)

\drawline[AHnb=0](6,.5)(7,1)

\drawline[AHnb=0](3.25,.5)(3.5,1)(3.5,3)(3,4)(3.5,5)(3.5,7)(3,8)(3.5,9)(4.5,9)(5,8)(4.5,7)(4.5,5)(5,4)(4.5,3)(4.5,1)(3.5,1)(3.75,.5)

\drawline[AHnb=0](4.75,.5)(4.5,1)(4.25,.5)

\drawline[AHnb=0](2.5,5.5)(2,4.5)

\drawline[AHnb=0](5.5,2.5)(6,3.5)
\drawline[AHnb=0](5.5,6.5)(6,7.5)

\drawline[AHnb=0](1.5,2)(3.5,3)(3.5,1)(1.5,2)
\drawline[AHnb=0](2.5,2.5)(2.5,1.5)(3.5,2)(2.5,2.5)(2,3.5)

\drawline[AHnb=0](1.5,6)(3.5,7)(3.5,5)(1.5,6)
\drawline[AHnb=0](2.5,6.5)(2.5,5.5)(3.5,6)(2.5,6.5)(2,7.5)

\drawline[AHnb=0](5,4)(7,5)(7,3)(5,4)
\drawline[AHnb=0](6,4.5)(6,3.5)(7,4)(6,4.5)(5.5,5.5)

\drawline[AHnb=0](5,8)(7,9)(7,7)(5,8)
\drawline[AHnb=0](6,8.5)(6,7.5)(7,8)(6,8.5)(5.5,9.5)

\drawline[AHnb=0](0,3)(1,3)(3,4)(1,5)(0,5)
\drawline[AHnb=0](1,4)(2,4.5)(2,3.5)(1,4)(0,4)

\drawline[AHnb=0](0,7)(1,7)(3,8)(1,9)(0,9)
\drawline[AHnb=0](1,8)(2,8.5)(2,7.5)(1,8)(0,8)

\drawline[AHnb=0](3.5,1)(4.5,1)(6.5,2)(4.5,3)(3.5,3)
\drawline[AHnb=0](4.5,2)(5.5,2.5)(5.5,1.5)(4.5,2)(3.5,2)

\drawline[AHnb=0](3.5,5)(4.5,5)(6.5,6)(4.5,7)(3.5,7)
\drawline[AHnb=0](4.5,6)(5.5,6.5)(5.5,5.5)(4.5,6)(3.5,6)

\drawline[AHnb=0](7,8)(8,8)
\drawline[AHnb=0](7,4)(8,4)
\drawline[AHnb=0](-1,.5)(0,1)(1,1)(2,0.5)
\drawline[AHnb=0](-1,1.75)(-.5,2)(-1,2.25)
\drawline[AHnb=0](-1,5.75)(-.5,6)(-1,6.25)
\drawline[AHnb=0](0,3)(-1,3.5)(0,4)(-1,4.5)(0,5)
\drawline[AHnb=0](0,7)(-1,7.5)(0,8)(-1,8.5)(0,9)

\drawline[AHnb=0](-1,3.5)(-1,4.5)
\drawline[AHnb=0](-1,7.5)(-1,8.5)

\drawpolygon[fillcolor=black](-1,1.75)(-.5,2)(0,3)(1,3)(3,4)(3.5,5)(4.5,5)(6.5,6)(7,7)(8,7)(8,6.85)(7.1,6.85)(6.6,5.85)(4.5,4.85)(3.6,4.85)(3.1,3.9)(1,2.85)(0.1,2.85)(-.4,1.85)(-1,1.65)

\drawpolygon[fillcolor=black](-1,5.75)(-.5,6)(0,7)(1,7)(3,8)(3.5,9)(4.5,9)(4.5,8.85)(3.6,8.85)(3.1,7.9)(1,6.85)(0.1,6.85)(-.4,5.85)(-1,5.65)

\drawpolygon[fillcolor=black](-.5,.75)(0,1)(1,1)(1.5,2)(3.5,3)(4.5,3)(5,4)(7,5)(8,5)(8,5.15)(7,5.15)(7.1,4.85)(5,3.85)(4.5,2.85)(3.6,2.9)(1.5,1.85)(1,.9)(.1,.85)(-.5,.65)

\drawpolygon[fillcolor=black](-1,4.55)(0,5)(1,5)(1.5,6)(3.5,7)(4.5,7)(5,8)(7,9)(8,9)(8,8.9)(7,8.9)(7.1,8.85)(5,7.85)(4.5,6.85)(3.6,6.9)(1.5,5.85)(1,4.9)(.1,4.85)(-1,4.4)


\drawpolygon[fillcolor=black](2,.5)(2.5,1.5)(3.5,2)(4.5,2)(5.5,2.5)(6,3.5)(7,4)(8,4)(8,3.9)(7,3.9)(6,3.4)(5.5,2.4)(4.5,1.9)(3.5,1.9)(2.5,1.4)(2,.4)

\drawpolygon[fillcolor=black](-1,3.5)(0,4)(1,4)(2,4.5)(2.5,5.5)(3.5,6)(4.5,6)(5.5,6.5)(6,7.5)(7,8)(8,8)(8,7.85)(7,7.85)(6.1,7.35)(5.6,6.35)(4.5,5.85)(3.5,5.85)(2.6,5.35)(2.1,4.35)(1,3.85)(0,3.85)(-1,3.35)

\drawpolygon[fillcolor=green](4.5,1)(4.5,3)(5,4)(4.5,5)(4.5,7)(5,8)(4.5,9)(4.65,9)(5.15,8)(4.65,7)(4.65,5)(5.15,4)(4.65,3)(4.65,1)

\drawpolygon[fillcolor=green](3.5,1)(3.5,3)(3,4)(3.5,5)(3.5,7)(3,8)(3.5,9)(3.65,9)(3.15,8)(3.65,7)(3.65,5)(3.15,4)(3.65,3)(3.65,1)

\drawpolygon[fillcolor=green](1,0)(1,1)(1.5,2)(1,3)(1,5)(1.5,6)(1,7)(1,9)(1.15,9)(1.15,7)(1.65,6)(1.15,5)(1.15,3)(1.65,2)(1.15,1)(1.15,0)

\drawpolygon[fillcolor=green](0,0)(0,1)(-.5,2)(0,3)(0,5)(-.5,6)(0,7)(0,9)(.15,9)(.15,7)(-.35,6)(.15,5)(.15,3)(-.35,2)(.15,1)(.15,0)

\drawpolygon[fillcolor=green](7,0)(7,1)(6.5,2)(7,3)(7,5)(6.5,6)(7,7)(7,9)(7.15,9)(7.15,7)(6.65,6)(7.15,5)(7.15,3)(6.65,2)(7.15,1)(7.15,0)

\drawpolygon[fillcolor=green](6,.5)(5.5,1.5)(5.5,2.5)(6,3.5)(6,4.5)(5.5,5.5)(5.5,6.5)(6,7.5)(6,8.5)(5.5,9.5)(5.65,9.5)(6.15,8.5)(6.15,7.4)(5.65,6.5)(5.65,5.5)(6.15,4.5)(6.15,3.5)(5.65,2.5)(5.65,1.5)(6.15,.5)

\drawpolygon[fillcolor=green](2,.5)(2.5,1.5)(2.5,2.5)(2,3.5)(2,4.5)(2.5,5.5)(2.5,6.5)(2,7.5)(2,8.5)(2.5,9.5)(2.65,9.5)(2.15,8.5)(2.15,7.4)(2.65,6.5)(2.65,5.5)(2.15,4.5)(2.15,3.5)(2.65,2.5)(2.65,1.5)(2.15,.5)

\put(-2,-1.5){\scriptsize {\tiny {\bf Figure 4.17.1:} Paths of type $S_1$  }} 

\end{picture}

\begin{picture}(0,0)(-25,48)
\setlength{\unitlength}{3.5mm}

\drawpolygon(0,0)(2,0)(3,-.5)(4,0)(6,0)(7,-.5)(8,0)(10,0)(11,-.5)(12,0)(12,1)(11,1.5)(10,1)(8,1)(7,1.5)(6,1)(4,1)(3,1.5)(2,1)(0,1)

\drawpolygon(1,3)(2,3.5)(4,3.5)(5,3)(6,3.5)(8,3.5)(9,3)(10,3.5)(12,3.5)(13,3)(14,3.5)(14,4.5)(13,5)(12,4.5)(10,4.5)(9,5)(8,4.5)(6,4.5)(5,5)(4,4.5)(2,4.5)(2,3.5)

\drawpolygon(3,6.5)(4,7)(6,7)(7,6.5)(8,7)(10,7)(11,6.5)(12,7)(14,7)(15,6.5)(16,7)(16,8)(15,8.5)(14,8)(12,8)(11,8.5)(10,8)(8,8)(7,8.5)(6,8)(4,8)(4,7)


\drawline[AHnb=0](0,1)(1,3)(2,1)(1,1)(.5,2)(1.5,2)(1,1)
\drawline[AHnb=0](4,1)(5,3)(6,1)(5,1)(4.5,2)(5.5,2)(5,1)
\drawline[AHnb=0](8,1)(9,3)(10,1)(9,1)(8.5,2)(9.5,2)(9,1)

\drawline[AHnb=0](2,4.5)(3,6.5)(4,4.5)(3,4.5)(2.5,5.5)(3.5,5.5)(3,4.5)
\drawline[AHnb=0](6,4.5)(7,6.5)(8,4.5)(7,4.5)(6.5,5.5)(7.5,5.5)(7,4.5)
\drawline[AHnb=0](10,4.5)(11,6.5)(12,4.5)(11,4.5)(10.5,5.5)(11.5,5.5)(11,4.5)

\drawline[AHnb=0](4,8)(5,10)(6,8)(5,8)(4.5,9)(5.5,9)(5,8)
\drawline[AHnb=0](8,8)(9,10)(10,8)(9,8)(8.5,9)(9.5,9)(9,8)
\drawline[AHnb=0](12,8)(13,10)(14,8)(13,8)(12.5,9)(13.5,9)(13,8)


\drawline[AHnb=0](3,1.5)(2,3.5)(4,3.5)(3,1.5)(2.5,2.5)(3.5,2.5)(3,3.5)(2.5,2.5)(1.5,2)

\drawline[AHnb=0](7,1.5)(6,3.5)(8,3.5)(7,1.5)(6.5,2.5)(7.5,2.5)(7,3.5)(6.5,2.5)(5.5,2)
\drawline[AHnb=0](11,1.5)(10,3.5)(12,3.5)(11,1.5)(10.5,2.5)(11.5,2.5)(11,3.5)(10.5,2.5)(9.5,2)

\drawline[AHnb=0](5,5)(4,7)(6,7)(5,5)(4.5,6)(5.5,6)(5,7)(4.5,6)(3.5,5.5)
\drawline[AHnb=0](9,5)(8,7)(10,7)(9,5)(8.5,6)(9.5,6)(9,7)(8.5,6)(7.5,5.5)
\drawline[AHnb=0](13,5)(12,7)(14,7)(13,5)(12.5,6)(13.5,6)(13,7)(12.5,6)(11.5,5.5)

\drawline[AHnb=0](7,8.5)(6,10.5)(8,10.5)(7,8.5)(6.5,9.5)(7.5,9.5)(7,10.5)(6.5,9.5)(5.5,9)
\drawline[AHnb=0](11,8.5)(10,10.5)(12,10.5)(11,8.5)(10.5,9.5)(11.5,9.5)(11,10.5)(10.5,9.5)(9.5,9)
\drawline[AHnb=0](15,8.5)(14,10.5)(16,10.5)(15,8.5)(14.5,9.5)(15.5,9.5)(15,10.5)(14.5,9.5)(13.5,9)


\drawline[AHnb=0](8,1)(9,3)
\drawline[AHnb=0](10,4.5)(11,6.5)
\drawline[AHnb=0](12,8)(13,10)(12,10.5)

\drawline[AHnb=0](12,1)(13,3)
\drawline[AHnb=0](14,4.5)(15,6.5)
\drawline[AHnb=0](16,8)(17,10)(16,10.5)

\drawline[AHnb=0](3.5,2.5)(4.5,2)
\drawline[AHnb=0](5.5,6)(6.5,5.5)
\drawline[AHnb=0](7.5,9.5)(8.5,9)

\drawline[AHnb=0](7.5,2.5)(8.5,2)
\drawline[AHnb=0](9.5,6)(10.5,5.5)
\drawline[AHnb=0](11.5,9.5)(12.5,9)

\drawline[AHnb=0](11.5,2.5)(12.5,2)
\drawline[AHnb=0](13.5,6)(14.5,5.5)
\drawline[AHnb=0](15.5,9.5)(16.5,9)

\drawline[AHnb=0](5,10)(6,10.5)
\drawline[AHnb=0](8,10.5)(9,10)(10,10.5)
\drawline[AHnb=0](13,10)(14,10.5)

\drawpolygon[fillcolor=black](6.5,9.6)(5.5,9)(4.5,9)(4,8)(4,7)(3,6.5)(2,4.5)(2,3.5)(1,3)(1.5,2)(.5,2)(0,1)(0,0)(0.1,0)(.1,1)(.6,1.9)(1.65,1.9)(1.2,2.9)(2.1,3.4)(2.1,4.4)(3.1,6.4)(4.1,6.9)(4.1,7.9)(4.6,8.9)(5.6,8.9)(6.5,9.5)

\drawpolygon[fillcolor=black](12,0)(11,-.5)(10,0)(8,0)(7,-.5)(6,0)(4,0)(3,-.5)(2,0)(1,0)(1,1)(2,1)(3,1.5)(4,1)(6,1)(7,1.5)(8,1)(10,1)(11,1.5)(11.5,2.5)(10.5,2.5)(9.5,2)(8.5,2)(7.5,2.5)(6.5,2.5)(5.5,2)(4.5,2)(3.5,2.5)(2.5,2.5)(3,3.5)(4,3.5)(5,3)(6,3.5)(8,3.5)(9,3)(10,3.5)(12,3.5)(12,4.5)(10,4.5)(9,5)(8,4.5)(6,4.5)(5,5)(4,4.5)(3,4.5)(3.5,5.5)(4.5,6)(5.5,6)(6.5,5.5)(7.5,5.5)(8.5,6)(9.5,6)(10.5,5.5)(11.5,5.5)(12.5,6)(13,5)(13.5,6)(14,7)(12,7)(11,6.5)(10,7)(8,7)(7,6.5)(6,7)(5,7)(5,8)(6,8)(7,8.5)(8,8)(10,8)(11,8.5)(12,8)(14,8)(15,8.5)(15.5,9.5)(14.5,9.5)(13.5,9)(12.5,9)(11.5,9.5)(10.5,9.5)(9.5,9)(8.5,9)(7.5,9.5)(6.5,9.5)(6.6,9.6)(7.6,9.6)(8.6,9.1)(9.5,9.1)(10.5,9.6)(11.5,9.6)(12.5,9.1)(13.5,9.1)(14.5,9.6)(15.6,9.6)(15.1,8.4)(14.1,7.9)(12,7.9)(11,8.4)(10,7.9)(8,7.9)(7,8.4)(6,7.9)(5.1,7.9)(5.1,7.1)(6.1,7.1)(7,6.6)(8,7.1)(10,7.1)(11,6.6)(12,7.1)(14.12,7.1)(13,4.85)(12.45,5.9)(11.5,5.4)(10.5,5.4)(9.5,5.9)(8.6,5.9)(7.5,5.4)(6.5,5.4)(5.5,5.9)(4.5,5.9)(3.6,5.4)(3.1,4.6)(4.1,4.6)(5,5.1)(6,4.6)(8,4.6)(9,5.1)(10,4.6)(12.1,4.6)(12.1,3.4)(10,3.4)(9,2.9)(8,3.4)(6,3.4)(5,2.9)(4,3.4)(3.1,3.4)(2.6,2.6)(3.6,2.6)(4.6,2.1)(5.5,2.1)(6.5,2.6)(7.5,2.6)(8.5,2.1)(9.5,2.1)(10.5,2.6)(11.6,2.6)(11.1,1.4)(10,.9)(8,.9)(7,1.4)(6,.9)(4,.9)(3,1.4)(2,.9)(1.1,.9)(1.1,.1)(2.1,.1)(3,-.4)(4,.1)(6,.1)(7,-.4)(8,.1)(10,.1)(11,-.4)(12,.1)

\drawline[AHnb=0](1,0)(1,1)
\drawline[AHnb=0](2,0)(2,1)
\drawline[AHnb=0](4,0)(4,1)
\drawline[AHnb=0](5,0)(5,1)
\drawline[AHnb=0](6,0)(6,1)

\drawline[AHnb=0](8,0)(8,1)
\drawline[AHnb=0](9,0)(9,1)
\drawline[AHnb=0](10,0)(10,1)

\drawline[AHnb=0](3,3.5)(3,4.5)
\drawline[AHnb=0](4,3.5)(4,4.5)
\drawline[AHnb=0](6,3.5)(6,4.5)
\drawline[AHnb=0](7,3.5)(7,4.5)
\drawline[AHnb=0](8,3.5)(8,4.5)
\drawline[AHnb=0](10,3.5)(10,4.5)
\drawline[AHnb=0](11,3.5)(11,4.5)
\drawline[AHnb=0](12,3.5)(12,4.5)

\drawline[AHnb=0](5,7)(5,8)
\drawline[AHnb=0](6,7)(6,8)
\drawline[AHnb=0](8,7)(8,8)
\drawline[AHnb=0](9,7)(9,8)
\drawline[AHnb=0](10,7)(10,8)
\drawline[AHnb=0](12,7)(12,8)
\drawline[AHnb=0](13,7)(13,8)
\drawline[AHnb=0](14,7)(14,8)

\put(-.8,-.6){\scriptsize {\tiny $a_{11}$}}
\put(.35,-.6){\scriptsize {\tiny $a_{12}$}}
\put(10.6,-1){\scriptsize {\tiny $a_{1i}$}}
\put(11.8,-.6){\scriptsize {\tiny $a_{11}$}}

\put(-1.45,.75){\scriptsize {\tiny $a_{21}$}}
\put(12.35,.7){\scriptsize {\tiny $a_{21}$}}

\put(-.95,1.9){\scriptsize {\tiny$a_{31}$}}
\put(12.75,1.7){\scriptsize {\tiny $a_{31}$}}

\put(-.45,2.9){\scriptsize {\tiny $a_{41}$}}
\put(.5,3.7){\scriptsize {\tiny $a_{42}$}}
\put(13.1,2.65){\scriptsize {\tiny $a_{41}$}}
\put(14,3.1){\scriptsize {\tiny $a_{42}$}}

\put(.5,4.4){\scriptsize {\tiny $a_{51}$}}
\put(14.2,4.4){\scriptsize {\tiny $a_{51}$}}

\put(1,5.25){\scriptsize {\tiny $a_{61}$}}
\put(14.7,5.2){\scriptsize {\tiny $a_{61}$}}

\put(2.85,9){\scriptsize {\tiny $a_{j1}$}}
\put(16.85,9){\scriptsize {\tiny $a_{j1}$}}

\put(3,10.5){\scriptsize {\tiny $a_{1(k+1)}$}}
\put(5.35,11){\scriptsize {\tiny $a_{1(k+2)}$}}
\put(15.75,11.15){\scriptsize {\tiny $a_{1k}$}}
\put(17,10.55){\scriptsize {\tiny $a_{1(k+1)}$}}

\put(-3,-.1){\scriptsize {\tiny $Q_{1}$}}
\put(-2.75,.85){\scriptsize {\tiny $Q_{2}$}}

\put(1.5,9.15){\scriptsize {\tiny $Q_{j}$}}
\put(1.8,10.35){\scriptsize {\tiny $Q_{1}$}}

\put(4.8,-1){\scriptsize $\ldots$}
\put(7.2,-1){\scriptsize $\ldots$} 
\put(12.5,11.25){\scriptsize $\ldots$}
\put(9,11.25){\scriptsize $\ldots$} 

\put(-.3,5.2){ $.$}
\put(-.5,4.8){ $.$}
\put(-.8,4.3){ $.$}

\put(-.5,-2.25){\scriptsize {\tiny {\bf Figure 4.17.2:} $M(i,j=6m+3,k)$ }} 

\end{picture}

\begin{picture}(0,0)(-83.5,43)
\setlength{\unitlength}{3.5mm}

\drawpolygon(0,0)(2,0)(3,-.5)(4,0)(6,0)(7,-.5)(8,0)(10,0)(11,-.5)(12,0)(12,1)(11,1.5)(10,1)(8,1)(7,1.5)(6,1)(4,1)(3,1.5)(2,1)(0,1)

\drawpolygon(1,3)(2,3.5)(4,3.5)(5,3)(6,3.5)(8,3.5)(9,3)(10,3.5)(12,3.5)(13,3)(14,3.5)(14,4.5)(13,5)(12,4.5)(10,4.5)(9,5)(8,4.5)(6,4.5)(5,5)(4,4.5)(2,4.5)(2,3.5)

\drawpolygon(3,6.5)(4,7)(6,7)(7,6.5)(8,7)(10,7)(11,6.5)(12,7)(14,7)(15,6.5)(16,7)(16,8)(15,8.5)(14,8)(12,8)(11,8.5)(10,8)(8,8)(7,8.5)(6,8)(4,8)(4,7)

\drawline[AHnb=0](8,10.5)(9,10)(10,10.5)(12,10.5)(13,10)(14,10.5)(16,10.5)(17,10)(18,10.5)(18,11.5)(17,12)(16,11.5)(14,11.5)(13,12)(12,11.5)(10,11.5)(9,12)(8,11.5)(6,11.5)(6,10.5)


\drawline[AHnb=0](0,1)(1,3)(2,1)(1,1)(.5,2)(1.5,2)(1,1)
\drawline[AHnb=0](4,1)(5,3)(6,1)(5,1)(4.5,2)(5.5,2)(5,1)
\drawline[AHnb=0](8,1)(9,3)(10,1)(9,1)(8.5,2)(9.5,2)(9,1)

\drawline[AHnb=0](2,4.5)(3,6.5)(4,4.5)(3,4.5)(2.5,5.5)(3.5,5.5)(3,4.5)
\drawline[AHnb=0](6,4.5)(7,6.5)(8,4.5)(7,4.5)(6.5,5.5)(7.5,5.5)(7,4.5)
\drawline[AHnb=0](10,4.5)(11,6.5)(12,4.5)(11,4.5)(10.5,5.5)(11.5,5.5)(11,4.5)

\drawline[AHnb=0](4,8)(5,10)(6,8)(5,8)(4.5,9)(5.5,9)(5,8)
\drawline[AHnb=0](8,8)(9,10)(10,8)(9,8)(8.5,9)(9.5,9)(9,8)
\drawline[AHnb=0](12,8)(13,10)(14,8)(13,8)(12.5,9)(13.5,9)(13,8)

\drawline[AHnb=0](6,11.5)(7,13.5)(8,11.5)(7,11.5)(6.5,12.5)(7.5,12.5)(7,11.5)
\drawline[AHnb=0](10,11.5)(11,13.5)(12,11.5)(11,11.5)(10.5,12.5)(11.5,12.5)(11,11.5)
\drawline[AHnb=0](14,11.5)(15,13.5)(16,11.5)(15,11.5)(14.5,12.5)(15.5,12.5)(15,11.5)


\drawline[AHnb=0](3,1.5)(2,3.5)(4,3.5)(3,1.5)(2.5,2.5)(3.5,2.5)(3,3.5)(2.5,2.5)(1.5,2)

\drawline[AHnb=0](7,1.5)(6,3.5)(8,3.5)(7,1.5)(6.5,2.5)(7.5,2.5)(7,3.5)(6.5,2.5)(5.5,2)
\drawline[AHnb=0](11,1.5)(10,3.5)(12,3.5)(11,1.5)(10.5,2.5)(11.5,2.5)(11,3.5)(10.5,2.5)(9.5,2)

\drawline[AHnb=0](5,5)(4,7)(6,7)(5,5)(4.5,6)(5.5,6)(5,7)(4.5,6)(3.5,5.5)
\drawline[AHnb=0](9,5)(8,7)(10,7)(9,5)(8.5,6)(9.5,6)(9,7)(8.5,6)(7.5,5.5)
\drawline[AHnb=0](13,5)(12,7)(14,7)(13,5)(12.5,6)(13.5,6)(13,7)(12.5,6)(11.5,5.5)

\drawline[AHnb=0](7,8.5)(6,10.5)(8,10.5)(7,8.5)(6.5,9.5)(7.5,9.5)(7,10.5)(6.5,9.5)(5.5,9)
\drawline[AHnb=0](11,8.5)(10,10.5)(12,10.5)(11,8.5)(10.5,9.5)(11.5,9.5)(11,10.5)(10.5,9.5)(9.5,9)
\drawline[AHnb=0](15,8.5)(14,10.5)(16,10.5)(15,8.5)(14.5,9.5)(15.5,9.5)(15,10.5)(14.5,9.5)(13.5,9)

\drawline[AHnb=0](9,12)(8,14)(10,14)(9,12)(8.5,13)(9.5,13)(9,14)(8.5,13)(7.5,12.5)
\drawline[AHnb=0](13,12)(12,14)(14,14)(13,12)(12.5,13)(13.5,13)(13,14)(12.5,13)(11.5,12.5)
\drawline[AHnb=0](17,12)(16,14)(18,14)(17,12)(16.5,13)(17.5,13)(17,14)(16.5,13)(15.5,12.5)


\drawline[AHnb=0](8,1)(9,3)
\drawline[AHnb=0](10,4.5)(11,6.5)
\drawline[AHnb=0](12,8)(13,10)(12,10.5)

\drawline[AHnb=0](12,1)(13,3)
\drawline[AHnb=0](14,4.5)(15,6.5)
\drawline[AHnb=0](16,8)(17,10)(16,10.5)
\drawline[AHnb=0](18,11.5)(19,13.5)(18,14)

\drawline[AHnb=0](3.5,2.5)(4.5,2)
\drawline[AHnb=0](5.5,6)(6.5,5.5)
\drawline[AHnb=0](7.5,9.5)(8.5,9)
\drawline[AHnb=0](9.5,13)(10.5,12.5)
\drawline[AHnb=0](13.5,13)(14.5,12.5)
\drawline[AHnb=0](17.5,13)(18.5,12.5)

\drawline[AHnb=0](7.5,2.5)(8.5,2)
\drawline[AHnb=0](9.5,6)(10.5,5.5)
\drawline[AHnb=0](11.5,9.5)(12.5,9)

\drawline[AHnb=0](11.5,2.5)(12.5,2)
\drawline[AHnb=0](13.5,6)(14.5,5.5)
\drawline[AHnb=0](15.5,9.5)(16.5,9)

\drawline[AHnb=0](5,10)(6,10.5)
\drawline[AHnb=0](8,10.5)(9,10)(10,10.5)

\drawline[AHnb=0](7,13.5)(8,14)

\drawpolygon[fillcolor=black](6.5,12.5)(6,11.5)(6,10.5)(5,10)(5.5,9)(4.5,9)(4,8)(4,7)(3,6.5)(2,4.5)(2,3.5)(1,3)(1.5,2)(.5,2)(0,1)(0,0)(0.1,0)(.1,1)(.6,1.9)(1.65,1.9)(1.2,2.9)(2.1,3.4)(2.1,4.4)(3.1,6.4)(4.1,6.9)(4.1,7.9)(4.65,8.9)(5.65,8.9)(5.15,9.9)(6.1,10.4)(6.1,11.4)(6.6,12.4)

\drawpolygon[fillcolor=black](12,0)(11,-.5)(10,0)(8,0)(7,-.5)(6,0)(4,0)(3,-.5)(2,0)(1,0)(1,1)(2,1)(3,1.5)(4,1)(6,1)(7,1.5)(8,1)(10,1)(11,1.5)(11.5,2.5)(10.5,2.5)(9.5,2)(8.5,2)(7.5,2.5)(6.5,2.5)(5.5,2)(4.5,2)(3.5,2.5)(2.5,2.5)(3,3.5)(4,3.5)(5,3)(6,3.5)(8,3.5)(9,3)(10,3.5)(12,3.5)(12,4.5)(10,4.5)(9,5)(8,4.5)(6,4.5)(5,5)(4,4.5)(3,4.5)(3.5,5.5)(4.5,6)(5.5,6)(6.5,5.5)(7.5,5.5)(8.5,6)(9.5,6)(10.5,5.5)(11.5,5.5)(12.5,6)(13,5)(13.5,6)(14,7)(12,7)(11,6.5)(10,7)(8,7)(7,6.5)(6,7)(5,7)(5,8)(6,8)(7,8.5)(8,8)(10,8)(11,8.5)(12,8)(14,8)(15,8.5)(15.5,9.5)(14.5,9.5)(13.5,9)(12.5,9)(11.5,9.5)(10.5,9.5)(9.5,9)(8.5,9)(7.5,9.5)(6.5,9.5)(7,10.5)(8,10.5)(9,10)(10,10.5)(12,10.5)(13,10)(14,10.5)(16,10.5)(16,11.5)(14,11.5)(13,12)(12,11.5)(10,11.5)(9,12)(8,11.5)(7,11.5)(7.5,12.5)(8.5,13)(9.5,13)(10.5,12.5)(11.5,12.5)(12.5,13)(13.5,13)(14.5,12.5)(15.5,12.5)(16.5,13)(17,12)(17.45,13)(18.45,12.5)(18.45,12.4)(17.55,12.9)(17,11.8)(16.4,12.9)(15.5,12.4)(14.4,12.4)(13.5,12.9)(12.5,12.9)(11.5,12.4)(10.5,12.4)(9.5,12.9)(8.5,12.9)(7.6,12.4)(7.1,11.6)(8,11.6)(9,12.1)(10,11.6)(12,11.6)(13,12.1)(14,11.6)(16.1,11.6)(16.1,10.4)(14,10.4)(13,9.9)(12,10.4)(10,10.4)(9,9.9)(8,10.4)(7.1,10.4)(6.6,9.6)(7.6,9.6)(8.6,9.1)(9.5,9.1)(10.5,9.6)(11.5,9.6)(12.5,9.1)(13.5,9.1)(14.5,9.6)(15.6,9.6)(15.1,8.4)(14.1,7.9)(12,7.9)(11,8.4)(10,7.9)(8,7.9)(7,8.4)(6,7.9)(5.1,7.9)(5.1,7.1)(6.1,7.1)(7,6.6)(8,7.1)(10,7.1)(11,6.6)(12,7.1)(14.12,7.1)(13,4.85)(12.45,5.9)(11.5,5.4)(10.5,5.4)(9.5,5.9)(8.6,5.9)(7.5,5.4)(6.5,5.4)(5.5,5.9)(4.5,5.9)(3.6,5.4)(3.1,4.6)(4.1,4.6)(5,5.1)(6,4.6)(8,4.6)(9,5.1)(10,4.6)(12.1,4.6)(12.1,3.4)(10,3.4)(9,2.9)(8,3.4)(6,3.4)(5,2.9)(4,3.4)(3.1,3.4)(2.6,2.6)(3.6,2.6)(4.6,2.1)(5.5,2.1)(6.5,2.6)(7.5,2.6)(8.5,2.1)(9.5,2.1)(10.5,2.6)(11.6,2.6)(11.1,1.4)(10,.9)(8,.9)(7,1.4)(6,.9)(4,.9)(3,1.4)(2,.9)(1.1,.9)(1.1,.1)(2.1,.1)(3,-.4)(4,.1)(6,.1)(7,-.4)(8,.1)(10,.1)(11,-.4)(12,.1)

\drawline[AHnb=0](16,14)(15,13.5)(14,14)
\drawline[AHnb=0](12,14)(11,13.5)(10,14)

\drawline[AHnb=0](1,0)(1,1)
\drawline[AHnb=0](2,0)(2,1)
\drawline[AHnb=0](4,0)(4,1)
\drawline[AHnb=0](5,0)(5,1)
\drawline[AHnb=0](6,0)(6,1)

\drawline[AHnb=0](8,0)(8,1)
\drawline[AHnb=0](9,0)(9,1)
\drawline[AHnb=0](10,0)(10,1)

\drawline[AHnb=0](3,3.5)(3,4.5)
\drawline[AHnb=0](4,3.5)(4,4.5)
\drawline[AHnb=0](6,3.5)(6,4.5)
\drawline[AHnb=0](7,3.5)(7,4.5)
\drawline[AHnb=0](8,3.5)(8,4.5)
\drawline[AHnb=0](10,3.5)(10,4.5)
\drawline[AHnb=0](11,3.5)(11,4.5)
\drawline[AHnb=0](12,3.5)(12,4.5)

\drawline[AHnb=0](5,7)(5,8)
\drawline[AHnb=0](6,7)(6,8)
\drawline[AHnb=0](8,7)(8,8)
\drawline[AHnb=0](9,7)(9,8)
\drawline[AHnb=0](10,7)(10,8)
\drawline[AHnb=0](12,7)(12,8)
\drawline[AHnb=0](13,7)(13,8)
\drawline[AHnb=0](14,7)(14,8)

\drawline[AHnb=0](7,10.5)(7,11.5)

\drawline[AHnb=0](8,10.5)(8,11.5)

\drawline[AHnb=0](10,10.5)(10,11.5)

\drawline[AHnb=0](11,10.5)(11,11.5)

\drawline[AHnb=0](12,10.5)(12,11.5)

\drawline[AHnb=0](14,10.5)(14,11.5)

\drawline[AHnb=0](15,10.5)(15,11.5)

\drawline[AHnb=0](16,10.5)(16,11.5)

\put(-.8,-.6){\scriptsize {\tiny $a_{11}$}}
\put(.25,-.6){\scriptsize {\tiny $a_{12}$}}

\put(10.6,-1){\scriptsize {\tiny $a_{1i}$}}
\put(11.8,-.6){\scriptsize {\tiny $a_{11}$}}

\put(-1.5,.75){\scriptsize {\tiny $a_{21}$}}
\put(12.35,.7){\scriptsize {\tiny $a_{21}$}}

\put(-.95,1.9){\scriptsize {\tiny$a_{31}$}}
\put(12.75,1.7){\scriptsize {\tiny $a_{31}$}}

\put(-.5,2.9){\scriptsize {\tiny $a_{41}$}}
\put(.5,3.7){\scriptsize {\tiny $a_{42}$}}
\put(13.1,2.65){\scriptsize {\tiny $a_{41}$}}
\put(14,3.1){\scriptsize {\tiny $a_{42}$}}

\put(.5,4.4){\scriptsize {\tiny $a_{51}$}}
\put(14.2,4.4){\scriptsize {\tiny $a_{51}$}}

\put(1,5.25){\scriptsize {\tiny $a_{61}$}}
\put(14.7,5.2){\scriptsize {\tiny $a_{61}$}}

\put(1.5,6.4){\scriptsize {\tiny $a_{71}$}}
\put(2.5,7.2){\scriptsize {\tiny$a_{72}$}}
\put(15.4,6.2){\scriptsize {\tiny $a_{71}$}}
\put(16.2,7.1){\scriptsize {\tiny $a_{72}$}}

\put(2.5,8){\scriptsize {\tiny $a_{81}$}}
\put(16.2,8){\scriptsize {\tiny $a_{81}$}}

\put(4.65,12.35){\scriptsize {\tiny $a_{j1}$}}
\put(18.75,12.2){\scriptsize {\tiny $a_{j1}$}}

\put(-3,-.1){\scriptsize {\tiny $Q_{1}$}}
\put(-2.75,.85){\scriptsize {\tiny $Q_{2}$}}

\put(4,13.5){\scriptsize {\tiny $Q_{1}$}}
\put(3.5,12.45){\scriptsize {\tiny $Q_{j}$}}

\put(4.8,-1){\scriptsize $\ldots$}
\put(7.2,-1){\scriptsize $\ldots$} 
\put(15,14.65){\scriptsize $\ldots$}
\put(11.5,14.65){\scriptsize $\ldots$}

\put(2.7,11.2){\scriptsize $.$}
\put(2.6,11){\scriptsize $.$}
\put(2.5,10.8){\scriptsize $.$}

\put(5.15,14){\scriptsize {\tiny $a_{1(k+1)}$}}
\put(7.25,14.5){\scriptsize {\tiny $a_{1(k+2)}$}}
\put(17.5,14.5){\scriptsize {\tiny $a_{1k}$}}
\put(18.75,14){\scriptsize {\tiny $a_{1(k+1)}$}}

\put(-.5,-2.25){\scriptsize {\tiny {\bf Figure 4.17.3:} $M(i,j=6m+6,k)$ }} 

\end{picture}

\vspace{5.5cm}

An $M(i,j,k)$ representation of $M$ follows by cutting $M$ one by one along a black and green colored cycles of type $S_1$ through a vertex $v$ with face-sequence $(3,4, 6,4)$. Without loss of generality, let the beginning adjacent face to the base horizontal cycle is a 4-gon.

\begin{lem}\label{l3.17.1} The DSEM $M$ of type $[3^3.4^2:3.4.6.4]$ admits an $M(i,j,k)$-representation iff the following holds: $(i)$ $j =3m$, where $m \in \mathbb{N}$,  $(ii)$ $i \geq 4$ and $i=4m$, $m \in \mathbb{N} $, $(iii)$ number of vertices of $M(i,j,k) = ij \geq 12 $, $(iv)$  $ k \in \{4r+3: 0 \leq r < i/4\}$.
\end{lem}

\noindent{\bf Proof.} Let $M$ be a DSEM of type $[3^3.4^2:3.4.6.4]$ having $n$ vertices. Clearly, its representation $M(i,j,k)$ has $j$ disjoint horizontal cycles of $S_{1}$ type, say $Q_0, Q_1,\ldots, Q_{j-1}$, of length $i$. Let $Q_{0}=C(w_{0,0},w_{0,1},\ldots,  w_{0,i-1}),Q_{1}=C(w_{1,0},w_{1,1},\ldots,w_{1,i-1}),\ldots,Q_{j-1}=C(w_{j-1,0},w_{j-1,1},\ldots,w_{j-1,i-1})$ be the list of horizontal cycles. Since the vertices of $M$ are in these cycles, the vertices in $M$ is $n = ij$. If $j = 1,2$ then $M$ is not a map. So $j \geq 3$. If $j \geq 3$ and $j \neq 3m, m \in \mathbb{N}$, then no vertex in the base horizontal cycle follow the face-sequences $(3^3,4^2)$ and $(3,4,6,4)$ after identifying the boundaries of $M(i,j,k)$. So $j = 3m, \text{ where } m \in \mathbb{N}$. Clearly if $i \leq 3$, $M$ is not a map. So $i \geq 4$. If $j=3m$ and $i \neq 4m$, where $m \in \mathbb{N}$, then the ${\rm lk}(w_{0,0})$ is not of type $(3,4,6,4)$. Which is not possible. So, $i = 4m$, where $m \in \mathbb{N}$. Thus, $ n = ij \geq 12$.

If $ k \in  \{r:0 \leq r \leq i-1\} \setminus (\{4r+3: 0 \leq r < i/4\}) $, then link of some vertices in $M(i,j,k)$ can not be completed. So, $ k \in \{4r+3: 0 \leq r < i/4\}$. This completes the proof. \hfill $\Box$

%

\vspace{.2cm}

\noindent{\bf Proof of Theorem \ref{t1}} Let $M$ be a DSEM of type $Y_1$, where $Y_1 \in \{[3^6:3^3.4^2]_1, [3^6:3^3.4^2]_2, [3^3.4^2:4^4]_1, [3^3.4^2:4^4]_2, [3^3.4^2:3^2.4.3.4]_2\}$ with its $M(i,j,k)$ representation (shown in Figures 4.1.3, 4.1.4, 4.2.2, 4.2.3 and 4.4.2). Now, if we remove the diagonal edges (indicated by thick black edges), we get a map, say $M^{r'}$, on the torus of type $[4^4]$. Clearly, $V(M^r) = V(M^{r'})$ and $E(M^{r'}) \subset E(M^r)$. By \cite[Theorem 4]{altshuler(1972)}, the map $M^{r'}$ contains a Hamiltonian cycle, say $H$. The cycle $H$ is also a Hamiltonian cycle in $M^r$ as $E(M^{r'}) \subset E(M^r)$. Hence, $M$ is Hamiltonian. \hfill$\Box$

Let $M$ be a DSEM of the type $[3.4^2.6: 3.6.3.6]_2$, $[3^6:3^2.6^2]$, $[3^2.4.3.4:3.6.4.6]$, $[3^6:3^4.6]_2$, $[3^6:3^2.4.12]$, $[3.4.6.4: 4.6.12]$ or $[3.4^2.6:3.4.6.4]$. Then the thick black cycle in its $M(i,j,k)$, shown respectively in Figures 4.6.4, 4.8.2, 4.10.2, 4.12.2, 4.13.2, 4.15.2 or 4.16.2, is Hamiltonian. Thus $M$ is Hamiltonian. 

Similarly, if $M$ is one of the types, $[3^3.4^2:3^2.4.3.4]_1$, $[3^6:3^2.4.3.4]$, $[3.4^2.6: 3.6.3.6]_2$, $[3^2.6^2:3.6.3.6]$, $[3^4.6:3^2.6^2]$, $[3^6:3^4.6]_1$, $[3.4.3.12: 3.12^2]$ or $[3^3.4^2:3.4.6.4]$, then depending on $j$, we get its $M(i,j,k)$ representation, shown in respective subsection. Note that the black thick cycle drawn in such representation is Hamiltonian. Hence $M$ is Hamiltonian. \hfill$\Box$


\section{Conclusion}

In this article, we have discussed the connectivity of DSEMs on the torus corresponding to the twenty 2-uniform tilings of the plane and shown that every such DSEM is either 3-connected or 4-connected. Using this, we have established the Nash-Williams conjecture for such maps. Recall that, in a map $M$, the combinatorial curvature of a vertex $v$ with the face-sequence $(p_1^{n_1}, \ldots, p_k^{n_k})$ is given by $\phi(v) = 1 - (\sum_{i=1}^k n_i)/2 + (\sum_{i=1}^k n_i)/p_i$. Note that the DSEMs, discussed here, have $\phi(v)=0$ for all the vertices. However, one can construct many more doubly semi-equivelar maps which may not have such curvature restriction. For instance, if we stack (subdividing a face by introducing a new vertex inside the face and joining this new vertex to each vertex of the face by an edge) all the 4-gonal faces of any semi-equivelar map of type $(4.8^2)$, then the resulting map is a DSEM with two types face-sequences $(3^4)$ and $(3^2.8^2)$, that is other than the twenty types. However, the idea used for the connectivity of DSEMs can be applied to arbitrary DSEMs on the torus. One can also explore similarly the Hamiltonicity of a class of maps on the torus corresponding to the remaining types $k$-uniform tilings for $3 \leq k \leq 7$. Thus a natural question arises for the readers to determine the complete classification of DSEMs on the torus in terms of types and to check their Hamiltonicity.

\end{document}